\documentclass[11pt]{article}
\usepackage{fullpage}
\usepackage{mathrsfs}
\usepackage{float}
\usepackage{amsmath, amsfonts, amssymb, amsthm, amsbsy, amscd, bm,
  bbm, mathtools, multirow,multicol}
\usepackage{thmtools} 
\usepackage{array}
\usepackage{booktabs}
\usepackage{graphicx}
\usepackage{tikz}
\usepackage[small,bf]{caption}
\usepackage{subcaption}
\usepackage{color}
\usepackage{ifthen}
\usepackage{xspace}
\usepackage{algorithm, algorithmic}
\usepackage[colorlinks,citecolor={blue},urlcolor={blue},linkcolor={blue}]{hyperref}
\usepackage{url}
\usepackage{tocbibind}
\usepackage{authblk}
\usepackage{geometry}
\usepackage{multicol}
\usepackage{etoc} 
\usepackage{titletoc}
\usepackage[round, sort, authoryear]{natbib}

\newcommand\E{\mathbb{E}}
\renewcommand\P{\mathbb{P}}

\newcommand\wt{\widetilde}
\newcommand{\iid}{\overset{\scriptsize{\mbox{i.i.d.}}}{\sim}}

\newcommand\cond{{\hspace{1mm} | \hspace{1mm}}}

\DeclareMathOperator*{\argmax}{argmax}

\newcommand\T{\mathcal{T}}
\newcommand\sE{\mathcal{E}}
\newcommand\A{{\mathcal{A}}}
\newcommand\at{{\mathcal{A}_t}}

\newcommand\F{\mathcal{F}}
\newcommand\C{\mathcal{C}}
\newcommand{\Fat}{{\mathcal{F}_{\mathcal{A}_t}}}
\newcommand\fat{f_{\mathcal{A}_t}}

\newcommand\cE{\mathscr{B}}
\newcommand\D{{\mathfrak{D}}}
\newcommand\cH{{\mathcal{H}}}
\newcommand\Y{{\mathcal{Y}}}

\newcommand\1{{(1)}}

\newcommand\0{{(0)}}

\newcommand{\m}{{(m)}}

\newcommand\ksup{{(k)}}
\newcommand\dsup{{(\mathfrak{D})}}
\newcommand\knew{{(k+1)}}

\newcommand\bias{{\mbox{\scriptsize{\normalfont{bias}}}}}

\newcommand{\appr}{{\mbox{\scriptsize{approx}}}}

\newcommand\tr{{\mbox{\scriptsize{tree}}}}
\newcommand\en{{\mbox{\scriptsize{forest}}}}
\newcommand\g{{\mbox{\scriptsize{\normalfont{Gibbs}}}}}

\newcommand{\Opt}{{\mbox{\tiny{OPT}}}}

\newcommand\btriplet{{\left(\D, \bm \upsilon, \{ \T_b^{h_b} \}_{b=1}^B
    \right)}}

\newcommand\forest{{\{ \T_b \}_{b=1}^B}}

\newcommand\newforest{{\mathcal{E}}}
\newcommand\newbforest{{\left( \mathcal{E} \right)}}

\newcommand\newforestp{{\mathcal{E}'}}
\newcommand\newbforestp{{\left( \mathcal{E}' \right)}}

\newcommand\newstarf{{ \mathcal{E}^\star }}
\newcommand\newbstarf{{\left( \mathcal{E}^\star \right)}}

\newcommand{\treek}{\{\mathcal{T}^{(k)}: 0 \leq k \leq \D \}}

\newcommand\kl{\mbox{KL}}

\newcommand\bentro{{H^{[\hspace{.3mm}]}}}
\newcommand\bnum{{N^{[\hspace{.3mm}]}}}

\renewcommand{\l}{\left(}
\renewcommand{\r}{\right)}

\newtheorem{theorem}{Theorem}[section]
\newtheorem{lemma}[theorem]{Lemma}
\newtheorem{corollary}[theorem]{Corollary}
\newtheorem{proposition}[theorem]{Proposition}
\theoremstyle{definition}
\newtheorem{definition}[theorem]{Definition}
\newtheorem{remark}{Remark}
\newtheorem{example}{Example}

\newtheorem{condition}{Condition}
\newtheorem{conditionalt}{Condition}[condition]
\newenvironment{conditionp}[1]{
  \renewcommand\theconditionalt{#1}
  \conditionalt
}{\endconditionalt}

\newcommand{\pref}[1]{\textcolor{black}{#1}}
\newcommand{\writing}[1]{{\color{black}{#1}}}
\newcommand{\tech}[1]{\textcolor{black}{#1}}

\newcommand{\rev}[1]{\textcolor{black}{#1}}

\makeatletter
\let\oldeqref\eqref
\renewcommand{\eqref}{%
  \@ifstar{\oldeqref}{\oldeqref}%
}
\makeatother

\makeatletter
\def\namedlabel#1#2{\begingroup
    #2%
    \def\@currentlabel{#2}%
    \phantomsection\label{#1}\endgroup
}
\makeatother

\title{On forest-type tree ensemble approaches to density
  learning under a generalized Bayesian framework}
\author[1]{Linxi Liu \thanks{linxi\_liu@pitt.edu}}
\affil[1]{Department of Statistics, University of Pittsburgh}
\author[2]{Li Ma \thanks{li.ma@uchicago.edu}}%
\affil[2]{Department of Statistics and Data Science Institute,
 University of Chicago}
\date{\today}

\begin{document}
\maketitle

\begin{abstract}
We study density forests (DFs), a class of forest-type tree ensembles designed for multivariate density learning. A DF is an aggregation of tree base learners for density functions, each growing under a small perturbation of the original data. We investigate the bias-variance trade-off of DFs 
and introduce Bayesian and generalized Bayesian strategies to enhance the data adaptivity of DFs, called Bayesian density forests (BDFs) and Gibbs density forests (GDFs), respectively. The BDF improves the data adaptivity of DFs by aggregating forests of varying depths with respect to the marginal posterior distribution on tree depth.  
The GDF, on the other hand, employs an empirical risk function for each forest topology and makes risk-based generalized Bayesian inference through the resulting Gibbs posterior to learn the forest topology from the data, thereby overcoming limitations of the likelihood in effectively quantifying the bias-variance trade-off for tree ensembles. For both BDF and GDF, we analyze their asymptotic properties, demonstrating that these ensembles can achieve faster convergence compared to single tree-based approaches for estimating smooth density functions. The
GDF, in particular, exhibits superior data adaptivity, achieving the minimax convergence rate up to a logarithmic term over 
anisotropic H\"older spaces \rev{with varying smoothness levels along
  different dimensions.} 
We identify a set of conditions
ensuring fast posterior contraction of GDFs in more general settings. Moreover, we provide empirical risk bounds that clarify the bias-variance decomposition for individual trees and forest-type ensembles in the density learning context. 
We propose efficient computational algorithms based on sequential Monte Carlo for computing the density estimates under BDF and GDF, and illustrate the practical performance of these methods through simulations in both low-dimensional and moderately high-dimensional settings.
\end{abstract}

\section{Introduction}

Tree-based methods are well-known for their effectiveness in nonparametric modeling and estimation. 
In supervised learning, for example, tree ensembles, such as random forests (RF, \citep{amit1997,ho1998,breiman2001random}), tree boosting
\citep{freund1997,friedman2001}, and Bayesian additive regression
trees (BART, \citep{chipman2010}), have become increasingly popular given their effectiveness in capturing nonlinear and interactive effects and addressing the high variance issue of single-tree learners. 

In contrast to the extensive theoretical analysis of tree ensembles in supervised contexts, tree ensembles are much less studied for unsupervised problems. Density learning---which includes both density estimation and generative modeling---is the quintessential unsupervised learning objective that can serve as the basis for a variety of downstream data analytical tasks such as synthetic data generation, sketching variations of data, profiling information content, and making risk-optimal decisions. In recent years, there has been a series of works that characterize the large-sample risk properties of single-tree density learners. Our main goal is to understand the gain of using tree ensembles over these single-tree models from a theoretical perspective. 

Although nonparametric regression and density learning are closely related problems, popular tree ensemble algorithms designed for supervised learning and their theoretical properties typically do not directly extend to the unsupervised context. One reason is that one would expect different types of spatial properties for the conditional mean versus those for the full distribution. Widely adopted assumptions for the conditional mean, such as additivity and {sparsity} in their usual sense, are typically not appropriate for the full distribution. Another reason is that probability densities are subject to additional constraints such as non-negativity and integrating-to-one. Therefore, it is of interest to investigate whether the theoretical properties of tree ensembles in the supervised setting, such as their asymptotic risk rates for function estimation, are still valid in unsupervised contexts. To this end, new frameworks for theoretical analyses are needed to characterize tree ensembles for density learning.

In this paper, we formulate such a framework from a \writing{generalized Bayesian} perspective. Our approach combines ideas from the popular supervised ensemble method random forest (RF) with Bayesian and generalized Bayesian strategies to achieve adaptive inference in density learning. We start by introducing a random-forest-type model averaging estimator for densities called a \writing{{\em density forest} (DF)}, and then enhance it with a layer of Bayesian and generalized Bayesian hierarchical modeling to fine-tune the complexity of trees in a \writing{DF}.

A first theoretical device we employ is to generalize the notion of bias-variance trade-off beyond mean squared errors and characterize it in terms of a decomposition of the likelihood function. 
Furthermore, by imposing a prior on tree structures and adapt the tree complexity via the marginal posterior on tree depth, we construct a density learner based on (Bayesian model) averaging on a collection of single-tree density estimators, which we call a \writing{{\em Bayesian density forest} (BDF)}. 
We show that the \writing{BDF} 
can achieve a faster convergence rate compared to single-tree density estimators over the isotropic H\"older class $\mathcal{C}^\alpha$ with $\alpha \in (1, 2]$. However, the rate is still \emph{suboptimal}, and we show that this is mainly due to the fact that the likelihood is \emph{not} an adequate measure of fit for tree ensembles. 

Motivated by the likelihood-based bias-variance trade-off, we address this limitation by introducing an empirical risk-based generalized Bayesian improvement. 
We define a novel risk function associated with each DF and make explicit bias-variance trade-off via \emph{information risk minimization} \citep{zhang2006a, zhang2006b}. We show, both theoretically and empirically, that this generalized Bayesian approach adapts the DF complexity to the underlying density, achieving {\em optimal} asymptotic rates over the aforementioned isotropic H\"older class without requiring \rev{oracle knowledge} of the smoothness.
We call this type of DFs as {\em Gibbs density forests} (GDFs).

Given the scarcity of theoretical analyses on tree ensembles outside supervised learning, we aim to acquire a deeper understanding of the driving forces behind good (or poor) performance of tree-based methods in density learning, which we believe will help guide the design of future methods to realize the full potential of tree ensembles in unsupervised learning. Specifically, achieving effective spatial adaptation to the potentially heterogeneous smoothness of the underlying density is a key determining factor of performance for density learners, including tree ensembles. 
For tree-based methods, splits should be made sparingly when the true density only has mild variation while being more intensively applied in regions with sharp changes. 
The spatial adaptation property of tree-based methods has drawn more attention in supervised learning \citep{veronika2020,Linero2018BayesianRT,klusowski2024,chi2022}, given its connection to variable selection. The common wisdom is that the more frequently a variable is split upon,
\rev{the more likely the response has a higher level of variation with respect to that variable.}
In the absence of a response, spatial adaptation in density learning is needed to accommodate the possibility that the underlying density function exhibits different levels of variation with respect to different variables, which can be modeled by a class of anisotropic H\"older continuous functions $\mathcal{C}^{\bm \alpha}$, where $\bm \alpha = (\alpha_1, \ldots, \alpha_d)$. 
Under our risk-based generalized Bayesian framework, we will show that \writing{GDF} can make an efficient allocation of the splits when spatial inhomogeneity is present, thereby achieving minimax convergence (up to a logarithmic term) over anisotropic H\"older classes $\mathcal{C}^{\bm \alpha}$ when $0 < \alpha_l \leq 2$.  
The rate is faster than that for single trees under certain conditions. 
In addition, we identify a set of sufficient conditions ensuring the good performance of \writing{GDF} in general.



Due to the complexity of the tree space, searching for a tree or forest topology that ``best'' fits the data can be computationally challenging. The (generalized) Bayesian framework 
naturally facilitates a stochastic search scheme. Instead of making a deterministic split at each step, one can sample a split and make adjustments with respect to the random sampling rigorously using sequential Monte Carlo (SMC). This type of SMC-based stochastic search scheme can bring significant improvements in moderately high-dimensional cases \citep{LJW}. In this work, for both types of forest estimators, we introduce such an algorithm and illustrate its performance through numerical experiments. 

We close this introduction by reviewing recent developments in tree-based approaches and related theory for both supervised and unsupervised learning tasks. Various tree-based methods have been developed for density estimation, including both Bayesian and non-Bayesian approaches. Examples of non-Bayesian techniques include multivariate histograms and their variations \citep{scott1979,scott1985ASH,pmlrcui21}, unsupervised tree boosting \citep{awaya2021}, and 
wavelet density estimation, particularly when using the Haar basis \citep{tribouley1995,donoho1996}.  
Within the Bayesian framework, several tree-type priors have been proposed to efficiently adapt to spatial and smoothness features of the underlying distribution.  These include the P{\'o}lya tree \citep{ferguson1974,lavine1992} 
and its more recent extensions \citep{wong2010,ma2017,randrianarisoa2021smoothing}. Extensive work has been done on understanding posterior concentration rates for Bayesian tree models, uncertainty quantification, and their spatial adaptation properties \citep{castillo2017,castillo2022spikeslab,castillo2022,liu2023,liu2024, randrianarisoa2021smoothing}.

In the context of supervised learning, the analysis of tree-based methods, especially the analysis of tree ensembles, has become a very active research area in recent years. Several important works have been done to reveal when tree-based methods can achieve good performance and why tree ensembles outperform individual trees. Key contributions include examining the consistency of RFs \citep{biau2012analysis,scornet2016,biau2008,denil2014,scornet2015,chi2022, syrgkanis2020}, 
making statistical inference for RFs \citep{mentch2016,athey2019,shi2024}, and establishing oracle inequalities for both classification and regression trees (CART, \citep{breiman1984classification}) and RFs \citep{klusowski2024,gey2005}. Many of these works rely on assumptions common in supervised learning, such as {\em additivity} \citep{scornet2015}. For additive models, spatial adaptation typically focuses on addressing {\em sparsity} \citep{klusowski2024,syrgkanis2020,chi2022}. Several works also rely on a local accuracy gain condition \citep{syrgkanis2020,chi2022,athey2019}, a condition often met by additive models but challenging to extend to unsupervised settings. The work by \cite{klusowski2024}, which leverages the greedy nature of tree construction algorithms, provides general empirical risk bounds, though they may not always be tight.  Additionally, non-adaptive tree-based methods, which use random splits instead of data-guided ones, have gained attention \citep{arlot2014analysis,mourtada2020,cattaneo2023inference}. However, these methods typically sacrifice spatial adaptivity.

Under the Bayesian framework,  there has been rigorous analysis of posterior contraction for Bayesian CART \citep{chipman1997} and uncertainty quantification \citep{veronika2020,castillo2021uncertainty}. A recent work by \cite{Linero2018BayesianRT} has examined the spatial properties of BART, with a focus on adaptation to sparsity. Spatial properties of BART have also been characterized over the anisotropic H\"older or region-wise anisotropic H\"older spaces \citep{jeong2020art,rockova2021ideal}, particularly when the regression function is less smooth, being at most Lipschitz continuous.

The rest of the paper is organized as the following. We review some single tree-based Bayesian models on probability densities, whose posterior mean we will adopt as the base learner in our tree ensembes,  and different types of regularizations on individual trees in Section \ref{sec:notation} to attain reliable single-tree-based estimates. There we also introduce assumptions on the true density function in Section \ref{sec:fassump}. Then we provide our first Bayesian model of forests for density learning in Section \ref{sec:BMA}. We characterize its convergence rate over the anisotropic H\"older space, and present a concentration inequality for the tree-based likelihood as a tool for decomposing the bias and the variance. In Section \ref{sec:gibbs}, we introduce our risk-based generalized Bayesian strategy for building a DF, and analyze the posterior contraction properties of the resulting ensemble. We provide a specific example of the generalized Bayesian forest, which can achieve near-optimal rate over the anisotropic H\"older space in an adaptive way. We demonstrate the performance of the two methods on simulated datasets in Section \ref{sec:simulations}.

\section{Elements for Building a Forest for Density Learning}
\label{sec:notation}
We study the density estimation problem in the $d$-dimensional Euclidean space. The dimension $d$ is viewed as fixed in this paper. Let $\Y_n = \{ Y_1, Y_2, \cdots, Y_n \}$ be a simple random sample of size $n$ from a distribution with a density $f_0$ with respect to a
$\sigma$-finite measure $\mu$ on a measurable space $(\Omega,
\mathcal{B})$. We focus on the case when $\Omega$ is a
bounded rectangle in $\mathbb{R}^d$, $\mathcal{B}$ is the Borel algebra, and $\mu$ is the Lebesgue measure. After translation and scaling, we may assume that the
sample space is the unit cube in $\mathbb{R}^d$.
Let $\mathcal{F}$ be the collection of all the density
functions on $\Omega$ with respect to $\mu$. 
We define the metric on $\mathcal{F}$ to be the Hellinger distance, denoted as $\rho$, and also consider the Kullback-Leibler (KL) divergence, which is denoted as $\kl (f_0, f) = \mathbb{E}_{f_0} \left( \log (f_0/f)(Y) \right)$.
It is well known that the KL divergence is a stronger
pseudo-metric compared to the Hellinger distance, in the sense that for any $f, g \in \mathcal{F}$, $
\rho^2(f, g) \leq \kl (f,g)$. 
We may also use the $\| \cdot \|_\infty$ metric, where $\| f-
g\|_\infty$ is defined to be $\sup_{y \in \Omega} |f(y) - g(y)|$.

We consider a class of forest-type density estimators as a weighted aggregation over base learners, each of which is a single-tree-based density learner: 
\begin{eqnarray}
  \label{eq:forestdef}
f_\en = \sum_{m} w_m f_m, \quad \mbox{where }w_m \geq 0, \sum_{m} w_m =1.
\end{eqnarray}
In this section, we will introduce the necessary elements for building such a forest-type estimator under the Bayesian
framework.

\subsection{An Overview of the Tree Aggregation Schemes}
\label{subsec:overview}
Throughout the paper, we index the complexity of a tree by its depth $\D$. Each tree within a DF can grow under a small perturbation, which will be formally characterized later in Section~\ref{subsec:perturb}. \rev{The perturbation we consider can be interpreted as making splits at dyadic points jittered towards a pre-specified direction, or equivalently, growing dyadic trees on slightly shifted data points.} The potential form of the perturbation is usually depth dependent; thus, we tentatively denote it as $h_\D$. Here we provide an overview of tree aggregation schemes of the form \eqref*{eq:forestdef} that we will consider in this paper, with step-by-step derivations to be presented in Sections~\ref{sec:BMA} and~\ref{sec:gibbs}. 

Under the fully Bayesian framework, we may view trees of different depths and under different forms of perturbations as candidate models for data, and impose a prior distribution $\Pi$ on these models. Then, the Bayesian estimator as the posterior mean can be written as a weighted average of single-tree-based estimators,  
\begin{eqnarray*}
  E( f | \Y_n ) = \sum_{\D} \sum_{h_\D} E( f | \Y_n, \D, h_\D) \cdot \Pi( h_\D |
  \Y_n, \D) \cdot \Pi( \D | \Y_n),
\end{eqnarray*}
where each Bayesian estimator $ E( f | \Y_n, \D, h_\D)$ can be viewed as a base learner (i.e., $f_m$ in \eqref*{eq:forestdef}) with $(\D,h_{\D})$ being its index ( ``$m$''), and the associated weight (``$w_m$'') 
is $\Pi( h_\D |\Y_n, \D)  \Pi( \D | \Y_n)=\Pi(\D, h_\D |\Y_n)$, which is the marginal probability of a tree model indexed by $(\D,h_\D)$. We may define a Bayesian density forest simply in this form. 

However, in practice, the posterior probability of a perturbed tree model $\Pi(h_\D |\Y_n, \D)$ usually suffers from high variance, in contrast to the fast concentration of the marginal posterior distribution on the depth of the tree $\Pi(\D |\Y_n)$. To avoid condensed weights on a very small subset of perturbed trees, we find that an effective strategy is to simply replace the data-dependent average using the posterior probability over perturbation $\Pi( h_\D |\Y_n, \D)$ by a pre-determined data-independent weight, e.g., a uniform weight as in random forest and bagging. This change results in the following hybrid between Bayesian and non-Bayesian model averaging which for simplicity we still refer to as a  {\em Bayesian density forest} (BDF)  
\begin{eqnarray*}
 \sum_{\D} \sum_{h_\D} E( f | \Y_n, \D, h_\D) \cdot u_\D \cdot \Pi( \D | \Y_n), \quad{ u_\D >0},
\end{eqnarray*}
where $u_\D$ is the depth-dependent weight that replaces $\Pi( h_\D |\Y_n, \D)$. This hybrid aggregation strategy allows us to learn an appropriate tree depth from data through the average with respect to $\Pi( \D | \Y_n)$, while still attaining diversity of trees via the unweighted average.

It turns out that this (hybrid) BDF adapts to the data, but can do so only through the likelihood associated with the single trees in the base learner, as the Bayesian part of the construction is induced only by the likelihood associated with trees. The tree-based likelihood, however, is inadequate in measuring the goodness-of-fit for a forest-type ensemble, resulting in suboptimal convergence rates of the BDF estimate. 

To address this limitation, we further introduce a novel risk function associated with each forest-type ensemble and make risk-based generalized Bayesian aggregation through the corresponding Gibbs posterior, resulting in what we call a {\em Gibbs density forest} (GDF). More specifically, a forest topology $\sE= \forest$ consisting of $B = |\sE|$ trees can be viewed as a candidate sampling model for the data, where the distribution from this model is assumed to be a mixture $f_\sE = (f_{\T_1}+ \cdots +f_{\T_B})/B$ of $\T_b$-supported density functions. We will rigorously justify that the risk function associated with every forest topology can adequately measure its goodness-of-fit. When endowed with a prior $q$ on the forest topology, the risk-based generalized Bayesian framework allows us to learn 
forest topologies from data through the Gibbs posterior, denoted as $\lambda^\g$. We consider two specific density estimators that can be readily constructed under this framework: (i) a {\em maximum a posteriori}, or ``MAP GDF'' estimator, defined on the forest topology that maximizes the Gibbs posterior,
$\widehat{\sE}_{\mbox{\scriptsize{MAP}}}$, and (ii) a ``predictive GDF''  obtained by averaging over all of the forest topologies in the support of the prior $q$ with respect to the Gibbs posterior. 


\subsection{A Base Learner for Density Estimation: the Optional P\'olya Tree}
\label{subsec:opt}
An important element for building a forest is the tree-supported base learner. In this paper, the base learner is defined as a Bayesian estimator under a tree-like prior for probability measures. 
There are various choices available for this type of priors, and our
ensuing analysis applies more broadly to these different
choices. Among them, a particularly attractive candidate is the optional P\'olya
tree (OPT, \citep{wong2010}), an unsupervised counterpart of the
CART model that has been extensively studied in the context of
single-tree density estimation. \rev{The posterior behavior under the
  OPT prior is well understood theoretically. It achieves near-minimax
  sup-norm contraction over the H\"older space \citep{castillo2022},
  allow for near optimal
  adaptation to spatial heterogeneity in Hellinger distance
  \citep{liu2024}, and yields Bayesian
  credible sets with adaptive, near-optimal width and valid
  frequentist coverage for H\"older-continuous densities \citep{castillo2022}. 
}
In the following of this subsection, we provide a brief overview of this prior under the context of planned analysis. 


\tech{We start our description from \emph{dyadic} trees.} To grow the
tree, the root $\Omega = [0, 1]^d$ can be split into two equally sized
level-1 nodes along the {\em midpoint} of a selected
dimension. We continue to make such type of splits recursively.
In general, for a level-$k$ node $A$, 
there are $d$ ways to partition it. For some $l \in [d]$, 
  $A$ is split into two level-$(k+1)$ nodes, denoted as $A_{l, 0}$ and $A_{l, 1}$, along the midpoint of the $l$-th dimension.
The collection of all possible level-$k$ nodes is denoted
as $\mathcal{R}^\ksup$, and the recursive partition of $\Omega$ obtained by making recursive dyadic splits is called a dyadic tree. Leaf nodes refer to the terminal level nodes. A tree is of depth $\D$ if all leaf nodes are of level at most $\D$ and at least one leaf node is of level $\D$. 

The specification of the OPT prior on probability measures can be finished in two steps: first, we impose a prior on random recursive partitions of $\Omega$, or equivalently we sample a random dyadic tree from the prior; second, given the tree, we generate a random probability measure $F$ which is uniformly distributed within each leaf node of this tree by recursively splitting the probability mass according to an independent beta random variable. 

More specifically, to generate a random tree structure, we start from the root $\Omega$ and run a recursive procedure with respect to the level. Suppose after $k$ steps of the recursion, we have obtained a random recursive partition of $\Omega$,
\begin{eqnarray*}
  \A = \{\Omega_j \}_{j=1}^t,  \quad \mbox{with } 
\Omega_j \in \cup_{k'\leq k } \mathcal{R}^{(k')}. 
\end{eqnarray*}
For any node $\Omega_j$ of level $k' < k$, the partition has already been stopped, and will remain intact in later steps. In next step, to sample a potential refinement of $\A$, for each level-$k$ node $A$, 
we introduce a \emph{stopping variable} as a Bernoulli
random variable, 
{\normalsize
\begin{eqnarray*}
  S(A) \sim \mbox{Bernoulli} \left (\tau (A) \right), \quad
  \mbox{stopping variable for a level-$k$ region $A$},
\end{eqnarray*}}
where the parameter $\tau (A)$ 
is a hyper-parameter that can be node-specific. If $S=0$, stop further partitioning of $A$. 
If $S=1$, sample the {\em split variable} $\Lambda(A)$ according to a discrete distribution $\zeta(A)$ on $[d]$ and split $A$ 
along the midpoint of the selected dimension. 
The notation $\zeta(A)$ emphasizes that the
prior on the split variable can also be node-specific. The recursion will continue until all regions have been stopped. 

To generate the random probability measure $F$, assume after $k$ steps, in addition to the random partition $\A$, we have obtained a random probability measure $F^\ksup$ which is uniformly distributed within each $\Omega_j
\in \A$. 
In the $(k+1)$-th step, the random probability measure $F^\knew = F^\ksup$ on stopped nodes. For a level-$k$ node $A$, if we choose to further partition it, then we split the probability mass $F^\ksup(A)$ according to an independent beta random variable with potentially node-specific hyper-parameters $a(A)$ and $b(A)$.  
The random measure $F$ can be viewed as the terminal level $F^\ksup$ when the recursion ends.

The collection of hyper-parameters for the OPT prior is {\small $\left \{ 
\tau(A), \zeta(A), a(A), b(A)\right\}_{A \in \cup_{k \geq 0}
\mathcal{R}^\ksup}$.} \tech{In the paper by \cite{wong2010}, the
authors have shown that the OPT is a conjugate prior for tree models.} 
We refer to their paper for more details of Bayesian inference with the OPT prior. For the rest of the paper, we use the notation $\mathcal{A}
\mbox{ or } \T= \{\Omega_j \}_{j=1}^t$ to denote a 
tree topology that can be generated under the OPT prior. 
Each rectangular region $\Omega_j$ is a leaf node of the tree.

\subsection{Regularization on Tree Density Learners}
\label{subsec:reg}
When fitting forests for supervised learning problems, \rev{regularization of fully grown trees helps control model complexity and reduce the variance of the forest estimator. We extend this idea to the unsupervised setting by introducing flexible control over model complexity.} 
In the following, we provide two directions for incorporating tree-level regularization: adjusting the prior on stopping variables and imposing shape constraints, both of which can be achieved by adjusting hyper-parameters in tree priors that control the complexity of the trees on which such priors are supported.


\subsubsection{Adjusting the Prior on Stopping Variables}
\label{subsubsec:depth}
The first approach is to make the prior probability on deep trees small, and this can be accomplished through specifying fast decaying prior to stopping variables. 
The hyper-parameter $\tau(A)$ denotes the prior probability of making a further split of $A$. Throughout the paper, 
we assume that it is completely determined by the level of the node. Following prior work (e.g., \cite{castillo2022}), for a level-$k$ node $A$, we consider an exponentially decaying prior probability for further partitioning, 
\begin{eqnarray}
  \label{eq:stoppingrule}
  \tau(A) = c_\tau 2^{-\nu k}, \quad \quad 0 < c_\tau \leq 1/2, \quad \nu \geq 1. 
\end{eqnarray}
Furthermore, we also consider trees of \rev{a pre-specified} depth $\D$. 
That is, $\tau(A)=0$ whenever $k=\D$.


\subsubsection{Regularization on the Shape of Tree Topologies}
\label{subsubsec:shape}
When we grow trees within a forest by making splits along observed data points or bootstrapped ones, although each tree is shallow or even of depth one, the induced partition by aggregating multiple trees can still be refined due to the randomness of data and is sufficient to capture local variation. However, with dyadic splits, the diversity of partitions induced by shallow tree ensembles is limited. We need to allow tree depth to grow at a slow rate to overcome such limits. On the other hand, with increasing depth, the complexity of the tree and forest spaces increases extremely fast. To reduce complexity, we consider imposing shape constraints on individual trees. Under the OPT prior, this is done by adjusting the node-specific prior on stopping variables $S(A)$ and split variables $\Lambda(A)$. Marginally, we assume that $\Lambda(A)$ follows a uniform distribution on $[d]$. Without constraints, both $S(A)$ and  $\Lambda(A)$ are assumed to be \emph{independent a priori}, implying that any potential further split of $A$ is independent of the other nodes in the current random recursive partition. 

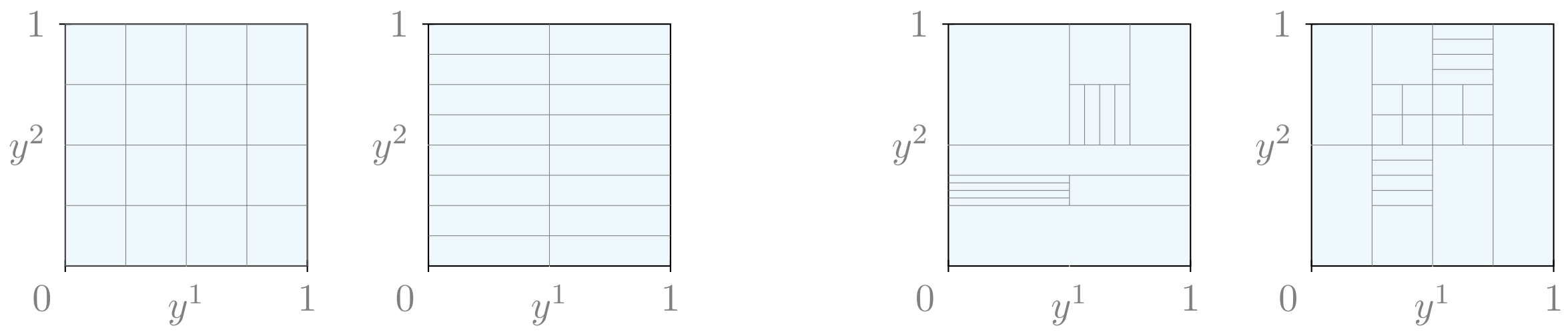
\begin{figure}
  \centering
\begin{multicols}{2}
 \begin{center}
\vspace{5mm}
\begin{tikzpicture}[xscale=0.6, yscale=0.6]
  \filldraw[color=black!, thin, fill= cyan!6] (1,1) rectangle (5,5);
\draw ( 1,1+0.1) -- (1 ,1-0.1) node[anchor=north east,gray] {0};
\draw ( 5,1+0.1) -- (5 ,1-0.1) node[anchor=north,gray] {1};
\draw[cyan!0,very thin] ( 3,1+0.1) -- (3 ,1-0.1) node[anchor=north,gray] {$y^1$};
\draw (1+0.1,5) -- (1-0.1,5) node[anchor=east,gray] {1};
\draw[cyan!0,very thin] (1+0.1,3) -- (1-0.1,3) node[anchor=east,gray] {$y^2$};

\draw[step=1cm, gray, very thin] (1,1) grid (5,5);



 \filldraw[color=black, thin, fill= cyan!6] (7,1) rectangle (11,5);
\draw ( 7,1+0.1) -- (7 ,1-0.1) node[anchor=north east,gray] {0};
\draw ( 11,1+0.1) -- (11 ,1-0.1) node[anchor=north,gray] {1};
\draw[cyan!0,very thin]( 9,1+0.1) -- (9 ,1-0.1) node[anchor=north,gray] {$y^1$};
\draw (7+0.1,5) -- (7-0.1,5) node[anchor=east,gray] {1};
\draw[cyan!0,very thin] (7+0.1,3) -- (7-0.1,3) node[anchor=east,gray] {$y^2$};




 \draw[gray, very thin] (9, 1) -- (9,5);

\draw[gray, very thin] (7, 3) -- ( 11, 3);
\draw[gray, very thin] (7, 1.5) -- ( 11, 1.5);
\draw[gray, very thin] (7, 2) -- ( 11, 2);
\draw[gray, very thin] (7, 2.5) -- ( 11, 2.5);
\draw[gray, very thin] (7, 3.5) -- ( 11, 3.5);
\draw[gray, very thin] (7, 4) -- ( 11, 4);
\draw[gray, very thin] (7, 4.5) -- ( 11, 4.5);
\end{tikzpicture}
\end{center}
 \begin{center}
\vspace{5mm}
\begin{tikzpicture}[xscale=0.6, yscale=0.6]
  \filldraw[color=black, thin,  fill= cyan!6] (1,1) rectangle (5,5);

\draw ( 1,1+0.1) -- (1 ,1-0.1) node[anchor=north east,gray] {0};
\draw ( 5,1+0.1) -- (5 ,1-0.1) node[anchor=north,gray] {1};
\draw[cyan!0,very thin] ( 3,1+0.1) -- (3 ,1-0.1) node[anchor=north,gray] {$y^1$};
\draw (1+0.1,5) -- (1-0.1,5) node[anchor=east,gray] {1};
\draw[cyan!0,very thin] (1+0.1,3) -- (1-0.1,3) node[anchor=east,gray] {$y^2$};


\draw[gray, very thin] (1, 3) -- ( 5, 3);
\draw[gray, very thin] (3, 3) -- ( 3, 5);
\draw[gray, very thin] (1, 2) -- ( 5, 2);


\draw[gray, very thin] (1, 2.5) -- ( 5, 2.5);
\draw[gray,very  thin] (1,2.25) -- ( 3, 2.25);
\draw[gray,very  thin] (1,2.125) -- ( 3, 2.125);
\draw[gray, very thin] (1,2.375) -- ( 3, 2.375);
\draw[gray, very thin] (3,2) -- ( 3, 2.5);

\draw[gray, very thin] (4, 3) -- ( 4, 5);
\draw[gray, very thin] (3.5, 3) -- ( 3.5, 4);
\draw[gray, very thin] (3.25, 3) -- ( 3.25, 4);
\draw[gray, very thin] (3.75, 3) -- ( 3.75, 4);
\draw[gray, very thin] (3, 4) -- (4,4);

 \filldraw[color=black, thin, fill= cyan!6] (7,1) rectangle (11,5);
\draw ( 7,1+0.1) -- (7 ,1-0.1) node[anchor=north east,gray] {0};
\draw ( 11,1+0.1) -- (11 ,1-0.1) node[anchor=north,gray] {1};
\draw[cyan!0,very thin]( 9,1+0.1) -- (9 ,1-0.1) node[anchor=north,gray] {$y^1$};
\draw (7+0.1,5) -- (7-0.1,5) node[anchor=east,gray] {1};
\draw[cyan!0,very thin] (7+0.1,3) -- (7-0.1,3) node[anchor=east,gray] {$y^2$};


\draw[gray, very thin] (9, 1) -- ( 9, 5);
\draw[gray, very thin] (7, 3) -- ( 11, 3);
\draw[gray, very thin] (8, 1) -- ( 8, 5);
\draw[gray, very thin] (10, 1) -- ( 10, 5);
\draw[gray, very thin] (8, 4) -- ( 10, 4);
\draw[gray, very thin] (8, 2) -- ( 9, 2);
 \draw[gray, very thin] (8, 3.5) -- ( 10, 3.5);
\draw[gray, very thin] (8.5, 3) -- ( 8.5, 4);
\draw[gray, very thin] (9.5, 3) -- ( 9.5, 4);

\draw[gray, very thin] (9, 4.5) -- ( 10, 4.5);
\draw[gray, very thin] (9, 4.25) -- ( 10, 4.25);
\draw[gray, very thin] (9, 4.75) -- ( 10, 4.75);

\draw[gray, very thin] (8, 2.5) -- ( 9, 2.5);
\draw[gray, very thin] (8, 2.25) -- ( 9, 2.25);
\draw[gray, very thin] (8, 2.75) -- ( 9, 2.75);

\end{tikzpicture}
\end{center}
\end{multicols}
\vspace{-5mm}
  \caption{{ \small Illustration of shape constraints in the two-dimensional case. The left panel shows two tree topologies obtained under the strong shape constraint. Both trees are of depth 4. The right panel illustrates the weak shape constraint. Both trees first grow without
    any constraint up to level $k_0 = 4$. Then the growth is coupled
    region-wise. }}
  \label{fig:constraint}
\vspace{-3mm}
\end{figure}

Under our first type of shape constraints, 
we assume that starting from the root $\Omega$, at each step of the recursion, the stopping variable and the split variable of all level-$k$ nodes are \emph{coupled} together: one either chooses to stop at level $k$, or selects a dimension and splits all level-$k$ nodes along the same dimension. 
We call the corresponding tree-like prior the OPT prior under the \emph{strong shape constraint}. In the left panel of Figure \ref{fig:constraint}, we provide two examples of tree topologies obtained under this constraint. 

The strong shape constraint is very restrictive and not commonly adopted in practice when applying tree-based density learners because it heavily relies on the homogeneous smoothness property to produce reasonable estimates. The inclusion of this strong shape constraint in our theoretical analysis is to show that single-tree methods can only achieve suboptimal rates of convergence for continuously differentiable densities even under such strong restrictions to reduce the variance, whereas a forest-type ensemble can achieve almost optimal rates when such constraints are relaxed, as the tree aggregation effectively removes the excessive variance associated with single-tree estimates.

For this purpose, we also consider the following \emph{weak shape constraint}, which is more relaxed and practically relevant. Under this constraint, a tree is first grown for  $\D_0$ levels 
without any constraints, obtaining a partition of depth $\D_0$, denoted as $\A_0 = \{ \Omega_j^0 \}_{j=1}^{t_0}$. Starting from the next level, within every node $\Omega_j^0$, further splits within $\Omega_{j}^0$ of the same level are coupled together as in the strong shape constraint. In the right panel of Figure \ref{fig:constraint}, we provide two examples to illustrate tree topologies generated under the weak shape constraint. 


In the context of the OPT prior, both the strong and weak shape constraints only limit the types of tree topologies that can be sampled under the prior. They do not affect the generation of a random probability measure given the tree, which still follows a distribution induced by the product of independent beta random variables. 

Under the (unconstrained) OPT prior, it has been shown that the posterior distribution \tech{conditional} on \emph{any} node $A$ is still an OPT with hyper-parameters fully determined by the data points within $A$ as well as hyper-parameters of the prior \citep{wong2010}. With shape constraints, if under the prior further partitioning of $A$ is coupled with any other node $A'$, then the posterior distribution conditional on $A$ is also informed by data in $A'$. It is these shape constraints 
that help avoid overfitting in deep levels of the tree when the underlying density is smooth, i.e., showing relatively stable local variation. This is indeed appropriate for the scope of this paper, as we mainly focus on the case when $f_0$ is H\"older continuous (we will detail the conditions on $f_0$ in Section~\ref{sec:fassump}).


\subsection{Decorrelation Through Probabilistic Cycle Spinning}
\label{subsec:perturb}
In supervised problems, trees within a forest are often fit on perturbed data, like bootstrap samples, in order to decorrelate the individual tree estimates, thereby reducing the variance of the aggregated forest estimate. In the unsupervised context, 
a major source of variance in tree-based density estimation arises from data points near partition boundaries: small changes of their locations can substantially alter the counts within leaf nodes and, in turn, to large fluctuations in, for example, the likelihood, which then propagate into the posterior when a Bayesian tree model is employed. We shall refer to this phenomenon as the ``boundary effect''. This issue can be more severe for dyadic trees, and has been noted both in previous works on P\'olya trees \citep{nieto2012,hanson2006,paddock2003}, and even earlier in the literature of wavelets, which are in essence a tree-based sampling model for functions. To alleviate the boundary effect, we consider a decorrelation strategy inspired by the \writing{``cycle-spinning'' technique in wavelet analysis} \citep{coifman1995}---through slightly shifting the recursive partition of $\Omega$ along pre-specified directions. Within the Bayesian paradigm, one can incorporate cycle-spinning in a probabilistic fashion, treating the random shifting $h$ as a parameter endowed with a hyper-prior.

 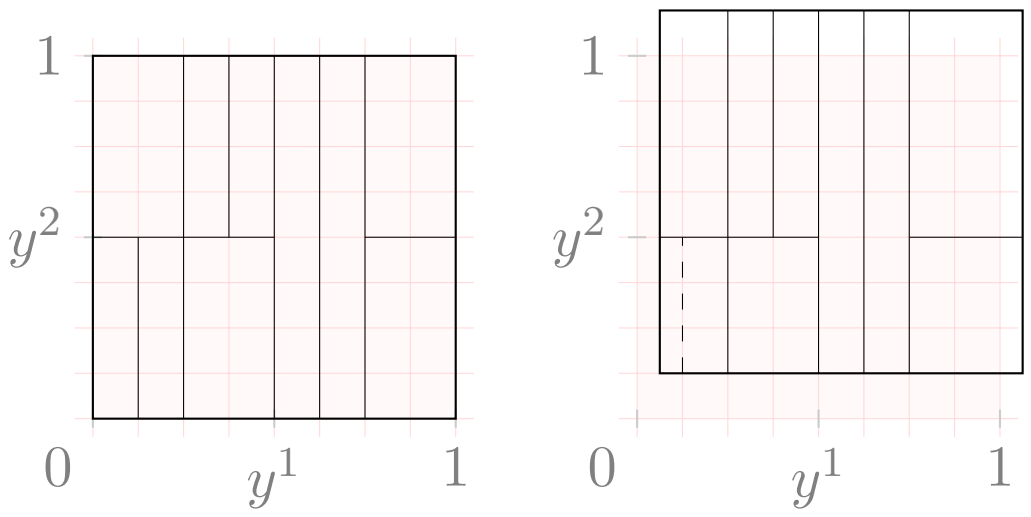
\begin{figure}[t]
  \vspace{-10mm}
  \centering 
\begin{center}
\vspace{5mm}
\begin{tikzpicture}[xscale=0.6, yscale=0.6]
  \filldraw[color=pink!0, fill= pink!10] (1,1) rectangle (5,5);
  \draw[step=.5cm,pink!60,very thin] (0.8,0.8) grid (5.2,5.2);
\draw[gray!40] ( 1,1+0.1) -- (1 ,1-0.1) node[anchor=north east,gray] {0};
\draw[gray!40] ( 5,1+0.1) -- (5 ,1-0.1) node[anchor=north,gray] {1};
\draw[gray!40] ( 3,1+0.1) -- (3 ,1-0.1) node[anchor=north,gray] {$y^1$};
\draw[gray!40] (1+0.1,5) -- (1-0.1,5) node[anchor=east,gray] {1};
\draw[gray!40] (1+0.1,3) -- (1-0.1,3) node[anchor=east,gray]{$y^2$};

  \draw[black] (1,1) rectangle (5,5);

\draw[black, very thin] (1, 3) -- ( 3, 3);
\draw[black, very thin] (4, 3) -- ( 5, 3);

\draw[black, very thin] (3, 1) -- (3,5);

\draw[black, very thin] (2, 1) -- (2,5);
\draw[black, very thin] (4, 1) -- (4,5);

\draw[black, very thin] (3.5, 1) -- (3.5,5);

\draw[black, very thin] (1.5, 1) -- (1.5,3);
\draw[black, very thin] (2.5, 3) -- (2.5,5);

\filldraw[color=pink!0, fill= pink!10] (7,1) rectangle (11,5);

\draw[step=.5cm,pink!60,very thin] (6.8,0.8) grid (11.2,5.2);
\draw[gray!40] ( 7,1+0.1) -- (7 ,1-0.1) node[anchor=north east,gray] {0};
\draw[gray!40] ( 11,1+0.1) -- (11 ,1-0.1) node[anchor=north,gray] {1};
\draw[gray!40] ( 9,1+0.1) -- (9 ,1-0.1) node[anchor=north,gray] {$y^1$};
\draw[gray!40] (7+0.1,5) -- (7-0.1,5) node[anchor=east,gray] {1};
\draw[gray!40] (7+0.1,3) -- (7-0.1,3) node[anchor=east,gray] {$y^2$};

\draw[black] (7+0.25,1+0.5) rectangle (11+0.25,5+0.5);

\draw[black, very thin] (7+0.25, 3) -- ( 9, 3);
\draw[black, very thin] (10, 3) -- ( 11+0.25, 3);

\draw[black, very thin] (9, 1+0.5) -- (9,5+0.5);

\draw[black, very thin] (8, 1+0.5) -- (8,5+0.5);
\draw[black, very thin] (10, 1+0.5) -- (10,5+0.5);

\draw[black, very thin] (9.5, 1+0.5) -- (9.5,5+0.5);

\draw[black,dashed, very thin] (7.5, 1+0.5) -- (7.5,3);
\draw[black, very thin] (8.5, 3) -- (8.5,5+0.5);

\end{tikzpicture}
\end{center}
\vspace{-6mm}
\caption{{\small An example of shifted trees in the two-dimensional case. Left: a depth-4 tree on $[0,1]^2$. Right: a depth-4 tree on the shifted sample space $\Omega^h$ with $h=(\frac{1}{16}, \frac{1}{8})$. The tree is obtained by making splits along the same point at every step as those for the tree in the left panel. The tree pruning scheme prevents a split shown in dashed line.}
}
\label{fig:shifted_partition}
\end{figure}


The perturbation can be incorporated equivalently by growing trees on shifted sample space. The shifted sample space $\Omega^h$ is defined as $\Omega^h = \{ y: y-h \in \Omega \}$. Under the shift $h$, we grow a tree on $\Omega^h$ with only splits along {\em dyadic} points allowed, and then shift the sample space back to obtain a partition of the original sample space $\Omega$. To avoid obtaining a leaf node of very small volume, we further introduce a pruning scheme for boundary regions. A shifted tree topology is usually denoted as $\A^h$ or $\T^h$. More details of tree growth under a shift can be found in Section~E.1 of the Supplementary Material \citep{liu2026b}.  In Figure~\ref{fig:shifted_partition}, we provide an example of the shifted tree in the two-dimensional case.

To carry out Bayesian inference, we can still introduce the OPT prior when considering shifted tree topologies. The major difference is to replace level-$k$ nodes by those obtained under a shift. The prior on stopping variables is specified in the same manner as that in \eqref*{eq:stoppingrule}. We set hyper-parameters for the beta distribution as $a(A)= b(A) = a_0$, where $0 < a_0 < 1$, for any node $A$, either with or without shifts. 
For the rest of the paper, we  denote the OPT prior under either the strong or the weak shape constraints as $\Pi(\cdot)$, the OPT of a pre-specified depth $\D$ as $\Pi (\cdot | \D)$, and the prior under a perturbation along the direction $h$ as $\Pi(\cdot | h)$. 


\section{Assumptions on $f_0$}
\label{sec:fassump}
We are mainly investigating the case when the true
density enjoys good smoothness properties. In specific, we assume that $f_0$ is anisotropic
H\"older continuous. Our definition of the function class follows that by \cite{barron1999}. Generally, for $k = (k_1, \ldots, k_d) \in \mathbb{N}^d$ and a multivariate function $f$, $\nabla^k f$ denotes the mixed partial derivatives of order $(k_1, \ldots, k_d)$ of $f$.
\begin{definition}
\label{def:holder}
For a function $f \in \mathcal{C}([0, 1]^d)$, $y \in [0,1]^d$, and
$1\leq l \leq d$, let $f_l (\cdot | y)$ denote the univariate function
$x \mapsto f( y_1, \ldots, y_{l-1}, x, y_{l+1}, \ldots, y_d)$.  The
anisotropic H\"older space with parameter $\bm \alpha \in
(0,\infty)^d$ and of norm $L$, denoted as $\mathcal{C}^{\bm \alpha}
(\Omega, L)$, is the set of functions $f$ which satisfies the following
two conditions, 
{\small
\begin{eqnarray}
\label{eq:holder1}
    \sup_{y \in [0, 1]^d} 
\sum_{k=0}^{\lfloor \alpha_l  \rfloor }
\left \| 
\nabla^k  f_l (\cdot |y) 
\right \|_\infty 
\leq L
\quad
\mbox{for all }
1 \leq l \leq d,
\end{eqnarray}
}
  \vspace{-2mm}
  {\small 
\begin{eqnarray}
  \label{eq:holder2}
  \hspace{3mm}
\sup_{y \in [0, 1]^d} 
\left \| 
\nabla^{\lfloor \alpha_l  \rfloor }
f_l ( x |y) 
-
\nabla^{\lfloor \alpha_l  \rfloor }
f_l (x' | y)
\right \|_\infty 
\leq L |x' - x|^{\alpha_l -\lfloor \alpha_l  \rfloor }
  \quad
  \text{for any } x, x' \in [0, 1], 
1 \leq l \leq d.
\end{eqnarray}
}
\end{definition}
\noindent When $\alpha_1 = \cdots =\alpha_d= \alpha$,
$\mathcal{C}^\alpha (\Omega, L)$ is the class of isotropic H\"older continuous functions. The anisotropic space allows the function to have different levels of smoothness along different directions, and is a commonly used function class to model spatial heterogeneity. In the analysis of
tree-based methods, this model for spatial heterogeneity has been previously considered by \cite{jeong2020art,liu2024}.

In certain cases, we also need to introduce a boundary
condition. 
Let $\partial \Omega = \mbox{closure}(\Omega) \cap \mbox{closure}(\mathbb{R}^d \setminus \Omega)$ denote the boundary of $\Omega$. For $\tau \geq 0$, let $\overline{\Omega}_\tau$ denote the set of points in $\Omega$ with distance no more than $\tau$ to the boundary
    $\overline{\Omega}_\tau =\{ y \in \Omega: d(y, \partial \Omega)
    \leq \tau \}$,
where $d(y, \partial \Omega)$ denotes the infimum Euclidean distance between $y$ and any point in $\partial \Omega$. We may consider the following boundary condition in later parts of the paper. 
\begin{conditionp}{(BOUND)}
 \label{cond:boundary}
The gradient of $f_0$ decays fast near the boundary. More explicitly,
when $f_0$ is \tech{continuously differentiable}, there exists a
$\tau_0 >0$, such that for any $\tau \leq \tau_0$ and any $1 \leq l
\leq d$, \rev{$\sup_{y: d(y, \partial \Omega) \leq \tau} \left |
  \frac{\partial f_0}{ \partial y^l}  \big |_y  \right |\leq L \tau $} , where $L$ is a positive constant.
\end{conditionp}

This boundary condition requires that derivatives of $f_0$ decay at a
polynomial rate at the boundary. In typical applications, there are
usually a limited amount of data in the boundary region, and
tree-based methods enforce strong regularization on the tail behavior
of the resulting density. Hence, the irregular behavior of $f_0$ at
the boundary often cannot be effectively recovered by tree-based
methods, especially forests. If $f_0$ only varies mildly, as implied
by decaying derivatives, we expect forest-based methods to continue
performing well, even with a limited number of data points for
estimating the density near the boundary.

\rev{Technically, Condition~\ref{cond:boundary} ensures that the
  approximation error of the forest decays to zero at a
  sufficiently fast rate in boundary regions. Instead of imposing
  a condition on the gradient, one may consider employing other
  methods to improve
  the approximation accuracy of forests at the boundary. For
  instance, in \cite{randrianarisoa2021smoothing}, the author proposed a stochastic
  weighting scheme. We focus on Condition~\ref{cond:boundary}
here as it is easier to work with in the multivariate case.}


\section{Bayesian Density Forests}
\label{sec:BMA}


We introduce our first model for density learning that involves the
aggregation of multiple single-tree learners in the form of
\eqref*{eq:forestdef}. We start with a fully Bayesian approach and view
trees of different depths and perturbations as
potential competing generative models for data. Averaging over the resulting single-tree density estimators leads to our
first ``ensemble'' learner. This learner is the posterior predictive density, involving a fully Bayesian model averaging (BMA, \citep{hoeting1999,bernardo1999bayesian}) on these competing
models. 
It turns out, however, that for reasons we will describe later, such a
BMA aggregation may not be the most effective approach to variance reduction
due to ``boundary
effects'' under single-tree generative models. As such, we next
introduce a modification to the BMA estimator, aiming to reduce
the variance by employing a hybrid aggregation strategy.


\subsection{A Fully Bayesian Aggregation}
\label{subsec:fullybma}

%
We consider the following hierarchical prior on tree models. First, we impose a prior on tree depth, in order to achieve a more exact control of the tree complexity. The prior takes the form,
\begin{eqnarray*}
  \Pi(\D) \propto \exp ( - \gamma \D 2^{\beta \D}), \quad \gamma \geq 1, \beta \geq 1,
\end{eqnarray*}
where both hyper-parameters $\beta$ and $\gamma$ can be adjusted to
indicate how strongly we are in favor of shallow trees. For individual
tree models, we set $\beta=1$, \rev{whereas for forests we use 
  $\beta >1$ to enforce a stronger penalty on tree depth}.

Given the depth, trees can grow after a global perturbation on the
underlying potential split values. Practically, such a perturbation
can be equivalently achieved by shifting the sample space along a
randomly chosen direction, as described in Section
\ref{subsec:perturb}. 
We assume that possible directions of the shift depend on the tree depth. This is based on the observation that
the shifts are introduced to reduce the variance due to the boundary
effect, as explained in Section \ref{subsec:perturb}, and the variance
level is closely related to the depth of the tree. More specifically,
at depth $\D$, we define the unit for shift as \tech{$\iota= 2^{-2
    \lfloor \D/d \rfloor}$}. \rev{The rationale behind this choice of unit will be discussed later.} Potential shifting directions are
\begin{eqnarray*}
    \mathcal{H}_\D = \{ (i_l \cdot \iota)_{1 \leq l \leq d}: i_l
    \in \mathbb{N}, 0 \leq i_l < \tech{2^{ \lfloor \D/d \rfloor}}
    \}.
\end{eqnarray*}  
Every direction $h$ also indexes a collection of tree models for data that grows under this shift. For simplicity, we may abuse the notation $h$ for both a shifting direction and the collection of tree models under the perturbation $h$ when there is no confusion. 
Given the depth $\D$, the prior on $\mathcal{H}_\D$ is \emph{uniform}. 

Given the depth $\D$ and the perturbation $h$, the base learner on the
space of probability measures is \tech{the OPT under a shape
  constraint}, denoted as $\Pi(\cdot | \D, h)$. The shape constraints
are mainly introduced for the theoretical analysis of the tree
ensemble. In practice, it is feasible to consider an OPT prior without
any constraints. \rev{Under the strong or weak shape constraint,
  although split variables 
$\Lambda(A)$'s for several intermediate nodes may be coupled, the
prior distribution on these variables remains uniform, meaning that each dimension has the
same prior probability of being selected for a further split.}

 Under this hierarchical prior, the posterior distribution is the random measure given by 
\begin{eqnarray}
  \label{eq:posteriordist}
  \Pi(f \in B |\Y_n)
 &=& \frac{\sum_{(\D,h) }
  \left( \int_B \prod_{i=1}^n f (Y_i) d
  \Pi(f | \D, h) \right)
  \Pi (\D, h)}
  {\sum_{(\D', h')} \left( \int_{\F} \prod_{i=1}^n f' (Y_i) d\Pi(f' | \D', h') \right) \Pi(\D', h')} \\
\nonumber  &=& \sum_{\D} \sum_{h \in \mathcal{H}_\D }\Pi( f \in B | \Y_n, \D, h)  \cdot \Pi( h | \Y_n, \D ) \cdot \Pi(\D | \Y_n),
\end{eqnarray}
where $\Pi(\D | \Y_n)$ is the posterior probability on a depth-$\D$ tree model 
under any suitable form of perturbation, $\Pi(h| \Y_n, \D)$ is given the depth $\D$ the posterior probability of a tree model under perturbation $h$, and $\Pi(\cdot | \Y_n, \D, h)$ is the posterior distribution on probability measures under the shape constrained shifted OPT prior. The full posterior distribution is a mixture of posterior distributions conditional on different tree models, and the tree model is indexed by both tree depth and form of the perturbation. 
The posterior mean as the BMA is,
\begin{eqnarray*}
  E( f | \Y_n ) = \sum_{\D} \sum_{h \in \mathcal{H}_\D} E( f | \Y_n, \D, h) \cdot \Pi( h |
  \Y_n, \D) \cdot \Pi( \D | \Y_n).
\end{eqnarray*}

\subsection{A Hybrid Aggregation Scheme}
\label{subsec:hybridbma}
 Due to the boundary effect described in Section
\ref{subsec:perturb}, the term $\Pi( h |\Y_n, \D)$ usually has a large
variance. As such, this fully Bayesian estimator typically incurs
large variance in practice. This phenomenon has been reported in many
previous works on single-tree-based Bayesian density estimators
\citep{nieto2012,hanson2006,paddock2003}. 
A simple strategy to avoid such overfitting is to replace the fully Bayesian estimate, which involves BMA with respect
to perturbations $h$ given $\D$ by an aggregation scheme with fixed weights $1 / |\mathcal{H}_\D|$ for all depth-$\D$ trees, and define the corresponding version of the ``posterior mean'' as 
\begin{eqnarray}
  \label{eq:pseudoposteriormean}
  E^s( f | \Y_n ) = \sum_{\D} \sum_{h \in
  \mathcal{H}_\D} E( f | \Y_n, \D, h) \cdot \frac{1}{|\mathcal{H}_\D|}
  \cdot \Pi( \D | \Y_n),
\end{eqnarray}
Every $E( f |\Y_n, \D,
h)$ is obtained under a tree base learner of maximum depth $\D$ 
with perturbation $h$. The DF is an aggregation of perturbed trees, weighted by deterministic (equal) weights on perturbations, and data-dependent weights on tree depths.  
The unweighted averaging over different perturbations $h$ makes the
above estimator a Bayesian/non-Bayesian hybrid one in that it weights
over different tree depths using their posterior probabilities but
simply employs an unweighted average over perturbations in a similar vein as in bagging and RFs. 
The hybrid aggregation is a simple and effective strategy to address
the boundary-related high-variance issue of individual tree density
learners, as will be confirmed in our later theoretical analysis. We
call the estimator in Equation \eqref*{eq:pseudoposteriormean} a
\writing{{\em Bayesian density forest} (BDF)}, despite its hybrid nature.

The implementation of the forest estimator
\eqref*{eq:pseudoposteriormean} is based on drawing multiple
tree-supported density samples from the posterior distribution
\eqref*{eq:posteriordist} under the hierarchical prior $\Pi$. We will introduce an
efficient algorithm for this purpose in Section~H.1 of the 
Supplementary Material \citep{liu2026b},  
and provide several numerical results later in \pref{Section~\ref{sec:simulations}} to illustrate the performance of \writing{BDF}.

To characterize the frequentist property of the \writing{BDF estimator}, we derive the convergence rate $\epsilon_n$ defined as 
\begin{eqnarray}
  \label{eq:postmeanconv}
  \rho \left(E^s( f | \Y_n ), f_0 \right) \leq \xi \epsilon_n   \quad \mbox{with
  probability tending to 1 under } \P_{f_0}^n, 
\end{eqnarray}
where $\xi$ increases at most as polynomials of $\log n$. 
Our characterization is different from
the classical frequentist characterization of a Bayesian approach,
which typically involves deriving the posterior concentration rate
\citep{ghosal00,shen2001}. We adopt this alternative route given the hybrid nature of our decision rule and our focus on estimation. 

Given a random sample of size $n$, 
we set 
\begin{eqnarray}
  \label{eq:maxdepth}
  \D_{\max} = 
\left \lfloor
 \log_2 (n / \log n )
\right \rfloor.
\end{eqnarray}
and do not allow the tree depth to be larger than
$\D_{\max}$. \rev{This sample-size-based maximum tree depth guarantees
  that, in high- to moderately high-density regions, the leaf nodes
  contain sufficiently many data points on average.} 
Due to the sample-size dependence of the tree prior, we denote the depth
truncated prior as $\Pi_n (\cdot)$. 

Under the strong shape
constraint, we expect that variation of the true distribution can guide ``selection'' of tree topologies and
the method can pick up the ``best''
candidates with high probability . To solidify such an intuition, we
need to introduce the following condition on $f_0$, which can be
understood as a bias gap between the best and the second best tree
approximations. \rev{For any tree topology $\A$, we use $f_\A$ to denote the
  best tree approximation to $f_0$ supported by $\A$ under the KL divergence.} 
\begin{conditionp}{(UNIQUE)}
 \label{cond:uniqueness}
Consider trees grown under the strong shape constraint without
perturbation. At any depth $ \D \leq \D_{\max}$, among all depth-$\D$
tree topologies, there is a unique depth-$\D$ tree $\A_\D^\star$
satisfying that, for any other depth-$\D$ tree topology $\A$, $\kl(f_0,
f_\A ) \geq \phi \kl(f_0, f_{\A_\D^\star})$, where $\phi>1$ is a constant.
\end{conditionp}

\noindent Let $\theta_l = \int  \frac{(\partial f_0/ \partial y^l)^2}
{f_0} d \mu= \E_{X \sim f_0} (\frac{\partial \log f_0}{\partial y^l} (X) )^2$. This quantity reflects the variation of $f_0$ along
the $l$-th dimension. When studying approximation errors to $f_0 \in
\mathcal{C}^{\bm \alpha}(\Omega, L)$ in Section~F of the
Supplementary Material \citep{liu2026b}, 
we will see that $\{\theta_l\}_{1 \leq l \leq d}$ plays an important role in characterizing the approximation error when $f_0$ is continuously differentiable and is closely related to the validity of Condition \ref{cond:uniqueness}.
In the following, our first theoretical result is on the convergence rate of the \writing{BDF estimator}
\eqref*{eq:pseudoposteriormean}
when the true density function $f_0$ is
smooth. 

\begin{theorem}
\label{thm:BMArateanisotropic}  
Assume $f_0$ lies in the anisotropic H\"older space $\mathcal{C}^{\bm \alpha}(\Omega, L)$ where $\bm \alpha =(\alpha_1, \ldots, \alpha_d)$ and $\alpha_l \in (1, 2]$, and both Condition \ref{cond:boundary} and Condition \ref{cond:uniqueness} hold. In addition, we assume that $f_0 \geq \kappa_L $ for some $\kappa_L >0$. Let $\alpha_{\min}
= \min_{1 \leq l \leq d} \alpha_l$. 
Then under the depth truncated prior $\Pi_n$, the convergence rate for
the \writing{BDF estimator} \eqref*{eq:pseudoposteriormean} is, 
\begin{eqnarray*}
  \epsilon_n =
  \begin{cases}
    (dL/\sqrt{\kappa_L})
^{\frac{\beta+2/d -\alpha_{\min}/d}{\beta+2/d}} 
n^{-\frac{\beta + 2/d -1}{2\beta+ 4 /d}}, 
& \mbox{when } 1\leq \beta \leq 1+ 2(\alpha_{\min} -1)/d;\\
(dL/\sqrt{\kappa_L})^{\frac{1}{\beta+2/d}}
n^{-\frac{ \alpha_{\min} /d}{\beta +2 /d}},
& \mbox{when } \beta \geq  1+ 2(\alpha_{\min}-1) /d.
  \end{cases}
\end{eqnarray*}
\end{theorem}


A special case of Theorem \ref{thm:BMArateanisotropic} is
$\alpha_1= \cdots \alpha_d= \alpha$, corresponding to an isotropic
H\"older space. When $\beta =1+ 2(\alpha-1)/d$, the rate is minimax up
to a logarithmic term over $\mathcal{C}^\alpha(\Omega, L)$, $\alpha
\in (1, 2]$. Except for this case, for the anisotropic H\"older space,
no matter how we adjust $\beta$, the rate for \writing{BDF}
\eqref*{eq:pseudoposteriormean} is \emph{slower} than the minimax
rate.


\rev{To justify the improvement of BDF over single-trees, we
  investigate a single-tree prior as the special case of the forest
  prior introduced above when the number of trees within a forest is
  equal to one. In the following, for spatially
  heterogeneous $f_0$ as considered in Theorem~\ref{thm:BMArateanisotropic}, we establish a lower bound on the
posterior contraction rate of this single-tree method. Our
characterization of the lower bound follows the work by \cite{castillo2021uncertainty}.}

\vspace{-4mm}
\rev{\begin{theorem}
    \label{thm:treelowerbound}
 We consider the single-tree prior as the special case of the forest
 prior introduced in Section~\ref{subsec:fullybma} by setting the
 number of trees within a forest to be one and assuming that the tree
 grows without any perturbation under the strong shape constraint. In particular, we assume that the prior on tree depth takes the form $\Pi_n(\D) \propto \exp(- \gamma \D 2^{\beta_{\text{tree}} \D})$ when $\D \leq \D_{\max}$ and $\Pi_n(\D) = 0$ when $\D > \D_{\max}$. The hyper-parameters $\gamma$ and $\beta_{\text{tree}}$ satisfy $\gamma \geq 1$ and $\beta_{\text{tree}} \geq 1$, respectively. Let $\epsilon_n = (\log n / n)^{\frac{1}{2+ \beta_{\text{tree}} d }}$ with $\beta_{\text{tree}} \geq 1$. Then, over the anisotropic H\"older space $\mathcal{C}^{\bm \alpha}(\Omega, L)$ with $\alpha_l \in (1, 2]$, there exists some $\underline{\eta} >0$, such that
 \begin{eqnarray}
   \label{eq:lowerbound}
  \inf_{f_0 \in \mathcal{C}^{\bm \alpha}(\Omega, L), f_0 \geq
   \kappa_L, 
   \text{ and Condition~\ref{cond:boundary} holds}} \E_{f_0}^n \Pi_n \left( \{ f: \rho(f_0, f) < \underline{\eta} \epsilon_n  \cond \Y_n \right) \rightarrow 0.   
  \end{eqnarray}
\end{theorem}}

\rev{The proof for Theorem~\ref{thm:treelowerbound} will be provided
  in Section~A.3 of the Supplementary
  Material \citep{liu2026b}. 
  The theorem reveals that 
 the posterior contraction rate for
 this single-tree method cannot be faster than $n^{-\frac{1}{d+2}}$ up to a logarithmic term. After simple calculation, we see that when
 $1 < \beta < \alpha_{\min} + 2(\alpha_{\min}-1)/d$, the rate for BDF
 is faster than that of the single tree. Moving beyond this specific
 single-tree model, it is more challenging to establish the lower bound
 under the OPT prior introduced in Section~\ref{subsec:opt}. We will
 leave this challenge for future work.}


 \rev{
   In Theorem~\ref{thm:BMArateanisotropic}, the faster rate for BDF
 over the single-tree method is mainly due to the better approximation
 accuracy of forests. From a theoretical perspective, the improved 
 approximation accuracy of forests is first rigorously established by
 \cite{arlot2014analysis}. This result has been later extended by a
 number of works, such as those by \cite{mourtada2020,
   cattaneo2023inference,pmlrcui21,randrianarisoa2021smoothing}. In
 particular, \cite{randrianarisoa2021smoothing} has shown that with
 carefully designed unequal weights, shifted tree ensembles can approximate
 ubitrarily smooth density functions well. 
 As we consider equal weights for tree aggregation rather than the specific
 weighting scheme introduced by \cite{randrianarisoa2021smoothing}, we limit $\alpha_l
 \leq 2$ in Theorem~\ref{thm:BMArateanisotropic}.
}

If we remove the boundary condition, the following corollary shows that insufficient approximation accuracy at the boundary may lead to a slower rate, especially for smoother density functions.
\begin{corollary}
  \label{cor:BMArateanisotropic}
Under the same conditions as those for Theorem
\ref{thm:BMArateanisotropic}. If we do not impose any boundary
assumptions, then with $\alpha_{\min}^* = \alpha_{\min} \wedge
\frac{3}{2} $, the convergence rate for \writing{BDF} over the anisotropic
H\"older space $\mathcal{C}^{\bm \alpha} (\Omega, L)$ is 
\begin{eqnarray*}
  \epsilon_n =
  \begin{cases}
   (d^{3/2}L/\sqrt{\kappa_L}) 
   ^{\frac{\beta+2/d-\alpha_{\min}^*/d}{\beta+2/d}}  
   n^{-\frac{\beta + 2/d-1}{2\beta+ 4 /d}}, 
& \mbox{when } 1 \leq \beta \leq 1 + 2(\alpha^*_{\min} -1 )/d;\\
(d^{3/2}L/\sqrt{\kappa_L})^{\frac{1}{\beta+2/d}} 
n^{-\frac{ \alpha_{\min}^* /d}{\beta +2 /d}},
& \mbox{when } \beta \geq 1 + 2(\alpha^*_{\min} -1 )/d. 
  \end{cases}
\end{eqnarray*}
\end{corollary}

The proofs for 
Theorem~\ref{thm:BMArateanisotropic} and Corollary~\ref{cor:BMArateanisotropic} will be provided
in Section~A 
of the Supplementary Material \citep{liu2026b}. In these results, the prior $\Pi$ or its truncated counterpart $\Pi_n$ is supported on a collection of trees of different depths and perturbations. We have introduced the hyper-parameter $\beta$ to adjust the
strength of penalty on tree depth, which plays a central role in the
bias-variance trade-off. As our base learner is a single tree, the
trade-off is achieved at the individual tree level. \rev{When $\beta =1$,
 rates of the BDF estimator are similar to that of a single-tree, and the proof is easier. A more challenging case
is when $\beta > 1$.}

To reveal the bias-variance trade-off, we first check the
approximation error under the KL divergence when $f_0$ lies in  the isotropic
or anisotropic H\"older space. 
Specifically, \rev{as shown
  in Lemma F.2 in the Supplementary Material \citep{liu2026b},} 
there exist $r >0$ and $c_{\appr} >0$, such that under any perturbation $h$ and at any depth $\D$, there is a depth-$\D$ tree-supported density function $f_\D^h$ such that
$\kl(f_0, f_\D^h) \leq (c_{\appr})^2 2^{-2r \D}$. For any $\beta \geq
1$, we define 
\begin{eqnarray}
  \label{eq:depthbeta}
 \D_\beta = 
\left \lfloor
\log_2 ( (c_\appr)^2 n / \log n)^
{\frac{1}{\beta + 2/d}}
\right \rfloor.
\end{eqnarray}
Our proof strategy is to first show that the posterior probability
mass on trees with depth much larger than $\D_\beta$ is exponentially
small. Then, when $\beta >1$, over the trees with depth around $\D_\beta$, an
aggregation of perturbed trees can bring the improvement.

\rev{With this proof strategy, we introduce the strong shape constraint and
Condition~\ref{cond:uniqueness}, as well as assumptions on the
smoothness level ($\alpha_l >1$) and on the prior ($\beta > 1$), to ensure that 
 at any given depth $\D$ around $\D_{\beta}$, posterior distributions under shifted OPT priors, denoted as
 $\Pi_n(\cdot | \Y_n, \D, h)$, concentrate on shifts of
 the same tree topology. This concentration result is required due to a technical artifact that effective smoothing can only be achieved through perturbations of the same tree topology. For less
 smooth densities ($\alpha_l \leq 1$), showing tree
 structure recovery may require additional conditions, which we do not
 pursue here.}

\rev{For BDFs with $\beta >1$, the tree-level bias-variance trade-off and tree topology
  recovery occur in the regime where the
  bias dominates the variance. Hence, at depths $\D$ around $\D_{\beta}$, the tree topology that best
  approximates $f_0$ is recovered with high probability. Because we
  assume the smoothness level $\alpha_l >1$, 
  while a single tree can only approximate Lipschitz-continuous functions
  well, the depth-$\D$ tree topology with smallest bias under the strong shape
  constraint is expected to have leaf
  nodes with edge length of order or smaller than $2^{- \lfloor \D / d \rfloor}$. The choice of shifting unit in BDF is consistent with this expectation.}

\subsection{Data Adaptivity of Individual Trees}
\label{subsec:likrisk}
Under the classical Bayesian framework, the likelihood function is employed to update the distribution of model
parameters after data come in, and plays a key role in making
data-adaptive inference. This is also the case for tree-based
methods. To better illustrate how well the likelihood function can be employed to adapt to data under a tree model, we provide a large-deviation type inequality displaying a bias-variance decomposition for the
log-likelihood (ratio).  

When we study the order of the variance for tree-based log-likelihood, we find that after appropriate scaling it usually depends on both the sample size $n$ and the number of leaf nodes $t$, and roughly takes a form $(t\log t) /n$ when both $n$ and $t$ are large enough. To avoid the difficulty of deriving sharp concentration results for shallow trees, we introduce the following quantities, 
\begin{eqnarray}
  \label{eq:depthlowerbound}
\D_{\min} = \lfloor \log_2 (n/ \log n)^{\frac{1}{10}} \rfloor, 
\quad
\mbox{and} \quad
t_{\min} = 2^{\D_{\min}}.
\end{eqnarray} 
For any depth-$\D$ tree with $\D< \D_{\min}$, we only show
concentration results as strong as those for a depth-$\D_{\min}$
tree. The exponent $\frac{1}{10}$ of $(n/\log n)$
in \eqref*{eq:depthlowerbound} can be replaced by any fraction
$\frac{1}{2r+1}$ with $r >2$. We also define 
\begin{eqnarray}
  \label{eq:deltadef}
  \delta(n, t) =
  \left( \frac{(t  \vee t_{\min} )\log (t \vee t_{\min}) }{n } \right)^{1/2}.
\end{eqnarray}
This quantity describes the order of the variance for the
scaled log-likelihood ratio based on a tree  with $t$ leaf nodes, and will appear multiple times later. 
The following theorem shows a bias-variance decomposition for the
likelihood function, and its proof will be provided in Section~D 
of the Supplementary Material \citep{liu2026b}. 
\begin{theorem}
 \label{thm:likbiasvar}
Given a tree with $t$ (\tech{$\leq \lfloor n / \log n \rfloor$}) leaf nodes $\mathcal{A}_t = \{\Omega_j\}_{j=1}^t$, 
we define the ``projection'' of the true density function onto the tree as
 \begin{eqnarray*}
  f_{\mathcal{A}_t}= \sum_{j=1}^t \frac{\mu_{f_0}(\Omega_j)}{\mu(\Omega_j)}
    \cdot \mathbbm{1}_{\Omega_j}, \quad \mbox{where } \mu_{f_0}(\Omega_j) =
    \int_{\Omega_j}f_0 d \mu.
  \end{eqnarray*}
$\mathcal{F}_{\mathcal{A}_t}$ denotes the collection of tree-like densities supported by $\mathcal{A}_t$. Let \tech{$\delta(n, t)$} be the quantity defined in \eqref*{eq:deltadef}. For $\eta>1$ and some constants $c_1, c_2, c_2' \in (0, 1)$, we have the following inequality when $n$ is sufficiently large, 
  \begin{eqnarray}
    \label{eq:OneSideIneqUpper}
  && \hspace{3mm}  \P_{f_0}^n \bigg( 
\sup_{f \in \mathcal{F}_{\mathcal{A}_t}} 
\prod_{i=1}^n \frac{f(Y_i)}{f_0(Y_i)} 
\geq 
\exp \Big( \underbrace{-n (1-c_1) \rho^2(f_0, \fat)}_{\mbox{bias}}
+
\underbrace{\eta^2 n \delta^2(n, t)}_{\mbox{variance}} \Big)
    \bigg)\\
\nonumber  && \hspace{1.8in}  \leq
\exp\left( -c_2' n \rho^2(f_0, \fat) \right)
+
    8\exp\left(-c_2 \eta^2 n \delta^2(n, t) \right),
  \end{eqnarray}
where $\P_{f_0}^n$ is the probability measure induced by $\Y_n$. Specific choices for $c_2$ and $c_2'$ will be provided in the
proof. 
\end{theorem}
\noindent Note that the theorem applies to any tree topology, with or without shape constraints. The purpose of this theorem is to reveal the asymptotic behavior of the
tree-based likelihood, which is independent of what type of tree
priors we impose. We also do not require $f_0$ to be bounded above or
bounded away from zero. In Section~F 
of the Supplementary Material \citep{liu2026b}, we will show that $\fat$ is indeed the best $\A_t$-supported tree
approximation to $f_0$ under the KL divergence. Based on inequality \eqref*{eq:OneSideIneqUpper}, it is clear that the highest likelihood that can be achieved by a tree topology
is driven by both bias and variance: the better approximation ability a tree topology has, the larger the likelihood is; while the fluctuation of the likelihood ratio due to sampling variability is bounded by $\eta^2 n \delta^2(n, t)$ with high probability under $\P_{f_0}^n$. 

In summary, the results in this section show that in certain cases forests can achieve faster convergence over trees for the density estimation problem. When the prior is appropriately imposed, the rate can be minimax for the isotropic class $\mathcal{C}^\alpha(\Omega, L)$, $\alpha \in (1, 2]$. But the rate is not adaptive, in the sense that in order to achieve the optimal rate, prior specification requires knowledge of the smoothness level $\alpha$.
For the anisotropic class, the \writing{BDF} generally cannot achieve
minimax convergence rate. 
This motivates us to study
an alternative model for forests, which will be described in the next section.

\section{Gibbs Density Forests}
\label{sec:gibbs}

The reason why the \writing{BDF}
estimator~\eqref*{eq:pseudoposteriormean} cannot achieve the optimal
rate especially over the anisotropic class is that it only allows for 
bias-variance trade-off at the individual tree level, which has an
explicit likelihood that can be decomposed following
Theorem~\ref{thm:likbiasvar}.
\rev{For forests, it is generally more challenging to formulate the
  phenomenon of bias-variance trade-off.} \rev{Our next objective is to overcome this limitation by incorporating
loss-based generalization to the likelihood decomposition so that new
bias-variance decomposition result becomes available at the forest
level. The decomposition further motivates the design of a new
forest algorithm that achieves minimax rates adaptively.}

\subsection{Risk-Based Generalized Bayesian Aggregation of Density Forests}
\label{subsec:gibbsframework}
For density estimation, a natural ``loss function'' is $-\log(f)$. The
risk associated with such a loss corresponds directly to the KL
divergence $\kl(f_0, f)$. This observation motivates us to consider a
generalized Bayesian inference framework that makes bias-variance trade-off via \emph{information risk minimization} \citep{zhang2006a,zhang2006b}, instead of the usual likelihood-based way. 
Specifically, following the idea of fitting ``honest forests'' \citep{biau2012analysis,denil2014,athey2019}, 
we randomly split the dataset $\Y_n$ into two subsets of sizes $n_0$
and $n_1$, respectively: {$\Y^\0= \{Y_1, \ldots, Y_{n_0}\}$} and {
  $\Y^\1 = \{ Y_{n_0+1}, \ldots, Y_{n}\}$}, with $n_0 \asymp n_1$ as
$n \rightarrow \infty$. We use the first subset of data for density estimation on leaf nodes given the trees through vanilla Bayesian estimation, and the second subset for learning of the trees in the ensemble through risk-based generalized Bayes. 

We use $\newforest = \forest$ to denote the collection of trees in an ensemble, where each $\T_b$ represents a tree as introduced in \pref{Section \ref{subsec:opt}}, with or without
perturbation. In distinguishing the discussion on a single
tree versus the whole ensemble, we refer to $\T_b$ as a {\em tree}
topology, and to { $\newforest= \forest$} a {\em forest} topology. The
forest topology is of depth $\D$ if each tree within the ensemble is
of depth at most $\D$. Given the forest topology, we estimate density
for different trees \emph{separately}, under the OPT prior. For each tree $\T_b$, since the tree topology is given, the (conditional) prior is only on random probability splits along the tree
nodes, which we denote as $\Pi (\cdot | \T_b)$. We use the
(conditional) posterior mean given the tree $\T_b$ and the first subset of data $\Y_n^\0$ for the estimation of densities on leaf nodes, and denote the predictive density given $\T_b$ by {$\tilde{f}_{\T_b} = E (f|\Y_n^\0, \T_b )$}. 
We still construct an ensemble estimate as an unweighted average over trees, analogous to RFs in the supervised setting, 
\begin{eqnarray}
\label{eq:rf}
\tilde{f}_\newforest = \tilde{f}_\forest =
 \frac{1}{B} 
 \left( \tilde{f}_{\T_1} + \tilde{f}_{\T_2} + \cdots + \tilde{f}_{\T_B}\right).
\end{eqnarray}
Compared to a fully Bayesian approach which would take a weighted average with respect to the posterior probability of each tree topology, the unweighted average is much easier
to compute in practice. It also makes the theoretical analysis
easier. Note that as long as each weak learner {$\tilde{f}_{\T_b}$} is a probability density, the forest is also a valid density. 

We measure the goodness of each forest topology using the following ``empirical risk'' 
\begin{eqnarray}
  \label{eq:riskdef}
  R_n \newbforest = -\frac{1}{n_1} \sum_{i=n_0+1}^n
  \log \left(  
\tilde{f}_\newforest
  (Y_i) \right).
\end{eqnarray}
We also define a ``centered'' version of risk function as 
\begin{eqnarray}
  \label{eq:modifiedriskdef}
  \wt{R}_n \newbforest = R_n \newbforest
  + \frac{1}{n_1} \sum_{i=n_0+1}^n \log f_0 (Y_i). 
\end{eqnarray}
The centered risk $\wt{R}_n$ is typically impossible to
calculate as $f_0$ is unknown and it serves merely as a
theoretical device. Later it will be clear that the analysis of the
centered risk is sufficient for deriving convergence rate for
forests.

\rev{When analyzing the likelihood
  associated with forests, an upper
bound for the variance is usually obtained by quantifying the complexity of the space of forest-supported functions with its metric entropy
or Vapnik-Chervonenkis (VC) dimension. However, the bound for the
entropy or the VC dimension of the
space of {\em all} forest-supported functions may significantly
overestimate the complexity, as
most forest algorithms
only search over specific types of ensembles, which is a very small
subset of the full space. As a consequence, the variance bound for the
likelihood can be quite loose due
  to the necessity to account for the complexity of the much finer
  forest-level partition. 
 Moreover,
the fully likelihood-based decision framework may limit the
possibility of accommodating commonly used
training techniques, such as tree-level regularization and sample-splitting, in the theoretical analysis.}

This new risk \eqref{eq:riskdef} enjoys several important properties that are essential
in overcoming the difficulty imposed by the model complexity of the
ensemble. First, the risk only involves leaf nodes fit at the
individual tree level. This implies that the variance for leaf node
estimation is \emph{not} related to the number of leaf nodes in the
induced forest-level partition, which can involve exponentially more
nodes than the individual trees. Second, one can readily control the number of possible forest topologies by imposing shape constraints on {\em individual} trees without complicating the evaluation of this
risk. The shape constraints, as described in Section
\ref{subsubsec:shape}, serve as regularization that imposes the
necessary penalty on the model complexity to ensure that, \rev{with these
constrained tree base learners, the number of
possible forest topologies does not grow too fast as the depth
increases.}

We consider generalized Bayesian inference on forest-type tree ensembles via information-risk minimization based on the
empirical risk \eqref*{eq:riskdef}, or equivalently the centered risk
\eqref*{eq:modifiedriskdef}. The corresponding {\em Gibbs posterior}, denoted by $\lambda$, with respect to a prior $q$ on forest topologies is given as 
\begin{eqnarray}
  \label{eq:gibbsposterior}
  \lambda^\g  = e^{-n \psi R_n} q \bigg /
  \bigg( 
\sum_{\newforest}
  e^{-n \psi R_n \newbforest}
  q \newbforest \bigg).
\end{eqnarray}
The specification of $q$ will be discussed later. Previous work on
risk-based generalized Bayesian inference \citep{zhang2006a, zhang2006b} have shown that among all measures $\lambda$ on the space of forest topologies with respect to the prior $q$, the Gibbs posterior is the minimizer for
\begin{eqnarray*}
\sum_{\newforest }
  \lambda  \newbforest
  n R_n \newbforest 
  +
  \psi^{-1} \kl (\lambda  , q),
\end{eqnarray*}
where 
$\psi^{-1}$ is a tuning parameter called the
``temperature''. \rev{Lower temperatures lead the Gibbs posterior to
  depend more heavily on the risk.} The choice of $\psi$ relies on
properties of the empirical risk and can be critical for making
generalized Bayesian inference \citep{syring2023,ommen2017}. We will have a
further discussion of the choice later. Analogous to the likelihood-based Bayesian inference,
under the generalized Bayesian framework,
we may consider defining the forest estimator
either as the ``posterior mode'' or as the ``posterior mean''. We 
call this type of \writing{DFs} as the
\writing{{\em Gibbs density forest (GDF)}}. More specifically, we define 
\writing{{\em MAP GDF}} as 
\begin{eqnarray}
  \label{eq:gibbsMAP}
  \hspace{0.2in}      \hat{f}_{\mbox{\tiny{MAP}}}^\g =
  \tilde{f}_{\widehat{\newforest}_{\mbox{\tiny{MAP}}}}, 
  \quad \mbox{where } 
\widehat{\newforest}_{\mbox{\tiny{MAP}}}= 
\argmax_{\newforest} \hspace{1mm} \lambda^\g (\newforest)
\quad 
\mbox{\writing{(MAP GDF)}},
\end{eqnarray}
and \writing{{\em predictive GDF}} as 
\begin{eqnarray}
  \label{eq:gibbsBMA}
  E^\g( f | \Y_n ) =
  \sum_{\newforest}
  \lambda^\g  \newbforest
\tilde{f}_{\newforest} \quad 
\mbox{\writing{(predictive GDF)}.}
\end{eqnarray}
If one can draw forest samples of high quality from the Gibbs posterior, then by further aggregating over multiple forests,
\writing{predictive GDF} \eqref*{eq:gibbsBMA} can achieve enhanced predictions through fully accounting for model uncertainty under the Bayesian framework. 
However, if a small fraction of forest
samples is of poor quality, the \writing{MAP GDF}
\eqref*{eq:gibbsMAP} generally produces more robust estimates.
\rev{In Section~\ref{subsec:gibbsrate}, we investigate the asymptotic
properties of the predictive GDF and establish its consistency
rates. In contrast, deriving rates for the MAP GDF is more
challenging. 
We only evaluate its performance in the simulation studies.}

\subsection{Properties of the Risk}
\label{subsec:riskproperty}
Under the generalized Bayesian framework, we are interested in whether
the risk function for forests behaves similarly to the negative
log-likelihood for trees such that a ``bias-variance trade-off'' can
be formulated and accordingly inference based on the generalized
likelihood can effectively adapt to the data without overfitting. To
answer this question, we provide a lower bound for the centered risk similar to the bound \eqref*{eq:OneSideIneqUpper} for likelihood ratios. 

Given a forest topology $\newforest$, we define the forest approximation to $f_0$ as
\begin{eqnarray}
  \label{eq:forestapprox}
  f_{\newforest} = \frac{1}{B} \sum_{b=1}^B f_{\T_b},
  \quad
  \mbox{with }
  f_{\T_b} = \sum_{j:\Omega_j \in \T_b}
  \frac{\mu_{f_0}(\Omega_j)}{\mu(\Omega_j)} \cdot
  \mathbbm{1}_{\Omega_j}.
\end{eqnarray}
Note that this density function is \emph{not} data-dependent once the forest topology is given. 
We still assume that the forest depth $\D$ cannot be greater than $\D_{\max}$, to ensure that on average there are enough data points in the leaf nodes, where $\D_{\max}$ is the same as that defined in \eqref*{eq:maxdepth}.

Recall that we split the data for empirical risk evaluation. For each
tree $\T_b$ within a forest, density learning is based on the first
subset of data and is finished via vanilla Bayesian inference under
the OPT prior. Due to certain technical issues, in the remainder of
this section we assume that the Bayesian inference for this step is
under an {\em upper truncated} OPT prior. That is, the density
functions sampled under a depth-$\D$ OPT prior cannot take values
greater than $ c_{\text{opt}} (\D \vee \D_{\min})$, where $\D_{\min}$
is defined in \eqref*{eq:depthlowerbound} and $c_{\text{opt}}\geq 1$ is
a constant. Such a modification will not affect the inference for
bounded density functions too much and is not implemented in
simulation studies either. We may avoid this technical modification if
we can show {\em uniform} Bayesian sup-norm convergence for density estimation on
leaf nodes given the tree topologies \citep{castillo2022}.  

When $f_0$ satisfies $\kappa_L \leq f_0 \leq \kappa_U$ for some
constants $\kappa_U \geq \kappa_L >0$, for a depth-$\D$ forest topology $\newforest$, we establish a lower bound for
the centered risk $\wt{R}_n \newbforest$ of the following form,
\begin{eqnarray}
\label{eq:risklowerpaper}
&& \wt{R}_n \newbforest 
\geq
(1-c_1-c_3) \kl (f_0, f_{\newforest}) \\
\nonumber && \hspace{1.2in}
             - \frac{c_{\text{opt}} \kappa_U^2 \D  }{c_3 \kappa_L^2} 
M_\tr^2 \cdot \eta_0^2 \delta^2(n_0, 2^\D) 
-\eta_1^2 \xi_n \delta^2(n_1,2^\D),
  \end{eqnarray} 
where $\eta_0, \eta_1, M_\tr \geq 1$ are some constants, $c_1$ and
$c_3$ take values between 0 and $\frac{1}{2}$, $\delta(
\cdot,\cdot)$ is the function defined in \eqref*{eq:deltadef}, and
$\xi_n \asymp (\kappa_U / \kappa_L) \log n $. We will provide a rigorous characterization
of the risk lower bound in Section~C 
of the Supplementary Material \citep{liu2026b}, where it becomes clear that the bound
\eqref*{eq:risklowerpaper} holds with probability tending to 1 under
$\P_{f_0}^n$. Moreover, when the number of depth-$\D$ forest topologies does
not increase too fast as $\D$ increases, this bound also holds {\em
  uniformly} for all forest topologies.  

The bound \eqref*{eq:risklowerpaper} suggests that the centered risk
\eqref*{eq:modifiedriskdef} behaves analogously to the log-likelihood
ratio, in the sense that we can derive a lower bound for the centered
risk, which shows a similar form of bias-variance decomposition. The
inequality \eqref*{eq:risklowerpaper} also implies that the variance
associated with the risk \eqref*{eq:modifiedriskdef} comes from two sources: density learning in leaf nodes of trees based on the first subset of data and risk evaluation based on the second subset. With high probability, the former source of variance is bounded by $\frac{c_{\text{opt}} \kappa_U^2 \D  }{c_3 \kappa_L^2} M_\tr^2 \cdot \eta_0^2 \delta^2(n_0, \rev{2^\D}) $, and the latter part is bounded by $2\eta_1^2 \xi_n \delta^2(n_1,\rev{2^\D})$. The bounds for these two parts are of the same order up to logarithmic terms. This explains why we set $n_0$ and $n_1$ of the same order when splitting the data. The ``optimal" sample splitting ratio can be obtained by balancing these two sources of variance. However, such an optimal ratio also requires knowledge of $f_0$, which is inaccessible in practice. Later in Section \ref{sec:simulations}, we will use simulated data to illustrate the effect of the sample-splitting ratio on the performance of \writing{GDFs}. 

\subsection{Contraction Rates for the Gibbs Posterior on Forest
Topologies}
\label{subsec:gibbsrate}
Our next objective is to investigate how fast the Gibbs posterior on forest topologies concentrates around the ``best'' ones. Rigorously, we analyze the contraction of the Gibbs posterior \eqref*{eq:gibbsposterior}. 
With a focus on forest topology learning, we also clarify the meaning of ``posterior contraction'', when making risk-based generalized Bayesian inference for density forests. 

Although the standard machinery for deriving posterior contraction 
rates \citep{ghosal00,shen2001,syring2023} cannot be applied directly,
we consider using the same proof strategy with some technical
modification. Specifically, we identify three conditions explicitly on
the design of tree ensembles and prior specification that ensure the
proper behavior of the Gibbs posterior. These conditions are similar
in flavor to those well-known ones first formulated by \cite{ghosal00}
to establish posterior contraction rates for likelihood-based Bayesian inference in general cases, although there may not be an exact one-to-one correspondence. 

Our first condition concerns the forest approximation error of the true density. 
\begin{conditionp}{(APPR)}
\label{cond:approx} 
The true density function $f_0$ can be well approximated by forests at
rate $r$. More specifically, for any depth $\D$ large enough, one can
always find a depth-$\D$ forest topology $\newforest$, such that
$\rho(f_0, f_{\newforest}) \leq c_{\appr} \cdot 2^{- r \D}$ for some
constant $c_{\appr} >0$ that does not depend on $r$ or $\D$, where $f_{\newforest}$ is as defined in Equation \eqref*{eq:forestapprox}. 
\end{conditionp}
\noindent This condition is satisfied by density functions lying in the H\"older space, isotropic or anisotropic. As an example, for the space $\mathcal{C}^\alpha(\Omega, L), \alpha \in (0, 1]$, Condition \ref{cond:approx} is satisfied with $r = \alpha /d$ and $c_{\appr} = dL$. In addition to the H\"older space, this condition also covers several other frequently considered function classes, such as piecewise or region-wise H\"older density functions and a mixture model with each component being H\"older continuous. 

Our second condition is about the prior distribution $q$ on forest topologies. 
\begin{conditionp}{(PRIOR)}
\label{cond:prior} 
There exist two constants $\overline{\gamma} > \underline{\gamma} >1$, such that for any depth-$\D$ forest topology $\newforest$, 
\begin{eqnarray*}
  q \newbforest \gtrsim \exp(-\overline{\gamma} \D \cdot 2^\D),
\end{eqnarray*} 
and the prior on the collection of all depth-$\D$ forest topologies satisfies 
\begin{eqnarray*}
  \sum_{\newforest \mbox{\scriptsize{ of depth $\D$}}} q \newbforest \lesssim
  \exp(-\underline{\gamma} \D \cdot 2^\D),
\end{eqnarray*}
where the notation ``$\lesssim$'' and ``$\gtrsim$'' mean the inequality holds up to a constant.  
\end{conditionp}
\noindent The first half of Condition \ref{cond:prior} in combination with Condition \ref{cond:approx} can lead to a result similar to the prior thickness, which is a standard condition formulated to study posterior contraction under the likelihood-based Bayesian framework. And the second half of Condition \ref{cond:prior} ensures that the prior decays sufficiently fast as the forest depth increases.   

Our last condition is imposed on the complexity of the forest space.
\begin{conditionp}{(COMP)}
  \label{cond:complexity} 
The number of depth-$\D$ forest topologies is bounded by \\$\exp(c_{\tiny \en} \D \cdot 2^\D)$ for some constant $c_{\tiny \en}>0$. And the number of trees in a depth-$\D$ forest is of order $o(\exp(2^\D))$.
\end{conditionp}
\noindent This condition ensures that the forest space is not too
complex, so that we can obtain a uniform empirical risk bound of the
form \eqref*{eq:risklowerpaper} for all forest topologies. It is
analogous to the entropy condition considered by \cite{ghosal00,
  shen2001}. In the proof, we will see that this condition holds under the weak shape constraint. 

We define
\begin{eqnarray*}
  \C_\en \coloneqq \left\{ \newforest: \rho(f_0, f_{\newforest}) \geq M_n
  \left(c_{\appr} \right)^{\frac{1}{2r+1}}
\left( n / \log n \right)^{-\frac{r}{2 r+1}}\right\}, 
\end{eqnarray*}
where $M_n$ increases at most as polynomials of $\log n$. 
We have the following result on posterior contraction of the Gibbs posterior on forests. 
\begin{theorem}
  \label{thm:gibbsrate}
Assume that Condition \ref{cond:approx}, Condition \ref{cond:prior} and Condition \ref{cond:complexity} hold. In addition, we assume that there exist constants $\kappa_U \geq \kappa_L >0$, such that $\kappa_L \leq f_0 \leq \kappa_U$. 
Then, for the Gibbs posterior defined in \eqref*{eq:gibbsposterior}
\tech{with an appropriately chosen $\psi$ satisfying $\psi =O \left(
    1/ \log n \right) $}, 
 \begin{eqnarray*}
   \lambda^\g 
\left( \C_\en \right)  
\tech{ \lesssim \exp 
\left(- \frac{1}{2} (\underline{\gamma} -1) (n / \log
   n)^{\frac{1}{2r+1}} \right) } 
\quad 
\mbox{in } \P_{f_0}^n\mbox{-probability}. 
 \end{eqnarray*}
\end{theorem}
\noindent The theorem characterizes how fast the Gibbs posterior
concentrates on promising forest topologies. The goodness of a forest
topology is evaluated by the Hellinger distance between $f_0$ and the
approximation of the form \eqref*{eq:forestapprox} supported by this topology. When density learning on leaf nodes also concentrates fast around the forest
approximation, which is the case as we argued in Section~\ref{subsec:riskproperty}, we expect that the \writing{predictive GDF} 
\eqref*{eq:gibbsBMA} also converges to $f_0$ at a fast rate. In the theorem,
we require that $f_0$ be bounded away from zero. This condition is
commonly imposed when studying posterior contraction for the density
estimation problem \citep{randrianarisoa2021smoothing, ghosal2001Bern,
  castillo2022}. In our proof, it is needed to derive \rev{the
  lower bound \eqref*{eq:risklowerpaper} for the empirical risk.}
We will provide the proof for Theorem~\ref{thm:gibbsrate}
in Section~B of the Supplementary
Material \citep{liu2026b}. 

\rev{The choice of $\psi$ relies on the empirical risk bound
\eqref*{eq:risklowerpaper}. In specific, it is determined by the
bound for the variance component of the form $ \frac{c_{\text{opt}} \kappa_U^2 \D  }{c_3 \kappa_L^2} 
M_\tr^2 \cdot \eta_0^2 \delta^2(n_0, 2^\D) 
+\eta_1^2 \xi_n \delta^2(n_1,2^\D)$ and the prior. Due to the logarithmic
terms $\D$ and $\xi_n$ both of order $ \log n$ and several constants, appearing 
in the empirical risk bound, we require $\psi = O(1 / \log n)$
to offset these extra logarithmic factors and unknown constants, in addition to the penalty
imposed by the prior distribution, which accounts for the terms
$ \delta^2(n_0, 2^\D) $ and $\delta^2(n_1,2^\D)$.}

\subsection{Gibbs Density Forests with Perturbations}
\label{subsec:gibbsholderrate}
\rev{Next, we provide an example of GDF in which 
trees are grown under perturbations. In this example, we also investigate the effects
of decorrelation with perturbed trees and} 
show that forests with perturbations can better adapt to spatial features of the underlying distribution under the generalized Bayesian framework. 

As the smoothness level of the true density may vary spatially, the perturbation scheme also needs to be adjusted to better adapt to data. 
To design the perturbation, a useful observation is that for
tree-based methods, the average size of the leaf nodes within a region
is usually a reflection of the smoothness level therein. Following
this observation, to design perturbations that can accommodate
spatially varying smoothness levels, we pair each design with a set of possible average grid sizes along different dimensions. Such type of pairing is quite straightforward under shape constraints, as with those constraints the grid size of the partition, or equivalently the shape of leaf nodes, within a region remains the same.  

More specifically, at depth $\D$, each nonnegative integer solution
$\bm k = (k_1, \ldots, k_d)$ to the equation $\sum_{l=1}^d k_l = \D$
can uniquely characterize the shape of a level-$\D$ leaf node without
any perturbation---the edge length of the rectangular leaf node along
the direction $l$ is exactly $2^{-k_l}$. Based on this shape, we can
design a perturbation scheme, where the shift unit along the $l$-th dimension is  $\upsilon_l = 2^{-2k_l}$. Potential shifting directions under this design are,  
\begin{eqnarray*}
    \cH_{\bm \upsilon } = \left \{ h = (i_l \upsilon_l)_{1 \leq l \leq d}: i_l \in
  \mathbb{N}, 0 \leq i_l < \upsilon_l^{-1/2} \right\},
  \quad
  \mbox{where }
  \bm \upsilon = (\upsilon_1, \ldots, \upsilon{_d}). 
\end{eqnarray*}
We use the vector $\bm \upsilon = (\upsilon{_1}, \ldots, \upsilon{_d})$ to denote a perturbation design scheme, which implicitly relies on the depth $\D$. 
The collection of perturbation design schemes at depth $\D$ is 
\begin{eqnarray*}
    \Upsilon_\D =\bigg \{ \bm \upsilon: \upsilon_l = 2^{-2k_l} \mbox{ for
     some } \bm k = (k_l)_{1 \leq l \leq d} \mbox{ satisfying } k_l \in
     \mathbb{N} \mbox{ and } k_1 + \cdots k_d = \D \bigg \}.
\end{eqnarray*}
Note that the cardinality of $\mathcal{H}_{\bm \upsilon}$ is exactly $2^\D$, independent of the specific design scheme $\bm \upsilon$. 


To grow a depth-$\D$ forest with perturbation, we grow multiple trees of the same depth separately, each under a small perturbation. The number of trees within a forest can be pre-determined or random with a prior. 
Given the depth, we randomly choose a perturbation design scheme $\bm
\upsilon$ from the set $\Upsilon_\D$. 
Then, 
the potential shifting directions are provided by
$\mathcal{H}_{\bm \upsilon}$. To grow the $b$th tree,
\tech{we first sample a shifting direction $h_b$ from
  $\mathcal{H}_{\bm \upsilon}$,} then we grow the tree under the shift
$h_b$ with the \emph{weak shape constraint} introduced in Section
\ref{subsubsec:shape}. We denote tree density learners within a depth-$\D$ forest as $\{f_{\T_b ^{h_b}}, 1 \leq b \leq B\}$, or simplified as $\{f_{\T_b}, 1 \leq b \leq B\}$, where the shifting directions $h_b$ are omitted whenever there is no confusion. 

We consider the following prior distribution on forest
topologies. First, the prior on forest depth $q(\D) \propto \exp ( -
\gamma \D 2^{\D})$ for some constant $\gamma > 1$, and is upper
truncated at $\D_{\max}$. Given the depth, the number of trees in a
depth-$\D$ forest is expected to satisfy $B \asymp O(2^\D)$. If $B$ is
viewed as random, we may consider a prior of form $q(B | \D) \propto
\exp (- \gamma_B \D \cdot B)$ and $B \leq c_B 2^\D$ for some
constants $\gamma_B, c_B > 1$. The perturbation design scheme is
chosen randomly from $\Upsilon_\D$, with a uniform prior on the
set. With a selected scheme $\bm \upsilon$, \tech{shifting directions
  for trees are uniformly sampled from the set $\mathcal{H}_{\bm
    \upsilon}$.} For the $b$-th tree within the depth-$\D$ forest, if
the sampled shifting direction is denoted as $h_b$, then the prior on
$f_{\T_b}$ is the shifted OPT $\Pi_n(\cdot | h_b)$ under the
\emph{weak shape constraint}.


Under this prior distribution, we have the following theorem on the convergence rate of \writing{GDFs} for densities lying in the anisotropic H\"older space. 
\begin{theorem}
\label{thm:gibbsholderrate}
Assume $f_0$ lies in the anisotropic H\"older space $\mathcal{C}^{\bm
  \alpha}(\Omega, L)$, $\alpha_l \in (0,2]$. $\overline{\alpha}$
satisfies $\overline{\alpha}^{-1} = \frac{1}{d} \sum_{l=1}^d
\alpha_l^{-1}$. In addition, we assume that $f_0$ is bounded away from
0 and Condition \ref{cond:boundary} holds. Under the prior
distribution $q$ on forest topologies specified above \rev{and with
  $\psi = O(1/ \log n)$}, for some sequence $\{M_n\}_{n\geq 1}$ that increases at most as polynomials of $\log n$, the Gibbs posterior on forest topologies contracts at the rate $ M_n (dL)^{\frac{d}{2\overline{\alpha}+d}}\left( n / \log n \right)^{-\frac{\overline{\alpha}}{2\overline{\alpha}+d}}$, in the sense described in Theorem \ref{thm:gibbsrate}. 
For the \writing{predictive GDF}, $E^\g(f | \Y_n)$, there exists a sequence $\{M_n'\}_{n \geq 1}$ that increases at most as polynomials of $\log n$, such that with probability tending to 1 under $\P_{f_0}^n$, 
    $\rho(f_0, E^\g(f|\Y_n)) \leq
    M_n'
    (dL)^{\frac{d}{2\overline{\alpha}+d}}
    \left( n / \log n \right)
    ^{-\frac{\overline{\alpha}}{2\overline{\alpha}+d}}$.
\end{theorem}
\noindent The theorem is a direct result of Theorem \ref{thm:gibbsrate}. We prove this by showing that under the prior distribution introduced in this section, Condition \ref{cond:prior} and Condition \ref{cond:complexity} hold. Condition \ref{cond:approx} is satisfied given the H\"older continuity and the boundary condition. If we remove the boundary condition, we may obtain a slower rate, summarized as the following:

\begin{corollary}
\label{cor:giibsholderrate}
Under the same conditions as those for Theorem
\ref{thm:gibbsholderrate}, if we remove Condition \ref{cond:boundary}
on $f_0$, then the contraction rate of the Gibbs posterior is
$\epsilon_n =  M_n (\rev{d^2
  L})^{\frac{d}{2\overline{\alpha^*}+d}}\left( n / \log n
\right)^{-\frac{\overline{\alpha^*}}{2\overline{\alpha^*}+d}}$, where
$M_n$ increases at most as polynomials of $\log n$,
$(\overline{\alpha^*})^{-1} = \frac{1}{d}
\sum_{l=1}^d(\alpha_l^*)^{-1}$ and $\alpha_l^* = \min\{\alpha_l,
3/2\}$. For the \writing{predictive GDF}, $E^\g(f | \Y_n)$, there
exists a sequence $\{M_n'\}_{n \geq 1}$ that increases at most as
polynomials of $\log n$, such that with probability tending to 1 under
$\P_{f_0}^n$, $\rho(f_0, E^\g(f|\Y_n)) \leq M_n'(\rev{d^2 L})^{\frac{d}{2\overline{\alpha^*}+d}}
\left( n / \log n \right)
^{-\frac{\overline{\alpha^*}}{2\overline{\alpha^*}+d}}$.
\end{corollary}
\noindent The proofs of Theorem \ref{thm:gibbsholderrate} and
Corollary \ref{cor:giibsholderrate} will be provided in Section~B 
of the Supplementary Material \citep{liu2026b}. 

To make inference for GDF with perturbations, \rev{the main 
  difficulty lies in efficiently exploring the space of forest
  topologies. To avoid directly searching over the entire forest space,} we
consider a proposal distribution induced by the following risk
function
\begin{eqnarray}
  \label{eq:proposalrisk}
   R_n' \newbforest = -\sum_{b=1}^B \frac{1}{n_1} \sum_{i=n_0+1}^n
  \log \left(  
\tilde{f}_{\T_b}
  (Y_i) \right).
\end{eqnarray}
When the prior $q$ on trees within a forest is independent OPT, this
proposal distribution allows us to sample trees one at a time.
We will
provide more details on the Bayesian inference for GDF in
Section~H.2 
of the Supplementary Material \citep{liu2026b}. 

\enlargethispage{\baselineskip}
\section{Numerical Results}
\label{sec:simulations}

We apply Algorithm~1 
from Section~H 
of the Supplementary Material \citep{liu2026b} to fit \writing{BDFs}
and Algorithm~2 
presented in the same section for \writing{GDFs}. 
For both methods, we build forests under
perturbations.

\rev{For the implementation of both BDF and GDF, a key step is
  sampling tree density
  learners. We employ a sequential
  importance sampling algorithm for this purpose. The general strategy is to
construct a tree-level proposal distribution that looks one-step ahead
and then adjust
the weights associated with the sampled trees accordingly. Both the
posterior distribution under the OPT prior for BDF and the proposal
distribution induced by $R_n'$ in \eqref{eq:proposalrisk} for GDF enable weight updates in closed
forms. In addition, to prevent rapid weight decay for less promising
tree structures, we implement a
resampling and pruning scheme.}
\rev{Once tree density learners are obtained under perturbations, BDF
  outputs forests as an unweighted average of trees, whereas GDF
  further learns a posterior distribution over aggregated forest
  topologies by 
  evaluating the risk \eqref*{eq:riskdef}, accounting for the proposal
  distribution induced by $R_n'$.}

\rev{Although our theoretical results are established with respect to
  the Hellinger distance, in simulation studies we evaluate
  the performance using the KL divergence, which is a stronger
  pseudo-metric.}
For each simulation setting, we compare the following types of
tree-based methods: the optional P\'olya tree (posterior mode), \rev{a
random-forest-type method for density estimation (unsupervised RF,
\cite{pmlr-wen22c}) for one-dimensional distributions,} the
\writing{BDF} with different choices of shifting unit $\iota_\beta$, and the
\writing{GDF} (\writing{including predictive GDF  and MAP GDF}). We
also compare tree-based methods with the kernel density estimator
(KDE). 
 \rev{More implementation details for these methods are provided in Section~H.4 of the
   Supplementary Material \citep{liu2026b}.}

We consider \rev{three} one-dimensional examples. The first
distribution (mixture of betas) has sharp local variation, while the
second density function (trigonometric distribution) is smooth and
shows relatively mild variation throughout the domain. \rev{The third
example (semi-circle distribution) has zero density on the interval
$[0,0.5]$.} The numerical results are summarized in
Section~H.5 of the Supplementary Material \citep{liu2026b}. 

The next example is of moderately-high dimensions. \rev{For this example, implementing unsupervised RF
  is more challenging, and therefore it is not included in the
  comparison.}
\begin{example}
  \label{eg:mixNorm}(Mixture of Gaussian distributions---5-30 dimensions)
The data are generated from the following location mixture of Gaussian
distributions:

\begin{eqnarray*}
  \begin{pmatrix} Y_1 \\ Y_2 \\ Y_3 \\ \end{pmatrix} &\sim& \frac{1}{2}
  \mathcal{N} \left( \begin{pmatrix} 0.25 \\ 0.25 \\ 0.25
    \\ \end{pmatrix}, \begin{pmatrix} 0.05^2 & 0.03^2 & 0 \\ 0.03^2 &
    0.05^2 & 0 \\ 0 & 0 & 0.05^2 \\ \end{pmatrix} \right) + \frac{1}{2}
  \mathcal{N} \left( \begin{pmatrix} 0.75 \\ 0.75 \\ 0.75 \\ \end{pmatrix},
  0.05^2 I_{3 \times 3} \right), \\
Y_4, Y_5 &\overset{\mbox{\scriptsize{i.i.d.}}}{\sim} &
                                                       \mathcal{N}(0.5, 0.01),
                                                      \quad 
Y_6, \ldots, Y_d \overset{\mbox{\scriptsize{i.i.d.}}}{\sim}  
\frac{1}{2} \mathcal{N}(0.35, 0.01) + \frac{1}{2} \mathcal{N}(0.6,0.05^2),\quad\mbox{for } d >5.
\end{eqnarray*}
\end{example}

\begin{table}[h]
  \centering
  \rev{
  \begin{tabular}{l l c c c}
    \toprule
   & & $d =5$ & $d = 10$ & $d = 30$
    \\
    \midrule
    {\footnotesize Optional P\'olya tree} 
   & & 0.055(0.009)
              &0.081(0.010)
                         &0.185(0.014) \\
    \midrule
    \multirow{5}{*}{{\footnotesize Bayesian density forest}}
    &{\footnotesize $\iota_\beta =2^{-6}$}
     &0.050 (0.009) 
              & 0.072 (0.009) 
                         &0.20 (0.025)\\
    &{\footnotesize $\iota_\beta =2^{-10}$}
     & 0.037 (0.008)
              & 0.053 (0.009)
                         &0.15 (0.024)\\
    &{\footnotesize $\iota_\beta =2^{-14}$} 
     &0.039 (0.009) 
              & 0.056 (0.009)
                         &0.155 (0.024)\\
     &{\footnotesize $\iota_\beta =0$} 
     &0.051 (0.009) 
              & 0.074 (0.009)
                         &0.21 (0.027)\\                
    \midrule
    {\footnotesize {Predictive GDF}} 
        &\multirow{2}{*}{\footnotesize $n_0/n = 0.5$} 
            &0.048 (0.014)
                & 0.061 0.012)
                    & 0.15 (0.020) \\
    {\footnotesize {MAP GDF}} 
        &
            &0.044 (0.012) 
                & 0.0575 (0.011)
                    & 0.145 (0.020)\\
    \addlinespace
    {\footnotesize {Predictive GDF}} 
        &\multirow{2}{*}{\footnotesize $n_0/n = 0.6$} 
            &0.045 (0.009)
                & 0.057 (0.010)
                    & 0.14 (0.021) \\
    {\footnotesize {MAP GDF}} 
        &
            &0.044 (0.0085) 
                & 0.055 (0.009)
                    & 0.14 (0.021)\\
    \addlinespace
    {\footnotesize {Predictive GDF}} 
        &\multirow{2}{*}{\footnotesize $n_0/n = 0.7$} 
            &0.045 (0.006)
                & 0.056 (0.010)
                    & 0.14 (0.021) \\
    {\footnotesize {MAP GDF}} 
        &
            &0.044 (0.0055) 
                & 0.054 (0.010)
                    & 0.14 (0.022)\\
     \addlinespace
    {\footnotesize {Predictive GDF}} 
       &\multirow{2}{*}{\footnotesize $n_0/n = 0.8$} 
            &0.048 (0.008)
                & 0.061 (0.013)
                    & 0.14 (0.026) \\
    {\footnotesize {MAP GDF}} 
        &
            &0.047 (0.007) 
                & 0.060 (0.012)
                    & 0.145 (0.026)\\
    \midrule
    {\footnotesize KDE }&
     & 0.75 (0.0020)
              & 1.76 (0.0041) 
                         &8.53 (0.012)
    \\
    \bottomrule
  \end{tabular}
  }
  \caption[5dnorm]{Numerical results for the mixture-of-Gaussian 
    example. For each case, the sample size is $1 \times 10^4$ and
    the average KL is calculated based on \rev{20}
    replicates. The standard deviation is reported in 
    parentheses.}
  \vspace{-3mm}
  \label{tab:5dnorm}
\end{table}
For one-dimensional examples, KDE performs well when the density is
smooth with mild variation, while can be
less satisfying in the presence of sharp changes. For the 
moderately-high dimensional example, the performance
of KDE worsens as the dimension increases, as it becomes more
and more challenging to choose an appropriate bandwidth.

When examining different types of tree-based methods, the performance
trends across these examples remain consistent. \rev{First, the performance of
  unsupervised RF depends on the choice of tree depth and may
  deteriorate for distributions with sharp local variations.} Second,
BDF and GDF with perturbations can 
outperform the single-tree method in most cases. Third, the
performance of \writing{BDFs} heavily relies on the choice of
$\iota_\beta$, which aligns with our theoretical result from Section~\ref{sec:BMA}. Specifically, \writing{BDFs}
with perturbations can achieve faster convergence when \rev{$\iota_\beta$} is appropriately chosen, but the rate is \emph{nonadaptive}. In our
simulations,  with access to \rev{oracle information}, we can identify an ``optimal'' choice of $\iota_\beta$, corresponding to the ``optimal'' tree depth. However,
such an optimal value \rev{may \emph{not} be} available when the 
true distribution is unknown. Last but not least, the strength 
of the \writing{GDF} approach lies in its adaptivity, by introducing a
novel risk function for the forest and making inference with the
posterior on forest topologies, \writing{GDF} 
adapts to spatial heterogeneity, effectively concentrating on the
most promising forest topologies. Although a data-splitting strategy
is employed in the generalized Bayesian inference, 
\writing{GDF's} performance remains comparable to or even better than
the best \writing{BDF} performance, which uses the full dataset and benefits from oracle information. The performance of GDF is
not sensitive to the sampling-splitting ratio. In all of our
simulations, its performance remains stable when the ratio $n_0/n$ ranges from 0.5
to 0.8. We recommend setting $n_0/n$ between 0.6 and 0.7, as GDF
performs best in this range. Given the fact that the \writing{GDF} is
fully data-adaptive while the best-performing \writing{BDF} requires
knowledge of the underlying true distribution, one may prefer the
former approach in real data analysis. For \writing{GDF}, the MAP-type
estimate performs slightly better in most simulation settings, as the
risk-based Gibbs posterior can be unstable and may inflate the
variance of \writing{predictive GDF}.

\rev{Under the Bayesian framework, we can naturally quantify the uncertainty level
  based on the full posterior distribution. The built-in uncertainty
  quantification (UQ) is also available under
  the generalized Bayesian framework. For GDF, we implement a
  sampling-based approach to UQ and obtain credible sets for the three
  one-dimensional examples. More details of the UQ method and 
  visualizations of the credible sets are provided in Section~H.6 of the
  Supplementary Material \citep{liu2026b}. These results show that empirically UQ under the
  generalized Bayesian framework enjoys reasonably good coverage
  properties. }

\section{Discussion}
\label{sec:discussion}

We have introduced and analyzed forest-type tree ensemble methods for density estimation, specifically focusing on two variants, \writing{BDF} and \writing{GDF}, based on Bayesian and generalized Bayesian tree aggregation, respectively. 
We showed that both \writing{BDF} and \writing{GDF} can achieve faster convergence than single-tree approaches over anisotropic H\"older spaces. In particular, \writing{GDF} stands out in its ability to adapt to spatial heterogeneity. Furthermore, we identified a set of sufficient conditions that ensure the success of \writing{GDFs}, offering valuable insights for the design of forest-type ensembles.

The DF introduced in this paper can alternatively be viewed as a mixture model, where each component corresponds to a tree-supported distribution. This perspective suggests that DFs can potentially be applied to mixture distribution learning and clustering problems. Furthermore, the definition of a GDF can be generalized: rather than aggregating trees with uniform weights $1/B$ as defined in \eqref*{eq:rf}, we may consider an aggregation scheme with unequal weights $\{v_b\}_{1 \leq b \leq B}$, where $\sum_{b=1}^B v_b =1$ and $v_b \geq 0$. As long as the weights are either {\em fixed} or derived from an {\em independent} dataset, empirical risk bounds \eqref*{eq:risklowerpaper} and posterior concentration results in the general case (Theorem \ref{thm:gibbsrate}) still hold. \rev{With unequal weights, we can build other types of
  GDFs for different purposes. For example, if we hope to better
  estimate smooth densities, we may adapt the tree aggregation scheme proposed by
  \cite{randrianarisoa2021smoothing} to our framework, achieving a
  faster rate, even near minimax, over the isotropic H\"older space $\mathcal{C}^\alpha$
  with any smoothness level $\alpha > 0$.} Another important implication of viewing DF as a mixture model is that they can be interpreted as a generative model. For each tree-supported distribution, samples can be drawn through a simple two-step procedure: first, select a leaf node according to the probability mass associated with each node and second, sample a point uniformly within the sampled leaf node. At the forest level, sampling from a DF is straightforward under this mixture model view. 

Under the generalized Bayesian framework, we derived concentration results for the Gibbs posterior on forest topologies. These results not only contribute to understanding convergence rates of the \writing{Gibbs MAP DF} and the \writing{Gibbs predictive DF}, but also provide information of the uncertainty level associated with the full Gibbs posterior. Under the classical Bayesian framework, uncertainty quantification (UQ) for Bayesian CART and optional P\'olya tree models have been carefully investigated by \citep{castillo2021uncertainty,castillo2022}. However, UQ for Bayesian tree ensembles is yet to be rigorously investigated. Under the generalized Bayesian framework, \rev{we provided preliminary numerical results for UQ in the Supplementary Material \citep{liu2026b}, while leaving a rigorous investigation of UQ for GDFs for future work}. 

\enlargethispage{\baselineskip}
\section*{Acknowledgement}
The BDF and GDF code, available at
\href{https://github.com/liulx06/Gibbs_density_forest}{https://github.com/liulx06/Gibbs\_density\_forest}, 
builds upon the Bayesian sequential partitioning (BSP, \citealt{LJW}) code
developed and maintained by Prof.\ Wing Hung Wong's Lab. The authors 
 gratefully acknowledge Prof.\ Wong for sharing the BSP code. The authors also thank
two anonymous reviewers, the Associate Editor, and the Editor for
their constructive comments, which greatly improved the quality of
this paper. Part of the work was completed when the second author was
at Duke University. This research was supported in part by the University of Pittsburgh
Center for Research Computing and Data, RRID:SCR\_022735, through the
resources provided. Specifically, this work used the H2P cluster,
which is supported by NSF award number OAC-2117681.
LM was supported in part by NSF award number DMS-2152999. 

\clearpage

\begin{appendix}

\renewcommand{\thefigure}{S\arabic{figure}}
\renewcommand{\thetable}{S\arabic{table}}
\renewcommand{\theexample}{S\arabic{example}}
\renewcommand{\theequation}{S\arabic{equation}}


\section{Proofs for the Results on BDF in Section~\ref{subsec:hybridbma}} 
\label{sec:proofbma}
Under the hybrid aggregation scheme introduced in Section
\ref{sec:BMA}, the sampling model for data is single tree based. Consequently,
the analysis of asymptotic properties of the \writing{BDF} estimator
\emph{only} \rev{relies} on the prior distribution on
individual trees and the likelihood function associated with the tree model. To better study
the posterior distribution \eqref*{eq:posteriordist}, especially the
marginal 
posterior distribution on tree depth, we introduce the
following quantify for any measurable set $B$ of probability measures, 
\begin{eqnarray}
\label{eq:omegadef}
 \omega_n(B; \D, h, \Y_n ) = \int_B \prod_{i=1}^n \frac{f(Y_i)}{f_0(Y_i)} d
  \Pi_n (f | \D, h).
\end{eqnarray}
This quantity plays an important role in the proofs of
Theorem~\ref{thm:BMArateanisotropic} and
Corollary~\ref{cor:BMArateanisotropic}. To characterize the decay of the posterior
distribution \eqref*{eq:posteriordist}, we need a lower bound for the
denominator and an upper bound for the numerator.

To provide a \emph{lower} bound for the denominator in
\eqref*{eq:posteriordist}, we study $\omega_n(B_\D^h(f_0); \D,
h, \Y_n)$, where $B_\D^h(f_0)$ is a collection of density functions
supported by a non-stopped depth-$\D$ tree topology. To be more
specific, we grow a tree under the shift $h$ from the root $\T^\0 = \{\Omega^h \}$, and for
all intermediate trees $\T^\0, \T^\1, \ldots, \T^{(\D-1)}$ of depths $0, 1,
\ldots, \D-1$, we do not \emph{stop} any intermediate node. The
terminal depth-$\D$ tree topology is denoted as $\T^{(\D)}$, and the
density functions supported by this tree topology grown in this
specific way is denoted as $\mathcal{F}_{\treek}^h$. For any density
function $f$ and measurable set $A$ with respect to the Lebesgue
measure on $\Omega$, $\mu_f(A) = \int_A f d \mu$,
where $\mu$ is the Lebesgue measure. The set $B_\D^h$
is defined as 
\begin{eqnarray*}
 && B_{\D}^h(f_0) \coloneqq \Big \{ 
f \in \mathcal{F}_{\treek}^h: | \mu_f(A) - \mu_{f_0}(A) | \leq 
(\eta')^2 \D /n,   \\
&& \hspace{6cm}
\mbox{ and } 
\mu_f(A) \geq (\frac{1}{2} (\eta')^2 /n)^{\D} \mbox{ for any } A \in \T^{(\D)}.
\Big \}
\end{eqnarray*}
We examine prior probability mass on the set $ B_{\D}^h(f_0)$
under the OPT prior in Section
\ref{sec:optanalysis}, and derive a lower bound for the
likelihood ratio in Section~\ref{subsec:concentrationlower}. Usually,
$\mathcal{F}_{\treek}^h$ is chosen in a way that $\T^{(\D)}$ is the tree
topology that can provide the \emph{best} depth-$\D$ approximation to
$f_0$ under the KL divergence. These two
pieces of results lead to a lower bound for $\omega_n(B_\D^h(f_0); \D,
h, \Y_n)$, which is formally stated in Section~\ref{subsec:bounds}.

To obtain an \emph{upper} bound for the numerator in
\eqref*{eq:posteriordist}, we check the integral $\omega_n(\Theta_\D^h; \D,
h, \Y_n)$, where $\Theta_\D^h$ is the collection of density functions
supported by any tree topology under the shift $h$ and of depth
\emph{at most} $\D$. We provide an upper bound for this quantity
in Section \ref{subsec:bounds}.

With these lower and upper bounds, the outline for proving Theorem~\ref{thm:BMArateanisotropic} is to first show
that the marginal posterior distribution on tree depth concentrates fast around
$\D_\beta$ defined through equation~\eqref*{eq:depthbeta} \rev{ in Lemma~\ref{lemma:BMAdepth}}, and then
to carefully analyze bias and variance of the \writing{BDF} estimator when the tree depth is
around $\D_\beta$ \rev{in 
  Lemma~\ref{lemma:bmabiasvariance}}. \rev{We bound the bias
  component using the approximation results from
  Section~\ref{sec:approx}, and the variance component using the posterior
  contraction results under the OPT prior from
  Section~\ref{sec:optpost}.}

\subsection{Contraction of the Posterior Distribution on Tree Depth}
\label{subsec:bmadepth}

We first characterize contraction of the posterior distribution on tree
depth when $f_0$ lies in the anisotropic H\"older space $\mathcal{C}^{\bm \alpha}(\Omega,
L)$. A specific focus is to reveal how the contraction is
affected by the penalty on tree complexity (reflected by the hyper-parameter
$\beta$) and the smoothness property of $f_0$ (indexed by $\alpha_l$'s). We present the result in a slightly more
general setting: instead of requiring parameters of the anisotropic
H\"older space to lie in $(1, 2]$, we only assume that they are positive. 

\begin{lemma}
\label{lemma:BMAdepth}
  Assume $f_0 \in \mathcal{C}^{\bm \alpha} (\Omega, L)$ with
  $\alpha_l >0$ and $ \kappa_L \leq f_0 \leq \kappa_U$. Let
  $\alpha_l^* = \min \{1, \alpha_l\}$ and $\overline{\alpha^*}$ is the
  value that satisfies $(\overline{\alpha^*})^{-1} = \frac{1}{d}\sum_{l=1}^d
  \alpha_l^{-1}$. We also assume that there exist constants
  $\overline{c}_{\text{KL}}$ and $\underline{c}_{\text{KL}}$, such
  that the best depth-$\D$ approximation under any shift $h$, denoted
  as $f_\D^h$, 
  satisfies $ \underline{c}_{\text{KL}} 2^{-2
    \overline{\alpha^*} /d } \leq \kl (f_0, f_\D^h ) \leq \overline{c}_{\text{KL}} 2^{-2
    \overline{\alpha^*} /d }$. Define 
  \begin{eqnarray}
\label{eq:depthbetadef}
    \D_\beta^*  = \left \lfloor 
\log_2 ( \underline{c}_{\scriptsize{\kl}} n / \log n)^
{\frac{1}{\beta + 2 \overline{\alpha^*}/d}}
\right \rfloor.
  \end{eqnarray}
  Then there exist constants
  $\underline{\Delta}$ and $\overline{\Delta}$, such that with
  probability tending to 1 under $\P_{f_0}^n$, 
  \begin{eqnarray*}
    \Pi_n ( \D \notin
    [\D_\beta^* - \underline{\Delta}, \D_\beta^* + \overline{\Delta}]
    | \Y_n )
    \leq \exp \left(- \gamma_{\text{depth}} \cdot 2^{\D_\beta^*} \right),
  \end{eqnarray*}
  where the constants $\underline{\Delta}$ and $\overline{\Delta}$
  are specified in the proof. 
\end{lemma}
\begin{remark}
  In this lemma, we assume that the approximation error has both upper
  and lower bounds. Later in Lemma \ref{lemma:approxtree} we 
  provide an example where such an assumption is satisfied. 
\end{remark}

\begin{proof}
  The posterior distribution on tree depth can be calculated as 
\begin{eqnarray*}
  \Pi_n ( \D | \Y_n) = \frac{\sum_{h}
  \omega_n (\Theta_\D^h;\D, h, \Y_n)
  \Pi_n(\D, h) }
  {\sum_{(\D',h')}
  \omega_n(\Theta_{\D'}^{h'}; \D', h', \Y_n) \Pi_n(\D', h')}.
\end{eqnarray*}
Define $\omega_n (\D; \Y_n) = \sum_{h} \omega_n (\Theta_\D^h; \D, h, 
\Y_n) \Pi_n (\D, h) $.

For the numerator, by applying inequality \eqref*{eq:intunifupper} from Lemma \ref{lemma:inttheta}, with
probability tending to 1 under $\P_{f_0}^n$, the following upper  
bound holds for all $\D$,
\begin{eqnarray*}
  \omega_n (\D ;\Y_n)
\leq  \Pi_n (\D) \exp
\left( 
-n (1 - \overline{c}_1) \underline{c}_{\text{KL}} 2^{- 2 \overline{\alpha^*}\D/d}
+2\eta^2 \D \cdot 2^\D 
\right).
\end{eqnarray*}

While the denominator can be bounded by applying 
Corollary \ref{cor:intbound}.
By
plugging in $\D = \D_\beta^*$, 
when $\beta>1$,
\begin{eqnarray*}
&&  \omega_n ( \D_\beta^*; \Y_n )  \\
&& \hspace{0.2in} 
\geq  
\exp\left( 
- \left(\gamma \cdot \frac{1}{\beta +2 \overline{\alpha^*}/d}
+(1+ \underline{c}_1) 
\frac{\overline{c}_{\scriptsize{\kl}}}
{\underline{c}_{\scriptsize{\kl}}}
        +o(1) \right)
(\underline{c}_{\scriptsize{\kl}})
^{\frac{\beta}{\beta +2 \overline{\alpha^*}/d }}
n
^{\frac{\beta}{\beta +2 \overline{\alpha^*}/d }} 
(\log n)^{\frac{2 \overline{\alpha^*}/d}{\beta +2 \overline{\alpha^*}/d }}
\right).
\end{eqnarray*}
And when $\beta = 1$, the lower bound becomes
\begin{eqnarray*}
  \omega_n ( \D_\beta^*; \Y_n )  
&\geq&  
\exp\bigg( 
- \bigg(\frac{\gamma + \nu \log 2 + 8 (\eta')^2(1 + 3 \log (\kappa_U/ \kappa_L)) }{1+2 \overline{\alpha^*}/d}
+(1+ \underline{c}_1) 
\frac{\overline{c}_{\scriptsize{\kl}}}
{\underline{c}_{\scriptsize{\kl}}}
        +1 +o(1) \bigg) \\
&& \hspace{1cm} 
\cdot
(\underline{c}_{\scriptsize{\kl}})
^{\frac{1}{1+2 \overline{\alpha^*}/d }}
n
^{\frac{1}{1+2 \overline{\alpha^*}/d }} 
(\log n)^{\frac{2 \overline{\alpha^*}/d}{1 +2 \overline{\alpha^*}/d }}
\bigg).
\end{eqnarray*}

Then there exists two constants $\overline{\Delta}$ and
$\underline{\Delta}$ (may depending on the constants $\gamma, \beta,
\underline{c}_1, \overline{c}_1, \underline{c}_{\scriptsize{\kl}}$,
$\overline{c}_{\scriptsize{\kl}},\overline{\alpha^*} $ and $d$, and also $\eta$,
$\eta'$ and $\nu$ in the case $\beta=1$), such that for any $\D
< \D_\beta^* - \underline{\Delta}$,
\begin{eqnarray*}
  \Pi(\D | \Y_n) &\leq& \frac{\omega_n(\D; \Y_n)}
  {\omega_n (\D_\beta^*; \Y_n)}\\
&\leq & 
\exp\bigg( 
- \left( (1-\overline{c}_1) 
2^{2 (\underline{\Delta} +1)\overline{\alpha^*}/d } 
-\gamma \cdot 
\frac{1}{\beta +2 \overline{\alpha^*}/d}
-(1+ \underline{c}_1) 
\frac{\overline{c}_{\scriptsize{\kl}}}
{\underline{c}_{\scriptsize{\kl}}} 
+ o(1) \right)\\
&& \hspace{3cm}
\cdot (\underline{c}_{\scriptsize{\kl}})
^{\frac{\beta}{1+2 \overline{\alpha^*}/d }}
n
^{\frac{\beta}{1+2 \overline{\alpha^*}/d }} 
(\log n)^{\frac{2 \overline{\alpha^*}/d}{1 +2 \overline{\alpha^*}/d }}
\bigg).
\end{eqnarray*} 
The equation above shows how the bound of posterior probability mass
look like when $\beta>1$, when $\beta=1$, the bound also involves
$\eta$, $\eta'$, and $\nu$.   
Similarly, we can also show that the posterior probablity on any $\D \geq \D_\beta^*
+ \overline{\Delta}$ is exponetially small with probability
tending to 1 under $\P_{f_0}^n$. 
\end{proof}


\subsection{A Bias-Variance Trade-off  (Proofs of Theorem~\ref{thm:BMArateanisotropic} and Corollary~\ref{cor:BMArateanisotropic})}
\label{subsec:bmabv}
 In the following, we provide an explicit bias-variance
decomposition for the \writing{BDF estimator} \eqref*{eq:pseudoposteriormean}
defined through the hybrid aggregation scheme. We consider the case when $f_0 \in
\mathcal{C}^{\bm \alpha} ( \Omega, L)$ with all $\alpha_l >1$. In this
case, 
the tree depth $\D_\beta^*$ defined in \eqref*{eq:depthbetadef} is the
same as $\D_\beta$ defined in \eqref*{eq:depthbeta}. 

\begin{lemma}
  \label{lemma:bmabiasvariance}
  Assume $f_0 \in \mathcal{C}^{\bm \alpha} ( \Omega, L)$ with
  $\alpha_l >1$, and $ \kappa_L \leq f_0 \leq \kappa_U$. Both 
  Condition \ref{cond:boundary} and  Condition
  \ref{cond:uniqueness} hold. $\D_\beta$ is
  the tree depth defined in \eqref*{eq:depthbeta} and $\eta>1$ is
 a constant. When $\beta > 1$, with
  probability tending to 1 under $\P_{f_0}^n$, the following
  inequality holds for sufficiently large $n$, 
  \begin{eqnarray}
    \label{eq:bmabiasvariance}
     \rho(E^s(f |\Y_n), f_0 ) 
        \lesssim 
        \sum_{l=1}^d (L/\sqrt{\kappa_L}) 2^{-\alpha_l (\D_{\beta}- \underline{\Delta})/d}
        +M_{\text{tree}} \cdot \eta \delta(n, 2^{(\D_\beta+\overline{\Delta})}),
  \end{eqnarray}
where $M_{\text{tree}}$ is a constant specified by Corollary
\ref{cor:posteriorratetreemean}. 
\end{lemma}
\begin{proof}
 Without loss of generality, we may assume $\theta_l = \int  \frac{
      (\partial f_0/ \partial y^l)^2}  {f_0}> 0$ for $l = 1,
 \ldots, d$, as $\theta_l$ being zero only leads to \emph{faster}
 convergence. Fix some $\beta >1$, for a forest estimator with depth
 around $\D_{\beta}$, we are in the regime that the bias dominates the
 variance. To simplify the notation, we denote $\frac{1}{|\mathcal{H}_\D|} \sum_{h \in \mathcal{H}_\D} E(f | \Y_n, \D, h)$ as $E^s(f| \Y_n, \D)$. 

We first check at the individual tree level, which tree topology can
provide the best tree-supported approximation to $f_0$. Under the strong shape constraint, all
leaf nodes in a tree topology without perturbation are 
$d$-dimensional rectangles of the same shape. We denote the edge
length along dimension $l$ as $\epsilon_l$, and $\epsilon = \max_{1
  \leq l \leq d} \epsilon_l$. Our first observation is that, when all
$\theta_l >0$, the tree topology providing the best depth-$\D$
approximation to $f_0$ should be among those with edge length of
leaf nodes satisfying $\epsilon_l = O( 2^{-
\D_\beta /d})$. Among this collection of tree topologies, 
based on Lemma \ref{lemma:approxtree} and Remark
\ref{remark:approxtree}, for any $f_0 \in \mathcal{C}^{\bm \alpha} (\Omega, L)$
with $\alpha_l >1$, the approximation accuracy can be written as 
\begin{eqnarray*}
    \int f_0 \log (f_0/f_\A)
  =\frac{1}{24} \sum_{l=1}^d
  \theta_l \epsilon_l^2
      + O \left(\sum_{l=1}^d d L^2 \epsilon^{1+\alpha_l} / \kappa_L\right), 
      \quad \mbox{as }
      \epsilon \rightarrow 0.
\end{eqnarray*}
Then, for $\D$ large enough, finding the depth-$\D$ tree that provides
the best the approximation under the KL divergence is equivalent to solving the optimization
problem 
\begin{eqnarray*}
  \min_{k_1, \ldots, k_d}  &&\sum_{l=1}^d \theta_l 2^{-2 k_l}, \\
\text{subject to } && \sum_{l=1}^d k_l = \D, k_l \in \mathbb{N}. 
\end{eqnarray*}
Condition \ref{cond:uniqueness} ensures the solution to this
problem is unique. As tree topologies with or without perturbations
share the same form of expansion (see Remark \ref{remark:approxtree}), we can also 
conclude that if $\A^\star_\D$ is the tree topology \emph{without}
perturbation that provides the best approximation to $f_0$, then
$(\A^\star_\D)^h$ should be the best tree topology among all trees
grown under shift $h$.

For the depth-$\D$ forest estimator $E^s(f| \Y_n, \D)$, 
\begin{eqnarray*}
 &&   \rho\left(
    E^s(f| \Y_n, \D), 
    f_0\right)\\
 && \hspace{0.5in}   \leq \rho\left(E^s(f| \Y_n, \D), 
    \frac{1}{|\mathcal{H}_\D|} 
    \sum_{h \in \mathcal{H}_\D} f_{(\A_\D^\star)^h} \right)
    + \rho\left( \frac{1}{|\mathcal{H}_\D|} 
    \sum_{h \in \mathcal{H}_\D} f_{(\A_\D^\star)^h}, f_0 \right).
\end{eqnarray*}
The second term is pure bias. According to Corollary \ref{cor:biasforestbma},
\begin{eqnarray*}
    \rho\left( \frac{1}{|\mathcal{H}_\D|} 
    \sum_{h \in \mathcal{H}_\D} f_{(\A_\D^\star)^h, f_0} \right)
    \lesssim \sum_{l=1}^d( L/\sqrt{\kappa_L}) 2^{-\alpha_l \D /d}. 
\end{eqnarray*}
The first term corresponds to estimation error. When the Condition
\ref{cond:uniqueness} holds and $n$ is sufficiently large, according to
  Proposition \ref{thm:posteriorratetree}, with probability at least
  $1- 16 |\mathcal{H}_\D| \exp (-\frac{1}{2} \tilde{c}_2\eta^2 \delta^2(n,2^\D))$, the
  posterior distributions concentrate on tree topologies
  $(\A_\D^\star)^h$ simultaneously over all trees. The posterior
  probability mass
  on $(\A_\D^\star)^h$ is at least $1- \exp (-
  \frac{1}{2}\underline{c}_1 M_{\text{bias}}^2 \cdot \eta^2
  \delta^2(n,2^\D))$. Then, 
    \begin{eqnarray*}
      && \rho^2(E(f|\Y_n, \D, h), f_{(\A_\D^\star)^h})  
      \\ 
      && \hspace{0.5in} \leq
     E(\rho^2(f, f_{(\A_\D^\star)^h}) | \Y_n, (\A_\D^\star)^h)
       \Pi(\mathcal{F}_{(\A_\D^\star)^h} | \Y_n, \D, h)
       + 2 \Pi(\mathcal{F}_{(\A_\D^\star)^h}^\complement | \Y_n, \D, h).
    \end{eqnarray*}
    The first term on the right hand side can be bounded by applying
    Corollary \ref{cor:posteriorratetreemean}, and the second term is
    bounded by the concentration result. Combining these bounds
    together, we have that uniformly for all trees with in the forest, with
    probability at least $1- 16 |\mathcal{H}_\D| \exp (-\frac{1}{2} \tilde{c}_2\eta^2 \delta^2(n,2^\D))$,
    \begin{eqnarray*}
         \rho^2(E(f|\Y_n, \D, h), f_{(\A_\D^\star)^h}) 
         \leq  (M_{\text{tree}} \cdot \eta)^2 \delta^2(n,2^\D).
    \end{eqnarray*} 
Therefore, with
  probability tending to 1 under $\mathbb{P}_{f_0}^n$, for
  sufficiently large $n$,
\begin{eqnarray*}
   && \rho^2\left(\frac{1}{|\mathcal{H}_\D|} 
    \sum_{h \in \mathcal{H}_\D} E(f | \Y_n, \D, h), 
    \frac{1}{|\mathcal{H}_\D|} 
    \sum_{h \in \mathcal{H}_\D} f_{(\A_\D^\star)^h} \right)\\
   && \hspace{0.5cm} \leq
   \frac{1}{|\mathcal{H}_\D|} 
   \sum_{h \in \mathcal{H}_\D}
   \rho^2\left(
     E(f | \Y_n, \D, h),  f_{(\A_\D^\star)^h} \right) 
     \leq  (M_{\text{tree}} \cdot \eta)^2 \delta^2(n, 2^{\D}). 
\end{eqnarray*}

The \writing{BDF} estimator $\eqref*{eq:pseudoposteriormean}$ takes another layer of averaging with respect to
the depth. Note that the posterior on depth concentrates on $[\D_{\beta}-
\underline{\Delta}, \D_{\beta}+ \overline{\Delta}]$ exponentially
fast as shown in Lemma \ref{lemma:BMAdepth}, and for any $\D$ within $[\D_{\beta}-
\underline{\Delta}, \D_{\beta}+ \overline{\Delta}]$ the same upper bound
for $\rho\left( E^s(f| \Y_n, \D), f_0\right)$ holds. This finishes the proof. 
\end{proof}

Theorem~\ref{thm:BMArateanisotropic} is a result of the bias-variance
trade-off.
\begin{proof}[Proof of Theorem~\ref{thm:BMArateanisotropic}]
We have provided a bias-variance decomposition for the \writing{BDF} estimator in equation \eqref*{eq:bmabiasvariance}. 
Let $\alpha_{\min} = \min_{1 \leq l \leq
  d} \alpha_l$. After some calculation, we can show that when $1 < \beta \leq 1 + 2(\alpha_{\min} -1) /d$,
the estimation error dominates and the convergence rate is $(dL/\sqrt{\kappa_L})^{\frac{1}{\beta+2/d}}n^{-\frac{2/d
    +\beta-1 }{\beta+2/d}}$; when $\beta < 1 + 2(\alpha_{\min} -1)
/d$, the bias is dominating and the rate is
$(dL/\sqrt{\kappa_L})^{\frac{\beta+2/d - \alpha_{\min}/d}{\beta+2/d}}
  n^{-\frac{\alpha_{\min}/d}{\beta + 2 /d}}$. When $\beta < \alpha_{\min} + 2(\alpha_{\min} -1)/d$, the rate of the \writing{BDF} estimator is faster
  than that based on a single tree model.
The proof above is mainly for the case when $\beta >1$. When $\beta=1$, the proof is still based on checking the bias term and the
variance term and is skipped. 
\end{proof}

If we remove the boundary condition, we can show
Corollary~\ref{cor:BMArateanisotropic} using a similar argument, only
modifying the bias term in the bound
\eqref*{eq:bmabiasvariance}.
\begin{proof}[Proof of  Corollary~\ref{cor:BMArateanisotropic}]
  \rev{When the boundary condition is not satisfied, the approximation
 error for forests may become larger due to the boundary effect. In
 particular, the approximation error can be bounded by using
 Corollary~\ref{cor:biasforestbma}. With a similar argument as that
 for Lemma~\ref{lemma:bmabiasvariance}, we can obtain the following
 bias-variance trade-off result without the boundary condition, 
  \begin{eqnarray*}
     \rho(E^s(f |\Y_n), f_0 ) 
        \lesssim 
        \sum_{l=1}^d (d^{3/2} L/\sqrt{\kappa_L}) 2^{-\min \{3/2,
    \alpha_l \} (\D_{\beta}- \underline{\Delta})/d}
        +M_{\text{tree}} \cdot \eta \delta(n, 2^{(\D_\beta+\overline{\Delta})}),
  \end{eqnarray*}
where $M_{\text{tree}}$ is a constant specified by
Corollary~\ref{cor:posteriorratetreemean}. Intuitively, we may view
$\alpha^* = \min \{\frac{3}{2}, \alpha \}$ as the effective smoothness
along each direction. With the argument for
Theorem~\ref{thm:BMArateanisotropic}, we can obtain the desired
result. }
\end{proof}

\subsection{A Lower Bound for the Posterior Contraction Rate of the Single-Tree Method (Proof of Theorem~\ref{thm:treelowerbound})}
\label{subsec:treelowerbound}
\rev{In this subsection, we provide the proof of Theorem~\ref{thm:treelowerbound}, which states a lower bound on the posterior contraction rate under a variation of the OPT prior. The prior distribution can be viewed as a special case of the forest prior introduced in Section~\ref{sec:BMA}, with the number of trees within the forest being one. This lower bound reveals the limit of the single-tree method. 
 Our proof strategy in inspired by that for Theorem 2 in the paper by \cite{castillo2021uncertainty}. }

\rev{Although this lower bound is established for the single-tree model rather than the BDF estimator, we present the proof in this section because part of the argument relies on analyzing the posterior distribution on tree depth, using an approach similar to that in Lemma~\ref{lemma:BMAdepth}.}



\begin{proof}[Proof of Theorem~\ref{thm:treelowerbound}]
  \rev{
    Our proof strategy is to find a sequence of $f_0$'s lying in the space $\mathcal{C}^{\bm \alpha}(\Omega, L)$ that satisfies Condition~\ref{cond:boundary} and is bounded away from zero, denoted as $f_0^n$, such that $\mathbb{E}_{f_0^n}^n \Pi_n \left( \{ f: \rho(f_0, f) < \underline{\eta}  \epsilon_n  \cond \Y_n \right) \rightarrow 0 $. In fact, we choose a same $f_0$ (not $n$-dependent), and show that  $\Pi_n \left( \{ f: \rho(f_0, f) < \underline{\eta}  \epsilon_n  \cond \Y_n \right) \rightarrow 0$ in $\mathbb{P}_{f_0}^n$-probability.
    }

    \rev{
      We first discuss the choice of $f_0$. To ensure that the space $\mathcal{C}^{\bm \alpha}(\Omega, L)$ contains a density function defined on $\Omega$, we assume that $L >1$. In the following, we first introduce a marginal density function $g_0$ on $[0,1]$, and define $f_0 (y) $ as $\prod_{l=1}^d g_0(y^l)$. We set $g_0 (x) = 1 - \phi \cos (2 \pi x)$ for some $\phi <1$.  We can adjust the value of $\phi$ such that $f_0 \leq L$ and the partial derivatives are Lipschitz continuous. We can also choose an appropriate $\phi$ such that the $g_0$ is bounded away from zero. For simplicity, we assume that $g_0 \geq 1/2$. Consequently, $f_0 \geq (1/2)^d$. Therefore, $f_0$ lies in the space $\mathcal{C}^{\bm \alpha}(\Omega, L)$. In addition, $f_0$ satisfies Condition \ref{cond:boundary}.
      }

      \rev{
    We denote the best depth-$\D$ strong-shape-constrained tree-supported approximation to $f_0$ under the Hellinger distance as $f_\D^{\text{H}}$. To study the approximation error, our first observation is that, with $\D' = d \lceil \D/d \rceil $,  $\rho(f_0 , f_\D^{\text{H}}) \geq \rho(f_0, f_{\D'}^{\text{H}}) $. This is because when we refine a partition, the approximation error alway decreases. For this distribution, based on the symmetry of the density function, the best depth-$\D'$ strong-shape-constrained single-tree approximation under the Hellinger distance is provided by the tree partition with splits equally spaced along different directions. This is to say, the best depth-$\D'$ approximation is provided by a partition $\A = \{ \Omega_j \}_{j=1}^t$ with the edge length of terminal regions along dimension $l$, denoted as $e_l$,  equal to $2^{- \lceil \D/d \rceil }$. By Lemma~\ref{lemma:bestHellinger}, the best depth-$\D'$ approximation $f_{\D'}^{\text{H}} = \sum_{j=1}^t \frac{w_j}{\mu(\Omega_j)} \mathbbm{1}_{\Omega_j}$ should satisfy
  \begin{eqnarray*}
     w_j =\frac{ u_j^2} {
   \lambda^2 \mu(\Omega_j)}
  \quad
  \mbox{with}
    \quad
    u_j = \int_{\Omega_j} \sqrt{f_0} d \mu
    \quad 
    \mbox{and}
    \quad
  \lambda^2 =\sum_{j=1}^t \frac{u_j^2}{ \mu( \Omega_j)}.
  \end{eqnarray*}
  }

  \rev{
  Next, we obtain a lower bound for the approximation error using the approximation result under the KL divergence provided in Lemma~\ref{lemma:approxtree}. Based on the form of the best approximation $f_{\D'}^{\text{H}}$ presented above, we can conclude that $f_{\D'}^{\text{H}}$ is bounded from below and above. Then, based on Lemma~\ref{lemma:kltoH}, there exists a constant $\eta^{\text{H}}$ that only depend on $f_0$, such that $\rho^2(f_0, f_{\D'}^{\text{H}}) \geq  \eta^{\text{H}} \kl (f_0, f_{\D'}^{\text{H}}) $. For the approximation error under the KL divergence, we can calculate the quantity $\theta_l$ defined in Section~\ref{subsec:hybridbma} and apply Lemma~\ref{lemma:approxtree} from Section~\ref{sec:approx} to derive a lower bound. For $f_0$ defined above, to calculate $\theta_l$'s, we only need to calculate $\theta_1$ given the symmetry. 
  \begin{eqnarray*}
    \frac{\partial{f_0}}{\partial{y^1}} =\frac{ \partial g_0(y^1)}{\partial y^1} \prod_{l=2}^d g_0(y^l).
  \end{eqnarray*}
  And
  \begin{eqnarray*}
    \theta_1 &=&\int_{\Omega} \frac{ \left( \partial f_0 / \partial y^1\right)^2 }{ f_0 } d\mu
    = \int_{\Omega} \frac{ \left( \partial g_0 (y^1)/ \partial y^1\right)^2 }{ g_0 (y^1)} \prod_{l=2}^d g_0(y^l) d\mu (y) \\
             &=& \int_0^1 \frac{ \left( \partial g_0 (y^1)/ \partial y^1\right)^2 }{ g_0 (y^1)} d y^1 \\
    & \geq& \int_0^1 \frac{ \left( 2 \pi \phi \sin(2 \pi x) \right)^2 }{ 2} d x = \pi^2 \phi^2,
  \end{eqnarray*}
  where the inequality is based on the fact that $g_0 \leq 2$.
Note that $f_{\D'}^{\text{H}}$ is supported by a balanced partition with edge lengths of terminal nodes along different dimensions all equal to $e =2^{- \lceil \D /d \rceil}$. For such a depth-$\D'$ partition $\A$, an application of Lemma~\ref{lemma:approxtree} implies that
  \begin{eqnarray*}
    \kl (f_0, f_{\D'}^{\text{H}}) \geq \frac{1}{24} \sum_{l=1}^d \theta_l e^2 + O \left(\sum_{l=1}^d d L^2 e^3 /(1/2)^d \right), \quad \mbox{as } \D \rightarrow \infty. 
  \end{eqnarray*}
 For $\D$ large enough, 
  \begin{eqnarray*}
    \kl(f_0, f_{\D'}^{\text{H}}) \geq \frac{\pi^2 \phi^2 d}{48} 2^{- 2 \lceil \D/d \rceil}.
  \end{eqnarray*}
  Therefore, for depth-$\D$ single-tree-supported density $f$ under the strong shape constraint, 
  \begin{eqnarray*}
    \rho^2(f_0, f) \geq \rho^2(f_0, f_\D^{\text{H}}) \geq
    \rho^2(f_0, f_{\D'}^{\text{H}})
    \geq
    \eta^{\text{H}} \kl(f_0, f_{\D'}^{\text{H}})
    \geq \frac{ \eta^{\text{H}} \pi^2 \phi^2 d}{48} 2^{- 2 \lceil \D/d \rceil}.
  \end{eqnarray*}
  The lower bound implies that, for $\D$ (and $n$ equivalently) large enough, only when $\D \geq \\ \lfloor \log_2 \left( \frac{\eta^{\text{H}} \pi^2 \phi^2 d }{196 \underline{\eta} ^2}\right)^{d/2} \left( \frac{n}{\log n}\right)^{\frac{1}{\beta_{\text{tree}} + 2/d}} \rfloor$, it is possible to find a tree-supported density $f$, such that $\rho (f_0, f ) < \underline{\eta} \epsilon_n$. For simplicity, we denote $\lfloor \log_2 \left( \frac{\eta^{\text{H}} \pi^2 \phi^2 d }{196 \underline{\eta} ^2}\right)^{d/2} \left( \frac{n}{\log n}\right)^{\frac{1}{\beta_{\text{tree}} + 2/d}} \rfloor$ as $\D_{\underline{\eta}}$.
  }

  \rev{
  Next, we show that we can choose some $ \underline{\eta}$, such that the posterior probability on depth-$\D$ trees with $\D \geq \D_{\underline{\eta}}$ decays to zero, which implies the lower bound result of the form \eqref*{eq:lowerbound}. The proof is very similar to that for Lemma~\ref{lemma:BMAdepth}. In particular, with notations introduced at the beginning of Section~\ref{sec:proofbma}, the posterior probability on depth-$\D$ trees can be calculated as
  \begin{eqnarray*}
    \Pi_n( \D | \Y_n) = \frac{
  \omega_n (\Theta_\D;\D, \Y_n)
  \Pi_n(\D) }
  {\sum_{\D'}
  \omega_n(\Theta_{\D'}; \D', \Y_n) \Pi_n(\D')}.
  \end{eqnarray*}
  Similar to the proof for Lemma~\ref{lemma:BMAdepth}, we can lower bound the denominator by applying Lemma~\ref{lemma:intbound} and upper bound the numerator by applying Lemma~\ref{lemma:inttheta}. With $\D_{\beta_{\text{tree}}}^*  = \left \lfloor 
\log_2 ( \underline{c}_{\scriptsize{\kl}} n / \log n)^
{\frac{1}{\beta_{\text{tree}} + 2 /d}}
\right \rfloor$, using a same argument, we can show that there exist $\underline{\Delta}$ and $\overline{\Delta}$, such that with
  probability tending to 1 under $\P_{f_0}^n$, 
  \begin{eqnarray*}
    \Pi_n ( \D \notin
    [\D_{\beta_{\text{tree}}}^* - \underline{\Delta}, \D_{\beta_{\text{tree}}}^* + \overline{\Delta}]
    | \Y_n )
    \leq \exp \left(- \gamma_{\text{depth}} \cdot 2^{\D_{\beta_{\text{tree}}}^*} \right).
  \end{eqnarray*}
 When $\beta_{\text{tree}} >1$, the choice of $\overline{\Delta}$ and
$\underline{\Delta}$ may depend on the constants $\gamma, \beta_{\text{tree}},
\underline{c}_1, \overline{c}_1, \underline{c}_{\scriptsize{\kl}}$,
$\overline{c}_{\scriptsize{\kl}} $ and $d$.  When $\beta_{\text{tree}} =1$, the choice of $\overline{\Delta}$ and
$\underline{\Delta}$ may additionally rely on $\eta$,
$\eta'$ and $\nu$. These contants are specified in the proof of Lemma~\ref{lemma:BMAdepth}. Therefore, we may make $\underline{\eta}$ smaller, such that $\D_{\underline{\eta}} > \D_{\beta_{\text{tree}}}^* + \overline{\Delta}$. This leads to the desired result.
}
\end{proof}



\section{Proofs for the Results on GDF in Section \ref{sec:gibbs}}
\label{sec:proofgibbs}
In this section, we provide detailed proofs for the results in Section \ref{sec:gibbs}. 
The posterior contraction result for the Gibbs posterior can be established 
by studying the risk function \eqref*{eq:riskdef}, 
prior thickness, and tail of the prior.

\rev{As presented in \eqref{eq:gibbsposterior}, the Gibbs posterior on forest
topologies can be
written as 
\begin{eqnarray*}
  \lambda^\g  = e^{-n \psi R_n} q \bigg /
  \bigg( 
\sum_{\newforest}
  e^{-n \psi R_n \newbforest}
  q \newbforest \bigg).
\end{eqnarray*}
To analyze the posterior contraction rates and prove
Theorem~\ref{thm:gibbsrate}, we need to derive an upper bound
for 
the numerator and a 
lower bound for the denominator, respectively, in order to show that the
Gibbs posterior on the set $\mathcal{C}_{\text{forest}}$ decays to zero.}

\rev{The upper bound of the numerator is derived mainly based on the
  lower bound of the risk function, shown in
  Section~\ref{subsec:proofrisk}, and the decay of the prior assumed
  by Condition~\ref{cond:prior}. Condition~\ref{cond:complexity}
  ensures that risk lower bounds hold uniformly for all forests.  }

\rev{To establish the lower bound of the denominator, when Condition \ref{cond:approx} 
holds, the following forest depth plays an important role,  
 \begin{eqnarray}
\label{eq:optimalforestdepth}
    \D^\star = \left \lfloor 
\log_2 
\left( (c_\appr)^2 n/ \log n \right)^{\frac{1}{2r+1}} 
\right \rfloor.
  \end{eqnarray}
Among all depth-$\D^\star$ forest topologies, we use $\newstarf$ to
denote the one that can provide the best approximation to $f_0$ in the
sense that $\kl(f_0, f_{\newstarf})$ is the smallest, where
$f_\newstarf$ is the forest-supported approximation defined in
\eqref*{eq:forestapprox}.} \rev{With the forest-supported approximation
  $f_{\newstarf}$ and the associated approximation error $\kl(f_0,
  f_{\newstarf})$, we can show an upper bound
  of the risk function associated with $\newstarf$, which, in
  combination with Condition~\ref{cond:prior}, 
  implies the lower bound of the denominator. We
  present the risk upper bound in Section~\ref{subsec:riskupper}.}

\rev{The proof of Theorem~\ref{thm:gibbsrate} proceeds in two steps:
  Lemma~\ref{lemma:shallowforesets} shows that the Gibbs posterior on
  ``deep'' forest topologies decays exponentially fast, and 
  Lemma~\ref{lemma:shapeselection} shows that, among ``shallow''
  forests, the Gibbs posterior assigns negligible weight to those with poor approximation ability.}

\rev{The proofs of Theorem~\ref{thm:gibbsholderrate} and
  Corollary~\ref{cor:giibsholderrate} are finished by checking the
  three conditions: Condition~\ref{cond:approx},
  Condition~\ref{cond:prior}, and Condition~\ref{cond:complexity}. In
  particular, we check Condition~\ref{cond:prior} and
  Condition~\ref{cond:complexity} in the proof, and verify Condition~\ref{cond:approx} using the
  approximation results presented in Section~\ref{sec:approx}. }


\subsection{Gibbs Posterior on Forest Depth}
To study the Gibbs posterior on forest depth, we use the risk bounds
derived in Section \ref{sec:risk}. The posterior concentration of
forest depth exactly reflects a bias-variance trade-off. We show that
under the Gibbs posterior, the variance is under control in the sense
that forests cannot be too deep.
  \begin{lemma}
    \label{lemma:shallowforesets}
Assume Conditions 
\ref{cond:approx}, \ref{cond:prior} and \ref{cond:complexity} 
hold. \writing{In addition, we assume $\kappa_L \leq f_0 \leq \kappa_U$ for some constants $\kappa_U \geq \kappa_L >0$.} $\D^\star$ is the forest depth defined in
\eqref*{eq:optimalforestdepth}. 
Let $\{\zeta_n\}_{n \geq 1}$ be a sequence that increases at most as
$\log \log n$.
For the Gibbs posterior \eqref*{eq:gibbsposterior} on forest topologies defined in Section \ref{sec:gibbs}, with probability tending to 1 under $\P_{f_0}^n$, 
  \begin{eqnarray*}
   \lambda^\g 
\left(\D \geq \D^\star 
    + \lceil \zeta_n \rceil 
\right) 
\lesssim\exp \left(- \frac{1}{2}(\underline{\gamma}-1) \D^\star \cdot
    2^{\D^\star} \right), 
  \end{eqnarray*}
where the choice of $\{ \zeta_n \}_{n \geq 1}$ and \rev{the inverse
  temperature parameter $\psi$} is specified in the proof. 
\end{lemma}
\begin{proof}
The Gibbs posterior on depth-$\D$ forest topologies can be calculated as
  \begin{eqnarray*}
    \lambda^\g ( \D) &=& 
\frac
    {
    \sum_{\newforest \text{ of depth } \D}
    \hspace{1mm} 
    e^{- n \psi \wt{R}_n \newbforest}
    q \newbforest }
    {\sum_{\D'=0}^{\D_{\max}} 
  \sum_{\newforestp \text{ of depth } \D'}
\hspace{1mm} 
    e^{- n \psi \wt{R}_n \newbforestp} 
    q \newbforestp} .
  \end{eqnarray*}
In the following, we provide bounds for the denominator and the numerator, respectively. 

To derive the lower bound for the denominator, we study the forest estimate $\tilde{f}_\en^\star$ supported by $\newstarf$. 
    \begin{eqnarray}
    \label{eq:gibbsthickness}
\hspace{0.2in} \sum_{\D'=0}^{\D_{\max}} 
  \sum_{\newforestp \text{ of depth } \D'}
\hspace{1mm} 
    e^{- n \psi \wt{R}_n \newbforestp} 
    q \newbforest
      \geq
       e^{- n \psi \wt{R}_n \newbstarf}
         q \newbstarf.
  \end{eqnarray}
  By applying the upper bound of the centered risk function $\wt{R}_n$ derived
  in Proposition \ref{prop:gibbsthickness}, we have with probability
  at least $1- \exp(-c_2^\star (\eta_1^\star)^2 n_1
    \delta^2(n_1, t^\star)) - 24 \exp (-2 c_2 \eta_0^2 t_{\min} \log
    t_{\min})$, 
  \begin{eqnarray*}
    \eqref*{eq:gibbsthickness} \geq 
    \exp \left( 
- \psi n (1 + c_1^\star)\kl (f_0, f_{\newstarf}) 
- c_{\text{den}}\psi \xi_n^\star \D^\star \cdot 2^{\D^\star}
   \right)
    \cdot
    q \newbstarf,
  \end{eqnarray*}
where \tech{$c_{\text{den}} \propto (c_1^\star)^{-1} M_\tr^2 \eta_0^2 \log 2
  +  (\eta_1^\star)^2 \log 2$}. In combination with the lower
bound on $q \newbstarf$ assumed by Condition \ref{cond:prior}, with the same
probability, 
\begin{eqnarray*}
  \eqref*{eq:gibbsthickness} &\gtrsim&
  \exp \left( - \psi n (1 + c_1^\star)\kl (f_0, f_{\newstarf}) 
- c_{\text{den}}\psi \xi_n^\star \D^\star \cdot 2^{\D^\star}
\right) \cdot 
\exp \left(- \overline{\gamma} \D^\star \cdot 2^{\D^\star} \right).
\end{eqnarray*}
For the bias term, under Condition \ref{cond:approx}, 
\begin{eqnarray*}
  n \kl (f_0,f_{\newstarf})
 & \leq&
  n \tilde{c}_{\bias} c_{\text{approx}}^2
  2^{-( \log_2 ((c_{\text{approx}})^2 n/\log n)^{1/(2r+1)}-1)2 r}\\
  &\leq& 
         2^{2 r }
         \tilde{c}_{\bias} (c_{\text{approx}})^{2/(2r+1)} (\log n) 2^{\D^\star+1},
\end{eqnarray*}
where $\tilde{c}_{\bias}$ is the constant for bounding the KL
divergence by the Hellinger distance for bounded densities (see Lemma~\ref{lemma:suptoH}). 
Therefore, \tech{for some constants $c_{\bias}$},
\begin{eqnarray*}
  \eqref*{eq:gibbsthickness} \gtrsim 
  \exp\left (- (\psi( c_{\bias} (c_{\text{approx}})^{2/(2r+1)}+
  c_{\text{den}} \xi_n^\star ) + \overline{\gamma} ) \D^\star \cdot 2^{\D^\star} \right).
\end{eqnarray*}

For the numerator, when Condition~\ref{cond:complexity} holds, by Theorem~\ref{prop:risklower} and Remark~\ref{remark:uniformrisklower}, we can show that with probability as
least
{\small
$$1- 9 \exp(-(c_2/2) \eta_1^2 t_{\min} \log t_{\min}) 
    - 3\exp(-(c_2'/2) \eta_1^2 t_{\min} \log t_{\min})
    - 24 \exp(-2c_2\eta_0^2 t_{\min} \log t_{\min}),$$ }
for any forest topology of depth $\D > \D^\star$, 
\begin{eqnarray*}
 e^{-n \psi \wt{R}_n \newbforest}
  \leq
  \exp \left(
  -n \psi (1 - c_1-c_3) \kl ( f_0, f_\newforest )
  + \psi c_{\text{num}} \xi_n \D \cdot 2^\D
  \right),
\end{eqnarray*}
where \tech{$c_1$ and $c_3$ are constants that can be chosen in $(0,1)$}, 
$$ \tech{c_{\text{num}} \propto
(c_{\text{opt}}/c_3)(\kappa_U/ \kappa_L)\eta_0^2 M_\tr^2 \log 2
+ 2\eta_1^2 \log 2},$$ 
and $\xi_n$ is of the form $c_\xi (\kappa_U /\kappa_L) \log n$ for some $c_\xi >0$. 
Then for any depth $\D > \D^\star$, 
\begin{eqnarray}
  \label{eq:forestnumerator}
 \nonumber \sum_{\newforest \text{ of depth } \D}
  e^{-n \psi \wt{R}_n \newbforest}
  q \newbforest 
  &\leq&
     \exp \left(
  \psi  c_{\text{num}} \xi_n \D \cdot 2^\D
     \right)
    \sum_{\newforest \text{ of depth } \D}
     q \newbforest\\
  &\lesssim & 
      \exp \left(
  \psi  c_{\text{num}} \xi_n \D \cdot 2^\D
     \right)
   \exp (- \underline{\gamma} \D \cdot 2^\D),
\end{eqnarray}
where the last inequality is by Condition \ref{cond:prior}. \tech{When
  $\psi \leq (  c_\xi  c_{\text{num}}( \kappa_U / \kappa_L)  \log
  n)^{-1}$, we may take 
  \begin{eqnarray*}
    \zeta_n = c_{\zeta}+\log_2 \frac{\psi \left( c_{\bias} (c_{\text{approx}})^{2/(2r+1)}+
  c_{\text{den}} \xi_n^\star \right) + \overline{\gamma}
    }{\underline{\gamma}-1},
    \quad
    \mbox{when }
    \psi \leq  \
    \frac{1}{ c_\xi  c_{\text{num}}( \kappa_U / \kappa_L)  \log n},
  \end{eqnarray*}
} where $c_{\zeta} \geq 1$ is a constant. Then with probability tending to 1 under $\P_{f_0}^n$, for any $\D \geq \D^\star +
\lceil \zeta_n \rceil$, 
\begin{eqnarray*}
  \sum_{\newforest \text{ of depth } \D}
\lambda^\g \newbforest
  \lesssim 
\exp\left(- \tech{\frac{1}{2}(\underline{\gamma}-1)}
(\D -\lceil \zeta_n \rceil )\cdot 2^{\D -\lceil \zeta_n \rceil} \right).
\end{eqnarray*}
An application of Lemma \ref{lemma:tooltail} finishes the proof. 
\end{proof}


\subsection{Posterior Contraction on Promising Forest Topologies
  (Proof of Theorem~\ref{thm:gibbsrate})}
In the previous subsection, we have shown that the Gibbs posterior on deep forest
topologies decays to zero exponentially fast. In this subsection,
we show that among the forest topologies of appropriate depth,
the Gibbs posterior concentrates on those that better fit the data,
where a ``better fit'' is in the
sense that the Hellinger distance between the approximations supported
by these forest topologies and the true density is small. Before we
proceed, we first introduce the following subset of forest topologies:
\begin{eqnarray}
  \label{eq:cforestdepthD}
  \C_\en' \coloneqq \left\{ \newforest: \D < \D^\star + \lceil \zeta_n
  \rceil \mbox{ and } \rho(f_0, f_{\newforest}) \geq M_n
  (c_{\text{approx}})^{\frac{1}{2r+1}}\left( n / \log n \right)^{-\frac{r}{2 r+1}}\right\}, 
\end{eqnarray}
where $\{\zeta_n\}_{n \geq 1} $ is introduced in Lemma
\ref{lemma:shallowforesets}, $f_{\newforest}$ is as defined in
\eqref*{eq:forestapprox}, and $\{M_n\}_{n \geq 1}$
is the same as that in the definition of $\C_\en$ in Section \ref{subsec:gibbsrate}. We have the
following posterior concentration result for the Gibbs posterior: 

\begin{lemma}
\label{lemma:shapeselection}
Under the same assumptions as those for Theorem \ref{thm:gibbsrate}, 
  \tech{with probability tending to 1 under 
    $\P_{f_0}^n$},
  \begin{eqnarray*}
\lambda^\g  (\C_\en') 
\lesssim
    \tech{ \exp \left(- \frac{1}{2} (\underline{\gamma} -1)
\D^\star \cdot 2^{\D^\star} \right) }.
  \end{eqnarray*}
\end{lemma}
\begin{proof}
  The posterior probability on $\C_\en'$ can be calculated as
  \begin{eqnarray*}
\lambda^\g (\C_\en' )
    =
    \frac{\sum_{\newforest \in \C_\en' }
    e^{-n \psi \wt{R}_n \newbforest}
    q \newbforest }
    {\sum_{\D'=0}^{\D_{\max}}
\sum_{\newforestp \mbox{\scriptsize{ of depth $\D'$}}}
    e^{-n \psi \wt{R}_n \newbforestp}
    q \newbforestp
    }.
  \end{eqnarray*}
  Similar to the proof for Lemma \ref{lemma:shallowforesets}, we can
  show that under $\P_{f_0}^n$ with probability at least 
  $$1- \exp(-c_2^\star (\eta_1^\star)^2 n_1\delta^2(n_1, t^\star)) 
  - 24 \exp (-2 c_2 \eta_0^2 t_{\min} \log t_{\min}),$$ 
  the denominator is lower bounded by 
$$\tech{ \exp\left (- (\psi( c_{\bias} (c_{\text{approx}})^{2/(2r+1)}+
  c_{\text{den}} \xi_n^\star ) + \overline{\gamma} ) \D^\star \cdot
  2^{\D^\star} \right)},$$ 
up to some constant. 

 For the numerator, again similar to the proof of Lemma \ref{lemma:shallowforesets}, by Theorem \ref{prop:risklower} and Remark
\ref{remark:uniformrisklower}, under $\P_{f_0}^n$ with probability at least 
{\small
\begin{eqnarray*}
  1 -9 \exp(- (c_2/2) \eta_1^2 t_{\min} \log t _{\min}) - 3 \exp
   (-(c_2'/2) \eta_1^2 t_{\min} \log t_{\min}) - 24 \exp (-2 c_2 \eta_0^2 t_{\min} \log
    t_{\min}), 
\end{eqnarray*}
}
  we have the following upper
 bound of the numerator by using the fact that for any $\newforest \in
 \C_\en'$, $n \rho^2(f_0, f_\en ) \geq M_n^2 (\log n) 2^{\D^\star} $,
 \begin{eqnarray*}
  &&\sum_{\newforest \in \C_\en' }
    e^{-n \psi \wt{R}_n \newbforest}
     q \newbforest \\
   && \hspace{1cm} \leq
     \exp \left(
  - \psi (1 - c_1-c_3) M_n^2 (\log n) 2^{\D^\star} 
  + \psi c_{\text{num}} \xi_n \D \cdot 2^\D
  \right) 
\times
      \sum_{\newforest \in \C_\en' } q \newbforest,
 \end{eqnarray*}
where $c_1$, $c_3$ and $c_{\text{num}}$ are the same as those in
  the proof of Lemma \ref{lemma:shallowforesets} and $\xi_n$ is of the form $c_\xi (\kappa_U /\kappa_L) \log n$.
  When the inverse temperature parameter $\psi$ satisfies
  \begin{eqnarray*}
    \tech{\psi \leq \frac{1}{c_{\text{num}} \xi_n}
    =
     \frac{1}{ c_\xi  c_{\text{num}}( \kappa_U / \kappa_L)  \log n} },  
  \end{eqnarray*}
  we may allow $M_n$ increases as
polynomials of $\log n$ ($M_n^2$ is at least of the order $c_{\text{den}}\xi_n^\star$), such that the upper bound for the numerator
decays faster than the lower bound for the denominator. This finishes
the proof. 
\end{proof}

Theorem~\ref{thm:gibbsrate} is a direct result of Lemma~\ref{lemma:shallowforesets} and Lemma~\ref{lemma:shapeselection}: the
Gibbs posterior probability on deep forest topologies is exponentially
small by Lemma~\ref{lemma:shallowforesets}; and among shallow forests the Gibbs posterior on those with
poor approximation ability is also exponentially small by Lemma~\ref{lemma:shapeselection}.
\begin{proof}[Proof of Theorem~\ref{thm:gibbsrate}]
  \rev{
    With the sequence $\{ \zeta_n \}_{n \geq 1}$ specified in the
    proof of Lemma~\ref{lemma:shallowforesets}  and the set $\mathcal{C}'_\en$ defined in
    \eqref*{eq:cforestdepthD}, we have 
    \begin{eqnarray*}
      \mathcal{C}_\en \subset \{ \newforest:
      \newforest \mbox{ is of depth } \D \geq \D^\star + \lceil \zeta_n
      \rceil \} \cup \mathcal{C}'_\en. 
    \end{eqnarray*}
    Therefore,
    \begin{eqnarray*}
      \lambda^\g (\mathcal{C}_\en) \leq 
     \lambda^\g \left( \D \geq \D^\star + \lceil \zeta_n \rceil
      \right)
      +\lambda^\g\left(\mathcal{C}'_\en \right).
    \end{eqnarray*}
    Based on Lemma~\ref{lemma:shallowforesets},
    \begin{eqnarray*}
       \lambda^\g \left( \D \geq \D^\star + \lceil \zeta_n \rceil
      \right) \lesssim \exp \left( - \frac{1}{2} ( \underline{\gamma}
      -1) \D^\star \cdot2^{\D^\star}\right).
    \end{eqnarray*}
    Based on Lemma~\ref{lemma:shapeselection},
    \begin{eqnarray*}
     \lambda^\g  \left(\mathcal{C}'_\en \right)
      \lesssim
      \exp \left( - \frac{1}{2} ( \underline{\gamma}
      -1) \D^\star \cdot2^{\D^\star}\right).
    \end{eqnarray*}
   Note that $2^{\D^\star} \asymp \left( \frac{n}{\log n} \right)^{\frac{1}{2r+1}}$. Combining these two bounds together, we obtain the desired
    result. 
  }
\end{proof}

\subsection{Posterior Concentration Over the H{\"o}lder Class (Proofs
  of Theorem \ref{thm:gibbsholderrate} and Corollary~\ref{cor:giibsholderrate}) }

We apply Theorem \ref{thm:gibbsrate} to show posterior
contration results over the anisotropic H\"older class. The proof is
finished by checking those three conditions introduced in Section
\ref{subsec:gibbsrate}. 
\begin{proof}[Proof of Theorem \ref{thm:gibbsholderrate}]
First, we show a contraction result in the form of Theorem
\ref{thm:gibbsrate}, with $r = \frac{\overline{\alpha}}{ 2
  \overline{\alpha}+d}$. To obtain the contraction
result for the Gibbs posterior $\lambda^\g$, we check Condition
\ref{cond:approx}, Condition \ref{cond:prior} and Condition
\ref{cond:complexity} one by one. 

For Condition \ref{cond:complexity}, we first check at the individual
tree level how many different tree topologies there are under the
\emph{strong shape constraint}. Each depth-$\D$ constrained tree
topology corresponds to a nonnegative integer solution to the equation 
    $\sum_{l=1}^d k_l = \D$. Therefore, the number of depth-$\D$ tree
    topology under the \emph{strong shape constraint} is bounded by $ {\D+d \choose d-1} =
    \frac{1}{(d-1)!} (\D+d) \cdots (\D+2) \asymp \D^{d-1} $. This is
    also the bound for the number of shift design schemes at depth $\D$. Under the
    \emph{weak shape constraint}, the forest grows up to level
    $\D_0$ without any constraint, which leads to at most $(d+1)^{t_0}$ different numbers of tree topologies of depth $\D_0$, where
    $t_0 = 2^{\D_0}$. Then starting from depth $\D_0+1$, we grow trees
    with \emph{strong shape constraint} region-wise up to depth
    $\D > \D_0$. In total, there are $O ((\D- \D_0)^{t_0(d-1)} \cdot
    (d+1)^{t_0})$ different numbers of tree topologies. If we further
    account for $2^\D$ different directions for shifting the sample
    space, the number of possible tree topologies would be at order $O
    ( 2^\D \cdot (\D- \D_0)^{t_0(d-1)} \cdot
    (d+1)^{t_0} )$ 

 At the forest level, assume there are $B$ trees within a depth-$\D$
 forest, then under the \emph{strong shape constraint}, the total
 number of forest topologies of depth $\D$ is at the order: 
    \begin{eqnarray*}
      \underbrace{\D^d}_{\mbox{No. shifts design}} \times
      \left( \underbrace{2^\D \cdot \D^d}_{\mbox{No. shapes for each tree}}
      \right)^B \asymp 
(\D)^{d(B+1)} \cdot (2^\D)^B. 
    \end{eqnarray*}
  Under the \emph{weak shape constraint}, the complexity is
    \begin{eqnarray*}
      \underbrace{\D^d}_{\mbox{No. shifts design}} \times
      \left(\underbrace{2^\D \cdot (d+1)^{t_0} \cdot \D^{t_0 d }}_{\mbox{No. shapes for each tree}}
      \right)^B \asymp
(d+1)^{B t_0}\cdot  (\D)^{d (B t_0 +1)}\cdot (2^\D)^B. 
    \end{eqnarray*}
Therefore, when there are $B \leq c_B 2^\D$ trees with in a depth-$\D$ forest,
Condition \ref{cond:complexity} is satisfied.

For Condition \ref{cond:prior}, as the prior on shift design schemes is uniform given the depth, there exist constants $c_{\text{depth}}, \underline{c}_{\upsilon}$ and $\overline{c}_{\upsilon}$ such that
\begin{eqnarray*}
  q (\D)
  =
  c_{\text{depth}} \cdot
  \exp (- \gamma {\D}  \cdot 2^{\D} ),
  \quad
  \mbox{and}
  \quad
   \underline{c}_\upsilon (\D)^{-d+1} \leq  q(\bm \upsilon| \D) \leq  \overline{c}_\upsilon (\D)^{-d+1}.
\end{eqnarray*}
The prior on the number of trees satisfies 
\begin{eqnarray*}
  q(B | \D) \asymp \D \exp(- \gamma_B \D B)
\quad
\mbox{and}
\quad
B \leq c_B 2^\D. 
\end{eqnarray*}
In terms of the prior on individual trees, there exist constants $\overline{c}_{\T} \geq \underline{c}_{\T}$, such that for any depth-$\D$ tree topology $\T_1$,
\tech{
  \begin{eqnarray*}
    \underline{c}_\T (d)^{-t_0 \D} \cdot 2^{- t_0 \nu \D}
    \leq
    q(\T_1 | \D, B, \bm \upsilon)
    \leq
    \overline{c}_\T (\D)^{-d} 2^{-\D}.
  \end{eqnarray*}
  As the optional P\'olya trees with in a forest are grown independently under the
  prior,
  \begin{eqnarray*}
    \underline{c}_\T^B (d)^{-t_0 B \D} \cdot 2^{- t_0 \nu B \D}
    \leq
    q( \newforest | \D, B,\bm \upsilon)
    \leq
    \overline{c}_\T^B (\D)^{-Bd} \cdot 2^{-B \D}.
  \end{eqnarray*}
}
Therefore, Condition \ref{cond:prior} is also satisfied. 


Last, we need to check Condition \ref{cond:approx}. Assume that the true density function $f_0$ lies in the anisotropic space
  $\mathcal{C}^{\bm \alpha}(\Omega, L)$, where $\bm \alpha =
(\alpha_l)_{1 \leq l \leq d}$ and  $0 < \alpha_l \leq 2$. The ``best''
approximation by a forest at level $\D$ can be obtained as the
following. We set  
  \begin{eqnarray*}
  k_{l} =  \left \lfloor
  \frac{\overline{\alpha}}{\alpha_l}
  \cdot \frac{\D}{d}
  \right \rfloor
  \quad 
  \mbox{with }
  \left(\overline{\alpha}\right)^{-1} =\frac{1}{d} \sum_{j=1}^d \frac{1}{\alpha_j} .
  \end{eqnarray*}
  It is also feasible to add 1 to some of $k_{l}$'s to make the
  summation of them exactly equal to $\D$. Without any shift, denote the tree topology with edge length of
  each leaf node equal to $2^{-k_{l}}$ as $\T$. To achieve the
  ``best" approximation, trees within the depth-$\D$ forest
  $\newforest$ are assumed to be shifted versions of $\T$, and the shift design scheme is indexed by $\bm \upsilon = (2^{-2k_l})_{1 \leq
    l \leq d}$. We define $ f_{\newforest}
  = \frac{1}{|H_{\bm\upsilon}|}
  \sum_{h \in H_{\bm\upsilon}} f_{\T^h}$. By Lemma \ref{lemma:biasforest}, we know that by such an ensemble of trees,
  \tech{when the boundary condition \ref{cond:boundary} holds, $f_{\newforest}$ satisfies
    \begin{eqnarray*}
      \| f_{\newforest}- f_0\|_\infty \leq c_\bias (2^{\D})^{-\overline{\alpha}/d}. 
    \end{eqnarray*}}
 With this bound, we can obtain bounds for $\rho(f_0, f_{\newforest})$
 and $\kl(f_0, f_{\newforest})$ respectively.

 This way, the three conditions are satisfied and an application of
 Theorem \ref{thm:gibbsrate} establishes the posterior concentration
 rate on forest topologies.

 \rev{Next, we show the convergence rate for the predictive GDF \eqref*{eq:gibbsBMA}. Based
   on the proofs for Theorem~\ref{prop:risklower},
   Corollary~\ref{cor:uniformrisklower}, and
   Proposition~\ref{prop:gibbsthickness}, we know that there exists an
 event $\mathcal{W}_n$, such that for any $\Y^n \in
 \mathcal{W}_n$, the following results hold: (i) tree base learners uniformly contract in the
 sense that $ \mathcal{W}_n \subset \cE_n^\0 \otimes \Omega^{n_1}$,
 where $\cE_n^\0$ is defined in Section~\ref{subsec:uniformtreecontract};
 (ii) The lower bound of the centered risk shown in
 Corollary~\ref{cor:uniformrisklower} holds uniformly for all forests;
(iii) The upper bound of the centered risk shown in
Proposition~\ref{prop:gibbsthickness} holds. In addition, there exists
such a sequence of $\mathcal{W}_n$'s satisfying $\mathbb{P}_{f_0}^n
(\mathcal{W}_n) \rightarrow 1$ as $n \rightarrow \infty$. Then, based on the proof
of Theorem~\ref{thm:gibbsrate}, which relies on the proofs of 
Lemma~\ref{lemma:shallowforesets} and
Lemma~\ref{lemma:shapeselection}, and the verification of the conditions
presented above, for any $\Y^n \in
 \mathcal{W}_n$, the Gibbs posterior on
 the set $\mathcal{C}_\en' \cup \{\newforest: \newforest \mbox{ is of
   depth } \D \geq \D^\star + \lceil \zeta_n \rceil \}$ decays to
 zero exponentially fast.
 For any forest estimator
 $\tilde{f}_\newforest$, the triangle inequality implies that $\rho (f_0, \tilde{f}_\newforest) \leq \rho(f_0,
 f_\newforest) + \rho( f_\newforest, \tilde{f}_\newforest).$ For any $\newforest \in
 \left( \mathcal{C}_\en
   \cup
   \{ \D \geq \D^\star + \lceil \zeta_n \rceil \}
 \right)^\complement$, $\rho(f_0, f_\newforest)$ converges at the rate
 $M_n (dL)^{\frac{d}{2 \overline{\alpha} +d }}
 \left( n/ \log n
 \right)^{-\frac{\overline{\alpha}}{2\overline{\alpha} +d }}$. In
 addition, if $\Y^n \in
 \mathcal{W}_n$, $\rho^2( f_\newforest, \tilde{f}_\newforest) \leq
 \frac{1}{B} \sum_{b=1}^B \rho^2( f_b, \tilde{f}_b) \lesssim
 M_\tr^2\cdot \eta_0^2\delta^2(n_0,
 2^{\D^\star + \lceil \zeta_n\rceil })$, where the
 last inequality is based on the uniform contraction result for tree
 base learners shown in Section~\ref{subsec:uniformtreecontract} and
 the fact that $\newforest$ under consideration is a forest of depth
 no larger than $\D^\star + \lceil \zeta_n\rceil$.  Combining these
 two pieces together, we have for any $\Y^n \in \mathcal{W}_n$, the
 Gibbs posterior concentrates on the set $ \newforest \in
 \left( \mathcal{C}_\en
   \cup
   \{ \D \geq \D^\star + \lceil \zeta_n \rceil \}
 \right)^\complement$, and for any forest topology $\newforest$ within
this set, $\rho(f_0, \tilde{f}_\newforest) \lesssim (dL)^{\frac{d}{2 \overline{\alpha} +d }}
 \left( n/ \log n
 \right)^{-\frac{\overline{\alpha}}{2\overline{\alpha} +d }}$, up to a
 logarithmic term. Based
 on the convexity of $\rho^2(\cdot, \cdot)$, we can conclude that the
 predictive GDF converges at the same rate. }
\end{proof}

\rev{The proof for Corollary~\ref{cor:giibsholderrate} is the same as
  that for Theorem~\ref{thm:gibbsholderrate}. The only difference is
  that we apply Corollary~\ref{cor:biasforest} to verify
  Condition~\ref{cond:approx} instead of
  Lemma~\ref{lemma:biasforest}. The proof is skipped here. }

\section{Analysis of the Risk Function for Forests}
\label{sec:risk}
The purpose of this section is to study the properties of
empirical risk \eqref*{eq:riskdef}, 
and develop several technical tools that provide a
bias-variance decomposition of the risk. \rev{In particular, we derive a
lower bound of the forest risk in Section~\ref{subsec:proofrisk} and
an upper bound in Section~\ref{subsec:riskupper}. In
Section~\ref{subsec:proofrisk}, we also show that when the number of
forest topologies does not increase too fast as depth increases, as
assumed by Condition~\ref{cond:complexity}, the risk lower bounds hold
uniformly for all forest topologies with high probability.}

\subsection{Posterior Concentration for Tree-Based Weak Learners}
\label{subsec:uniformtreecontract}
In Section \ref{sec:optpost}, we
show that given a tree topology, the posterior under the constrained OPT prior
concentrates around a good tree approximation at a fast rate. In this section, we show that such type of contraction holds uniformly for all trees
within a forest and over all possible forest topologies, as long as the
number of trees within the forest does not increase too fast. 

With data splitting, we use $\mathscr{C}_\D^\0$ to denote the event that density learning on leaf nodes based on the first subset of data behaves \emph{unsatisfying} for some tree in a depth-$\D$ forest:
{\small
\begin{eqnarray*}
\mathscr{C}_\D^\0 = 
\bigcup_{\{\T_b\}_{b=1}^B \mbox{\scriptsize{of depth $\D$}}} 
\left \{ 
\Y_n^\0: 
\rho \left(f_{\T_b} , E(f | \Y_n^\0, \F_{\T_b}) \right)
\geq 
\writing{M_\tr}\cdot \eta_0 \delta(n_0, 2^\D) 
\mbox{ for any } 1 \leq b \leq B
\right  \},
\end{eqnarray*}
}
where both $\eta_0$ and $M_\tr$ are constants. $\eta_0$ can take any value greater than 1, and a possible choice of $M_\tr$ is \writing{25}. 
The superscript indicates that the event \emph{only} depends
on $\Y_n^\0$. We define the event $\cE_n^\0$ as the following complement, meaning leaf node estimation based on the first subset of data is reasonably good for \emph{all} trees in \emph{any} forest topology: 
\begin{eqnarray*}
  \cE_n^\0 = \left(\cup_{\D=0}^{\D_{\max}}\mathscr{C}_\D^\0 \right) ^\complement.
\end{eqnarray*}


Our first lemma shows {\em uniform} contraction for density learning
on leaf nodes given the trees within {\em any} depth-$\D$ forest. 

\begin{lemma}
  \label{lemma:treeconcentrationdepthuniform}
Consider depth-$\D$ forests $\forest$ in which each tree $\T_b$
  grows with or without shape constraints. Each $\T_b$ may grow
  with perturbations. Assume that Condition \ref{cond:complexity}
  holds. 
  Let $t = 2^\D$ and $\eta_0 \geq 1$ to be a constant. 
For the event $\mathscr{C}_\D^\0 $ defined above, 
 when $n$ is large enough, 
$
  \P_{f_0}^n (\mathscr{C}_\D^\0) \leq 8 \exp(- 2c_2\eta_0^2  n_0 \delta^2(n_0, t)).
$
\end{lemma}
\begin{proof}
  We can rewrite the event $\mathscr{C}_\D^\0$ as
  \begin{eqnarray*}
    \mathscr{C}_\D^\0 &=&    
   \bigcup_{\{\T_b\}_{b=1}^B \text{ of depth } \D}
   \bigcup_{b=1}^B
\left \{ 
\Y_n^\0: \rho(f_{\T_b} , E(f | \Y_n^\0,  \T_b) )
    \geq M_\tr \cdot \eta_0 \delta(n_0, t) \right \}. 
  \end{eqnarray*}
  Applying the Bonferroni inequality, we have
 \begin{eqnarray*}
    \P_{f_0}^n(\mathscr{C}_\D^\0) \leq
    \sum_{\forest \text{of depth } \D}
    \sum_{b=1 }^B 
    \P_{f_0}^n \left( 
\rho(f_{\T_b} , E(f | \Y_n^\0,  \T_b) )
    \geq M_\tr \cdot \eta_0 \delta(n_0, t) \right ),
  \end{eqnarray*}
where $\P_{f_0}^n$ can be understood as the outer measure on the product space $\Omega^n$ induced by i.i.d. random variables $Y_1, \ldots, Y_n$. 
Under Condition \ref{cond:complexity}, by Corollary \ref{cor:posteriorratetreemean} we have
\begin{eqnarray*}
\sum_{b=1}^B \P_{f_0}^n \left( 
\rho(f_{\T_b} , E(f | \Y_n^\0,  \T_b) )
\geq 
M_\tr \cdot \eta_0 \delta(n_0, t) \right )
\leq 8 \exp (- 4 c_2 \eta^2 n_0 \delta^2(n_0, t) +t).
\end{eqnarray*}
  When the number of depth-$\D$ forests is bounded by $\exp(c_\en \D
  \cdot 2^\D)$, for $n$ large enough and
  an appropriately chosen $\eta_0$,
  \begin{eqnarray*}
    \P_{f_0}^n(\mathscr{C}_\D^\0) &\leq& 
8 \exp (- 4c_2 \eta_0^2 n_0 \delta^2(n_0, t) 
+ c_\en t \log t / \log 2+ t) \\
&\leq& 8 \exp (-2c_2 \eta_0^2 n_0 \delta^2(n_0, t) ).
  \end{eqnarray*}
This finishes the proof. 
\end{proof}

Next, we extend the result of Lemma \ref{lemma:treeconcentrationdepthuniform}
to all forest topologies by showing that it holds uniformly over
all forests of different depths with probability tending to 1. 

\begin{lemma}
  \label{cor:treeconcentrationuniform} Under the same condition as
  those for Lemma \ref{lemma:treeconcentrationdepthuniform}, define
  the event
    $\mathscr{C}_n^\0 = \cup_{\D=1}^{\D_{\max}} \mathscr{C}_\D^\0$.
   When $n$ is large enough, 
    $\P_{f_0}^n \left(\mathscr{C}_n^\0\right) \leq
    24 \exp \left( -c_2 \eta_0^2 t_{\min} \log t_{\min} \right)$. 
\end{lemma}
\begin{proof}
Based on the Bonferroni inequality and Lemma \ref{lemma:treeconcentrationdepthuniform},
  \begin{eqnarray*}
    \P_{f_0}^n \left(\mathscr{C}_n^\0 \right)
    \leq
           \sum_{\D=0}^{\D_{\min}} 8 \exp(-2 c_2 \eta_0^2 t_{\min} \log t_{\min} )
           +
           \sum_{\D= \D_{\min}+1}^{\D_{\max}} 8 \exp\left (-2 c_2 \eta_0^2 2^\D \log (2^\D) \right).
  \end{eqnarray*}
In combination with Lemma \ref{lemma:tooltail}, we obtain the desired
upper bound whenever $c_2 \eta_0^2 \geq1$. 
\end{proof}

Lemma \ref{cor:treeconcentrationuniform} implies that uniformly for all
forest topologies, tree-based weak learners are always ``close'' to the best possible tree approximation. It is easy to see that $\P^n_{f_0} (\cE_n^\0 ) \rightarrow 1$
as $n_0 \rightarrow \infty$. In the next subsection, we study properties of the empirical risk conditional on $\cE_n^\0$, and show that 
the risk function defined through data splitting exhibits several enjoyable properties.

\subsection{A Lower Bound for the Risk}
\label{subsec:proofrisk}
The following theorem claims that the risk behaves similarly
to the negative log-likelihood for trees conditional on $\cE_n^\0$. The result
is a lower bound for the risk (indeed the centered risk), which presents a bias-variance decomposition.  

\begin{theorem}
\label{prop:risklower}
Assume that $\newforest= \forest$ is a depth-$\D$ forest with each tree
$\T_b$ growing with or without perturbation and $\D \leq \writing{\D_{\max}}$. 
The risk and centered risk functions associated with $\newforest$ are defined in \eqref*{eq:riskdef} and \eqref*{eq:modifiedriskdef}, 
respectively. Let $t = 2^\D$ and $\wt{\D} = \D \vee \D_{\min}$. $\delta(n, t)$ is the function defined in \eqref*{eq:deltadef}, 
and $f_{\newforest}$ is the forest approximation defined in \eqref*{eq:forestapprox}. 
Assume that $\kappa_L \leq f_0 \leq \kappa_U$ for some constants $\kappa_U \geq \kappa_L >0$. $\eta_1 >1$ is a constant. When $B = o (e^t)$, 
\begin{eqnarray}
\label{eq:risklower}
&& \mathbbm{1}_{\{\Y_n^\0 \in \cE_n^\0 \}}\P_{f_0}^n \bigg( \wt{R}_n \newbforest 
\leq
(1-c_1-c_3) \kl (f_0, f_{\newforest})\\
\nonumber&& \hspace{1.5in} - \frac{c_{\text{opt}} \wt{\D} \kappa_U^2}{c_3 \kappa_L^2} 
M_\tr^2 \cdot \eta_0^2 \delta^2(n_0, t) 
-\eta_1^2 \xi_n \delta^2(n_1,t)
\Big | \Y_n^\0\bigg)  \\
\nonumber&& \hspace{0.5in} \leq  
    \tech{3
    \exp(-(2c_2/3) \eta_1^2 n \delta^2(n_1, t))
+\exp(-c_2' \eta_1 ^2 n \delta^2(n_1, t))},
  \end{eqnarray}
where $c_1, c_3 \in (0, 1)$ and $c_{\text{opt}} >0$ are constants,
$\eta_0$ and $M_\tr$ appear in the definition of $\cE_n^\0$, $c_2$ are specified in Lemma \ref{lemma:uniuppervariance}, $c_2'$
depends on $c_1$ in the same way as that noted in Remark
\ref{remark:approxBernstein},
a possible choice of $\xi_n$ is $c_\xi (\kappa_U / \kappa_L) \log n$
for some constant $c_\xi$ large enough, 
and $\P_{f_0}^n$ is the probability measure induced by $\Y_n$.
\end{theorem}

In order to simplify the notation, for the remaining parts of this  
section, we skip explicitly specifying the tree of forest topology
  $\forest$ when writing a tree- or forest-based approximation or estimate. In
  particular, given a forest topology $\newforest= \forest$, when there is no
  confusion, we use $f_\newforest$ to denote the forest-based approximation
  $f_{\forest}$, $\tilde{f}_\newforest$ to denote the forest-based estimate
  $\tilde{f}_\forest$, $f_b$ for the tree-supported 
  approximation $f_{\T_b}$, and $\tilde{f}_b$ for the tree-based weak
  learner $\tilde{f}_{\T_b} = E( f | \Y^\0, \T_b)$. All these
  density functions depend on the forest topology implicitly.
  Whenever there is a ``$ \hspace{1mm}\tilde{\hspace{2mm} }
  \hspace{1mm}$'' symbol over a density function, it implies 
  that the density function is a Bayesian point estimate, and 
  in most cases it is the posterior
  mean. First, we consider the risk function associated with a given
  forest topology $\newforest$ and prove Theorem \ref{prop:risklower}.  
\begin{proof}[Proof of Theorem \ref{prop:risklower}]
For the centered risk associated with the forest topology $\newforest$,  
\begin{eqnarray*}
  \wt{R}_n \newbforest
  =\frac{1}{n_1}
      \sum_{i=n_0+1}^n
      \log \left (f_0(Y_i) / f_\newforest (Y_i) \right)
      -  \frac{1}{n_1}
      \sum_{i=n_0+1}^n
        \log \left (\tilde{f}_\newforest(Y_i) / f_\newforest (Y_i) \right).
\end{eqnarray*}
The right-hand side can be approximately viewed as a bias-variance
decomposition, where the first term is mainly driven by the bias and the second by the variance. 

For the first term, it is independent of the event $\cE_n^\0$. By
Remark \ref{remark:approxBernstein}, for some $\eta_1 > 1$ 
\begin{eqnarray*}
  &&\P_{f_0}^n \left( 
\frac{1}{n_1}
      \sum_{i=n_0+1}^n
      \log \left (f_0(Y_i) / f_\newforest (Y_i) \right)
\leq 
(1-c_1) \kl (f_0, f_\newforest) - c_1 \eta_1^2 \delta^2(n_1, t) \right)\\
&& \hspace{1cm} \leq  
\exp \left( - c_2' n \max\{ \kl (f_0, f_\newforest), \eta_1^2 \delta^2(n_1,
   t)\} \right) 
\leq \exp \left(- c_2' n_1 \eta_1^2 \delta^2(n_1,t) \right). 
\end{eqnarray*}

For the second term,
{\normalsize
\begin{eqnarray*}
- \frac{1}{n_1}
      \sum_{i=n_0+1}^n
        \log \left (\frac{\tilde{f}_\newforest(Y_i) }{ f_\newforest (Y_i)} \right)
  &=& - \frac{1}{n_1}
  \sum_{i = n_0 +1}^n
  \log \left( 1 +
  \frac{\frac{1}{B} \sum_{b=1}^B \tilde{f}_b (Y_i)}
  {f_\newforest (Y_i)} -1
  \right)\\
  & \geq &
   -  \frac{1}{n_1}\sum_{i = n_0 +1}^n
           \left(
           \frac{\frac{1}{B} \sum_{b=1}^B \tilde{f}_b (Y_i)}
           {f_\newforest (Y_i)}
           -1 \right),
\end{eqnarray*}}
where the first equality is obtained by plugging in the forest
estimate, and the second inequality is based on
the fact that $\log (1+x) \leq x$ for any $x >-1$. After
rearranging the summand and switching the order of
summations with respect to $b$ and $i$, we obtain
{\normalsize
\begin{eqnarray*}
 -\frac{1}{n_1}
      \sum_{i=n_0+1}^n
        \log \left (\frac{\tilde{f}_\newforest(Y_i) }{ f_\newforest (Y_i)} \right)
  &\geq&
   -  \frac{1}{n_1}\sum_{i = n_0 +1}^n
           \left(
  \frac{1}{B} \sum_{b=1}^B
  \left(
  \frac{ \tilde{f}_b (Y_i)}
           {f_\newforest (Y_i)}
  -1 \right)
  \right)\\
  &=&  - \frac{1}{B} \sum_{b=1}^B
      \left( \frac{1}{n_1} \sum_{i=n_0+1}^n
      \left(
      \tilde{f}_b(Y_i) - f_b(Y_i) \right)/ f_\newforest (Y_i)  
      \right).
\end{eqnarray*}
}
By Lemma \ref{lemma:uniuppervariance}, we can show that conditional on
any $\Y_n^\0 \in \cE_n^\0$, with probability at
least \tech{$1- 3 \exp(-(2c_2/3) \eta_1^2 n_1 \delta^2(n_1, t) )$},
the right hand side
is lower bounded 
\begin{eqnarray*}
  - c_3 \kl (f_0, f_\newforest)- (c_{\text{opt}} \D/ c_3)(
  \kappa_U/ \kappa_L)^2  M_\tr^2 \cdot \eta_0^2 \delta^2(n_0, t) -
 \eta_1^2 \xi_n \delta^2(n_1,t).
\end{eqnarray*}
Combining the two lower bounds together, we obtain the desired
result. 
\end{proof}

We introduce the following bound for the term
$\frac{1}{n_1} \sum_{i=n_0+1}^n (\tilde{f}_b(Y_i) - f_b(Y_i) )/ f_\newforest
(Y_i) $, which holds simultaneously for all trees within the
forest. The bound is obtained by studying the empirical process
indexed by $(f- f_b)/f_\newforest$ with $f \in \mathcal{F}_{\T_b}$. As we analyze the risk function conditional on the event $\cE_n^\0$, we may view the Bayesian point estimate
$\tilde{f}_b$ as deterministic, lying within a small Hellinger ball (or $L_2$ ball) around $f_{b}$. In the remainder of this subsection,
we use $\mathbb{G}_{n_1}(u)$ to denote the empirical process
$(1/\sqrt{n_1})\sum_{i=n_0+1}^n u(Y_i)$ and define $P_{n_1}(u) =
\frac{1}{n_1} \sum_{i=n_0+1}^n u(Y_i)$. 

\begin{lemma}
  \label{lemma:uniuppervariance}
Assume that $\newforest$ is a depth-$\D$ forest topology. Let $u_b^f = (f -f_b) / f_\newforest$ with $f \in \mathcal{F}_{\T_b}$, where both $f_b$ and $f_\newforest$ are approximating densities associated with $\forest$. Set $\wt{\D} = \D \vee \D_{\min}$. Under the same conditions as Theorem \ref{prop:risklower}, for some $\eta_1 >1$, $\xi_n \asymp (\kappa_U/ \kappa_L) \log n$, 
any $\Y_n^\0 \in \cE_n^\0$, 
  \begin{eqnarray*}
 && 
    \P_{f_0}^n \bigg(   
\max_{1 \leq b \leq B}
P_{n_1}\left( (\tilde{f}_b- f_b)/f_\newforest \right) 
     \geq 
    \tech{ c_3\kl (f_0, f_\newforest)}\\
    && \hspace{5.5cm}
+  \tech{\frac{c_{\text{opt}} \wt{\D} \kappa_U^2 }{c_3\kappa_L^2} M_\tr^2 \cdot \eta_0^2 \delta^2 (n_0, t)}
    + \eta_1^2 \xi_n \delta^2 (n_1, t)
\Big |
\Y_n^\0
\bigg) \\
&& \hspace{6cm} 
\leq 3 \exp\left(- (2c_2/3) \eta_1^2 n_1 \delta^2 (n_1, t) \right),
  \end{eqnarray*}
where $\eta_0$ and $M_\tr$ are the constants appeared in the definition of the event
$\cE_n^\0$, $c_3 \in (0, 1)$ is a constant, and value of $c_2$ is provided in the proof.  
\end{lemma}

\begin{remark}
The constant $c_2$ in the statement of the lemma can be different from that in Proposition \ref{thm:var_deviation}. We abuse the notation, as both of them are some constants between 0 and 1 and appear in the bounds for tail probabilities of empirical processes. 
\end{remark}

\begin{proof}
Conditional on any $ \Y_n^\0 \in \cE_n^\0$, $\tilde{f}_b$ lies within a small Hellinger ball around $f_b$ as shown by Lemma \ref{cor:treeconcentrationuniform}. \tech{Under the upper truncated prior, $\tilde{f}_b \leq  c_{\text{opt}} (\D \vee \D_{\min})$.} Then with $\mathcal{B}_b = \{f \in \mathcal{F}_{\T_b}: \rho(f_b, f) \leq
M_\tr \cdot \eta_0 \delta(n_0, t), f \leq c_{\text{opt}} \wt{\D}\}$, for any $\Y_n^\0 \in \cE_n^\0$, 
{\small
  \begin{eqnarray}
    \label{eq:empirical2bound}
  &&  \hspace{0.1in} \P_{f_0}^n \left(   
P_{n_1}\left ( \frac{\tilde{f}_b- f_b}{f_\newforest}  \right) 
     \geq 
     \tech{ c_3\kl (f_0, f_\newforest) 
+  \frac{ c_{\text{opt}} \wt{\D} \kappa_U^2 }{c_3\kappa_L^2} M_\tr^2 \cdot \eta_0^2 \delta^2 (n_0, t)}
    + \eta_1^2 \xi_n \delta^2 (n_1, t)
\Big | \Y_n^\0
\right)\\
\nonumber  && \leq \P_{f_0}^n \left ( 
     \sup_{f \in \mathcal{B}_b}
P_{n_1} (u_b^f)
\geq 
    \tech{ c_3\kl (f_0, f_\newforest) 
+  \frac{ c_{\text{opt}} \wt{\D}\kappa_U^2 }{c_3\kappa_L^2} M_\tr^2 \cdot \eta_0^2 \delta^2 (n_0, t)}
    + \eta_1^2 \xi_n \delta^2 (n_1, t)
\Big | \Y_n^\0
\right). 
  \end{eqnarray}
}
  As the event on the right-hand side of the inequality only depends on
  $\Y_n^\1$, being independent from $\cE_n^\0$, we may
simply work on the unconditional event for the rest of the proof. The
  probability above can be further bounded by
{\normalsize
  \begin{eqnarray*}
 &&    \P_{f_0}^n \bigg(
       \sup_{ f \in \mathcal{B}_b}
         \frac{1}{\sqrt{n_1}} \mathbb{G}_{n_1} ( u_b^f)
\geq
 \tech{ c_3\kl (f_0, f_\newforest) } 
 + \eta_1^2 \xi_n \delta^2 (n_1, t) \\
 && \hspace{1.7in}
+ \tech {\frac{c_{\text{opt}} \wt{\D} \kappa_U^2 }{c_3\kappa_L^2} M_\tr^2 \cdot \eta_0^2 \delta^2 (n_0, t)}   
 -
   \sup_{f \in \mathcal{B}_b} \E(u_b^f(Y_n))
\bigg) 
          ,
\end{eqnarray*}}
where $\mathbb{G}_{n_1} ( \cdot ) $ is the empirical process defined on 
$\Y^\1$. 
If we examine the term $\E(u_b^f (Y_n))$ carefully, 
\begin{eqnarray*}
   \E(u_b^f (Y_n)) 
              &=& \int f_\newforest (f-f_b)/f_\newforest +
                  \int (f_0 - f_\newforest )(f - f_b) /f_\newforest \\
              &=& 0 + \int (f_0/f_\newforest - 1) (f-f_b).
\end{eqnarray*}
It is not hard to see that, for any $f \in \mathcal{B}_B$, $ \E(u_b^f (Y_n)) \leq
 (4 \sqrt{c_{\text{opt}} \wt{\D}}) (\kappa_U/ \kappa_L) \rho(f_b, f) \kl^{1/2}(f_0, f_\newforest)$.
Therefore, to bound the probability in \eqref*{eq:empirical2bound}, it
suffices to bound
\begin{eqnarray*}
    \P_{f_0}^n \left(
       \sup_{f \in \mathcal{B}_b}
         \mathbb{G}_{n_1}( u_b^f) 
       \geq  
       \eta_1^2 n_1^{1/2} \xi_n \delta^2 (n_1, t) 
       \right).
\end{eqnarray*}
We use tools for empirical processes to bound this probability. In specific, Theorem~3 in
the paper by \cite{shen1994} can provide us with the desired result. To apply the
theorem, we need to check the conditions one-by-one. For any $f \in
\mathcal{B}_b $, $u_b^f \leq (\kappa_U/ \kappa_L) c_{\text{opt}} \wt{\D}$. 
$
\mbox{Var}(u_b^f (Y_n) ) \leq  \E \left( \left(u_b^f(Y_n) \right)^2\right) = \int f_0 (f - f_b)^2 / f_\newforest^2
  \leq (4 c_{\text{opt}} \wt{\D}) (\kappa_U /\kappa_L)^2 M_\tr^2 \cdot \eta_0^2 \delta^2(n_0, t).$ 
Note that the bracketing metric entropy of $\mathcal{B}_b$ under the $L^2$-norm is bounded by that under the supnorm up to some constant, as $f_\newforest$ is bounded from below. By setting $M$ in the
theorem as $n_1^{1/2} \eta_1^2 \xi_n \delta^2 (n_1, t)$, we can show that 
for some $c_2 \in (0, 1)$,  
\begin{eqnarray*}
 \eqref*{eq:empirical2bound} \leq 3 \exp (- c_2\eta_1^2 n_1 \delta^2
  (n_1, t) ). 
\end{eqnarray*}
A possible choice of $\xi_n$ is $c_\xi (\kappa_U / \kappa_L) \log n$
for some constant $c_\xi$ large enough. 

When the number of trees within the forest satisfies $B = o (\exp (t))$, an application of the Bonferroni inequality provides the uniform bound over all trees. 
\end{proof}

The next corollary shows that when forest complexity is under control, the result of Lemma \ref{lemma:uniuppervariance} holds uniformly for all forest topologies. 
\begin{corollary}
\label{cor:uniformrisklower}
Assume that the conditions of Theorem \ref{prop:risklower} hold. $\eta_1$, $c_2$, $c_3$ and $\{ \xi_n \}_{n\geq 1}$ are the same as those in Lemma \ref{lemma:uniuppervariance}. In addition, we assume that Condition \ref{cond:complexity} 
holds. Then conditional on any $\Y_n^\0 \in \cE_n^\0$, under $\P_{f_0}^n$ with probability at least 
$$1- 9 exp( -(c_2/2) \eta_1^2 t_{\min} \log t_{\min}) - 3 \exp (- (c_2'/2)
 \eta_1^2 t_{\min} \log t_{\min} ),$$ the following inequality holds uniformly for all forest topologies
 \begin{eqnarray*}
\wt{R}_n \newbforest 
    \geq 
    (1-c_1-c_3) \kl (f_0, f_\newforest)
- \tech{\frac{ c_{\text{opt}} \wt{\D} \kappa_U^2}{ c_3 \kappa_L^2} M_\tr^2 \cdot \eta_0^2 \delta^2(n_0, t)}
-\eta_1^2 \xi_n \delta^2(n_1,t).
  \end{eqnarray*}
\end{corollary}
\begin{proof}
To obtain the simultaneous bound, we need to carefully count the
number of times that the inequality in Theorem \ref{prop:risklower} needs to be applied. Note that at depth $\D$, we
need to consider at most $\exp (c_\en \D \cdot 2^\D)$ different forest topologies. In total, based on Bonferroni's inequality, for depth-$\D$ forests, the probability that the desired inequality does not hold is bounded by
$$ 3 \exp (-(c_2/2) \eta_1^2 n_1 \delta^2(n_1, 2^\D) ) + \exp (-(c_2'/2)
\eta_1^2 n_1 \delta^2(n_1,2^\D)),$$
where $\eta_1$ should be chosen large enough. With such a bound for
all depth-$\D$ forest topologies, we can obtain the uniform bound by using a similar argument as that for Corollary
\ref{cor:treeconcentrationuniform}. 
\end{proof}

\begin{remark}
\label{remark:uniformrisklower}
As the event $\cE_n^\0$ occurs with probability tending to 1 under $\P_{f_0}^n$, we can conclude that the lower bound for the risk $\wt{R}_n$ in Corollary
\ref{cor:uniformrisklower} also holds uniformly for all forest topologies with probability tending to 1 under the measure induced by $\Y_n$. 
\end{remark}

\subsection{An Upper Bound for the Risk}
\label{subsec:riskupper}
In the following proposition, we provide an upper bound for the
risk function associated with a specific forest topology, namely
$f_{\newstarf}$ defined in the beginning of Section \ref{sec:proofgibbs}. 

\begin{proposition}
\label{prop:gibbsthickness}
Assume Condition \ref{cond:approx} 
is satisfied. In addition, we assume $\kappa_L \leq f_0 \leq \kappa_U$ for some constants
$\kappa_U \geq \kappa_L >0$. 
For any $\Y_n^\0 \in \cE_n^\0$,  
\begin{eqnarray*}
  \nonumber && \P_{f_0}^n
  \bigg(
\wt{R}_n \newbstarf
   \geq ( 1 + c_1^\star) \kl (f_0, f_{\newstarf}) 
+ \frac{1}{c_1^\star} \xi_n^\star M_\tr^2 \cdot \eta_0^2 \delta^2(n_0, t^\star)\\ 
  && \hspace{1.6in}
+ \xi_n^\star (\eta_1^\star)^2 \delta^2(n_1, t^\star)
  \hspace{1mm} \big | \hspace{1mm}
    \Y_n^\0 \bigg) 
  \leq \tech{\exp(-c_2^\star(\eta_1^\star)^2 n_1 \delta^2(n_1, t^\star))},
\end{eqnarray*}
where $ \P_{f_0}^n ( \cdot | \Y_n^\0)$ is understood as the product probability measure induced by $\Y_n^\1$, $\{ \xi_n^\star \}_{n \geq 1}$ is a sequence increasing as polynomials of $\log n$, $t^\star = 2^{\D^\star}$, and the function $\delta(n, t)$ is
defined by \eqref*{eq:deltadef}. 
Consequently, 
\begin{eqnarray}
 \label{eq:riskupper}
   && \P_{f_0}^n
  \bigg(
\wt{R}_n \newbstarf
   \geq ( 1 + c_1^\star) \kl (f_0, f_{\newstarf}) 
+ \frac{1}{c_1^\star} \xi_n^\star M_\tr^2 \cdot \eta_0^2 \delta^2(n_0, t^\star)
\\ 
  \nonumber&& \hspace{0.5in}
  + \xi_n^\star (\eta_1^\star)^2 \delta^2(n_1, t^\star)
  \bigg) 
  \leq \tech{\exp(-c_2^\star(\eta_1^\star)^2 n_1 \delta^2(n_1,
     t^\star))} + \P_{f_0}^n \left( ( \cE_n^\0 )^\complement
     \right).
\end{eqnarray}
\end{proposition}
\begin{proof}
 Given the forest topology $\newstarf$ and conditional on any $ \Y_n^\0
  \in \cE_n^\0$, we may view the density forest
  estimator $\tilde{f}_{\newstarf}$ as deterministic. By triangle
  inequality for the Hellinger distance, we have the following
  bias-variance decomposition,
  \begin{eqnarray*}
    \rho(f_0, \tilde{f}_{\newstarf}) \leq \rho(f_0, f_{\newstarf}) 
+ \rho(f_{\newstarf}, \tilde{f}_{\newstarf}). 
  \end{eqnarray*}
By joint convexity of $\rho^2 (\cdot, \cdot)$ and definition
  of the event $\cE_n^\0$, for any $\Y_n^\0 \in \cE_n^\0$,
  \begin{eqnarray*}
    \rho^2(f_{\newstarf}, \tilde{f}_{\newstarf})
 \mathbbm{1}_{\{ \Y_n^\0 \in \cE_n^\0 \}}
\leq \frac{1}{B} \sum_{b=1}^B 
\rho^2(f_{\T_b^\star},\tilde{f}_{\T_b^\star} )
 \mathbbm{1}_{\{ \Y_n^\0 \in \cE_n^\0 \}}
\leq \writing{M_\tr^2} \cdot \eta_0^2 \delta^2 (n_0, t^\star). 
  \end{eqnarray*}
With the choice of $\D^\star$ specified in \eqref*{eq:optimalforestdepth},
$\rho(f_0, f_{\newstarf})$ is of the same or smaller order compared to
$\rho(f_{\newstarf},
  \tilde{f}_{\newstarf}) $. This way, we know the order of $\rho(f_0,
  \tilde{f}_{\newstarf})$ approximately. 

Keeping the bound for $\rho(f_0,
  \tilde{f}_{\newstarf})$ in minde, our goal is to show that for any $\Y_n^\0 \in \cE_n^\0$, 
{\normalsize
\begin{eqnarray*}
\mathbbm{1}_{\{ \Y_n^\0 \in \cE_n^\0\}} \P_{f_0}^n \left( 
- \wt{R}_n \newbstarf
+ \E_{f_0}^n \left( 
\wt{R}_n \newbstarf \cond \Y_n^\0 \right)  
\leq - \xi_n^\star (\eta_1^\star)^2 \delta^2(n_1, t^\star) 
\hspace{1mm} \Big | \hspace{1mm}
\Y_n^\0
\right) 
\rightarrow 0.
\end{eqnarray*}
}
The probability shown in the equation above is denoted as
$P^\star (\Y_n^\0) $ in the remaining part of the proof. 
The proof is based on the Bernstein's inequality.  

We first consider a simpler case when $ f_0 \leq \kappa_U $ and
$\tilde{f}^\star_\en \geq \kappa_L$. In this case, we can apply the
Bernstein's equality for upper bounded random variables, as
$-Z_{\tilde{f}_{\newstarf}} + \E_{f_0} (Z_{\tilde{f}_{\newstarf}} | \Y_n^\0) \leq \log
(\kappa_U / \kappa_L)$. By Lemma \ref{lemma:kltoH}, $\E_{f_0}^n \left( (Z_{\tilde{f}_{\newstarf}}
  (Y_n) )^2 | \Y_n^\0\right )\leq 4 (1 +\log (\kappa_U/ \kappa_L) )^2 \rho^2
(f_0, \tilde{f}_{\newstarf})$ when $\kappa_U/ \kappa_L \geq 2$. By
the Bernstein's inequality for bounded random variables,
{\small
\begin{eqnarray*}
  P^\star\left( \Y_n^\0 \right) \leq \exp \left(
- \frac{(n_1^{1/2}  \xi_n^\star (\eta_1^\star)^2 \delta^2(n_1, t^\star) )^2}
{8 (1 + \log (\kappa_U/ \kappa_L) )^2\rho^2(f_0, \tilde{f}_{\newstarf})
+ (2/3) \log(\kappa_U/ \kappa_L) \xi_n^\star  (\eta_1^\star)^2 \delta^2(n_1, t^\star)}
\right).
\end{eqnarray*}
}
By the definition of $\D^\star$,
$\rho(f_0, f_{\newstarf}) = O\left( \rho(f_{\newstarf},
  \tilde{f}_{\newstarf}) \right)$
\tech{when $\eta_0$ is appropriately chosen.} Thus the right hand side
of the inequality is bounded by $\exp(-c_2^\star n_1   (\eta_1^\star)^2
\delta^2(n_1, t^\star) )$ for some constant $c_2^\star$. 

In the case $\tilde{f}_{\newstarf}$ is not bounded from below, we can
find a decaying lower bound for $\tilde{f}_{\newstarf}$. As we use the
posterior mean for density learning given the tree topology, it is not hard to see that
each $\tilde{f}_b^\star$ satisfies $\tilde{f}_b^\star \geq (\alpha
/n_0)^{\D^\star}$. Only $\log (1/ \kappa_L)$ is involved in the bound
for $P^\star\left( \Y_n^\0 \right)$. Therefore, we can find a sequence
$\{ \xi_n^\star \}_{n \geq 1}$ that increases as polynomials of $\log n$, such that in the case when
$\tilde{f}_{\newstarf}$ is not bounded from below, the bound 
 $ P^\star( \Y_n^\0) \leq \exp\left(-c_2^\star n_1   (\eta_1^\star)^2
\delta^2(n_1, t^\star) \right)
$ 
still holds. 


Another observation is that, for $\Y_n^\0 \in \cE_n^\0$,
\begin{eqnarray*}
  \mathbbm{1}_{\{\Y_n^\0 \in \cE_n^\0\} }\E_{f_0}^n (\wt{R}_n \newbstarf | \Y_n^\0)
\leq (1 +c_1^\star) \kl (f_0, f_{\newstarf}) + \frac{1}{c_1^\star} \xi_n^\star 
M_\tr^2 \cdot \eta_0^2 \delta^2(n_0, t^\star),
\end{eqnarray*}
 where $\{ \xi_n^\star \}_{n\geq 1}$ is
still a sequence increasing as polynomials of $\log n$. This shows the
first half of our claim. 

\eqref*{eq:riskupper} can be obtained by applying the Fubini's
theorem. 
\end{proof}

\section{Analysis of the Likelihood Function for Individual Trees}
\label{sec:riskanalysis}
In this section, we analyze the
asymptotic behavior of the likelihood function associated with
individual trees. The upper and lower
bounds of the (log) likelihood present a bias-variance
decomposition. \rev{These results are necessary tools for establishing 
  posterior contraction results for individual tree models.} In order to make the results from this section
applicable in general cases, we do not consider data splitting and
derive bounds for the likelihood evaluated on the full dataset. 
\subsection{ An Upper Bound for the Likelihood Ratio Associated with Trees (Proof of
  Theorem \ref{thm:likbiasvar} )} 
\label{subsec:likelihoodupper}
We first focus on the upper bound for
the largest possible likelihood defined on a tree topology.
In order to simplify the notation in the proof, we define $Z_f (Y_i)= \log (f/f_0)(Y_i)$. 
 An obstacle here is that the log-likelihood ratio $Z_f$ does not always have absolute moment generating functions. To
handle this difficulty, we introduce the lower-truncated version of
$Z_f$: 
\begin{eqnarray}
   \label{eq:lowertruncation}
  \tilde{f} =
\begin{cases} f, & \text{if $ f> \exp(-\tau) f_0$, } \\
\exp(-\tau)f_0, & \text{if $f \leq \exp(-\tau) f_0$.}
\end{cases}
\quad 
\widetilde{Z}_f = Z_{\tilde{f}} =
\begin{cases} Z_f, &\text{if $Z_f > -\tau$,} \\
-\tau, &\text{if $Z_f \leq -\tau$.}
\end{cases} 
\end{eqnarray}
In this definition, $\tau$ is a truncation constant. An outline of the
proof for Theorem \ref{thm:likbiasvar} is as follows.
\begin{proof}[Proof of Theorem \ref{thm:likbiasvar}] 
We first decompose the left hand side of (\ref{eq:OneSideIneqUpper})
into the following two parts:  
\begin{eqnarray*}
  \P_{f_0}^n \left( \sup_{f \in \Fat}
\prod_{i=1}^n \frac{f(Y_i)}{f_0(Y_i)} \geq 
\exp\left( -n(1-c_1) \rho^2 \left( f_0,\fat) \right) +
  \eta^2 n \delta^2(n, t) \right)\right)  \leq P_1 + P_2, 
\end{eqnarray*}
where
\begin{eqnarray*}
P_1 &=& \P_{f_0}^n \left( 
\prod_{i=1}^n \frac{\fat(Y_i)}{f_0(Y_i)} \geq 
\exp\left( -(1-c_1)n \rho^2 \left( f_0,\fat) \right) \right)\right), \\
  P_2 &=& \P_{f_0}^n \left( \sup_{f \in  \Fat}
\prod_{i=1}^n \frac{f(Y_i)}{\fat (Y_i)} \geq 
\exp\left( 
          \eta^2 n \delta^2(n, t)  \right)\right).
\end{eqnarray*}
The proof is mainly based on the fact that $P_1$ cannot be too large given the approximation ability by $\Fat$ and $P_2$ can be controlled by studying empirical processes indexed by $Z_f$'s. We provide
bounds for these two terms in the following.

A bound for $P_1$ can be obtained by applying Lemma
\ref{lemma:approxBernstein}. For $P_2$, by Corollary
\ref{cor:lik_deviation_fat}, we have
\begin{eqnarray*}
  P_2 \leq 8 \exp\left( - c_2 \eta^2 n \delta^2(n, t) \right).
\end{eqnarray*}
The desired result follows.
\end{proof}

The following lemma is on bounding $P_1$. 
\begin{lemma}
\label{lemma:approxBernstein}
    For the density function $f_0$ on $\Omega$, let $g$ be an approximation to $f_0$. $\Y_n$ is a random sample from $f_0$. Then for some constant $c_1 \in (0, 1)$,
    \begin{eqnarray*}
 \P_{f_0}^n \left( 
\prod_{i=1}^n \frac{g(Y_i)}{f_0(Y_i)} \geq 
\exp\left( -(1-c_1)n \rho^2 \left( f_0,g \right) \right)
\right)
\leq \exp\left(-c_2'n \rho^2 (f_0, g) \right).
    \end{eqnarray*}
    The choice of $c_2'$ is provided in the proof.
\end{lemma}
\begin{proof}
The proof is based on the Bernstein's
inequality. Denote the probability of the event by $P_1$. First, we apply a lower truncation to the likelihood
ratio with respect to $f_0$ as defined in equation \eqref*{eq:lowertruncation}. 
With the truncated
likelihood ratio, 
\begin{eqnarray*}
  P_1 \leq \P_{f_0}^n \left( \sum_{i=1}^n \left( \widetilde{Z}_{g}(Y_i) 
  - \E \wt{Z}_{g}(Y_i)\right)
  \geq
  -(1-c_1) n \rho^2(f_0, g) - n \E \wt{Z}_{g} (Y_1)
  \right).
\end{eqnarray*}
By Lemma \ref{lemma:boundofexpectation}, $- \E \wt{Z}_{g} (Y_1)
\geq ( 1- \gamma) \rho^2(f_0, g)$, where $\gamma$ is a constant
depending on the truncation parameter $\tau$ as defined in the
lemma. Plugging in this bound for the truncated log-likelihood ratio, we have
\begin{eqnarray*}
  P_1 \leq \P_{f_0}^n \left( n^{-1/2}\sum_{i=1}^n \left( \widetilde{Z}_{g}(Y_i) 
  - \E \wt{Z}_{g}(Y_i)\right)
  \geq
  (c_1- \gamma) n^{1/2} \rho^2(f_0, g) 
  \right).
\end{eqnarray*}
With appropriately chosen $\gamma$ and $c_1$, we may assume $c_1 -
\gamma >0$. 
Then by the Bernstein's inequality under moment conditions, with $c_0 = (\exp(\tau/2) -1 -
\tau/2)/(1- \exp(-\tau/2))^2$, again $\tau$ is the truncation
parameter, we obtain
\begin{eqnarray*}
 \nonumber  P_1 &\leq& \exp\left( - \frac{(c_1- \gamma)^2 n \rho^4(f_0, g) }
  {8 \left( 8c_0 \rho^2(f_0, g)  + 2  (c_1- \gamma) \rho^2(f_0,
  g) \right)}\right)\\
  &\leq& \exp\left(- \frac{(c_1- \gamma)^2}{(2^6c_0 +2^4 (c_1 - \gamma))}
  \cdot n \rho^2(f_0, g)\right).
\end{eqnarray*}
Therefore, we may set $c_2' = \frac{(c_1- \gamma)^2}{(2^6c_0 +2^4 (c_1 - \gamma))}$. This finishes the proof. 
\end{proof}

In the proof of Lemma \ref{lemma:approxBernstein}, we apply a truncation to the
likelihood ratio and characterize the bias by using the Hellinger
distance, as the negative part of the log-likelihood is not always
bounded. When both $f_0$ and $g$ are bounded away from zero, we can obtain a tighter bound, presented in the following remark. 
\begin{remark}
\label{remark:approxBernstein}
In Lemma \ref{lemma:approxBernstein}, if both $f_0$ and $g$ are
bounded away from zero, and $f_0 \leq \kappa_U$ for some $\kappa_U>0$, 
then 
\begin{eqnarray*}
 \P_{f_0}^n \left( 
\prod_{i=1}^n \frac{g(Y_i)}{f_0(Y_i)} \geq 
\exp\left( -(1-c_1)n \kl \left( f_0,g \right) \right)
\right)
\leq \exp\left(-c_2'n \kl (f_0, g) \right).
    \end{eqnarray*}
    To show the result, the lower truncation is not needed,
  as the function $ g/f_0$ is bounded from below. Let $\tau =
  \inf_\Omega \log( g/f_0)$. By using a similar argument as that for
  Lemma \ref{lemma:approxBernstein}, we can show that the inequality
  hold when $c_2' =
  c_1^2 /(2^6c_0 + 2^4 c_1 ) $ and $c_0 = (\exp(\tau/2) -1 -
\tau/2)/(1- \exp(-\tau/2))^2$. Note that to simplify the notation, we
use the same notation for $c_2'$ as that in Lemma
\ref{lemma:approxBernstein}, although the value can be slightly different. 
    \end{remark}
When $f_0$ is bounded from above and below, so is $\fat$. By applying Remark \ref{remark:approxBernstein}, we obtain an alternative inequality for the likelihood ratio. 
\begin{corollary}
  \label{cor:likbiasvar}
Under the same conditions as those for Theorem \ref{thm:likbiasvar},
further assume that $f_0$ is bounded from below and above. Then for $n$ large enough, 
 \begin{eqnarray*}
 \nonumber &&   \P_{f_0}^n \bigg( 
\sup_{f \in \mathcal{F}_{\mathcal{A}_t}} 
\prod_{i=1}^n \frac{f(Y_i)}{f_0(Y_i)} 
\geq 
\exp \Big(-n (1-c_1) \kl (f_0, \fat)
+
\eta^2 n \delta^2(n, t) \Big)
    \bigg)\\
  && \hspace{5cm}  \leq
\exp\left( -c_2' n \kl (f_0, \fat) \right)
+
    8\exp\left(-c_2 \eta^2 n \delta^2(n, t) \right).
  \end{eqnarray*}
\end{corollary}
\begin{proof}
  The proof is similar to that for Theorem \ref{thm:likbiasvar},
  except for the bound for $P_1$.  
\end{proof}

The following proposition provides a bound for the tail probability
associated with empirical processes indexed by $Z_f$'s. The
proposition is useful for bounding $P_2$ in the proof of Theorem
\ref{thm:likbiasvar}, and for deriving posterior concentration results
on under the tree-like prior. 
\begin{proposition}
   \label{thm:var_deviation}
 Under the same assumptions as those for Theorem \ref{thm:likbiasvar}, for
$\eta>1$ and some constants $c_2 \in (0,1)$ (to be specified in the proof), 
\begin{eqnarray}
\label{eq:var_deviation_lik}
 && \P_{f_0}^n \left( 
\sup_{f \in  \Fat, \rho(f, \fat) \geq
            2 \eta \delta(n,t)}
\prod_{i=1}^n \frac{f(Y_i)}{\fat (Y_i)} \geq 
\exp\left( 
         - \frac{1}{4} n \eta^2 \delta^2(n,t) \right) \right) \\
\nonumber && \hspace{3in} \leq 
4\exp\left(-\frac{1}{16}c_2\eta^2 n \delta^2(n, t) \right).
\end{eqnarray}
\end{proposition}
The result can be obtained by applying the one-sided large deviation
inequality in Lemma \ref{lemma:largedeviation} iteratively. 
Before providing a detailed proof, we introduce an alternative version of lower-truncated likelihood ratio,
where the true density $f_0$ is replaced by the ``best approximation''
$f_{\mathcal{A}_t}$ within $\mathcal{F}_{\mathcal{A}_t}$ under the KL divergence.
Let $Z_{\mathcal{A}_t, f} (Y_i) = \log 
(f(Y_i) /f_{\mathcal{A}_t}(Y_i))$. The lower-truncated versions of $f$
and $Z_{\mathcal{A}_t, f}$ are defined similarly as before,
\begin{eqnarray*}
  \tilde{f} |_{\Fat} =
  \begin{cases}
    f, & \mbox{if } f > \exp(-\tau) f_{\mathcal{A}_t}, \\
    \exp(-\tau) f_{\mathcal{A}_t}, & \mbox{if } f \leq \exp(-\tau) \fat.
  \end{cases}
\quad     
\widetilde{Z}_{\at, f }=
  \begin{cases}
    Z_{\at, f}, & \mbox{if } Z_{\at, f} > - \tau, \\
    -\tau, & \mbox{if } Z_{\at, f} \leq -\tau. 
  \end{cases}
\end{eqnarray*}
When it does not lead to ambiguity, we still use the notation 
$\tilde{f}$ to denote the truncated likelihood with respect to $\fat$.
We can show that the truncated likelihood ratio with respect to
  $\fat$ still enjoys most desired properties which lead to expected 
  concentration results. In particular, Lemma \ref{lemma:boundofexpectation} and Lemma
  \ref{lemma:largedeviation} still hold if we replace $f_0$ by
  $\fat$. Without checking the proofs for those two lemmas, we can
  verify this point by the following ``change-of-measure'' type
  argument: given the tree topology, the associated parameter space, which is $\mathcal{F}_{\mathcal{A}_t}$, can be
  viewed as a multinomial model. The sufficient
  statistic for this model is the vector of counts $(\sum_{i=1}^n
  \mathbbm{1}_{X_i \in \Omega_j})_{1 \leq j \leq t}$. The distribution
of the random vector remains the same if we view the data as they were generated
from $\fat$. This implies that, in Lemma \ref{lemma:boundofexpectation} and Lemma
  \ref{lemma:largedeviation}, the probability of the events involved
  in the proof is the same as that measured under
  $\fat$. Therefore, we may replace $f_0$ by $\fat$ in these two
  lemmas and the results still hold.
\begin{proof}
In the proof, Lemma \ref{lemma:largedeviation} is applied iteratively to
obtain a ``sharp'' bound of the likelihood ratio. Let
$\delta_k(n, t) = \left(\frac{(t \vee t_{\min} )\log (t \vee t_{\min})}{n
    } \right)^{\omega_k}$, where $\omega_k$ is a
sequence defined by 
\begin{eqnarray}
  \label{eq:omega}
 \omega_1= 1/4 \quad \mbox{and} \quad \omega_{k+1} = \frac{1}{2}
\omega_k + \frac{1}{4}.
\end{eqnarray}
Define $\C_1:= \{ f \in \mathcal{F}_\at: \rho (f, \fat) \geq
\eta \delta_1(n,
t) \}$ and $\C_k := \{ f \in \mathcal{F}_\at: \eta \delta_k(n,t) \leq
\rho (f, \fat) < \eta \delta_{k-1}(n,t) \}$ for $k\geq 2$.

When $k=1$, applying the lower truncation of the log-likelihood ratio,
we have 
\begin{eqnarray}
\label{eq:toLowerTrunc}
&&\P_{f_0}^n \left( \sup_{f \in \C_1} \sum_{i=1}^n
Z_{\at, f} (Y_i) \geq  - \frac{1}{4}\eta^2 n \delta^2(n, t) \right)\\
\nonumber && \hspace{1in} \leq 
\P_{f_0}^n \left( \sup_{f \in \C_1} \sum_{i=1}^n
\widetilde{Z}_{\at, f} (Y_i) \geq  - \frac{1}{4}\eta^2 n \delta^2(n, t) \right).
\end{eqnarray}
Note that the modified version of Lemma \ref{lemma:boundofexpectation}
can be applied to $\widetilde{Z}_{\at, f}$, which implies that 
\begin{eqnarray*}
  \E  \widetilde{Z}_{\at, f} (Y_1 ) \leq - (1
  -\gamma) \rho^2( f, \fat), \quad \mbox{and} \quad
  \sup_{f \in \mathcal{C}_1}
  \left( \E  \widetilde{Z}_{\at, f} (Y_1 ) \right)
  \leq -(1-\gamma) \eta^2 \delta_1^2 (n, t).
\end{eqnarray*}
We can also choose the truncation parameter $\tau$ such that $\gamma$
in Lemma \ref{lemma:boundofexpectation} satisfies $\gamma \leq 2/3.$ 
Together with the previous upper bound (\ref{eq:toLowerTrunc}), we obtain 
{\small
\begin{eqnarray}
\label{eq:toEP}
   && \hspace{0.5cm}\P_{f_0}^n \left( \sup_{f \in \C_1} \sum_{i=1}^n
Z_{\at, f} (Y_i) \geq - \frac{1}{4} \eta^2 n \delta^2(n, t) \right) \\
\nonumber  && \hspace{0.5cm} \leq  
\P_{f_0}^n \left( \sup_{f \in \C_1} 
\sum_{i=1}^n 
\left(\widetilde{Z}_{\at, f} (Y_i) 
- \mathbb{E} \widetilde{Z}_{\at, f} (Y_i ) \right) 
\geq
- \sum_{i=1}^n 
\sup_{f \in \C_1} \left(\mathbb{E} \widetilde{Z}_{\at, f} (Y_i)
\right) -\frac{1}{4} n \eta^2 \delta^2(n, t) \right)\\
\nonumber && \hspace{0.5cm} \leq
\P_{f_0}^n \left( \sup_{f \in \C_1} n^{-1/2}
\sum_{i=1}^n \left(\widetilde{Z}_{\at, f} (Y_i ) - \mathbb{E}
   \widetilde{Z}_{\at, f} (Y_i ) \right) \geq
\frac{1}{4} n^{1/2} (1-\gamma ) \eta^2 \delta_1^2(n, t) \right).
\end{eqnarray}
}
The last inequality is due to $\delta_1(n, t) \geq \delta(n,t)$ and
$\gamma \leq 2/3$. In order to apply Lemma \ref{lemma:largedeviation} to bound the
probability on the right hand side of inequality (\ref{eq:toEP}), we
need to check the conditions
(\ref{eq:condition1}) and (\ref{eq:condition2}). Note that $ \| f^{1/2} - f_\at ^{1/2} \|_2 ^2 = 2-2 \int_{\Omega}
\sqrt{f(y) f_\at (y) } dy \leq 2$. We may set $\epsilon$ in Lemma
\ref{lemma:largedeviation} to be $\sqrt{2}$ here. When $ t \leq n/\log
n$, 
$\delta_1 (n, t) <1$. 
$L$ in Lemma \ref{lemma:largedeviation} is
$\frac{1}{4} n^{1/2}
(1-\gamma) \eta^2 \delta_1^2(n,t)$ here. First (\ref{eq:condition1}) is satisfied, since
\begin{equation*}
  L = \frac{1}{4} n^{1/2} (1-\gamma) \eta^2 \delta_1^2(n,t) = o\left(b n^{1/2}
 \epsilon^2/4 \right). 
\end{equation*}
For the function class $\mathcal{F}_\at$, we can use the bound for the bracketing
entropy in Lemma \ref{lemma:boundofH}. Let $\epsilon' =\frac{b
  (1-\gamma)}{32} \eta^2 \delta_1^2(n,t)$, and $c = -1/ \log
\epsilon'$.  
\begin{eqnarray*}
 \int_{\frac{b (1-\gamma)}{32} \eta^2 \delta_1^2(n,t)}^d \left( \bentro ( u/(2 \exp{\tau/2}),
  \mathcal{F}_\at, \rho) \right)^{1/2} du  
\lesssim (c+ \log c / \log \epsilon' )\sqrt{t \log t_{\min}}.
\end{eqnarray*}
Therefore, (\ref{eq:condition2}) is also satisfied when $n$ and $t$
are large enough. With Lemma \ref{lemma:largedeviation}, we have
\begin{eqnarray*}
 \nonumber \P_{f_0}^n \left( \sup_{f \in \C_1} \sum_{i=1}^n
Z_{\at, f} (Y_i) \geq  - \frac{1}{4} \eta^2 n \delta^2(n, t) \right)
&\leq& 3\exp
  \left(-(1-b) \varphi(L, \epsilon^2, n) \right) \\
&\leq&3 \exp\left( -\frac{(1-b) (1-\gamma)^2 \eta^2}{16(2^8 c_0 + 8)}
       \cdot  n \delta^2(n, t)\right).
\end{eqnarray*}

When $k\geq 2$, we still want to apply Lemma
\ref{lemma:largedeviation} to establish the inequality 
\begin{eqnarray}
  \label{eq:largedevCm}
   & &\P_{f_0}^n \left( \sup_{f \in \C_k} \sum_{i=1}^n
Z_{\at, f} (Y_i) \geq  - \frac{1}{4} \eta^2 n \delta^2(n, t) \right) \\
\nonumber&& \hspace{0.5in} \leq  3 \exp\left( -\frac{(1-b) (1-\gamma)^2 \eta^2}{16(2^8 c_0 + 8)}
       \cdot  n \delta^2(n, t) \right).
\end{eqnarray}
To do so, at the first step, by a similar argument as that for
inequalities (\ref{eq:toLowerTrunc}) and
(\ref{eq:toEP}), we have 
\begin{eqnarray*}
 \nonumber & &\P_{f_0}^n \left( \sup_{f \in \C_k} \sum_{i=1}^n
Z_{\at, f} (Y_i) \geq  - \frac{1}{4} \eta^2 n \delta^2(n, t) \right)  \\
&&\hspace{0.8cm} \leq \P_{f_0}^n \left( \sup_{f \in \C_k} n^{-1/2} \sum_{i=1}^n
\left(\widetilde{Z}_{\at, f}(Y_i) - \mathbb{E} \widetilde{Z}_{\at, f} (Y_j) \right) \geq
\frac{1}{4} n^{1/2}(1-\gamma) \eta^2 \delta_k^2(n,t) \right).
\end{eqnarray*}
Now $L$ in Lemma \ref{lemma:largedeviation} is $\frac{1}{4} n^{1/2}(1-\gamma) \eta^2
\delta_k^2(n,t)$ here and $\epsilon$ can be set as $\eta \delta_{k-1}(n,t)$. By our assumptions, $\delta_k(n,t) = o(\delta_{k-1}(n,t))$ for $k\geq
2$. Then condition (\ref{eq:condition1}) is satisfied.
Based on the same type of calculation,
\begin{eqnarray*}
  \int_{\frac{b (1-\gamma)}{128} \eta^2
    \delta_k^2(n,t)}^{\eta \delta_{k-1}(n,t)} \left( \bentro ( u/(2 \exp{\tau/2}),
  \mathcal{F}_\at, \rho) \right)^{1/2} du 
\lesssim t^{1/2} \delta_{k-1}(n,t) \sqrt{\log \frac{8}{\delta_{k-1}^2(n,t)}}.
\end{eqnarray*}
Plugging in $\delta_{k-1}(n,t) = (\frac{(t \vee t_{\min})
  \log (t \vee t_{\min})}{n})^{\omega_{k-1}}$, after some calculation, we can
verify that condition (\ref{eq:condition2}) is also satisfied locally when
$n$ and $t$ are large enough. 
Note that 
\begin{eqnarray*}
 & & \hspace{-1.5cm} \varphi(\frac{1}{4}n^{1/2}(1-\gamma) \eta^2
\delta_k^2(n,t), (\eta \delta_{k-1}(n,t))^2, n) \\
&\geq& \frac{n(1-\gamma)^2 \eta^2
\delta_k^4(n,t)/16} {8(8c_0 \delta_{k-1}^2(n,t) +(1-\gamma)
\delta^2_{k-1}(n,t))} \\
&\geq& \frac{(1-\gamma)^2 \eta^2 }{16(2^8 c_0 + 8)} \cdot n
       \delta^2(n, t).
\end{eqnarray*}
Applying Lemma \ref{lemma:largedeviation}, we get the inequality
(\ref{eq:largedevCm}).

Based on (\ref{eq:omega}), it is easy to see $\omega_k = \frac{1}{2}-
\left( \frac{1}{2} \right)^{k+1}$. 
If we set $K(n,t)$ to be $\frac{\log \log \frac{1}{\delta(n,t)} - \log
  \log 2}{\log 2}$, then $ \eta \delta(n, t) \leq \eta \delta_{K(n,t)} \leq 2 \eta
\delta(n,t)$. It implies that 
\begin{eqnarray*}
 & &\P_{f_0}^n \left( \sup_{f \in \mathcal{F}_{\mathcal{A}_t}: \rho(f, \fat) \geq 2 \eta \delta(n,t) } \prod_{i=1}^n
    \frac{f(Y_i)}{\fat(Y_i)} \geq \exp( -\frac{1}{4} \eta^2 n \delta^2(n, t))
    \right)\\
&&\hspace{1cm} \leq
\sum_{k=1}^{K(n,t)} \P_{f_0}^n \left( \sup_{f \in \C_k} \sum_{i=1}^n
Z_{\at, f} (Y_i) \geq  - \frac{1}{4} \eta^2 n \delta^2(n, t) \right) \\
&&\hspace{1cm} \leq 
3 K(n,t) \exp\left( -\frac{(1-b) (1-\gamma)^2
   }{16(2^8 c_0 + 8)} \cdot 
 \eta^2 n \delta^2(n, t) \right).
\end{eqnarray*}
Therefore, for some $\eta>1$ and appropriately chosen $c_2 =
(1-b)(1-\gamma^2)/(2^8 c_0 +8)$, the
inequality (\ref{eq:var_deviation_lik}) holds.
\end{proof}

With a similar argument, we can also obtain the following weaker
bound for the tail of the likelihood ratio. The proof is omitted, as it
follows exactly the same steps as those for Proposition 
\ref{thm:var_deviation}. 
\begin{corollary}
   \label{cor:var_deviation_weaker}
Under the same assumptions as Theorem \ref{thm:likbiasvar}, for
$\eta>1$ and some constants $c_2 \in (0,1)$ (same as the one in
Proposition \ref{thm:var_deviation}), 
\begin{eqnarray*}
  \P_{f_0}^n \left( 
\sup_{f \in  \Fat, \rho(f, \fat) \geq
             \eta \delta(n,t)}
\prod_{i=1}^n \frac{f(Y_i)}{\fat (Y_i)} \geq 
  \exp\left( 
         n \eta^2 \delta^2(n,t) \right) \right) 
\leq 
4\exp\left(-c_2\eta^2 n \delta^2(n, t) \right).
\end{eqnarray*}
\end{corollary}

Another useful corollary is an inequality for likelihood ratios
indexed by density functions defined on a specific partition or tree
topology. This corollary is applied to bound $P_2$ in the proof
of Theorem \ref{thm:likbiasvar}.  
\begin{corollary}
  \label{cor:lik_deviation_fat}
Under the same assumptions as Theorem \ref{thm:likbiasvar}, for
$\eta>1$ and some constants $c_2 \in (0,1)$ (same as the one in
Proposition \ref{thm:var_deviation}), 
\begin{eqnarray*}
  \P_{f_0}^n \left( 
\sup_{f \in  \Fat}
\prod_{i=1}^n \frac{f(Y_i)}{\fat (Y_i)} \geq 
 \exp\left( 
         n \eta^2 \delta^2(n,t) \right) \right) 
\leq 
  8\exp\left(-c_2\eta^2 n \delta^2(n, t) \right).
\end{eqnarray*}
\end{corollary}
\begin{proof}
  To show the inequality, we first make the following decomposition,
\begin{eqnarray}
 \nonumber  && \P_{f_0}^n \left( 
\sup_{f \in  \Fat}
\prod_{i=1}^n \frac{f(Y_i)}{\fat (Y_i)} \geq 
               \exp\left( 
  n \eta^2 \delta^2(n,t) \right) \right)\\
  && \hspace{1cm} \leq  \label{eq:P1part1}  \P_{f_0}^n \left( 
\sup_{f \in  \Fat, \rho(f, \fat) \geq
             \eta \delta(n,t)}
\prod_{i=1}^n \frac{f(Y_i)}{\fat (Y_i)} \geq 
\exp\left( 
          n \eta^2 \delta^2(n,t) \right) \right) \\
&& \hspace{1.5cm} + \label{eq:P1part2}
\P_{f_0}^n \left( 
\sup_{f \in  \Fat, \rho(f, \fat)
   <
              \eta \delta(n,t)}
\prod_{i=1}^n \frac{f(Y_i)}{\fat (Y_i)} \geq 
\exp\left( 
          n \eta^2 \delta^2(n,t) \right) \right). 
\end{eqnarray}

For the probability in (\ref{eq:P1part1}), in Corollary 
\ref{cor:var_deviation_weaker} we have shown that it is bounded by $4 \exp(- c_2 n \eta^2 \delta^2(n, t))$.

For the probability of the other event in (\ref{eq:P1part2}), we can apply Lemma \ref{lemma:largedeviation} to obtain an
upper bound. In specific, let $\mathcal{B} = \{f \in  \Fat:
  \rho(f, \fat) < \eta \delta(n,t)\}$, 
{\small
\begin{eqnarray*}
 && \P_{f_0}^n \left(
  \sup_{f \in  \mathcal{B}}
  \sum_{i=1}^n Z_{\mathcal{A}_t, f} (Y_i)
  \geq 
    n \eta^2 \delta^2(n, t)
    \right) \\
  &&\hspace{0.2in}
     \leq
     \P_{f_0}^n \left(
  \sup_{f \in  \mathcal{B}}
  \sum_{i=1}^n \wt{Z}_{\mathcal{A}_t, f} (Y_i)
  \geq 
     n \eta^2 \delta^2(n, t)
     \right) \\
  && \hspace{0.2in}
     \leq
     \P_{f_0}^n \left(
     \sup_{f \in  \mathcal{B}}
     n^{-1/2} \sum_{i=1}^n
     \left(\wt{Z}_{\mathcal{A}_t, f} (Y_i) -
     \E \wt{Z}_{\mathcal{A}_t, f} (Y_i) \right)
     \geq
     -\sup_{f \in  \mathcal{B}}
     n^{1/2}\E \wt{Z}_{\mathcal{A}_t, f} (Y_i) + n^{1/2} \eta^2 \delta^2(n, t)
     \right). 
\end{eqnarray*}
}
To apply Lemma \ref{lemma:largedeviation}, $\epsilon$ in the lemma is $ \eta
\delta(n,t)$ here 
and $L$ is $n^{1/2} \eta^2 \delta^2(n, t)$. By using the same argument
as that in the proof of Proposition \ref{thm:var_deviation}, we can
check both conditions \eqref*{eq:condition1} and
\eqref*{eq:condition2} are satisfied. 
Note that 
\begin{eqnarray*}
 \varphi( L, \epsilon^2  , n) 
&=& \frac{n ( \eta^2 \delta^2(n, t) )^2} {8(8c_0 (
    \eta \delta (n,
    t))^2 + \eta^2 \delta^2(n, t) )}\\
&=& \frac{n ( \eta^2 \delta^2(n, t) )^2 } 
{(2^6 c_0 +8 ) \eta^2 \delta(n,t) ^2} \\
&\geq& \frac{\eta^2 }{2^6 c_0 + 8} n \delta^2(n, t).
\end{eqnarray*}
By Lemma \ref{lemma:largedeviation}, we have
\begin{eqnarray*}
  \eqref*{eq:P1part2} \leq 4 \exp ( - c_2 \eta^2 n \delta^2(n, t)).
\end{eqnarray*}

Combining the upper bounds for \eqref*{eq:P1part1} and
\eqref*{eq:P1part2}, we get the desired result.
\end{proof}

\subsection{A Large Deviation Inequality as the Lower Bound}
\label{subsec:concentrationlower}
In order to show uniform convergence
of all tree-based weak leaners within a forest, we need a sharper lower bound
for the likelihood ratio. In previous works on posterior contraction
of Bayesian methods \citep{ghosal00,shen2001}, such a lower bound is usually obtained
by applying the Chebyshev's inequality. In this section, we provide an improved result by using the Bernstein's inequality when
$f_0$ is bounded away from 0. 
Let 
\begin{eqnarray*}
\wt{\mathcal{B}}_\infty= \left \{f \in  \Fat:
  \| f - \fat \|_\infty \leq \eta^2 \delta^2 (n, t), f \geq t
  (\eta^2/(2n))^{2 \lceil \log_2 t \rceil} \right \}.
\end{eqnarray*}
This is a ball around $\fat$. The definition of
$\wt{\mathcal{B}}_\infty$ ensures that any $f \in \wt{\mathcal{B}}_\infty$ is bounded from below
and consequently $\kl (\fat, f)$ is well
defined. We want to show that when $f_0$ is bounded away from zero,
the likelihood ratio for density functions lying in this ball can not
be too small either. Usually, we set $\A_t$
    as a non-stopped tree topology. 
  \begin{lemma}
    \label{lemma:empineqinf}
Under the same conditions as those for Theorem
  \ref{thm:likbiasvar}, if we further assume that $f_0$ is bounded away from
  zero, then
\begin{eqnarray*}
  \P_{f_0}^n \left(
  \inf_{f \in  \wt{\mathcal{B}}_\infty}
  \sum_{i=1}^n Z_{\mathcal{A}_t, f} (Y_i)
  \leq 
   - 4n \eta^2 \delta^2(n, t) - n\sup_{f \in \wt{\mathcal{B}}_\infty} \kl(\fat, f)
\right) 
\leq 3 \exp(-4 c_2 \eta^2 n \delta^2(n, t)).
\end{eqnarray*}
  \end{lemma}
  \begin{proof}
When $f_0$ is bounded away from zero, so is $\fat$. Note that 
{\small
\begin{eqnarray*}
 && \P_{f_0}^n \left(
  \inf_{f \in  \wt{\mathcal{B}}_\infty}
  \sum_{i=1}^n Z_{\mathcal{A}_t, f} (Y_i)
  \leq 
   - 4 n \eta^2 \delta^2(n, t) -  n^{1/2}\sup_{f \in \wt{\mathcal{B}}_\infty} \kl(\fat, f)
    \right) \\
  && \hspace{0.2in}
     \leq
     \P_{f_0}^n \bigg(
     \inf_{f \in  \wt{\mathcal{B}}_\infty}
     n^{-1/2} \sum_{i=1}^n
     \left(Z_{\mathcal{A}_t, f} (Y_i) -
     \E Z_{\mathcal{A}_t, f} (Y_i) \right)
     \leq
     - \inf_{f \in \wt{\mathcal{B}}_\infty} n^{1/2} \E Z_{\mathcal{A}_t, f} (Y_1) \\
     && \hspace{2.5in}
    - 4 n^{1/2} \eta^2 \delta^2(n, t) -n^{1/2}\sup_{f \in \wt{\mathcal{B}}_\infty} \kl(\fat, f)
     \bigg),
\end{eqnarray*}
}
where $\E Z_{\mathcal{A}_t, f} (Y_1) = - \kl (\fat, f)$ for any $f \in
\Fat$. $n^{1/2}\sup_{f \in \wt{\mathcal{B}}_\infty}
  \kl(\fat, f)$ cancels with the positive term $ - \inf_{f \in
  \wt{\mathcal{B}}_\infty} n^{1/2} \E Z_{\mathcal{A}_t, f} (Y_1) $. 
It suffices to bound the probability
 \begin{eqnarray*}
   && \P_{f_0}^n \left(
     \inf_{f \in  \wt{\mathcal{B}}_\infty}
     n^{-1/2} \sum_{i=1}^n
     \left(Z_{\mathcal{A}_t, f} (Y_i) -
     \E Z_{\mathcal{A}_t, f} (Y_i) \right)
     \leq
    -4 n^{1/2} \eta^2 \delta^2(n, t) 
      \right) \\
   && \hspace{0.5cm} =
      \P_{f_0}^n \left(
     \sup_{f \in  \wt{\mathcal{B}}_\infty}
     n^{-1/2} \sum_{i=1}^n
     \left(- Z_{\mathcal{A}_t, f} (Y_i) +
     \E Z_{\mathcal{A}_t, f} (Y_i) \right)
     \geq 
     4 n^{1/2} \eta^2 \delta^2(n, t) 
      \right).
 \end{eqnarray*}
 We can alternatively study the empirical process indexed by 
 $-Z_{\mathcal{A}_t, f}$. As we have assumed that the density function
 is bounded away from zero, on the set $\wt{\mathcal{B}}_\infty$
 log-likelihood ratio $Z_{\mathcal{A}_t, f}$ is bounded from above and
below. The two empirical processes indexed by either
$Z_{\mathcal{A}_t, f}$ or $-Z_{\mathcal{A}_t, f}$ can be
 studied by using similar techniques. In specific, Lemma
 \ref{lemma:largedeviation} can still be applied to the one associated
 with $-Z_{\mathcal{A}_t, f}$. Note that by Lemma
   \ref{lemma:suptoH} for any $f \in \wt{\mathcal{B}}_\infty$, $\rho(f,
 \fat) \leq 2 \eta \delta(n, t)$. Consequently, the probability above is bounded by $3 \exp(-4 c_2
 \eta^2 n\delta^2(n, t))$. This finishes the proof.
\end{proof}


\begin{remark}
\label{remark:empineqinf}
When $f_0$ is bounded from below and $\| f
  - \fat \|_\infty \leq \epsilon^2$ for some $\epsilon^2 \leq
  \frac{1}{2} \min\{
  \inf_{x \in \Omega} f_0(x), 1\}$,  we may claim that $1/2 \leq  f/\fat \leq 3/2$. By
  Lemma \ref{lemma:kltoH}, $\kl (\fat, f) \leq (3 \log 2)\rho^2(\fat,
  f)$. And by Lemma \ref{lemma:suptoH}, $\rho(\fat, f) \leq 2
  \epsilon$. Combining the two bounds together we have for any $f \in
  \wt{\mathcal{B}}_\infty$, $\kl( \fat, f) \leq 12 \eta^2 \delta^2(n, t)$. 
  Therefore, Lemma \ref{lemma:empineqinf} implies that, 
  \begin{eqnarray}
  \label{eq:empineqinfalt}
  \P_{f_0}^n \left(
  \inf_{f \in  \wt{\mathcal{B}}_\infty}
  \sum_{i=1}^n Z_{\mathcal{A}_t, f} (Y_i)
  \leq 
   - 16 n \eta^2 \delta^2(n, t) 
\right) 
\leq 3 \exp(-4 c_2 \eta^2 n \delta^2(n, t)).
\end{eqnarray}
\end{remark}

By using Lemma \ref{lemma:empineqinf}, we can derive an inequality as a lower bound for the likelihood ratios. 
\begin{theorem}
\label{lemma:likineqinf}
  Besides the same set of conditions imposed in Theorem
  \ref{thm:likbiasvar}, we further assume that $f_0$ is bounded
 from above and below. Then for some $\underline{c}_1 \in (0, 1), \underline{c}_2'
>0$ to be specified in the proof, and the same $c_2$ as that in Proposition \ref{thm:var_deviation}, 
\begin{eqnarray}
    \label{eq:likineqinf}
 &&  \hspace{0.2in} \P_{f_0}^n \bigg( 
\inf_{f \in \wt{\mathcal{B}}_\infty} 
\prod_{i=1}^n \frac{f(Y_i)}{f_0(Y_i)} 
\leq  
\exp \Big( 
-n (1+\underline{c}_1) \kl (f_0, \fat)
-4\eta^2 n \delta^2(n,t) \bigg)\\
  \nonumber && \hspace{0.6in}
- n \sup_{f \in \wt{\mathcal{B}}_\infty} \kl (\fat, f) 
\Big)
      \leq
\exp\left( - \underline{c}_2' n  \kl (f_0, \fat) \right)
+
    3\exp\left(-4 c_2 \eta^2 n \delta(n, t)^2 \right).
  \end{eqnarray}
\end{theorem}

\begin{proof}
The proof is similar to that for Theorem \ref{thm:likbiasvar}. 
We first decompose the left hand side of \eqref*{eq:likineqinf}
into the following two parts:  
{\small
\begin{eqnarray*}
 && \P_{f_0}^n \left( \inf_{f \in \wt{\mathcal{B}}_\infty}
\prod_{i=1}^n \frac{f(Y_i)}{f_0(Y_i)} \leq  
\exp\Big( -n (1+\underline{c}_1) \kl (f_0, \fat)
-4\eta^2 n \delta^2(n,t) 
- n \sup_{f \in \wt{\mathcal{B}}_\infty} \kl (\fat, f) 
 \Big)\right)  \\
 && \hspace{1in}
 \leq P_1 + P_2, 
\end{eqnarray*}
}
where
\begin{eqnarray*}
P_1 &=& \P_{f_0}^n \left( 
\prod_{i=1}^n \frac{\fat(Y_i)}{f_0(Y_i)}  \leq 
\exp\left( -n (1+\underline{c}_1) \kl (f_0, \fat)\right)\right), \\
  P_2 &=& \P_{f_0}^n \left( \inf_{f \in  \wt{\mathcal{B}}_\infty}
\prod_{i=1}^n \frac{f(Y_i)}{\fat (Y_i)} \leq  
\exp \Big( 
          -4\eta^2 n \delta^2(n,t) 
- n \sup_{f \in \wt{\mathcal{B}}_\infty} \kl (\fat, f)   \Big)\right).
\end{eqnarray*}
Similar to the proof for Theorem \ref{thm:likbiasvar}, we apply the
Bernstein's inequality to bound $P_1$, and apply Lemma
\ref{lemma:empineqinf} to bound $P_2$. As $f_{\mathcal{A}_t}$ is
  bounded from above and below, the function
  $f_{\mathcal{A}_t}/f_0$ is bounded from below. Therefore, we apply
  the Bernstein inequality with moment conditions on $\log (f_{\mathcal{A}_t}/f_0)$. The constant $\underline{c}_2'$ can be set as $\frac{(\underline{c}_1)^2}{2^6 c_0 + 2^4 \underline{c}_1}$, where $c_0 = (\exp(\tau/2) -1 -
\tau/2)/(1- \exp(-\tau/2))^2$ and $\tau = \inf_{y \in \Omega}
\log(f_0/\fat)(y)$.
\end{proof}

\section{Posterior Concentration Results under the OPT Prior}
\label{sec:optpost}
In this section, we study
posterior concentration properties under variations of the OPT
prior. Specifically, we study shape-constrained OPTs in Section~\ref{subsec:optpostshape} and early-stopped OPTs in Section~\ref{subsec:optposttree}. The results apply in a more general setting. Proposition~\ref{thm:posteriorrateshape} shows that, given a fixed tree topology, the posterior distribution concentrates around the best approximation supported by that tree structure. Proposition~\ref{thm:posteriorratetree} complements this result by characterizing the subset of tree topologies around which the posterior distribution concentrates when the OPT is early stopped. \rev{The proofs rely on the likelihood
  (ratio) bounds derived in Section~\ref{sec:riskanalysis}.} Before we study posterior concentration
properties, we provide a detailed description of how to generate
a recursive partition of $\Omega$ under a shift, and how to prune a
tree, especially at the boundary region, to avoid obtaining a leaf
node of very small volume. \rev{We also summarize several properties of the
OPT prior in Section~\ref{sec:optanalysis}.} 
\subsection{More Details for Growing a Tree Under a Shift}
\label{subsec:shiftedtree}
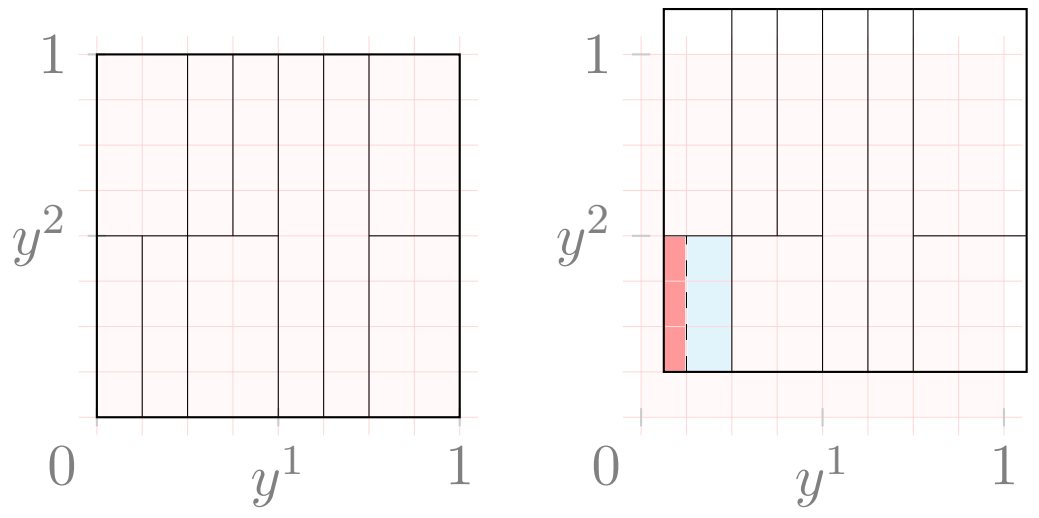
\begin{figure}[t]
\vspace{-5mm}
\begin{center}
\vspace{5mm}
\begin{tikzpicture}[xscale=0.6, yscale=0.6]
  \filldraw[color=pink!0, fill= pink!10] (1,1) rectangle (5,5);
  \draw[step=.5cm,pink!60,very thin] (0.8,0.8) grid (5.2,5.2);
\draw[gray!40] ( 1,1+0.1) -- (1 ,1-0.1) node[anchor=north east,gray] {0};
\draw[gray!40] ( 5,1+0.1) -- (5 ,1-0.1) node[anchor=north,gray] {1};
\draw[gray!40] ( 3,1+0.1) -- (3 ,1-0.1) node[anchor=north,gray] {$y^1$};
\draw[gray!40] (1+0.1,5) -- (1-0.1,5) node[anchor=east,gray] {1};
\draw[gray!40] (1+0.1,3) -- (1-0.1,3) node[anchor=east,gray]{$y^2$};

  \draw[black] (1,1) rectangle (5,5);

\draw[black, very thin] (1, 3) -- ( 3, 3);
\draw[black, very thin] (4, 3) -- ( 5, 3);

\draw[black, very thin] (3, 1) -- (3,5);

\draw[black, very thin] (2, 1) -- (2,5);
\draw[black, very thin] (4, 1) -- (4,5);

\draw[black, very thin] (3.5, 1) -- (3.5,5);

\draw[black, very thin] (1.5, 1) -- (1.5,3);
\draw[black, very thin] (2.5, 3) -- (2.5,5);

\filldraw[color=pink!0, fill= pink!10] (7,1) rectangle (11,5);

\filldraw[color=red!0, fill= red!40] (7+0.25,1+0.5) rectangle
(7.5,3);
\filldraw[color=red!0, fill= cyan!10] (7.5,1+0.5) rectangle (8,3);
\draw[step=.5cm,pink!60,very thin] (6.8,0.8) grid (11.2,5.2);
\draw[gray!40] ( 7,1+0.1) -- (7 ,1-0.1) node[anchor=north east,gray] {0};
\draw[gray!40] ( 11,1+0.1) -- (11 ,1-0.1) node[anchor=north,gray] {1};
\draw[gray!40] ( 9,1+0.1) -- (9 ,1-0.1) node[anchor=north,gray] {$y^1$};
\draw[gray!40] (7+0.1,5) -- (7-0.1,5) node[anchor=east,gray] {1};
\draw[gray!40] (7+0.1,3) -- (7-0.1,3) node[anchor=east,gray] {$y^2$};

\draw[black] (7+0.25,1+0.5) rectangle (11+0.25,5+0.5);

\draw[black, very thin] (7+0.25, 3) -- ( 9, 3);
\draw[black, very thin] (10, 3) -- ( 11+0.25, 3);

\draw[black, very thin] (9, 1+0.5) -- (9,5+0.5);

\draw[black, very thin] (8, 1+0.5) -- (8,5+0.5);
\draw[black, very thin] (10, 1+0.5) -- (10,5+0.5);

\draw[black, very thin] (9.5, 1+0.5) -- (9.5,5+0.5);

\draw[black,dashed, very thin] (7.5, 1+0.5) -- (7.5,3);
\draw[black, very thin] (8.5, 3) -- (8.5,5+0.5);

\end{tikzpicture}
\end{center}
\vspace{-5mm}
\caption{{\small An example of shifted tree in the two-dimensional case. Left: a depth-4 tree built on $[0,1]^2$. Right:
 a depth-4 tree built on the shifted sample space $\Omega^h$ with $h
 =(\frac{1}{16}, \frac{1}{8})$. The tree is obtained by making splits
 along the same point at every step as those for the
 tree in the left panel. The tree pruning scheme prevents a split like
 the one in the dashed line.}
}
\label{fig:shifted_partition_app}
\vspace{-3mm}
\end{figure}

We first provide more details for tree construction under a shift. The shifted sample space $\Omega^h$ is defined as $\Omega^h = \{ y: y-h \in \Omega \}$. Under the shift $h$, we grow a tree on $\Omega^h$ with only splits along {\em dyadic} points allowed. More specifically, if we choose to divide an elementary region $A =\otimes_{l=1}^d [a_l, b_l)$ along the $l$th dimension, we first identify the largest integer $\kappa$ such that there exists a unique number of the form $i \cdot 2^{-\kappa}$ satisfying $i \in \mathbb{N}$ and $a_l < i \cdot 2^{-\kappa} < b_l$. Then $A$ can be \textit{tentatively} divided into two child nodes $A_l$ and $A_r$ along the dimension $l$ as $A_l =\{y \in \mathbb{R}^d: y^l < i \cdot 2^{-\kappa} \} \cap A$ and $A_r = A \setminus A_l$. In the case when we grow a depth-constrained tree, after making the tentative splits, in order to avoid obtaining a leaf node of very small volume, especially near the boundary of $\Omega$, we check the volume of all leaf nodes before finalizing the tree structure. That is, when we grow a depth-$\D$ tree under the shift $h$, if any of the leaf nodes has a volume
smaller than $ \frac{1}{2} \cdot 2^{-\D}$, then the node will
be \emph{merged} with the neighboring leaf node, recovering their
parent node. In Figure \ref{fig:shifted_partition_app}, the left panel
shows a depth-4 dyadic tree built on the sample space $\Omega$, and
the right panel shows the shifted sample space along $h =
(\frac{1}{16},\frac{1}{8})$ and a depth-4 tree making splits along the
same points. Note that for the tree in the right panel, after the
initial tree growth, the red rectangular region in the bottom left corner is of
volume $\frac{3}{128}$, while the limit on leaf node volume for a
depth-4 tree is $\frac{1}{2} \cdot \frac{1}{16}$. Therefore, we merge
the red node with its neighboring blue node and keep their parent node
in the final tree structure. By shifting $\Omega^h$ back, we
obtain a partition of the original sample space $\Omega$.

When making Bayesian inference under the shifted OPT prior, we always
set hyper-parameters for the beta distribution as $a(A)= b(A) = a_0$,
$0 < a_0 < 1$ for all nodes under a shift. 
The merge operation in the tree punning step is assumed to happen with prior probability~1.

\subsection{Posterior Concentration Under the Shape Constrained OPT Prior}
\label{subsec:optpostshape}
First, we study given a tree topology,
how fast the posterior distribution contracts around
the ``best'' approximation supported by the constrained tree
topology. We first focus on the case when there is no perturbation. 
\begin{proposition}
  \label{thm:posteriorrateshape}
Assume $f_0 \in \mathcal{C}^{\bm
    \alpha}(\Omega, L)$, $\alpha_l >0$, and $f_0$ is bounded away from zero. 
  With a random sample $\Y_n = \{Y_1, \ldots, Y_n\}$ from $f_0$, we 
  make Bayesian inference under the OPT prior or OPT with shape constraints. 
  Assume $\A$ 
  is a tree topology of depth
  $\D \leq \D_{\max}$, 
  $\D_{\min}$ and $t_{\min}$ are as defined in
  \eqref*{eq:depthlowerbound}, and $\delta(n, t)$ is the quantity
  defined in \eqref*{eq:deltadef}. Let $t = 2^\D$ and
  $\eta \geq 1$ be some constant. 
  Then with probability at least
  $1 - 8\exp(- 4c_2 \eta^2 n \delta^2(n, t))$,  
  \begin{eqnarray*}
    \Pi \left(f: \rho (f, f_{\A}) \geq 24 \eta  \delta(n, t) 
\hspace{1mm}| \hspace{1mm}
    \Y_n, \mathcal{F}_\A \right)
    \leq
    \exp \left(- \eta^2 n\delta^2(n, t) \right).
\end{eqnarray*}
\end{proposition}
\begin{remark}
  The result also holds if we replace $\delta(n, t)$ by any 
  $\delta'(n, t) \geq \delta(n, t)$. 
The prior on unsupervised trees can be either OPT or OPT with any
form of shape constraints, as for these types of prior distributions,
conditional on the tree topology $\A$ the constrained prior $\Pi( \cdot | \mathcal{F}_\A)$ remains the same. 
\end{remark}
\begin{proof}
  Let $\mathcal{C} = \{ f \in \F_\A: \rho (f, f_{\A}) \geq 2 \eta' \delta(n, t)
  \}$. The posterior probability on $\mathcal{C}$ can be calculated as 
  \begin{eqnarray*}
    \Pi (\mathcal{C}  | \Y_n, \F_\A)
    = \frac{\int_{\mathcal{C}} \prod_{i=1}^n f(Y_i) d \Pi(f | {\F_\A})}
    {\int_{\mathcal{F}_{{\A}}}\prod_{i=1}^n f(Y_i) d \Pi(f
        | {\F_\A})}.
  \end{eqnarray*}
  Note that conditional on tree topology, the prior $\Pi(f | \F_\A)$ is
  simply a prior on probability masses on leaf nodes, which is
 induced by products of independent beta random variables under the OPT or shape constrained OPT prior. Dividing
  both the numerator and denominator by $\prod_{i=1}^n f_\A (Y_i)$,
  we
  have
  \begin{eqnarray*}
      \Pi (\mathcal{C}  | \Y_n, \F_{\A}) =
    \frac{\int_{\mathcal{C}} \prod_{i=1}^n \frac{f(Y_i)}{f_\A(Y_i)} d \Pi(f|\F_\A )}
    {\int_{\mathcal{F}_{{\A}}}\prod_{i=1}^n \frac{f(Y_i)}{f_\A (Y_i)} d \Pi(f|\F_\A)}.
  \end{eqnarray*}
  
  To bound the denominator, assume that $\A'$ is a non-stopped depth-$\D$ tree topology, and it is obtained by making refinement of $\A$. 
  Then the size of
  $\A'$ is $t = 2^{\D} $. 
  Define a ball abound $f_\A$ as
  \begin{eqnarray*}
    \mathcal{B} = \{ f \in \mathcal{F}_\A: \|f - f_{\A} \|_\infty \leq \eta^2  \delta^2(n, t), f \geq (\eta^2/(2n))^\D 2^\D
    \} \mbox{ for some } \eta \geq 1.
  \end{eqnarray*}
 As $\A'$ is a refinement of $\A$, both $f \in \mathcal{F}_\A$ and
 $f_\A$ can be viewed as a density supported by $\A'$, and
 \begin{eqnarray*}
  \mathcal{B} =  \mathcal{B}_\D = \{ f \in \mathcal{F}_{\A'}: \|f - f_{\A} \|_\infty \leq \eta^2  \delta^2(n, t), f \geq (\eta^2/(2n))^\D 2^\D
    \}.
 \end{eqnarray*}
 As making refinement only lowers the prior on $\mathcal{B}_\D$, by Lemma \ref{lemma:prioronball} from Section \ref{sec:optanalysis}, 
  \begin{eqnarray*}
    \Pi(\mathcal{B} |\F_\A) \geq \Pi(\mathcal{B}_\D | \F_{\A'})
    \geq
    \exp \left( -t \log n + t \log
    \left(
    \eta^2 /2 
    \right) 
    - t \log (2\Gamma(a_0)/ \Gamma(2a_0) ) \right).
  \end{eqnarray*}
  Note that $f_0$ is bounded away from zero. Based on Remark \ref{remark:empineqinf}, a lower bound for the likelihood ratio is 
  \begin{eqnarray*}
   \P_{f_0}^n \left( 
     \inf_{f \in \mathcal{B}_\D} \prod_{i=1}^n \frac{f(Y_i)}{f_\A (Y_i)}
    \leq
    \exp\left( - 16 \eta ^2 n \delta^2(n, t)\right)
    \right)
    \leq
    3\exp \left( 
    -4 c_2 \eta^2 n \delta^2(n, t)
    \right).
  \end{eqnarray*}

%

  For the numerator, the prior probability is upper bounded by 1. The
  likelihood ratio can also be bounded from above by
  applying the large deviation inequality from Proposition \ref{thm:var_deviation}.
  With probability at least $1- 4 \exp\left(-\frac{1}{16}c_2(\eta')^2 n
    \delta^2(n, t) \right)$,
  \begin{eqnarray*}
    \sup_{f \in  \F_\A, \rho(f, f_\A) \geq
            2 \eta' \delta(n, t)}
\prod_{i=1}^n \frac{f(Y_i)}{f_\A (Y_i)} \leq
    \exp\left( 
         - \frac{1}{4} n (\eta')^2 \delta^2(n, t) \right).
  \end{eqnarray*}
  
Combining
these bounds together, we obtain with probability at least $$1-
4 \exp\left(-\frac{1}{16}c_2(\eta')^2 n \delta^2(n, t) \right)
-3\exp \left( - 4c_2 \eta^2 n \delta^2(n, t)\right),$$
  \begin{eqnarray*}
    \Pi(\mathcal{C}| \Y_n,  \F_\A)
    \leq
       \scriptsize{
       \frac{
\exp\left( 
         - \frac{1}{4} n (\eta')^2 \delta^2(n, t) \right)
}
       {\exp \left( 
-t \log n
+t \log(\eta^2/2 )
- t \log (2\Gamma(a_0)/ \Gamma(2a_0) )
-16 \eta ^2 n \delta^2(n, t)
\right) }
       },
  \end{eqnarray*}
Therefore, when $\eta' \geq  12 \eta$,  and $\eta \geq 1$, $\Pi(\mathcal{C}|\Y_n,
       \F_\A) \leq \exp (-  n \eta ^2 \delta^2(n, t))$
       for 
        $n$ and $t$ are large enough.
      \end{proof}

      The following corollary shows that the posterior mean is also
      close the best tree approximation $f_\A$ with high probability. 
      \begin{corollary}
        \label{cor:posteriorratetreemean}
        Under the conditions of
       Proposition \ref{thm:posteriorrateshape},
        given a tree topology $\A$ of depth $\D$, with probability at least
    $1 - 8 \exp(- 4c_2 \eta^2 n \delta^2(n, t))$, for $n$ large enough,
  \begin{eqnarray*}
    \rho \left( f_\A, E(f | \Y_n, \F_\A) \right) \leq 25 \eta \delta(n, t).
  \end{eqnarray*}
\end{corollary}
\begin{proof}
  The proof is based on the convexity of $\rho^2(\cdot, \cdot)$. By
 the Jensen's inequality,
  \begin{eqnarray*}
    \rho^2( f_\A, E( f | \Y_n, \F_\A)) \leq
    E\left( \rho^2( f_\A, f) | \Y_n, \F_\A \right), 
  \end{eqnarray*}
  where $E( \cdot | \Y_n, \F_\A)$ is understood as the expectation with
  respect to the posterior distribution under shape constraint on
  $\F_\A$, which is $\Pi(\cdot | \Y_n, \F_\A)$. In Proposition 
  \ref{thm:posteriorrateshape}, we have shown that this posterior
  distribution concentrate around $f_\A$ with probability at least $1-
    8 \exp(- 4c_2 \eta^2 n \delta^2(n, t))$. And when the concentration result holds, let
  $\mathcal{C} = \{ f \in \F_\A: \rho(f_\A, f) \geq  24 \eta \delta(n,
  t) \}$, then 
  \begin{eqnarray*}
    E\left( \rho^2( f_\A, f) | \Y_n, \F_\A \right) &=&
    \int_{f \in \mathcal{C}^\complement}
    \rho^2(f_\A, f) d \Pi(f| \Y_n, \F_\A)
    +  \int_{f \in \mathcal{C}}
  \rho^2(f_\A, f) d \Pi(f| \Y_n, \F_\A) \\
 & \leq&
         (24 \eta)^2 \delta^2(n, t) \Pi(\mathcal{C}^\complement | \Y_n, \F_\A)
         + 2 \Pi(\mathcal{C} | \Y_n, \F_\A) \leq  (25 \eta)^2 \delta^2(n, t).
  \end{eqnarray*}
  This finishes the proof. 
\end{proof}

\begin{remark}
  \label{remark:posteriorrateshape}
In Proposition \ref{thm:posteriorrateshape} and Corollary
\ref{cor:posteriorratetreemean}, if we assume that $\A$ is a tree
topology grown under a shift $h$, with or without shape constraints,
both the posterior concentration rate and the convergence rate of the
posterior mean inflate by a constant. This is because to derive the lower bound of the likelihood ratio,
we need to define an $L^\infty$-ball $\mathcal{B}_\D$ of \emph{doubled} radius $2
\eta^2 \delta^2(n,t)$. 
\end{remark}

\subsection{Posterior Concentration Under the Early Stopped OPT Prior}
\label{subsec:optposttree}
In this subsection, we study posterior
contraction properties under the early stopped OPT prior. First, we briefly
discuss the exact meaning of an ``early stop'' when $f_0 \in
\mathcal{C}^{\bm \alpha} ( \Omega, L)$, $\alpha_l >0$. 
%
%
Over the anisotropic H\"older space, trees can only approximate functions with $\alpha_l \in (0, 1]$ well. Let $\alpha_l^* = \min \{1, \alpha_l \}$ and  $\left( \overline{\alpha^*}  \right)^{-1}= \frac{1}{d} \sum_{l=1}^d
  (\alpha_l^*)^{-1}$. 
The tree depth at which bias and variance are approximately of the
same order is 
  \begin{eqnarray*}
    \D_\tr =\left \lfloor \log_2 (n / \log n)^{\frac{d}{2
    \overline{\alpha^*}+d}} \right \rfloor.
  \end{eqnarray*}
If the prior distribution forces all trees to stop at a level $\D
  < \D_\tr$, then the
  bias term will be dominating, which is the case involved in the proof for
  the BMA scheme. The following proposition shows that in this case,
  the posterior distribution concentrates on those tree topologies
  that can provide better approximations to $f_0$. 
We study tree level posterior concentration result under the early stopped
OPT prior \emph{without} perturbation. 
\begin{proposition}
  \label{thm:posteriorratetree}
  Assume $f_0 \in \mathcal{C}^{\bm
    \alpha}(\Omega, L)$, $\alpha_l > 0$, and $\kappa_L \leq f_0 \leq \kappa_U$. 
  $\Y_n = \{Y_1, \ldots, Y_n\}$ is a random sample from $f_0$. Under
  the early stopped OPT prior at depth $\D$, with or without shape constraints, set 
$\delta (n, t) = \left( \frac{(t \vee t_{\min})
      \log (t \vee t_{\min}) }{n }\right)^{1/2} $ with $t = 2^\D$
  and $t_{\min} = 2^{\D_{\min}}$.
Assume $f_\D$ is the best depth-$\D$ approximation to $f_0$ under
the KL divergence. When $\D$ is at the order of $\log n$ and
$\kl(f_0, f_\D) \geq M_{\text{bias}} \eta^2 \delta^2(n, t)$ for some constant $M_{\text{bias}} >1$,  
  with probability at least
 $1 -
  16\exp(- \frac{1}{2} c_2'' \eta^2 n \delta^2 (n,t))$, for $n$
  large enough, 
 \begin{eqnarray*}
    \Pi_n \left(\mathcal{F}_{\A}: \kl (f_0, f_\A) \geq 
   \frac{1+ 2 \underline{c}_1}{1-\overline{c}_1}\kl (f_0, f_\D)
\cond
    \Y_n, \D \right)
    \leq
    \exp \left( - \frac{1}{2} n \underline{c}_1 
   \kl (f_0, f_\D)\right).
\end{eqnarray*}
$\overline{c}_1, \underline{c}_1 \in (0, 1)$ can be some small
constants, and $c_2'' = \min\{\overline{c}_2',
\underline{c}_2', c_2 \}$ with $\overline{c}_2'$, $\underline{c}_2'$, 
and $c_2$ being the same as
those in Theorem \ref{thm:likbiasvar} and Lemma
\ref{lemma:empineqinf}. A possible choice of $M_{\text{bias}}$ is $ 2\underline{c}_1^{-1}\eta^{-2} (2 \eta^2 + 4(\eta')^2 (1+ 3 \log(\kappa_U/\kappa_L)) + \nu + \log n/(\D \log 2))$.
\end{proposition}
The proof is similar to that for Proposition 
\ref{thm:posteriorrateshape}. We use the following notations in
the proof. Similar to the definition of $\omega_n(B; \D, h, \Y_n)$ in
Section \ref{sec:proofbma}, we
also define 
\begin{eqnarray}
  \label{eq:omega2def}
  \omega_n (B; \D, \Y_n ) = \int _B \prod_{i=1}^n
  \frac{f(Y_i)}{f_0(Y_i)} d \Pi_n(f | \D),
\end{eqnarray}
where $\Pi_n(\cdot |\D)$ is the early stopped OPT prior \emph{without}
shift. The set $B_\D(f_0)$ is the special case of $B_\D^h(f_0)$ in
Section \ref{sec:proofbma} when
$h=0$, and matches the definition in Section
\ref{sec:optanalysis}. $\Theta_\D$ is the collection of all
tree-supported density functions defined on tree topologies of depth
at most $\D$, where trees can grow with or without shape constraints. 
\begin{proof}
 We define $\mathcal{C} = \Big \{\F_\A:
   \kl (f_0, f_{\A} ) \geq \frac{1+ 2\underline{c}_1}{1- \overline{c}_1}
\kl (f_0, f_\D) \Big \}$. The posterior
probability mass on $\mathcal{C}$ can be calculated as 
\begin{eqnarray*}
    \Pi_n (\mathcal{C}  | \Y_n, \D ) = 
\frac{\omega_n (\mathcal{C}; \D,\Y_n )}
{ \omega_n (\Theta_\D; \D, \Y_n)}
\leq 
\frac{\omega_n (\mathcal{C}; \D,\Y_n )}
{ \omega_n (B_\D; \D, \Y_n)}
\end{eqnarray*}

  To bound the denominator, assume $\T^\dsup$ is a non-stopped
  depth-$\D$ tree structure that can provide the best
  approximation to $f_0$, and $f_\D = f_{\T^\dsup}$. 
Then, by Lemma \ref{lemma:intbound}, Lemma \ref{lemma:priorontree}, and Lemma \ref{lemma:prioronball}, 
  with probability at least $1- 3\exp(-4c_2
  (\eta')^2 n \delta^2(n, t))- 
\exp\left( -\underline{c}_2' n \kl (f_0, f_\D) \right)$,
  \begin{eqnarray*}
  \omega_n (B_\D; \D, \Y_n)  
&\geq& \exp(-n (1+\underline{c}_1) \kl (f_0, f_\D)
-4(\eta')^2 (1 + 3 \log (\kappa_U/ \kappa_L)) 2^\D \cdot \D \log 2
    \\
&& \hspace{-0.3in} 
\cdot \exp \left(
- \nu 2^\D \cdot \D \log 2 
+ 2^\D \log \left( ( c_\tau (\eta')^2/2) \Gamma(2a_0)/ (d 
     \cdot  2\Gamma(a_0)) \right) 
 -2^\D \log n
  \right),
  \end{eqnarray*}
where the constant $\underline{c}_2'$ depends on
 $\underline{c}_1$ in the way that is specified by Lemma \ref{lemma:empineqinf}.

  For the numerator, for any depth $\D$ partition $\A$ grown through the
  path $\A^\0, \A^\1, \ldots, \A^{(\D)} = \A$, the size of the partition
  is bounded by $t = 2^\D$. By Lemma \ref{lemma:intfa}, 
  with probability at least $1
- \exp\left( -\frac{2 \overline{c}_2'}{3} n \delta^2(n, t) \right)
- 8\exp\left(-\frac{2c_2}{3} \eta^2 n \delta^2(n, t)\right)$,
uniformly over all tree topologies of depth at most $\D$,  
\begin{eqnarray*}
\omega_n(\F_\A; \D, \Y_n)
\leq 
\Pi_n( \F_\A| \D) \exp \Big( -n (1-\overline{c}_1) \kl (f_0, f_\A)
+
2 \eta^2  2^\D \cdot \D \log 2\Big).
  \end{eqnarray*}
The constant $\overline{c}_2'$ depends on $\overline{c}_1$ in the way
specified by Theorem \ref{thm:likbiasvar}. 

For any depth-$\D$ tree structure $\A$, if $\kl(f_0, f_\A)  \geq \frac{1+ 2
\underline{c}_1}{1- \overline{c}_1} \kl (f_0, f_\D)$, then the
posterior probability on $\A$ is bounded by $\Pi_n( \F_\A| \D) \exp(- n \underline{c}_1
\kl (f_0, f_\D) + O( (2 \eta^2 +4(\eta')^2 (1 + 3 \log (\kappa_U/ \kappa_L)) + \nu + \frac{\log n}{\D
  \log 2} ) 2^\D \cdot \D \log 2  )
\leq \Pi_n( \F_\A| \D) \exp (- \frac{1}{2} n \underline{c}_1 \kl (f_0, f_\D)) $, where
the last inequality is based on the fact that the bias is
dominating. 
And this
event happens with probability at least 
\begin{eqnarray*}
 && 1
- \exp\left( -\frac{2 \overline{c}_2'}{3} \eta^2 n \delta^2(n, t) \right)
-  8\exp\left(- \frac{2 c_2}{3} \eta^2 n \delta^2(n, t)\right)\\
&& \hspace{5mm} - 3\exp(-4c_2
  (\eta')^2 n \delta^2(n, t))- 
\exp\left( -\underline{c}_2' n \kl (f_0, f_\D) \right)
\geq 1- 16\exp \left( - \frac{1}{2} c_2'' \eta^2 \delta^2(n, t)\right),
\end{eqnarray*}
where the last inequality is true when $\eta'$ is set to be the same as $\eta$. $c_2''$ can be set as $c_2'' = \min
\{c_2, \overline{c}_2', \underline{c}_2'\}$.
This finishes the proof.
      \end{proof}

      \begin{remark}
        \label{remark:posteriorratetree} 
 In Proposition \ref{thm:posteriorratetree}, if the early stopped OPT
is grown under a shift $h$, in order to obtain
the same type of concentration result, we need to require $\kl(f_0, f_\D^h)
\geq M_{\text{bias}} \cdot \eta^2  \delta^2(n, t)$ for a slightly
larger $M_{\text{bias}}>1$ 
This is mainly because under a
shift, the lower bound for $\omega_n(B_\D^h(f_0); \D, h, \Y_n)$ is in
the form provided by Corollary \ref{cor:intbound}. A possible choice
of $M_{\text{bias}}$ is $ 2\underline{c}_1^{-1}\eta^{-2} (2 \eta^2 + 8(\eta')^2 (1+ 3 \log(\kappa_U/\kappa_L)) + \nu + \log n/(\D \log 2))$.
      \end{remark}


\subsection{Several Bounds}
\label{subsec:bounds}
Assume that $\Y_n = \{ Y_1, \ldots, Y_n \}$ is a random sample from
$f_0$. In this section, we first provide bounds for $\omega_n (B_\D; \D,
\Y_n)$ and $\omega_n (\Theta_\D; \D,
\Y_n)$, where $\omega_n (\cdot ; \D, \Y_n)$ is defined in
\eqref*{eq:omega2def}, and then generalize these bounds to the case when 
the OPT grows under a shift. 
\begin{lemma}
  \label{lemma:intbound}
  Assume that $ \kappa_l \leq f_0  \leq
  \kappa_U$ and $\D \leq \D_{\max}$. $f_\T$ is the approximation to
  $f_0$ supported by a depth-$\D$ tree $\T$, grown as $\T^\0, \ldots,
  \T^{(\D)} = \T$. Let $t = 2^\D$. With probability at least $1- 3\exp(-4c_2
  (\eta')^2 n \delta^2(n, t))- 
\exp\left( -\underline{c}_2' n \kl (f_0, f_\T) \right)$,
\begin{eqnarray*}
  \nonumber \omega_n \left( B_\D (f_0); \D, \Y_n \right)
&\geq& \exp \left(-n (1+\underline{c}_1) \kl (f_0, f_{\T})
       -4(\eta')^2(1 + 3 \log(\kappa_U/ \kappa_L)) n \delta^2 (n, t) \right)
\\
&& \hspace{-0.6in}
\cdot \pi_{\Opt} \left(\mathcal{F}_{\treek} \right) \exp \left(
    2^\D \log \left( ((\eta')^2/2)
   \cdot \Gamma(2a_0)/ (2\Gamma(a_0))
   \right) 
 -2^\D \log n
  \right).
  \end{eqnarray*}
In this claim, the constant $\underline{c}_2'$ depends on
 $\underline{c}_1$ in the way specified by Theorem
 \ref{lemma:likineqinf}, and $\eta' \geq 1$. 
\end{lemma}
\begin{remark}
  In this lower bound, $\pi_{\Opt}$ denote the OPT prior
  without any constraint. Such a lower bound also holds in the case
  that the individual tree grows under shape constraints, as
  introducing shape constraints only increases the prior probability on
  $\F_{\treek}$. 
\end{remark}
\begin{proof}
 When $f_0$ is bounded from below and above, so is any density
  function in $B_\D$. As a result of Lemmas \ref{lemma:suptoH} and \ref{lemma:kltoH}, for any $f
\in B_\D$, $\kl (f_{\T}, f) \leq  12 \log (\kappa_U/ \kappa_L) (\eta')^2 \delta^2 (n, t)$. By Theorem \ref{lemma:likineqinf}, we
obtain the following inequality, 
{\small
\begin{eqnarray*}
  && \P_{f_0}^n \Bigg( \inf_{f \in B_\D}
\prod_{i=1}^n \frac{f(Y_1)}{f_0(Y_i)} \leq 
\exp\left( -n (1+\underline{c}_1) \kl \left( f_0, f_{\T} \right) -
 4 (\eta')^2n \delta^2 (n, t) 
- n \sup_{f \in B_\D} \kl (f_\T, f) \right) \Bigg)   \\
  \nonumber && \hspace{2.2in}
     \leq \exp \left(- \underline{c}_2' n 
\kl \left( f_0, f_{\T} \right) \right)
+ 3\exp\l- 4 c_2 (\eta')^2 n \delta^2 (n, t)\r.
\end{eqnarray*}
}
The desired inequality is a result of the above lower bound in
combination with Lemma \ref{lemma:priorontree}.
\end{proof}

When the tree is grown under a shift, for any 
$f \in B_\D^h(f_0)$, $\| f - f_\T^h\|_\infty$ may become larger, as
the volume of a node $A \in \T^h$ can be smaller. After the
\tech{pruning step described in Section \ref{sec:notation}}, the
volume of each $A \in \T^h$ at depth $\D$ is lower bounded by
$2^{-\D-1}$. Therefore, $\| f - f_\T^h\|_\infty \leq 2 \eta \D 2^\D/n$
for any $f \in B_\D^h(f_0)$. With this new bound for the
supnorm, we can derive a slightly looser upper bound for
$\omega(B_\D^h(f_0); \D, h, \Y_n)$, summarized in the following
corollary. The proof is skipped as it is very similar to that of Lemma
\ref{lemma:intbound}. 
\begin{corollary}
\label{cor:intbound}
Assume that $ \kappa_l \leq f_0  \leq
  \kappa_U$ and $\D \leq \D_{\max}$. $f_\T^h$ is the approximation to
  $f_0$ supported by a depth-$\D$ tree $\T^h$ under a shift $h$, grown as $\T^\0, \ldots,
  \T^{(\D)} = \T^h$. Let $t = 2^\D$. With probability at least $1- 3\exp(-4c_2
  (\eta')^2 n \delta^2(n, t))- 
\exp\left( -\underline{c}_2' n \kl (f_0, f_\T^h) \right)$, 
\begin{eqnarray*}
 \omega_n \left( B_\D^h (f_0); \D, h, \Y_n \right)
&\geq& \exp \left(-n (1+\underline{c}_1) \kl (f_0, f_{\T}^h)
       -8(\eta')^2(1 + 3 \log(\kappa_U/ \kappa_L)) n \delta^2 (n, t) \right)
\\
\nonumber&& \hspace{-0.6in}
\cdot \pi_{\Opt} \left(\mathcal{F}_{\treek}^h \right) \exp \left(
    2^\D \log \left( ((\eta')^2/2)
   \cdot \Gamma(2a_0)/ (2\Gamma(a_0))
   \right) 
 -2^\D \log n
  \right).
  \end{eqnarray*}
\end{corollary}

In the following, we provide upper bounds for $\omega_n(\cdot; \D,
\Y_n)$ and $\omega_n( \cdot; \D, h, \Y_n)$. 

\begin{lemma} 
\label{lemma:intfa}
Assume that $f_0$ is bounded away from zero and $\D \leq \D_{\max}$. For any depth-$\D$ tree topology $\A$ grown
as $\A^\0, \ldots, \A^{(\D)} = \A$,  the number of leaf nodes is
bounded by $t = 2^\D$. With probability at least $1
- \exp\left( -\overline{c}_2' \eta^2 n \delta^2 (n,
t)\right)
- 8\exp\left(-c_2 \eta^2 n \delta^2(n, t)\right)$, 
\begin{eqnarray}
  \label{eq:intupper}
&&\omega_n(\F_{\{ \A^\ksup: 0 \leq k \leq \D \}} ; \D, \Y_n) \\
\nonumber && \hspace{0.5in}\leq \Pi_n ( \F_{\{ \A^\ksup: 0 \leq k \leq \D \}} |\D )
     \cdot
    \exp\left( -n (1- \overline{c}_1 ) 
  \kl \left( f_0, f_\A \right)+
  2 \eta^2 n \delta^2 (n, t) \right),
\end{eqnarray}
for tree $\A^h$ grown under a shift $h$, with at least the same probability 
\begin{eqnarray}
  \label{eq:intuppershift}
&& \hspace{0.2in} \omega_n(\F_{\{ \A^\ksup: 0 \leq k \leq \D \}}^h ; \D, h, \Y_n) \\
\nonumber&& \hspace{0.6in} \leq \Pi_n ( \F_{\{ \A^\ksup: 0 \leq k \leq \D \}}^h |\D, h)
     \cdot
    \exp\left(-n (1- \overline{c}_1 ) 
   \kl \left( f_0, f_\A^h \right)+
  2 \eta^2 n \delta^2 (n, t) \right),
\end{eqnarray}
where $\overline{c}_1, \overline{c}_2'$, and $c_2 \in (0, 1)$ are as
specified in Theorem \ref{thm:likbiasvar}. 
\end{lemma}

\begin{proof}
The upper bound \eqref*{eq:intupper} is derived based on the following inequality derived
in Theorem \ref{thm:likbiasvar} and Remark \ref{remark:approxBernstein},
\begin{eqnarray*}
 \nonumber && \P_{f_0}^n \left( \sup_{f \in \F_{\A}}
\prod_{i=1}^n \frac{f(Y_i)}{f_0(Y_i)} \geq 
\exp\left( -n (1- \overline{c}_1 ) \kl \left( f_0, f_\A \right) +
  2 \eta^2 n \delta^2 (n, t) \right) \right)   \\
&& \hspace{3cm} \leq \exp \left(- \overline{c}_2' n 
   \max \{ 
   \kl \left(f_0,f_\A  \right),
   \eta^2 \delta^2(n, t) 
   \} \right)
+ 8\exp\left( - c_2 \eta^2 n \delta^2(n, t) \right) ,
\end{eqnarray*}
where $\eta \geq 1$, and $\overline{c}_1, c_2, \overline{c}_2' \in (0,1)$ are some constants
specified in the theorem and remark. The proof for
\eqref*{eq:intuppershift} is similar. 
\end{proof}

Recall that $\Theta_\D^h$ is the collection of density functions
supported by perturbed (along $h$) depth-$\D$ trees grown with weak or strong shape constraints. Then we
have the following bound for the integral over $\Theta_\D^h$. 

\begin{lemma}
\label{lemma:inttheta}
 Assume that $f_0$ is bounded away from zero and $\D \leq \D_{\max}$. Then with probability at least $1
- \exp\left( -\frac{2 \overline{c}_2'}{3} \eta^2 n \delta^2 (n,
t)\right)
- 8\exp\left(-\frac{2c_2}{3} \eta^2 n \delta^2(n, t)\right)$, 
\begin{eqnarray}
  \label{eq:intunifupper}
\omega_n(\Theta_\D^h ; \D, h, \Y_n) 
 \leq 
    \exp\left(  -n (1- \overline{c}_1 ) \kl \left( f_0,
  f_\D^h \right)+
  2 \eta^2 n \delta^2 (n, t) \right),
\end{eqnarray}
where $\overline{c}_1, \overline{c}_2'$, and $c_2 \in (0, 1)$ are as
specified in Theorem \ref{thm:likbiasvar}, and $f_\D^h$ is understood
as the best 
depth-$\D$ approximation to $f_0$ within $\Theta_\D^h$ under the KL divergence. 
\end{lemma}
\begin{proof}
  According to the definition of $f_\D^h$,
  \begin{eqnarray*}
    \omega_n (\F_\A; \D, h, \Y_n) \leq
    \Pi_n ( \F_\A | \D, h)  \cdot 
    \exp\left( -n (1- \overline{c}_1 ) \kl \left( f_0,
  f_\D^h \right)+
  2 \eta^2 n \delta^2 (n, t) \right),
  \end{eqnarray*}
  with  probability at least $1
- \exp\left( -\overline{c}_2' \eta^2 n \delta^2 (n,
t)\right)
- 8\exp\left(-c_2 \eta^2 n \delta^2(n, t)\right)$. The number
of depth-$\D$ tree is bounded by $(d+1)^{t}$ with $t=2^\D$ without any
shape constraints. 
Therefore, with probability at least $1
- \exp\left( -\frac{\overline{2c}_2'}{3} \eta^2 n \delta^2 (n,
t)\right)
- 8\exp\left(-\frac{2c_2}{3} \eta^2 n \delta^2(n, t)\right)$, the bound
holds uniformly. The inequality \eqref*{eq:intunifupper} is based on
the observation that $\omega_n
(\Theta_\D^h; \D, h, \Y_n) = \sum_{\A \text{ is depth-$\D$ tree} }
\omega_n (\F_\A; \D, h, \Y_n)$. 
\end{proof}

\subsection{Analysis of the Optional P\'olya Tree Prior}
\label{sec:optanalysis}
The purpose of this subsection is to obtain the prior thickness result
for $\Pi_n$. Most results presented in this subsection have been obtained in the
previous work \citep{liu2024}. We keep the results here for the sake of
completeness. In the following, we first analyze the OPT prior {\em without}
any constraints or perturbations. The hyper-parameters are specified
as the following.
\begin{itemize}
   \item For any level-$k$ node $A$, the stopping variable follows an independent Bernoulli
    distribution with parameter
     $\tau(A) =\tau_k = c_\tau 2^{-\nu k}$
for some $\tech{0 < c_\tau \leq 1/2}, \nu
\geq 1$. Recall that $\tau_k$ denotes the probability to make a split at level
  $k$, and $1- \tau_k$ is the probability to stop at level $k$.
  \item When imposing the prior on probability masses on leaf nodes, the
    beta random variables are independent and with parameters $(a_0,
    a_0)$ for some $0 < a_0< 1$. Throughout this paper, we set $a_0 =
    1/2$. 
\end{itemize}
For the OPT prior with such specification of 
  hyper-parameters and \emph{without} any shape constraints or perturbations, we denote
  the prior as $\pi_{\Opt}$. We establish the prior thickness result
  for $\pi_{\Opt}$. Such result is sufficient for showing prior thickness of
  $\Pi_n$, as imposing shape constraints or considering the data-dependent truncation only
\emph{increases} the prior probability on the set $B_\D(g) $ defined in
the following. 

We consider a special type of tree topologies of depth exactly equal to $\D$, generated under the OPT prior by growing the tree 
without stopping at \emph{all} intermediate nodes of level $k <\D$,
and stopping all level-$\D$ nodes. This implies that the volume of all terminal nodes is $2^{-\D}$.
 We denote such a depth-$\D$ tree as
$\T$, and assume it is generated through non-stopped intermediate trees $\T^\0,
\T^\1, \ldots, \T^{(\D-1)}$ of depth $0, 1, \ldots, \D-1$
respectively ($\T^\0 $ is simply $\Omega$). The terminal depth-$\D$ tree can also be denoted as
$\T^{(\D)}$ with a convention $\T = \T^{(\D)}$. 
The number of nodes of the intermediate tree $\T^\ksup$ is equal to $2^k$, as there is no stopping at lower levels.
We use $\mathcal{F}_{\treek}$ to denote the collection of density functions
defined on $\T$, which grows as $\T^\0, \T^\1, \ldots, \T^{(\D)}$. 
The following lemma provides bounds for prior distribution on $\mathcal{F}_{\treek}$.


\begin{lemma}
\label{lemma:priorontree}
Under the OPT prior, we have the following lower and upper bounds of
  the prior on a non-stopped tree,
 {\small
  \begin{eqnarray*}
 \log \pi_{\Opt}(\mathcal{F}_{\treek}) \geq \log (1- c_\tau) + (2^\D -2)
     \log ( c_\tau/d)- \nu \log 2 \left( (\D-1) 2^\D - 2^\D +2\right). 
  \end{eqnarray*}
  }
\end{lemma}
\begin{proof}
  See Lemma 8 in the paper by \cite{liu2024}.
\end{proof}

Assume $\T$ is a non-stopped tree topology as described in the previous
subsection. It grows as $\T^\0, \T^\1, \ldots, \T^\dsup$. Let 
$g^\dsup$ be the $L^2$-projection of an arbitrary density $g$
onto $\mathcal{F}_{\treek}$. In the following, we consider an
$L^\infty$-ball around $g^\dsup$, and show a lower bound for the
prior probability on this ball. For a density function $f$ and a set
$A$, we denote $\mu_f (A) = \int_A f d \mu$, where $\mu$ is the
Lebesgue measure. For some constant $\eta \geq 1$, we define a ball
around $g^\dsup$ to be the collection of tree-supported densities
that are uniformly close to $g^\dsup$ on all leaf nodes, and are bounded away from 0. 
\begin{eqnarray*}
  B_\D (g) &\coloneqq& \Big \{ 
f \in \mathcal{F}_{\treek}: |\mu_f (A ) - \mu_{g}(A)| \leq \eta \D/ n \\
&& \hspace{3.5cm}  \mbox { and } \mu_f(A) \geq \left( \frac{1}{2} \eta/n \right)^\D \mbox{ for all } A \in \T^\dsup.
\Big \}
\end{eqnarray*}
As $\T$ is non-stopped with each leaf node of volume $2^{-\D}$, the set
$B_\D (g)$ can be equivalently written as 
\begin{eqnarray*}
  B_\D (g) \coloneqq \left \{ 
f \in \mathcal{F}_{\treek}: \|f - g^\dsup \|_\infty \leq 
\eta \D 2^\D /n, \mbox{ and } f \geq  \left( \frac{1}{2} \eta/n \right)^\D 2^\D.
\right \}
\end{eqnarray*}
Then, we show the following lemma for the prior probability on
$B_\D (g)$.
\begin{lemma}
\label{lemma:prioronball}
Under the optional P\'olya tree prior, for the set $B_\D$ defined above
   \begin{eqnarray*}
   \log \pi_{\Opt}(B_\D (g)) \geq \log \pi_{\Opt}(\mathcal{F}_{\treek}) - 2^\D \log n +
    2^\D \log ((\eta/2) \cdot \Gamma(2 a_0) /( 2 \Gamma(a_0))).
  \end{eqnarray*}
\end{lemma}
\begin{proof}
  See Lemma 9 in the paper by \cite{liu2024}.
\end{proof}

\begin{remark}
\label{remark:boundedbelow}
It is worthy to note that as $B_\D$ is defined by using the supnorm,
  when $g$ is bounded from above and below, density functions
  in $B_\D$ are also bounded from above and below for appropriately
  chosen $n$ and $\D$. Typically, $\D$ is set to be $\lfloor \log_2
  (n/\log n)^{\frac{1}{2 r +1}} \rfloor$ for some $r >0$. So this is true when
  $n$ is large enough. 
\end{remark}

\section{Approximation Results}
\label{sec:approx}
In this section, we study how well tree- and forest-supported density functions can approximate
those in $\mathcal{C}^{\bm \alpha} (\Omega, L)$. The first lemma
characterizes the best tree approximation under the KL divergence. 
\begin{lemma}
  Assume that $f_0$ is a density function on $\Omega$ and $\A= \{ \Omega_j
  \}_{j=1}^t$ is a recursive partition of $\Omega$. Let $f_\A = \sum_{j=1}^{t}
\frac{\int_{\Omega_j} f_0}{\mu(\Omega_j)} \cdot
\mathbbm{1}_{\Omega_j}$. Then for any other piecewise constant density
function $g$ defined on $\A$, $\kl (f_0, f_\A) \leq \kl (f_0, g)$. 
\end{lemma}


\begin{proof}
  For any $\A$-supported piecewise constant density function $g$, 
  \begin{eqnarray*}
    \kl (f_0, g) = \sum_{j=1}^t \int_{\Omega_j} f_0 \log (f_0 / g)
    &=& \sum_{j=1}^t \left(  \int_{\Omega_j} f_0 \log f_0 - (\log g
    )\int_{\Omega_j} f_0 \right).
  \end{eqnarray*}
  Let $a_j = \int_{\Omega_j} f_0$. The value of $g$ on $\Omega_j$
  is denoted as $w_j$. Then to minimize the KL is
  equivalent to solving the following convex optimization problem,
  \begin{eqnarray*}
    \min  &&- \sum_{j=1}^t a_j \log w_j, \\
    \mbox{subject to} && \sum_{j=1}^t w_j \mu(\Omega_j) = 1.  
  \end{eqnarray*}
  The Lagrangian is $L (w_1, \ldots, w_t, \lambda) = -\sum_{j=1}^t a_j
  \log w_j + \lambda (\sum_{j=1}^t w_j \mu(\Omega_j) -1 )$. It is easy
  to see that
  \begin{eqnarray*}
    \frac{\partial L }{ \partial w_j } = - \frac{a_j}{w_j} + \lambda
    \mu(\Omega_j), \quad \mbox{ for } j = 1, \ldots, t. 
  \end{eqnarray*}
  By setting partial derivatives to 0, we have $w_j = a_j /
  (\lambda \mu(\Omega_j)) $. Plugging this into the constraint, we
  obtain $\lambda=1$. Therefore, when $w_j = a_j / \mu(\Omega_j)$ the
  objective function is minimized. This finishes the proof. 
\end{proof}

The following lemma shows that when $f_0$ is smooth enough and the
tree grows under the strong shape constraint, the
approximation error under the KL divergence has an analytical form as
the tree depth increases. We
define $\theta_l = \int \frac{(\partial f_0/ \partial y^l)^2 }{f_0}$,
which is a reflection of the variation of $f_0$ along the $l$-th
dimension. The quantify $\theta_l$ is well defined when $f_0 \in
\mathcal{C}^{\bm \alpha}(\Omega, L)$ for $\alpha_l >1$. 
\begin{lemma}
\label{lemma:approxtree}
  Assume $f_0 \in \mathcal{C}^{\bm \alpha} (\Omega, L)$ with $\alpha_l
  \in (1, 2]$, and $f_0$ is bounded away from
zero. $\A = \{ \Omega_j \}_{j=1}^t$ is a tree topology obtained under the strong
  shape constraint without perturbation. The edge length of each terminal region $\Omega_j$
  along
  dimension $l$ is denoted as $\epsilon_l$ with $\epsilon_l = O(2^{-
    \D/d})$ as $\D \rightarrow \infty$. Define $f_\A = \sum_{j=1}^{t}
\frac{\int_{\Omega_j} f_0}{\mu(\Omega_j)} \cdot
\mathbbm{1}_{\Omega_j}$. Let $\theta_l = \int \frac{(\partial f_0/ \partial
  y^l)^2 }{f_0}$ for $l = 1, \cdots, d$. 
Then, 
\begin{eqnarray*}
  \int f_0 \log (f_0/f_\A)
  =\frac{1}{24} \sum_{l=1}^d
  \theta_l \epsilon_l^2
      + O \left( \sum_{l=1}^dd L^2 \epsilon^{1+\alpha_l}/ \kappa_L \right),
\end{eqnarray*}
where $\epsilon = \max_{1 \leq l \leq d} \epsilon_l$ and $\kappa_L =
\inf_\Omega f_0$. 
\end{lemma}
\begin{proof}
 Without loss of generality, we also assume
that $\theta_l >0$ for $l = 1,\ldots, d$. As $f_0 \in \mathcal{C}^{\bm \alpha} (\Omega, L)$, for any
$x_{j0}, x_j \in \Omega_j$,
\begin{eqnarray*}
  \big |
  f_0(x_j) - f_0(x_{j0})
  - \nabla f_0(x_{j0})^\top (x_j - x_{j0})
  \big |
  \leq \sum_{l=1}^d L \epsilon_l^{\alpha_l}.
\end{eqnarray*}
Then
{\small
\begin{eqnarray*}
  \left( \int_{\Omega_j} f_0(x_j) d x_j \right) / \mu(\Omega_j)
  -f_0( x_{j0}) &=&
  \frac{1}{\mu(\Omega_j)}
  \int_{\Omega_j } \nabla f_0(x_{j0})^\top (x_j - x_{j0}) d x_j
  + O \left(\sum_{l=1}^d L \epsilon_l^{\alpha_l} \right)\\
  &=& \nabla f_0(x_{j0})^\top ((b_{j, U}+ b_{j, L})/2 -
      x_{j0})
      + O \left(\sum_{l=1}^d L \epsilon_l^{\alpha_l} \right)\\
  &=& O\left(\sum_{l=1}^d L \epsilon_l \right) = O ( d L \epsilon). 
\end{eqnarray*}
}
  This can be viewed as a bound for $\|f_0 - f_\A\|_\infty$. 

Keeping this point-wise bound for $f_0 - f_\A$, we proceed to bound
the KL divergence between $f_0$ and $f_\A$. 
{\small
\begin{eqnarray*}
  \kl(f_0, f_\A) = \int f_0 \log (f_0/f_\A) 
= \int \left(f_0 \log f_0 -  f_\A \log f_\A + f_\A \log f_\A - f_0
  \log f_\A \right).
\end{eqnarray*}
}
Based on the Taylor's expansion for function $x \log x$, we have 
for any $x \in\Omega$,
{\small
\begin{eqnarray*}
 && f_0(x) \log f_0(x) - f_\A (x)\log f_\A(x) \\
&& \hspace{0.2in} =
 (1+ \log f_\A(x) )(f_0(x) - f_\A(x) )
  + \frac{1}{2 f_\A(x)}(f_0(x) - f_\A(x) )^2 +
  O ((f_0(x) - f_\A(x) )^3/ f_\A^2(x)).
\end{eqnarray*}
}
The upper bound for $\|f_0 - f_\A\|_\infty$ ensures such an expansion is
  valid. 
When $f_0$ is bounded away from zero, 
based on the dominated convergence theorem, 
{\small
\begin{eqnarray*}
  \int f_0 \log (f_0/f_\A)
&=& \sum_{j=1}^t \int_{\Omega_j} \frac{1}{2 f_\A(x)}
\left( f_0(x) - f_\A(x) \right)^2 + O\left( d^3L^3 \epsilon^3/ \kappa_L^2
    \right).\\
\end{eqnarray*}
}
Within each $\Omega_j$, let $x_j
\in \Omega_j$ be the data point that satisfies $f_\A(x_j) = f_0
(x_j)$. By H\"older-continuity of $f_0$ and the mean value theorem,
such a data point exists. Given the fact that $f_0$ is bounded away
from zero, $\kl(f_0, f_\A)$ can be further written
as 
{\small
\begin{eqnarray*}
  \int f_0 \log (f_0/f_\A) &=& \sum_{j=1}^t\int_{\Omega_j} \frac{1}{2f_0(x_j)}
  \left( \nabla f_0(x)^\top (\frac{1}{2}(b_{j, U}+ b_{j, L}) -
  x )\right)^2 + O \left( \sum_{l=1}^d d L^2 \epsilon^{1+\alpha_l} / \kappa_L\right).
\end{eqnarray*}
}
For any $x \in \Omega_j$, $| (\partial f_0/ \partial x^l) (x) -
(\partial f_0/ \partial x^l) (x_j) | \leq
\sum_{l=1}^d L \epsilon^{\alpha_l-1}$. In combination with the fact
that $f_0$ is bounded away from zero, we have 
{\small
\begin{eqnarray*}
  \int f_0 \log (f_0/f_\A) &=&\sum_{j=1}^t\int_{\Omega_j} \frac{1}{2f_0(x_j)}
  \left( \nabla f_0(x_j)^\top (\frac{1}{2}(b_{j, U}+ b_{j, L}) -
  x )\right)^2 +O \left( \sum_{l=1}^d d L^2 \epsilon^{1+\alpha_l}/ \kappa_L \right),
\end{eqnarray*}
}
Based on the fact that $|(\partial f_0 /\partial x^l)^2 (x) /f_0(x) - (\partial f_0 /\partial x^l)^2 (x_j)
/f_0(x_j) | = O(d L^2 \epsilon_l^{\alpha_l-1}/\kappa_L)$, with some further
calculation,
{\small
\begin{eqnarray*}
  \int f_0 \log (f_0/f_\A)
  &=&\sum_{j=1}^t\int_{\Omega_j}
  \frac{\sum_{l=1}^d
  ( \partial f_0 /\partial x^l)^2(x_j)(\frac{1}{2}(b_{j, U}^l+ b_{j, L}^l) -
  x_l )^2}
  {2 f_0(x_j)}
  + O \left( \sum_{l=1}^d  d L^2 \epsilon^{1+\alpha_l} / \kappa_L\right)\\
  &=&\frac{1}{24} \sum_{l=1}^d
  \epsilon_l^2\int  \frac{
      (\partial f_0 /\partial x^l)^2}  {f_0}
      + O \left( \sum_{l=1}^d d L^2 \epsilon^{1+\alpha_l}/ \kappa_L \right).
\end{eqnarray*}
}
This finishes the proof. 
%
\end{proof}

\begin{remark}
  \label{remark:approxtree}
If $\A^h$ is a tree topology grown under the strong shape constraint
with perturbation $h$, as $\epsilon_l \rightarrow 0$, the expansion in
Lemma \ref{lemma:approxtree} still holds for H\"older continuous $f_0$, with $\epsilon_l$'s
representing edge lengths of $\A$---the tree topology without
perturbation. 
\end{remark}

\rev{
  The next lemma characterizes the best approximation to $f_0$ under
  the Hellinger distance when the tree topology is fixed. 
  \begin{lemma}
    \label{lemma:bestHellinger}
  Assume that $f_0$ is a density function on $\Omega$ and $\A= \{ \Omega_j
  \}_{j=1}^t$ is a recursive partition of $\Omega$. Let $u_j =
  \int_{\Omega_j} \sqrt{f_0}$ and $f_\A^{\text{H}} = \sum_{j=1}^{t}
\frac{w_j}{\mu(\Omega_j)} \cdot
\mathbbm{1}_{\Omega_j}$, where $w_j$'s satisfy
\begin{eqnarray*}
  w_j =\frac{ u_j^2} {
   \lambda^2 \mu(\Omega_j)}
  \quad
  \mbox{with}
  \quad
  \lambda^2 =\sum_{j=1}^t \frac{u_j^2}{\mu( \Omega_j)}.
\end{eqnarray*}
Then for any other piecewise constant density
function $g$ defined on $\A$, $\rho (f_0, f_\A^{\text{H}} ) \leq \rho (f_0, g)$. 
\end{lemma}
\begin{proof}
  For any $\A$-supported piecewise constant density function $g =
  \sum_{j=1}^t \frac{w_j} {\mu(\Omega_j)} \cdot \mathbbm{1}_{\Omega_j}$, 
  minimizing the Hellinger distance between $f_0$ and $g$ is
  equivalent to minimizing $- \int_\Omega \sqrt{g f_0} = -
  \sum_{j=1}^t \sqrt{\frac{w_j} {\mu(\Omega_j)}} \int_{\Omega_j}
  \sqrt{f_0}$.  Let $u_j = \int_{\Omega_j} \sqrt{f_0}$, then 
  to minimize the Hellinger distance is
  equivalent to solving the following convex optimization problem,
  \begin{eqnarray*}
    \min  &&- \sum_{j=1}^t  u_j \sqrt{\frac{w_j} {\mu(\Omega_j)}}, \\
    \mbox{subject to} && \sum_{j=1}^t w_j = 1.  
  \end{eqnarray*}
  The Lagrangian is $L (w_1, \ldots, w_t, \lambda) = -\sum_{j=1}^t u_j \sqrt{\frac{w_j} {\mu(\Omega_j)}}+ \lambda (\sum_{j=1}^t w_j -1 )$.
  \begin{eqnarray*}
    \frac{\partial L }{ \partial w_j } = -
    \frac{u_j}{\sqrt{\mu(\Omega_j)}} \cdot \frac{1}{2 \sqrt{w_j}} + \lambda, \quad \mbox{ for } j = 1, \ldots, t. 
  \end{eqnarray*}
  By setting partial derivatives to 0, we have $w_j = u_j^2 /
  ( 4 \lambda^2 \mu(\Omega_j)) $. Plugging this into the constraint, we
  obtain $\lambda^2 = \sum_{j=1}^t \frac{u_j^2}{4 \mu( \Omega_j)}$. This finishes the proof. 
\end{proof}
}

Our next lemma on approximation error of a type of specifically
designed forests shows that aggregating perturbed trees can lead to
improved approximation accuracy. 
\begin{lemma}
    \label{lemma:biasforest}
 Assume $f_0 \in \mathcal{C}^{\bm \alpha} (\Omega,L)$ with $\alpha_l
 \in (0, 2]$ and the
 Condition \ref{cond:boundary} is satisfied. $\overline{\alpha}$ is
calculated as $(\overline{\alpha})^{-1} = \frac{1}{d} \sum_{l=1}^d
\alpha_l^{-1}$. Assume that trees are grown with perturbations and there are $B =
O(2^\D)$ trees within any ensemble 
of depth $\D$. Then for all $\D >0$, there exist a perturbation design
scheme $\bm \upsilon$ and a corresponding forest-supported 
 density $f_\D^\en$ of the form 
 \begin{eqnarray*}
   f_\D^\en = \frac{1}{B} \sum_{m=1}^{B} f_\D^\m, \quad f_\D^\m \mbox{
   is a depth-$\D$ tree grown under the shift } h_m \in \mathcal{H}_{\D, \bm \upsilon}, 
 \end{eqnarray*}
  such that 
  \begin{eqnarray*}
    \| f_0 - f_\D^\en \|_\infty \lesssim 4 d L 2^{-\D \overline{\alpha}/d}.
  \end{eqnarray*}
  \end{lemma}
\begin{proof}
 Without loss of generality, we may assume $\alpha_l >1$, as
 single trees are already sufficient to approximate H\"older continuous
 functions well. We construct a depth-$\D$ forest approximation to $f_0$ as the
  following. Each density within the
  forest is supported by the tree topology $\A_\D^\star$ satisfying
  the strong shape constraint or shifted
  version of it. The edge length of $\A_\D^\star$ along direction $l$
  is $\epsilon_l $. We may adjust some of
  $\epsilon_l$'s such that $\sum_{l=1}^d - \log_2 \epsilon_l = \D$. The
  perturbation design scheme is $\bm \upsilon = (\upsilon_l)_{1
  \leq l \leq d}$ with $\upsilon_l = \epsilon_l^2$. Potential shifting
directions are $h =
(i_l \upsilon_l)_{1 \leq l \leq d}$, where each $0 \leq i_l < \epsilon_l/\upsilon_l$ is an integer. The total
number of possible shifting directions is exactly $B=2^\D$. Let $(h_m)_{1
  \leq m \leq B}$ be an enumeration of these directions. Then, at
depth $\D$ the tree approximation $f_\D^\m$ defined on the $m$-th tree
with $t_m$ leaf nodes is 
\begin{eqnarray*}
  f_{\D}^\m (x) = \sum_{j=1}^{t_m}\frac{1 }{\mu(\Omega_j^{h_m}
  \cap \Omega) }\left( \int_{\Omega_j^{h_m}}
  f_0 (x') dx' \right) \mathbbm{1}_{\Omega_j^{h_m}\cap \Omega} (x), 
\end{eqnarray*}
where $h_m = (h_m^1, \ldots, h_m^d )^\top$ is the vector representing the
shift, $\Omega_j^{h_m}$ is a leaf node obtained under a shift. We
define the forest approximation as $f_\D^\en =( \sum_{m=1}^{B}
f_\D^\m)/B $. In the following, we show
$f_\D^\en$ can achieve the desired approximation rate.

We divide $\Omega$ into an interior region and a boundary region
as
\begin{eqnarray*}
  \Omega = I_{\bm \epsilon} \cup J_{\bm \epsilon}, 
\end{eqnarray*}
where $I_{\bm \epsilon} = \{ x \in \Omega: x^l > 2\epsilon_l
\mbox{ and } ( 1-x^l)  > 2\epsilon_l \mbox{ for some } 1 \leq l \leq
d\}$ and $J_{\bm \epsilon} = \Omega \setminus
I_{\bm \epsilon}$.

For any
$x_0 \in I_{\bm \epsilon}$, 
\begin{eqnarray*}
  f_0(x_0) - f_\D^\en(x_0) =
  f(x_0) -
  \frac{1}{B} \sum_{m=1}^{B}
  \sum_{j: x_0 \in \Omega_j^{h_m}}
  \frac{1}{\mu(\Omega_j^{h_m}) } \left( \int_{\Omega_j^{h_m}} f_0(x) dx \right), \quad 0 \leq h_m^l< \epsilon_l.
\end{eqnarray*}
Assume that $x_0 \in \Omega_{j_0}$ for some $j_0$, where
$\Omega_{j_0}$ is a leaf node of $\A_\D^\star$, which is grown without
perturbation. If we allow $h_m^l$ to
take negative values, then the approximation at a point $x_0$ can be rewritten as
\begin{eqnarray}
   \label{eq:approxholder}
  f^\en_\D(x_0) = \frac{1}{B}\sum_{h_m: x_0 \in \Omega_{j_0}^{h_m} } 
   \frac{1}{\mu( \Omega_{j_0}^{h_m}) } \left( \int_{ \Omega_{j_0}^{h_m}} f_0(x) dx \right), \quad -\epsilon < h_m^l < \epsilon. 
\end{eqnarray}
Let $\Omega_{j_0} = \otimes_{l=1}^d [b_{L, j_0} ^l, b_{U,
  j_0}^l)$. The condition $x_0 \in \Omega_{j_0}^{h_m}$ is satisfied
when $ b_{L, j_0}^l +h_m^l \leq x_0^l <  b_{U, j_0}^l + h_m^l$ with
$|h_m^l| < \epsilon_l$. Although for the definition of density forests we
require $h_m^l \in [0, \epsilon_l)$, here an $h_m$ with negative entries
can be understood as the node is obtained under the shift 
$\tilde{h}_m$, where $\tilde{h}_m^l = h_m^l$ for $h_m^l \geq 0$ and
$\tilde{h}_m^l = \epsilon_l + h_m^l$ for $h_m^l < 0$. Under the shift $\tilde{h}_m$,
we can find an $\Omega_{j'}$ such that  $x_0 \in
\Omega_{j'}^{\tilde{h}_m}$, and the set $\Omega_{j'}^{\tilde{h}_m}$
can be equivalently viewed as $\Omega_{j_0}^{h_m}$ without affecting
the calculation involved in \eqref*{eq:approxholder}.

For any $x, x_0$, the difference between $f_0(x)$ and $f_0(x_0)$ can
be written as $f_0(x) - f_0(x_0) = (f_0(x^1,x^2, \ldots, x^d) -
f_0(x_0^1, x^2, \ldots, x^d)) + \cdots + (f_0(x_0^1, \ldots, x_0^{d-1},
x^d) - f_0(x_0^1,\ldots, x_0^{d-1},x_0^d))$. To simplify the notation, we denote
each of these $d$ terms as $\Delta_1(x, x_0), \ldots, \Delta_d(x,
x_0)$, respectively. Then,
\begin{eqnarray*}
  |f_\D^\en(x_0) - f_0(x_0) | 
\leq \sum_{l=1}^d \bigg| \frac{1}{B} \sum_{h_m: x_0 \in \Omega_{j_0}^{h_m} }
      \frac{1}{\mu(\Omega_{j_0}) } \left( \int_{\Omega_{j_0}^{h_m}}
  \Delta_l(x, x_0) dx \right) \ \bigg|.
\end{eqnarray*}
In the following, we provide a bound for each of the $d$ terms on the
right hand side. Without loss of generality, we take $l=1$ as an
example. The other terms can be bounded similarly.

As $f_0$ lies in the H\"older space $\mathcal{C}^{\bm \alpha} (\Omega,
L) $, applying the mean value theorem, we have that for some $t_1 \in (0, 1)$ and any $x, x_0$ lying in the interior region
of $\Omega$:
  $f_0(x^1, \ldots, x^d) -
f_0(x_0^1, x^2, \ldots, x^d) = \frac{\partial }{\partial x^1} f_0 \left(x_0^1 +
  t_1(x^1-x_0^1), x^2, \ldots, x^d \right) (x^1-x_0^1), 
$
Then,
\begin{eqnarray*}
  & &  \bigg |   f_0(x^1, \ldots, x^d) -f_0(x_0^1, x^2, \ldots, x^d) -  \frac{\partial }{\partial x^1} f_0 \left(x_0^1, x^2, \ldots, x^d \right) (x^1-x_0^1) \bigg|\\
  & & \hspace{5mm}= \left| \left(  \frac{\partial }{\partial x^1} f_0 \left(x_0^1 +
  t_1(x^1-x_0^1), x^2, \ldots, x^d \right) -  \frac{\partial }{\partial x^1} f_0 \left(x_0^1, x^2, \ldots, x^d \right)\right) ^\top (x^1-x_0^1)  \right|\\
  & & \hspace{5mm} \leq L \epsilon_1^{\alpha_1 }.
\end{eqnarray*}
Therefore, for any $x_0 \in I_{\bm \epsilon}$, by taking integral over
$\Omega_{j_0}^{h_m}$'s in (\ref{eq:approxholder}), we have
{\small
\begin{eqnarray*}
  & & \left| \frac{1}{B} \sum_{h_m: x_0 \in \Omega_{j_0}^{h_m} }
      \frac{1}{\mu(\Omega_{j_0}) } \left( \int_{\Omega_{j_0}^{h_m}}
  \Delta_1(x, x_0) dx \right) \ \right|\\
  & & \hspace{5mm} \leq \left| \frac{1}{B} \sum_{h_m: x_0 \in \Omega_{j_0}^{h_m} }
      \frac{1}{\mu(\Omega_{j_0}) } \left( \int_{\Omega_{j_0}^{h_m}} \frac{\partial }{\partial x^1} f_0 \left(x_0^1, x^2, \ldots, x^d \right) (x^1-x_0^1) dx \right)\ \right|
      +L \epsilon_1^{\alpha_1}.
\end{eqnarray*}}
Applying the mean value theorem for integrals, we can find a point lying in
$I_{\bm \epsilon}$ of the form $(x_0^1, \tilde{x}^2, \ldots,\tilde{x}^d)$, such
that the integral above can be calculated as 
{\small
  \begin{eqnarray*}
  & & \hspace{15mm}  \left| \frac{\partial }{\partial x^1} f_0(x_0^1, \tilde{x}^2, \ldots,\tilde{x}^d) \Bigg( \frac{1}{B}  \sum_{h_m: x_0 \in \Omega_{j_0}^{h_m} }
      \frac{1}{\mu(\Omega_{j_0}) } \int_{\Omega_{j_0} } (x^1+h_m^1-x_0^1) dx \Bigg)\ \right|
      +  L \epsilon_1^{\alpha_1}\\
  & & \hspace{15mm} =  \left| \frac{\partial }{\partial x^1} f_0(x_0^1, \tilde{x}^2, \ldots, \tilde{x}^d )\Bigg(
      \frac{1}{\mu(\Omega_{j_0}) } \int_{\Omega_{j_0} } (x^1-x_0^1) dx
      + \frac{1}{B}  \sum_{h_m: x_0 \in \Omega_{j_0}^{h_m} } h_m^1
      \Bigg)
     \right|
      +  L \epsilon_1^{\alpha_1}.
\end{eqnarray*}
}
We calculate the integral more explicitly as 
   $\frac{1}{\mu(\Omega_{j_0}) } \int_{\Omega_{j_0} } (x^1 - x_0^1) dx
   = \frac{1}{2}(b_{j_0, L}^1 + b_{j_0, U}^1) - x_0^1$. 
At the same time, all $h_m$'s with $x_0 \in \Omega_{j_0}^{h_m}$ satisfy $ b_{L, j_0}^1 + h_m^1 \leq x_0^1 <  b_{U, j_0}^1 + h_m^1$ and are equally spaced along each dimension. 
If $x_0 \in [b_{L, j_0}^1 + m_0\upsilon_1 , b_{L, j_0}^1 + (m_0+1)
\upsilon_1)$ for some integer $0 \leq m_0 < \lceil
\epsilon_1/\upsilon_1 \rceil $, then it is easy to show that,
\begin{eqnarray*}
 \bigg| \frac{1}{\mu(\Omega_{j_0}) } \int_{\Omega_{j_0} } (x^1 - x_0^1) dx + \frac{1}{B} \sum_{h_m: x_0 \in \Omega_{j_0}^{h_m} } h_m^1\bigg|
      \leq \upsilon_1/2.
\end{eqnarray*}

Repeating the same argument for the bias due to each $\Delta_l(x,
x_0)$, together with the boundedness of partial derivatives, we have
\begin{eqnarray*}
   \left|f_0 (x_0) - f_\D^\en (x_0) \right| \leq \sum_{l=1}^d L \upsilon_l + \sum_{l=1}^d L \epsilon_l^{\alpha_l}
  \quad \mbox{for all } x_0 \in I_{\bm \epsilon}.
\end{eqnarray*}
In this upper bound, $\sum_{l=1}^d L \epsilon_l^{\alpha_l}$ is the
dominating term, as $\upsilon_l = \epsilon_l^2$ by design and
$\alpha_l \leq 2$.

Meanwhile, for any $x_0
  \in J_{\bm \epsilon}$, we can find $l_0$, such that $x_0$ lies in the boundary
  region $J_{\bm \epsilon}^{l_0}= \{ x: x^{l_0} < 2\epsilon_{l_0} \mbox{ or
} 1-x^{l_0} < 2 \epsilon_{l_0}\}$. Under the Condition
  \ref{cond:boundary}, $ \left | \frac{\partial f_0}{\partial x^l}
    |_{x_0} \right| 
\leq L \epsilon_{l_0} $ for all $1 \leq l\leq d$. Therefore, 
  $|  f_0(x_0)- f_\D^\en (x_0) | \leq L \epsilon_{l_0} \sum_{l=1}^d
  \epsilon_l + \sum_{l=1}^d L \epsilon_l^{\alpha_l}.
$
The term $\sum_{l=1}^d L \epsilon_l^{\alpha_l}$ is still
dominating. 

Combining the approximation results for the interior part and for
  the boundary, the approximation error for a depth-$\D$ forest is
  mininized when $\epsilon_l = O\left(  2^{- \D \lfloor
    \overline{\alpha} / (d\alpha_l) \rfloor } \right)$. With such
choices of $\epsilon_l$'s, 
  \begin{eqnarray*}
    \| f_0 - f_\D^\en \|_\infty \lesssim 4d L 2^{- \D \overline{\alpha}/d}.
  \end{eqnarray*}
This finishes the proof. 
\end{proof}

\rev{Without the boundary condition, we may only achieve an accuracy level
as good as single-tree approximations at the boundary. Specifically,
we have the following corollary on the approximation accuracy without
the boundary condition.}
\begin{corollary}
  \label{cor:biasforest}
Under the same assumptions as those for Lemma \ref{lemma:biasforest},
if we remove Condition \ref{cond:boundary} while assuming $f_0 \geq \kappa_L$, then for any $f_0 \in
\mathcal{C}^{\bm \alpha} (\Omega, L)$, we can find a depth-$\D$ forest
approximation $f_\D^\en$, such that $\rho(f_0, f_\D^\en) \lesssim ( 4d^2L/
\sqrt{\kappa_L}) 2^{- \D \overline{\alpha^*} /d} $, where $(\overline{\alpha^*})^{-1} =
\frac{1}{d} \sum_{l=1}^d (\alpha_l^*)^{-1}$, $\alpha_l^* = \min\{
\alpha_l, 3/2 \}$ and $\kappa_L = \inf_\Omega f_0$.
\end{corollary}
\begin{proof} 
  \rev{We study the approximation error bound for the forest
    approximation considered in the proof of Lemma~\ref{lemma:biasforest}. We divide $\Omega$ into an interior region and a boundary region
as
\begin{eqnarray*}
  \Omega = I_{\bm \epsilon} \cup J_{\bm \epsilon}, 
\end{eqnarray*}
where $I_{\bm \epsilon} = \{ x \in \Omega: x^l > 2\epsilon_l
\mbox{ and } ( 1-x^l)  > 2\epsilon_l \mbox{ for all } 1 \leq l \leq
d\}$ and $J_{\bm \epsilon} = \Omega \setminus
I_{\bm \epsilon}$.}

\rev{At any boundary region $J_{\bm \epsilon}^{l}=\{ x: x^{l} < 2\epsilon_{l} \mbox{ or
} 1-x^{l} < 2 \epsilon_{l}\}$, the approximation error bound is
looser. More specifically, for $x \in J_{\bm \epsilon}^{l}$, $\|f_0 (x)- f_\D^\en(x) \|_\infty
\leq L \sum_{l'=1}^d \epsilon_{l'}$. }

\rev{If $f_0$ is bounded away from zero, when we calculate the
  Hellinger distance, in terms of the integral over the boundary
  region, 
\begin{eqnarray*}
\int_{J_{\bm \epsilon}} (\sqrt{f_0} - \sqrt{f_\D^\en } )^2 
\leq
\sum_{l'=1}^d \int_{J_{\bm \epsilon}^{l'}} (\sqrt{f_0} - \sqrt{f_\D^\en } )^2 
\leq 
\sum_{l'}^d 4\epsilon_{l'}( \frac{L}{\sqrt{\kappa_L}}\sum_{l=1}^d \epsilon_l)^2
  ,
\end{eqnarray*}
where $\kappa_L = \inf_\Omega f_0$.
A loose bound for $\int_{J_{\bm \epsilon}} (\sqrt{f_0} -
\sqrt{f_\D^\en } )^2$ is $\frac{ 4 d^3 L^2 }{\kappa_L} \max_{1 \leq l
  \leq d}
\epsilon_l^3$.
For the integral over the interior region,
\begin{eqnarray*}
  \int_{I_{\bm \epsilon}} (\sqrt{f_0} - \sqrt{f_\D^\en } )^2
  \lesssim ( \frac{L}{\sqrt{\kappa_L}}\sum_{l=1}^d
  \epsilon_l^{\alpha_l})^2 \lesssim \frac{dL^2}{\kappa_L}
  \sum_{l=1}^d \epsilon_l^{ 2\alpha_l}.
\end{eqnarray*}
Therefore,
\begin{eqnarray*}
  \rho^2(f_0, f_\D^\en) \lesssim \frac{8 d^3 L^2}{\kappa_L} \sum_{l=1}^d
  \max \{\epsilon_l^3, \epsilon_l^{2 \alpha_l} \}.
\end{eqnarray*}
Given any depth $\D$, the righthand side is minimized when $\epsilon_l = 2^{- \D \lfloor
    \overline{\alpha^*} / (d\alpha_l^*) \rfloor }$. The corresponding
  approximation error bound is $\rho(f_0, f_\D^\en)  \lesssim (4 d^2 L/
  \sqrt{\kappa_L} )2^{- \D \overline{\alpha^*} /d }$.}
\end{proof}

Using an argument similar to that for Lemma~\ref{lemma:biasforest} and
Corollary~\ref{cor:biasforest}, we can show the
following result. 
\begin{corollary}
  \label{cor:biasforestbma}
Assume $f_0 \in \mathcal{C}^{\bm \alpha} (\Omega, L)$ and Condition~\ref{cond:boundary} is satisfied. $\A$ is a depth-$\D$ tree
topology grown under the strong shape constraint without
perturbation. The edge length of all leaf nodes in $\A$ along the dimension $l$ is $\epsilon_l = O( 2^{- \D/d})$. $\mathcal{H}_\D$ is a set of potential shifting directions defined in Section~\ref{subsec:fullybma}. Define the forest approximation as 
\begin{eqnarray*}
  f_\D^\en = \frac{1}{|\mathcal{H}_\D|}
\sum_{h \in \mathcal{H}_\D} f_{\A^h}, 
\end{eqnarray*}
then, $\| f_0 - f_\D^\en \|_\infty
\lesssim 2 L \sum_{l=1}^d \epsilon_l^{\alpha_l}$. If we remove Condition~\ref{cond:boundary} but further assume that $f_0$ is bounded from
below, then $\rho(f_0 , f_\D^\en) \lesssim \sum_{l=1}^d \left( 2 d L
  \sqrt{2d} / \sqrt{\kappa_L} \right) \epsilon_l^{\alpha_l^*}$, where $\alpha_l^* = \min\{
\alpha_l, 3/2\}$.  
\end{corollary}
\begin{proof}
  \rev{First, note that the forest approximation considered in this corollary matches
    that described at the beginning of Lemma~\ref{lemma:biasforest}.}

  \rev{When Condition~\ref{cond:boundary} is satisfied, repeating the
    argument used in the proof of Lemma~\ref{lemma:biasforest}, we can
    show that for any $x \in \Omega$
    \begin{eqnarray*}
      | f_0 (x)  - f_\D^\en (x) | \lesssim 2L \sum_{l=1}^d
      \epsilon_l^{\alpha_l}. 
    \end{eqnarray*}
  }

  \rev{When Condition~\ref{cond:boundary} is {\em not} satisfied, for
    the tree approximation considered in this corollary, with
    the argument for Corollary~\ref{cor:biasforest}, we have
    \begin{eqnarray*}
  \rho^2(f_0, f_\D^\en) \lesssim \frac{8 d^3 L^2}{\kappa_L} \sum_{l=1}^d
      \max \{\epsilon_l^3, \epsilon_l^{2 \alpha_l} \}
      \lesssim \frac{8 d^3 L^2}{\kappa_L}
      \left( \sum_{l=1}^d \max \{\epsilon_l^{3/2}, \epsilon_l^{\alpha_l}
      \}
      \right)^2.
\end{eqnarray*}
  This finishes the proof.}
\end{proof}

\section{Auxiliary Lemmas}
\label{sec:tools}
Most results in Section \ref{subsec:entropy}
and Section \ref{subsec:truncation} have been obtained in the previous
work \citep{liu2023}. We still summarize them here for the sake of
completeness. We also make some improvements of the previous
results. 
\subsection{Bounds for the Bracketing Metric Entropy}
\label{subsec:entropy}
We first define the metric
entropy with bracketing of a space. 
\begin{definition}
  Let $ (\Theta, \rho)$ be a separable pseudo-metric
  space. $\Theta(\epsilon)$ is a finite set of pairs of functions $\{
  (u_j, v_j), j=1,\cdots, N\}$ satisfying
  \begin{eqnarray}
    \label{eq:entropyDefC1}
    \rho (u_j, v_j) \leq \epsilon\ for\ j=1,\cdots,N,
  \end{eqnarray}
and for any $f \in \Theta$, there is a $j$ such that
\begin{eqnarray}
  \label{eq:entropyDefC2}
  u_j \leq f \leq v_j.
\end{eqnarray}
Let 
\begin{eqnarray*}
  \bnum (\epsilon, \Theta, \rho) = min \{|\Theta(\epsilon)|:
  \Theta(\epsilon) \mbox{ is a set satisfying 
  (\ref{eq:entropyDefC1}) and
  (\ref{eq:entropyDefC2})} \}.
\end{eqnarray*}
The metric entropy with bracketing of $\Theta$ is defined to be 
\begin{eqnarray*}
  \bentro (\epsilon, \Theta, \rho) = \log \bnum (\epsilon, \Theta, \rho).
\end{eqnarray*}
\end{definition}

The following lemma provides an upper bound for the bracketing metric entropy of a tree model when the underlying partition is given.  
\begin{lemma}
  \label{lemma:boundofH}
Take $\rho$ to be the Hellinger distance. Let 
$\Theta_t^{\mathcal{A}, L} = \{ f:\ f\ is\ a\ density\ function\\ supported\ by\ the\ recursive\
partition\ \mathcal{A} = \{ \Omega_j\}_{j=1}^t,\ and\ \rho(f, f_{\A})
\leq L \}$, where $f_{\A}$ is any fixed piecewise constant density
function defined on $\A$. Then for all $\epsilon \leq L$, 
\begin{equation}
  \label{eq:boundofH}
  \bentro (\epsilon, \Theta_t^{\mathcal{A}, L}, \rho) \leq
  \frac{1}{2} t \log t + t
  \log(L/\epsilon) +b_{\text{entropy}}, 
\end{equation}
where $b_{\text{entropy}}$ is a constant. 
\end{lemma}
\begin{proof}
  See Lemma 14 in Section 7 of the paper by \cite{liu2023}.
\end{proof}

The number of different tree structures with $t$ leaf nodes is roughly
bounded by $d^{t-1} t!$, where $d$ is the dimension of the
problem. Such a bound is also used in the proof to characterize the complexity of the tree space.

\subsection{Lower Truncation of Log-Likelihood Ratios}
\label{subsec:truncation}
We cite two results by \cite{wong1995}
to demonstrate
that after the truncation, $\widetilde{Z}_f$ still maintains some key properties
of the log-likelihood ratio. This guarantees that the behavior of the
empirical process
indexed by the truncated log-likelihood ratios does not differ much
from that of the original one, and at the same time some existing
techniques which fail for the original process now can be
applied after the truncation. The proofs are omitted here.

The first lemma shows that after the truncation, the log-likelihood ratio
still has negative expected value. 
\begin{lemma}
\label{lemma:boundofexpectation}
  Let $ \gamma= 2\exp(-\tau/2) / (1-\exp(-\tau/2))^2$. Then 
  \begin{eqnarray}
    \label{eq:boundofexpectation}
    \mathbb{E} \widetilde{Z}_f \leq -(1-\gamma) \| f^{1/2} - f_0 ^{1/2} \|_2^2.
  \end{eqnarray}
\end{lemma}
\begin{proof}
  See the paper by \cite{wong1995}, Lemma 4 in Section 2.
\end{proof}

The second lemma provides a
one-sided large deviation inequality for the empirical process indexed by the
lower-truncated log-likelihood ratios. 
\begin{lemma}
\label{lemma:largedeviation} (One-sided large deviation inequality.)
 Let $\nu_n (\widetilde{Z}_f) =
 n^{-1/2} \sum_{i=1}^{n} (\widetilde{Z}_f(Y_i) -\mathbb{E} \widetilde{Z}_f(Y_i))$. Let
 $\mathcal{G}$ be a class of densities with bracketing Hellinger metric
entropy $\bentro (u, \mathcal{G}, \rho)$. For $\epsilon>0$, consider the empirical
 process 
 \begin{equation}
   \nonumber \{ \nu_n(\widetilde{Z}_f): f \in \mathcal{G}, \rho(f, f_0)
   \leq \epsilon \}
 \end{equation}
induced by the truncated log-likelihood ratios for $f \in \mathcal{G}$
inside a Hellinger ball around $f_0$.\\
For any $\epsilon>0$, $0<b<1$ and $L>0$, let 
\begin{equation}
  \nonumber \varphi(L, \epsilon^2, n) = L^2 / 8(8 c_0 \epsilon^2 + L/ n^{1/2}),
\end{equation}
where $c_0$ is set to be $(\exp(\tau/2) -1-
\tau/2)/(1-\exp(-\tau/2))^2$. Assume that 
\begin{equation}
  \label{eq:condition1}
  L \leq b n^{1/2} \epsilon^2 /4,
\end{equation}
and 
\begin{equation}
  \label{eq:condition2}
  \int_{b L/ (32n^{1/2})}^{\epsilon} \l  \bentro (u/(2 \exp(\tau/2)),
  \mathcal{G}, \rho) \r^{1/2}  du \leq L b^{3/2} / (2^{10}(c_0+1/8)).
\end{equation}
Then
\begin{equation}
\nonumber \mathbb{P}_{f_0} \left( \sup_{ \{ \|f^{1/2} - f_0^{1/2} \|_2
    \leq \epsilon, f \in \mathcal{G}
    \} } \nu_n( \widetilde{Z}_f) \geq L \right) \leq 3 \exp( -(1-b) \varphi
  (L, \epsilon^2, n)),
\end{equation}
where $\mathbb{P}_{f_0}$ is understood to be the outer probability
measure under $f_0$.
\end{lemma}
\begin{proof}
  See the paper by \cite{wong1995}, Lemma 7 in Section 2.
\end{proof}

\subsection{Bounds for Metrics}
\label{subsec:metricbound}
In this subsection, we study relations between different types of metrics.

\begin{lemma}
\label{lemma:suptoH}
For any two densities $f$ and $g$ defined on $\Omega$, if $\|f -g \|_\infty \leq \epsilon^2$, then $\rho(f, g) \leq 2 \epsilon$.
\end{lemma}
\begin{proof}
 According to the definition,
 \begin{eqnarray*}
   \rho^2(f, g) &=& \int_\Omega (\sqrt{f}- \sqrt{g} )^2
= \int_{f < \epsilon^2/2} (\sqrt{f}- \sqrt{g} )^2
+ \int_{f \geq \epsilon^2/2} (\sqrt{f}- \sqrt{g} )^2\\
& \leq & \epsilon^2/2 + \epsilon^2/2+ \epsilon^2 
+ \int_{f \geq \epsilon^2/2} (f-g)^2/ (\sqrt{f}+\sqrt{g} )^2,
 \end{eqnarray*}
where the inequality is based on the fact that when $f < \epsilon^2/2$, $g < \epsilon^2/2 + \epsilon^2$. Note that when $f \geq \epsilon^2/2$, $(\sqrt{f}+\sqrt{g} )^2 \geq \epsilon^2/2$, then 
  $ \rho^2(f, g) \leq 2 \epsilon^2 + \epsilon^4/(\epsilon^2/2) = 4 \epsilon^2$. This finishes the proof. 
\end{proof}

\begin{lemma}
\label{lemma:kltoH}
Assume $ L^2 \leq f/g \leq U^2$ for some
  constants $U \geq 1 \geq  L >0$, and
  $\rho(f, g) \leq \epsilon$. Then
  $\kl (f, g) \leq \max\{ 1, \log U, - \log L\} \cdot 6 \epsilon^2$ and $\int
  f \left( \log (f /g) \right)^2 \leq \max\{1, ( U/(U-1)
    \log U )^2 \} \cdot 4 \epsilon^2$.
\end{lemma}
\begin{proof}
  Note that
  \begin{eqnarray*}
    \kl (f, g) &=& \int f \log (f/g)
    = \int (f-g) \log \left(f/g\right) + \int g \log ( 1 + (f-g)/g)\\
               &\leq& \int (f-g) \log \left(f/g\right)
                      + \int g (f-g)/g  = \int (f-g) \log \left( f/g
                      \right),
  \end{eqnarray*}
  where the inequality is based on the fact that $\log (1+x) \leq x$
  and $g$ is nonnegative. For the first part, following the argument
  by \cite{wong1995}, let $y = (f/g)^{1/2}-1$ and write $a = (f^{1/2} -
  g^{1/2})^2 = g y^2$ and $b = (f-g) \log (f/g) = 2 g y (y+2) \log
  (1+y)$. For $0 < y < 1$, $\log(1+y) / y < 1$; and for $y \geq 1$,
  $(y+2)/y \leq 3$, then we have
  \begin{eqnarray*}
    b/a = 2( (y+2) \log (y+1))/y \leq
    \begin{cases}
      6, & \mbox{when } 0 < y < 1, \\
      6 \log (y+1), &\mbox{when } y \geq 1. 
    \end{cases}
  \end{eqnarray*}
  Hence, for any $K > 1$, we have
  \begin{eqnarray*}
    \int_{1 < f/g \leq K^2} b = \int_{0 < y \leq K-1} b \leq  \max\{
    1, \log K \} \cdot 6\epsilon^2.
  \end{eqnarray*}
  As the definition of $a$ and $b$ is symmetric in $f$ and $g$,
  similarly, $\int_{ L^{2} \leq f/g <1} b \leq  \max\{ 1, - \log L
  \} \cdot 6\epsilon^2.$ 
  Then
  \begin{eqnarray*}
    \kl (f, g) &\leq&
    \int_{L^{2}\leq f/g < 1 } (f-g) \log (f/g)
    + \int_{1 < f/g \leq U^2 } (f-g) \log (f/g)\\
    &\leq &\max\{ 1, \log U, - \log L\} \cdot 6\epsilon^2. 
  \end{eqnarray*}

  For the second inequality, we decompose the integral as
  \begin{eqnarray*}
    \int f ( \log (f/g))^2 &=&
    \int_{f/g >1} f ( 2 \log \sqrt{f/g})^2
    +
    \int_{f/g <1} f ( 2 \log \sqrt{g/f})^2\\
    &\leq&
           \int_{f/g >1} f ( \log f/g)^2
           +\int_{f/g < 1} 4 f (\sqrt{g/f}-1)^2.
  \end{eqnarray*}
For the first part of the integral, it is not hard to see that the
function $ \log x /(1 - 1/ \sqrt{x})$ is increasing in $x$ when $x >
1$. Therefore, $\log (f/g) /(1 - \sqrt{g/f}) \leq \log (U^2) / (1 - (1/U))$.
  Then the desired result follows. 
\end{proof}

\subsection{Other Auxiliary Lemmas}
When we analyze the Gibbs posterior, the following lemma is applied several times to bound how fast the prior or the Gibbs posterior decays as the forest depth increases. 
\begin{lemma}
  \label{lemma:tooltail}
  For any $\gamma >0$, and integers $1\leq \D_1 < \D_2$, 
  \begin{eqnarray*}
    \sum_{\D=\D_1}^{\D_2} \exp\left( - \gamma \D \cdot 2^\D \right)
    \leq
    (1+\frac{1}{\gamma}) \exp\left(-\gamma \D_1 \cdot 2^{\D_1} \right).
  \end{eqnarray*}
\end{lemma}
\begin{proof}
Note that 
  \begin{eqnarray*}
    \sum_{\D=\D_1}^{\D_2}  \exp\left( - \gamma \D \cdot 2^\D \right)
    &=& \exp\left(-\gamma \D_1 \cdot 2^{\D_1}\right)
        + \sum_{\D=\D_1+1}^{\D_2} \exp\left( - \gamma \D \cdot 2^\D \right)\\
    &\leq& \exp\left(-\gamma \D_1 \cdot 2^{\D_1}\right)
           + \int_{\D_1}^{\D_2} \exp\left( - \gamma u \cdot 2^u \right) du.
  \end{eqnarray*}
  For the second term, let $t= 2^u$, and $t_1 = 2^{\D_1}$ and $t_2= 2^{\D_2}$, then
  \begin{eqnarray*}
    \int_{\D_1}^{\D_2}  \exp( - \gamma u \cdot 2^u) du =
    \int_{t_1}^{t_2} \frac{1}{t \log 2} \exp( - \gamma t \log t /
    \log 2 ) dt.
  \end{eqnarray*}
  When $\D_1 \geq 1$, $1/t \leq (\log t +1)$. Therefore,
  \begin{eqnarray*}
    \int_{\D_1}^{\D_2} \exp( - \gamma u \cdot 2^u) du
    &\leq&
     \int_{t_1}^{t_2} (\log t +1) / \log 2 \cdot \exp( - \gamma t \log t /
    \log 2 ) dt\\
&\leq& - \frac{1}{\gamma} \exp( - \gamma t \log t /\log 2) \Big|_{t_1}^{t_2} 
\leq
       \frac{1}{\gamma}  \exp\left(-\gamma \D_1 \cdot 2^{\D_1}\right).
  \end{eqnarray*}
The desired result follows. 
\end{proof}



\section{Inference for {BDF} and GDF with Perturbations}
\label{sec:algs}
\subsection{Inference for {Bayesian Density Forest}}
\label{subsec:algbdf}
To make Bayesian inference for the forest estimator introduced in Section \ref{sec:BMA}, we need to draw multiple tree-supported density samples from the posterior distribution
\eqref*{eq:posteriordist}. A useful observation is that the posterior distribution on tree depth concentrates rapidly around $\D_{\beta}$ defined in \eqref*{eq:depthbeta} when $f_0$
enjoys good regularity properties. This has been rigorously proved in Lemma \ref{lemma:BMAdepth}. Consequently, the unit for shifting the sample space should pair with
the depth of the tree as $\iota_\beta = 2^{-2 \D_\beta /d}$. However, when implementing \writing{BDF}, we are unable to calculate the value of $\D_\beta$ exactly, as it depends not only on $n$, $d$, and $\beta$, but also on a constant $c_\appr$. The constant may change sharply for different density functions, leaving the value of $\D_\beta$ inaccessible. In turn, it is impossible to calculate the unit for the shift $\iota_\beta$ based on the one-to-one correspondence between the depth of the tree $\D_\beta$ and $\iota_\beta$. The tree growth strategy that we use in practice is to set a large upper bound for tree depth and allow the algorithm to learn when to stop based on the data. The value of $\iota_\beta$ is set using \eqref*{eq:depthbeta} and an estimate of $c_\appr$, like $c_\appr =5$. In particular, we apply the following Algorithm \ref{algo:bma} to draw tree samples from the posterior, and output the forest. 
\begin{algorithm}[t]
\caption{Bayesian Density Forest}
\label{algo:bma}
\begin{algorithmic}[1]
\REQUIRE A simple random sample $Y_i \iid f_0, i = 1, \ldots, n$. The
hyper-parameters $\beta$ and $\gamma$ for tree depth. Hyper-parameters
of the OPT prior $\Pi_n$. 
\STATE Set $\widehat{\D} = \left \lfloor \log_2 (25 d n/ \log
  n)^{\frac{1}{\beta + 2/d}} \right \rfloor$, $\iota_\beta = 2^{-2
  \widehat{\D}/d}$, and $B = O(2^{\widehat{\D}})$.
\FOR{$b= 1, 2, \cdots, B$}
      \STATE Randomly sample a shifting direction $h \in
      \mathcal{H}_{\widehat{\D}} = \{(i_l \cdot \iota_\beta )_{1 \leq i \leq d}:
      i_l \in \mathbb{N}, 0 \leq i_l < 2^{\lfloor \widehat{\D} /d \rfloor} \}$.
      \STATE Calculate the posterior mean $\tilde{f}_b$ under the OPT
      prior under perturbation $h$, with (or without) the strong shape constraint. 
\ENDFOR
\STATE Out put the forest estimate $\tilde{f}_\en = \frac{1}{B} \sum_{b=1}^B \tilde{f}_b$.
\end{algorithmic}
\end{algorithm}

When implementing the algorithm, to calculate the posterior mean $\tilde{f}_b$ under the OPT prior, we apply an algorithm introduced in Section \ref{sec:optsampling}.

\subsection{Generalized Bayesian Inference for Gibbs Density Forest}
\label{subsec:alggdf}
To sample from the Gibbs posterior on forest topologies, the main difficulty still lies in efficiently searching the forest topology space. To avoid directly searching over the over-complex forest space, we define an alternative risk function associated with a depth-$\D$ forest (with perturbations $\{h_b\}_{1 \leq b \leq B}$ from the perturbation scheme $\bm \upsilon$) as 
\begin{eqnarray}
  \label{eq:proposal}
   R_n' \btriplet = -\sum_{b=1}^B \frac{1}{n_1} \sum_{i=n_0+1}^n
  \log \left(  
\tilde{f}_{\T_b^{h_b}}
  (Y_i) \right).
\end{eqnarray}
We sample from a proposal distribution defined as the Gibbs posterior induced by $R_n'$. 
Given the forest depth $\D$ and the perturbation design scheme $\bm \upsilon$, when prior distributions on trees are independent, the Gibbs posterior induced by the risk $R_n'$ can be factorized as a product of the posterior distributions on trees, allowing us to sample one tree at a time. More specifically, we consider a proposal distribution for sampling as the Gibbs posterior associated with $R_n'$ under the prior $q'$, denoted as $\lambda'$. When $q'$ is independent for trees within the forest, the posterior on a tree $\T_b^{h_b}$ can be written as, 
\begin{eqnarray*}
  \lambda'  (\T_b^{h_b} | \D, \bm \upsilon, \Y_n)  
  \propto  
q'(h_b | \D, \bm \upsilon) 
q' (\T_b^{h_b} | \D, \bm \upsilon, h_b) 
\exp \left( \frac{ \psi_\tr \cdot n}{n_1} \sum_{i=n_0+1}^n \log \left(\tilde{f}_{\T_b^{h_b}} (Y_i) \right) \right),
\end{eqnarray*}
where $q'(\cdot | \D, \bm \upsilon, h_b)$ denotes the shifted OPT
prior with a set of potentially different hyper-parameters, $q'(h_b,
\D, \bm \upsilon)$ is usually the same as $q(h_b, \D, \bm \upsilon)$,
and $\psi_\tr$ is the ``temperature'' associated with the Gibbs
posterior $\lambda' $. After we sample $B$ trees from this
posterior, at the forest level, we obtain a sample of forests from the
target Gibbs posterior $\lambda^\g$ with weight
{\small
\begin{eqnarray}
\label{eq:weight}
&&  v \btriplet \\
\nonumber && \hspace{0.5cm} 
=\frac{
\prod_{b=1}^B q(h_b, \D, \bm \upsilon) q ( \T_b^{h_b} | \D, \bm \upsilon, h_b)
\cdot \exp \left( 
\frac{\psi_\en \cdot n}{n_1} \sum_{i=n_0+1}^n \log ( 
\tilde{f}_{\{ \T_b^{h_b}\}_{b=1}^B} (Y_i) )
\right)}
{\prod_{b=1}^B q'(h_b, \D, \bm \upsilon) q' (\T_b^{h_b} | \D, \bm \upsilon, h_b) 
\cdot 
\prod_{b=1}^B \exp \left( 
\frac{\psi_\tr \cdot n}{n_1} \sum_{i=n_0+1}^n \log ( 
\tilde{f}_{\T_b^{h_b}} (Y_i) )
\right)},
\end{eqnarray}
}
where $ \psi_\en$ is the temperature associated with the target Gibbs
posterior $\lambda^\g $. Each term involved in the proposal
distribution $\lambda' $ or the target Gibbs posterior
$\lambda^\g$ has an explicit form. We provide more details about
how to sample a tree from $\lambda'$ in Appendix~\ref{sec:optsampling}. 

When we make Bayesian inference, we also need to explore the landscape of the Gibbs posterior on forest depth $\D$ and perturbation design scheme $\bm \upsilon$. As the dimension $d$ increases, the computational cost of an exhaustive search over both $\D$ and $\bm \upsilon$ can be prohibitively high. To reduce the cost, we make the following two simplifications when implementing the generalized Bayesian forest. First, we do not aggregate forests of different depths. Instead, we implicitly set a large upper limit for forest depth and allow the tree growth algorithm (for sampling from $\lambda'$) to learn an ``optimal'' tree depth based on data. Second, we only consider perturbation design schemes $\bm \upsilon = (\upsilon_l)_{1 \leq l \leq d}$ with $\upsilon_1 = \cdots =\upsilon_d$. In summary, we use the following Algorithm \ref{algo:gibbs} to build \writing{GDFs}.

\begin{algorithm}[t]
\caption{Gibbs Density Forest}
\label{algo:gibbs}
\begin{algorithmic}[1]
\REQUIRE A simple random sample $Y_i \iid f_0, i = 1, \ldots,
n$. Prior $q$ on forest topologies. Hyper-parameters for the OPT prior $\Pi_n'$ associated with the
proposal distribution $\lambda'$. Temperatures $\psi_\tr$ and $\psi_\en$. Size of
the forest $B$. Number of repetitions $n_{\text{rep}}$. 
\STATE Set $ k_{\min} = \max \{ \lfloor 
\log_2 (16 d n)^{\frac{1}{d+ 4}}
\rfloor, 
1 \} $, and
$k_{\max} =  \max \{ \left \lfloor (\log_2 n )/d \right \rfloor,
k_{\min}+5 \}$.
\FOR{$i = 1, 2, \cdots, n_{\text{rep}}$}
\FOR{$ k = k_{\min}, k_{\min}+1, \cdots, k_{\max}$}
\STATE Set $\bm \upsilon = (\upsilon_l)_{1 \leq l \leq d}$ with
$\upsilon_1 = \cdots =\upsilon_d = 2^{-2k}$. 
\FOR{$b= 1, 2, \cdots, B$}
      \STATE Randomly sample a shifting direction $h \in
      \mathcal{H}_{\bm \upsilon} = \{ (i_l \upsilon_l)_{1 \leq l\leq
        d}: i_l \in \mathbb{N}, 0 \leq i_l < 2^k \}$.
      \STATE Sample a tree $\tilde{f}_b$ from $\lambda'$ under the
      shifted OPT prior with (or without) the weak shape constraint. 
\ENDFOR
\STATE Calculate the forest estimate $\tilde{f}_{\sE} =
\frac{1}{B} \sum_{b=1}^B \tilde{f}_b$ and the weight associated with the forest according to \eqref*{eq:weight}.
\ENDFOR
\STATE Output the \writing{predictive GDF} or the \writing{MAP GDF}.
\ENDFOR
\end{algorithmic}
\end{algorithm}

\subsection{Tree Growth Under the OPT Prior}
\label{sec:optsampling}
When we grow a tree under the OPT prior, we apply a sequential
algorithm. Each tree with $t$ leaf nodes $\mathcal{A}_t = \{ \Omega_j \}_{j=1}^t$ is obtained by
recursively partitioning the sample space. We can use a sequence of
trees 
$\mathcal{A}_1, \A_2, \cdots, \A_t$ to keep track of the path leading
to $\A_t$. For any intermediate tree $\A_k$, $1 \leq k \leq t$, we
assume that all leaf nodes have been stopped. Let $\Pi_n(\cdot)$ denote the posterior
distribution $\Pi_n ( \cdot |\Y_n)$ for simplicity. Then
$\Pi_n(\A_t)$ can be decomposed as 
\begin{eqnarray*}
  \Pi_n ( \A_t) = \Pi_n ( \A_1) \Pi_n ( \A_2 | \A_1) \cdots \Pi_n(
  \A_t | \A_{t-1}).
\end{eqnarray*}
The conditional distribution $\Pi_n ( \A_{k+1} | \A_k)$ can be
calculated by $\Pi_n( \A_{k+1}) / \Pi_n(\A_{k})$. However, the
computation of the marginal distribution $\Pi_n(\A_k)$ is sometimes
infeasible. In the following, we discuss how to choose a further
split at each step from a proposal distribution, and how to adjust the
weight associated with each tree sample. The general strategy is to
consider a proposal distribution that looks one-step ahead. 

 Assume $A$ is a leaf node in $\A_k$ ($1 \leq k \leq t$), and could be
 potentially split into two new leaf nodes  $A_{\Lambda(A), 0}$ and
 $A_{\Lambda(A), 1}$ along a randomly chosen dimension $\Lambda(A) \in
 [d]$. After this split, the tree grows into
 $\A_{k+1}$. Then the Bayes factor associated with the two tree models
 $\A_k$ and $\A_{k+1}$ is simply 
 \begin{eqnarray*}
   \text{BF}(\A_k, \A_{k+1}) = 2^{n(A)} \cdot \frac{B \left(n(A_{\Lambda(A), 0})+1/2,
n(A_{\Lambda(A), 1})+1/2 \right )}
{B(1/2, 1/2)},
 \end{eqnarray*}
where for any region $A$, $n(A)$ denotes the number of data points
therein, and $B(\cdot, \cdot)$ is the beta function. We sample
the split from a proposal distribution determined by the Bayes factor,
the split of a leaf node $A$ along dimension $l$ is sampled with
probability 
\begin{eqnarray*}
  \Pi_n^k (A, l) \propto 2^{n(A)} \cdot \frac{B \left(n(A_{l, 0})+1/2,
n(A_{l, 1})+1/2 \right)}
{B(1/2, 1/2)},
 \end{eqnarray*}
where $\Pi_n^k$ denotes the proposal for the tree growth from $\A_k$
to $\A_{k+1}$.  
We can also apply a smoothing of $\Pi_n^k$ to make the tree
structure diversified, similar to that applied in the paper by \cite{LJW}. 

To update the weight associated with $\A_k$'s, we denote the tree grow from $\A_k$ to
$\A_{k+1}$ as $\A_k \rightarrow \A_{k+1}$. Then for the 
tree $\A_t$, the associated weight accounting from sampling from the
proposal distribution should be
\begin{eqnarray*}
  w_{t} (\A_t) = 
\frac{{\Pi_n (\A_1)}
\Pi_n (\A_1 \rightarrow \A_2 | \A_1) \cdots
\Pi_n (\A_{t-1} \rightarrow \A_t |\A_{t-1} )  }
{\Pi_n^1 (\A_1)
\Pi_n^2 (\A_1 \rightarrow \A_2 | \A_1) \cdots
\Pi^\D_n  (\A_{t-1} \rightarrow \A_t |\A_{t-1} ) }.
\end{eqnarray*}
Then we can calculate the weight in a stepwise manner alternatively as 
\begin{eqnarray*}
  w_{k} (\A_k) =
  \begin{cases}
    1, & \mbox{when } k=1,\\
w_{k-1} ( \A_{k-1}) u_k (\A_{k-1} \rightarrow \A_k), & \mbox{when }
k > 1.
  \end{cases}
 \end{eqnarray*}
%
And the adjustment of weight $u_k(\cdot)$ can be calculated as 
\begin{eqnarray*}
  u_k (\A_{k-1} \rightarrow \A_k)
\propto 
\frac{\Pi_n (\A_k)}
{\Pi_n (\A_{k-1})},
\end{eqnarray*}
where the constant does not depend on a specific tree topology. After
some simple calculation, we can show that if node $A$ is picked up to
be split into $A_{\Lambda(A), 0}$ and $A_{\Lambda(A), 1}$ in this
step, 
\begin{eqnarray*}
&&  u_k (\A_{k-1} \rightarrow \A_k)
\propto 
\frac{\tau(A) (1-\tau(A_{\Lambda(A), 0})) (1- \tau(A_{\Lambda(A), 1}))}
   {d(1- \tau(A))}\\
  && \hspace{2in}\times
\frac{2^{n(A)} \cdot \frac{B \left(n(A_{\Lambda(A), 0})+1/2,
n(A_{\Lambda(A), 1})+1/2 \right)}
{B(1/2, 1/2)}}
{\Pi_n^{k-1}(A, \Lambda(A))}.
\end{eqnarray*}
\noindent Alternatively, we can assign each piece in the weight update a name,
and write the weight update as, 
\begin{eqnarray*}
  u_k (\A_{k-1} \rightarrow \A_k)
\propto 
\text{prior ratio} 
\cdot 
\frac{\text{BF}(\A_{k-1}, \A_{k})}
{\text{proposal}}.
\end{eqnarray*}

To sample from the posterior distribution on trees under the OPT
prior, we usually grow multiple trees simultaneously. When
implementing the algorithm, we also employ a resampling and pruning
scheme. We refer to the paper by \cite{LJW} for more details.

For generalized Bayesian inference, our algorithm is almost the same as that for likelihood-based Bayes, except for the way we calculate the ``Bayes factor''. For risk based Bayesian inference, as the risk function associated with each tree has an explicit form, the risk-associated "Bayes factor" is very easy to calculate. In specific, if a leaf node $A$ in $\A_k$ is further split into $A_{\Lambda(A),0}$ and $A_{\Lambda(A),1}$, leading to a new tree structure $\A_{k+1}$, then the following ratio 
\begin{eqnarray*}
\frac{
\exp \left(
\frac{n \psi_\tr}{n_1} \sum_{i=n_0+1}^n 
\log ( \tilde{f}_{\A_{k+1}} (Y_i) ) 
\right)
}
{\exp \left(
\frac{n \psi_\tr}{n_1} \sum_{i=n_0+1}^n 
\log ( \tilde{f}_{\A_k} (Y_i) ) 
\right)} 
\end{eqnarray*}
has an explicit form. Our sequential algorithm is the same as that for likelihood-based Bayesian inference, except for that we replace the Bayes factor the the ratio presented above. 

\vspace{3mm}
In order to make the data points as uniform as possible, we apply a copula transformation to each variable in advance whenever the dimension exceeds 3. More specifically, we
estimate the marginal distribution of each variable $Y_j$ by our
approach, denoted as $\hat{f}_j$ (we use $\hat{F}_j$ to denote the
c.d.f. of $Y_j$), and transform each point $(y^1, \cdots, y^d)$ to
$(\hat{F}_1(y^1), \cdots, \hat{F}_d(y^d))$. Another advantage of this
transformation is that after the transformation the sample space
naturally becomes $[0, 1]^d$.

\subsection{Implementation Details for Simulation Studies}
\label{subsec:kdeRF}
Next we provide implementation details for different
methods considered in simulation studies. 

\subsubsection{BDF and GDF}
For \writing{BDF}, as suggested by Algorithm~\ref{algo:bma} from the
Supplementary Material, we need to fix the value of $\beta$ in advance, estimate $\D_{\beta}$ (defined in Equation
\eqref*{eq:depthbeta}) as $\left \lfloor \log_2 (25 d n /\log
n)^{\frac{1}{\beta + 2/d}} \right \rfloor$, and obtain $\iota_\beta$ via $\iota_\beta = 2^{-2 \lfloor \D_\beta/d \rfloor }$. According to the definition of $\D_\beta$, it also relies on the constant $c_{\text{approx}}$ that appears in the approximation error. 
However, in practice, the constant $c_{\text{approx}}$ is {\em unknown}. Therefore, in our simulations, we set \rev{$\beta= 1$}, $\gamma=1$, and try different values of $\D_\beta$ (and correspondingly $\iota_\beta$). 

For \writing{GDF}, as we consider another layer of prior on forest
depth and the perturbation design scheme, and as Bayesian inference
has shown to be adaptive, Algorithm~\ref{algo:gibbs} from the Supplementary Material is expected to learn an ``optimal'' design based on data. When applying the algorithm, we set the prior $q$ on the tree depth
proportional to $\exp(- \D 2^\D)$. We also study the effect of sample splitting ratio. When evaluating the empirical risk \eqref*{eq:riskdef}, the proportion of data for density learning given tree topologies ($n_0/n$) ranges from 0.5 to 0.8. 

For both \writing{BDF} and \writing{GDF}, our theoretical results
require individual trees to grow under some shape constraint. In
simulations, we found that unconstrained trees perform quite well. In
the following, we choose to grow trees within a forest without any
shape constraint, as such algorithms are closer to how tree-based
methods are used in practice and are easier to implement. In terms of
specifications of hyper-parameters of the OPT prior, to make a fair
comparison, we set them the same in all scenarios as \rev{$c_\tau =
0.5$, $\nu = 0.3$ when $d=1$, $\nu=0.5$ when $d > 1$}, and the pseudo counts for beta distributions to be
0.5.

For Example~\ref{eg:mixNorm}, the sample size is $1 \times 10^4$. With $d$ ranging from 5, 10 to 30,
we set the forest size at 200, \rev{ and $\psi = 0.5, 0.3$, and $0.1$,
respectively}. For each case,
we run \rev{20} replicates to calculate the KL divergence. The
specification of hyper-parameters for the three one-dimensional
examples is provided together with the numerical results in
Section~\ref{subsec:1Dsimulation}.

\subsubsection{KDE}
\rev{In our simulation studies, we compare our method with the kernel
density estimation (KDE) method.
Let $Y_1,\dots,Y_n \in \mathbb{R}^d$ denote the training sample.}

\rev{In the univariate case ($d=1$), the estimator takes the form
\[
\widehat f_h(x)
=
\frac{1}{n h}
\sum_{i=1}^{n}
\phi\!\left(\frac{x - Y_i}{h}\right),
\]
where $\phi(\cdot)$ denotes the standard normal density and $h>0$ is the bandwidth parameter. 
The bandwidth is selected via $K$-fold cross-validation (with $K=2$ in our implementation). 
Specifically, for each candidate bandwidth $h$, we partition the data into $K$ folds. 
For each fold $k$, we construct the estimator $\widehat f_h^{(-k)}$ using the observations excluding fold $k$, and evaluate it on the held-out fold. 
The cross-validation criterion is the average negative log-likelihood,
\[
\mathrm{CV}(h)
=
-\frac{1}{n}
\sum_{k=1}^{K}
\sum_{Y_i \in \mathcal{I}_k}
\log \widehat f_h^{(-k)}(Y_i),
\]
where $\mathcal{I}_k$ denotes the index set of fold $k$. 
We select the bandwidth minimizing $\mathrm{CV}(h)$ over a grid of candidate values of the form 
$h = n^{-\lambda/(d+4)}$,
where $\lambda$ varies over a prespecified range.}

\rev{In the multivariate case ($d \ge 2$), we employ a Gaussian kernel with a full bandwidth matrix $H$. 
The estimator is defined as
\[
\widehat f_H(x)
=
\frac{1}{n (2\pi)^{d/2} |H|^{1/2}}
\sum_{i=1}^{n}
\exp\!\left\{
-\frac{1}{2}
(x - Y_i)^\top
H^{-1}
(x - Y_i)
\right\}.
\]
Following the Scott's rule \citep{Scott1992}, the bandwidth matrix is chosen as
\[
H = n^{-2/(d+4)} \widehat{\Sigma},
\]
where $\widehat{\Sigma}$ is the sample covariance matrix of the training data. 
This choice yields an elliptically shaped kernel that adapts to
differences in scale and correlations across dimensions.}

\subsubsection{Unsupervised RF}

\rev{In addition, we also compare our method with a tree-based method for
density learning \citep{pmlr-wen22c}, which we refer to as
the unsupervised random forest (unsupervised RF) in the paper.}

\rev{For this method, each tree is built by recursively splitting the
sample space along 
a randomly selected dimension at the midpoint of the current 
leaf node, until the pre-specified tree depth is reached. Given independent 
observations $Y_1,\ldots,Y_n \in \mathbb{R}^d$, the density estimate from 
a single tree is defined as
\begin{eqnarray*}
  \widehat f^{(m)}(x)
=
\frac{N^{(m)}(x)}{n\, \mu \left( \Omega^{(m)}(x) \right)},
\end{eqnarray*}
where $\Omega^{(m)}(x)$ denotes the leaf node containing $x$ in the $m$-th tree, 
$\mu \left( \Omega^{(m)}(x) \right)$ is its volume, and $N^{(m)}(x)$ is the number of observations 
among $\{Y_i\}_{i=1}^n$ that fall in that node. The forest estimator averages 
over $M$ independently generated trees:
\begin{eqnarray*}
  \widehat f(x)
=
\frac{1}{M}\sum_{m=1}^M \widehat f^{(m)}(x).
\end{eqnarray*}
For unsupervised RF, the pre-specified tree depth is a tuning
parameter and reflects a tradeoff between the bias and the variance.}

\rev{In order to make a fair comparison, we also add perturbations to
individual trees within the forest, in the same way as that for GDF
with perturbations described in Section~\ref{subsec:gibbsholderrate}. Specifically, we shift the sample space or
equivalently the split points along a randomly selected direction
$h$. The unit for shift is depth-dependent. In the one-dimensional
case, for depth-$\D$ forest, the
unit is $2^{-2\D}$.}

\subsection{Numerical Results for the One-Dimensional Examples}
\label{subsec:1Dsimulation}

\begin{example}(Mixture-of-beta distribution)
  \label{eg:beta}
 Assume the distribution is
  \begin{eqnarray*}
    Y \sim 0.7 \mbox{Beta}(3, 12) + 0.3 \mbox{Beta}(300, 10).
  \end{eqnarray*}
\end{example}

For this example, we generate data of sizes $5 \times 10^2$, $1\times
10^3$, and $1 \times 10^4$. The forest size is set to be 50, 100 and
200, respectively. \rev{The inverse temperature parameter $\psi$ is set to be 0.5}. In each case, we run \rev{50} replicates. We measure the performance of different methods by calculating the average KL divergence of the estimated distribution from the true one. The results are shown in Table \ref{tab:mixBeta}.


\begin{table}[b]
  \centering
  \rev{
  \begin{tabular}{l l c c c}
    \toprule
   & & $n = 5 \times 10^2$ & $n = 1 \times 10^3$ & $n = 1 \times 10^4$
    \\
    \midrule
    {\footnotesize Optional P\'olya tree} 
        & 
            & 0.0495(0.016)
                &0.0335 (0.0090)
                    &0.0069 (0.0011) \\
    \midrule
    \multirow{5}{*}{{\footnotesize Bayesian density forest}}
        &{\footnotesize $\iota_\beta =2^{-6}$}
            &0.041 (0.012) 
                & 0.025 (0.0071) 
                    &0.00565 (0.0009)\\
        &{\footnotesize $\iota_\beta =2^{-10}$}
            & 0.028 (0.0088)
                & 0.016 (0.0045)
                    &0.0032 (0.0006)\\
         &{\footnotesize $\iota_\beta =2^{-14}$} 
            & 0.030 (0.0085)
                & 0.018 (0.0045)
                    & 0.0032 (0.0006) \\
        &{\footnotesize $\iota_\beta =0$} 
            &0.044 (0.012) 
                & 0.027 (0.0074)
                    &0.0058 (0.0010)\\
    \midrule
    {\footnotesize {Predictive GDF}}  
        &\multirow{2}*{\footnotesize $n_0/n =0.5$}
            &0.044 (0.022)
                & 0.026 (0.011)
                    & 0.0045 (0.0028) \\
    {\footnotesize {MAP GDF}} 
        &
            &0.034 (0.0095) 
                & 0.019 (0.0060)
                    &0.0037 (0.0007)\\
    \addlinespace
    {\footnotesize {Predictive GDF}}  
        &\multirow{2}*{\footnotesize $n_0/n =0.6$}
            &0.033 (0.015)
                & 0.024 (0.0145)
                    & 0.0046 (0.00385) \\
    {\footnotesize {MAP GDF}} 
        &
            &0.030 (0.0092) 
                & 0.018 (0.0061)
                    &0.0038 (0.0009)\\
     \addlinespace
    {\footnotesize {Predictive GDF}}  
        &\multirow{2}*{\footnotesize $n_0/n =0.7$}
            &0.030 (0.015)
                & 0.022 (0.013)
                    & 0.0036 (0.00065) \\
    {\footnotesize {MAP GDF}} 
        &
            &0.029 (0.0086) 
                & 0.018 (0.0053)
                    &0.0039 (0.0008)\\
     \addlinespace
    {\footnotesize {Predictive GDF}}  
        &\multirow{2}*{\footnotesize $n_0/n =0.8$}
            &0.028 (0.0069)
                & 0.020 (0.0125)
                    & 0.00415 (0.0007) \\
    {\footnotesize {MAP GDF}} 
        &
            &0.032 (0.0090) 
                & 0.019 (0.0064)
                    &0.00435 (0.0008)\\
    \midrule
    {\footnotesize KDE }&& 0.115 (0.059)& 0.059 (0.038) &0.018 (0.010)
    \\
    \midrule
    \multirow{5}{*}{{\footnotesize Unsupervised RF}}
        &{\footnotesize depth 4}
            &0.41 (0.00525) 
                & 0.40 (0.0030) 
                    &0.40 (0.0002)\\
        &{\footnotesize depth 5}
            & 0.17 (0.0070)
                & 0.16 (0.0030)
                    &0.16 (0.0003)\\
         &{\footnotesize depth 6} 
            & 0.16 (0.0086)
                & 0.15 (0.00415)
                    & 0.15 (0.0004) \\
        &{\footnotesize depth 8} 
            &0.049 (0.013) 
                & 0.0345 (0.00755)
                    &0.013(0.0009)\\
    \bottomrule
  \end{tabular}
  }
  \caption[mixBeta]{Simulation results for the mixture-of-beta
    example. For each case, \rev{50} replicates are performed to
    calculate the average KL divergence. The standard deviation over
    the \rev{50} replicates is reported in parentheses.}
  \label{tab:mixBeta}
\end{table}

\begin{figure}[b]
\centering
\includegraphics[width=0.9\linewidth]{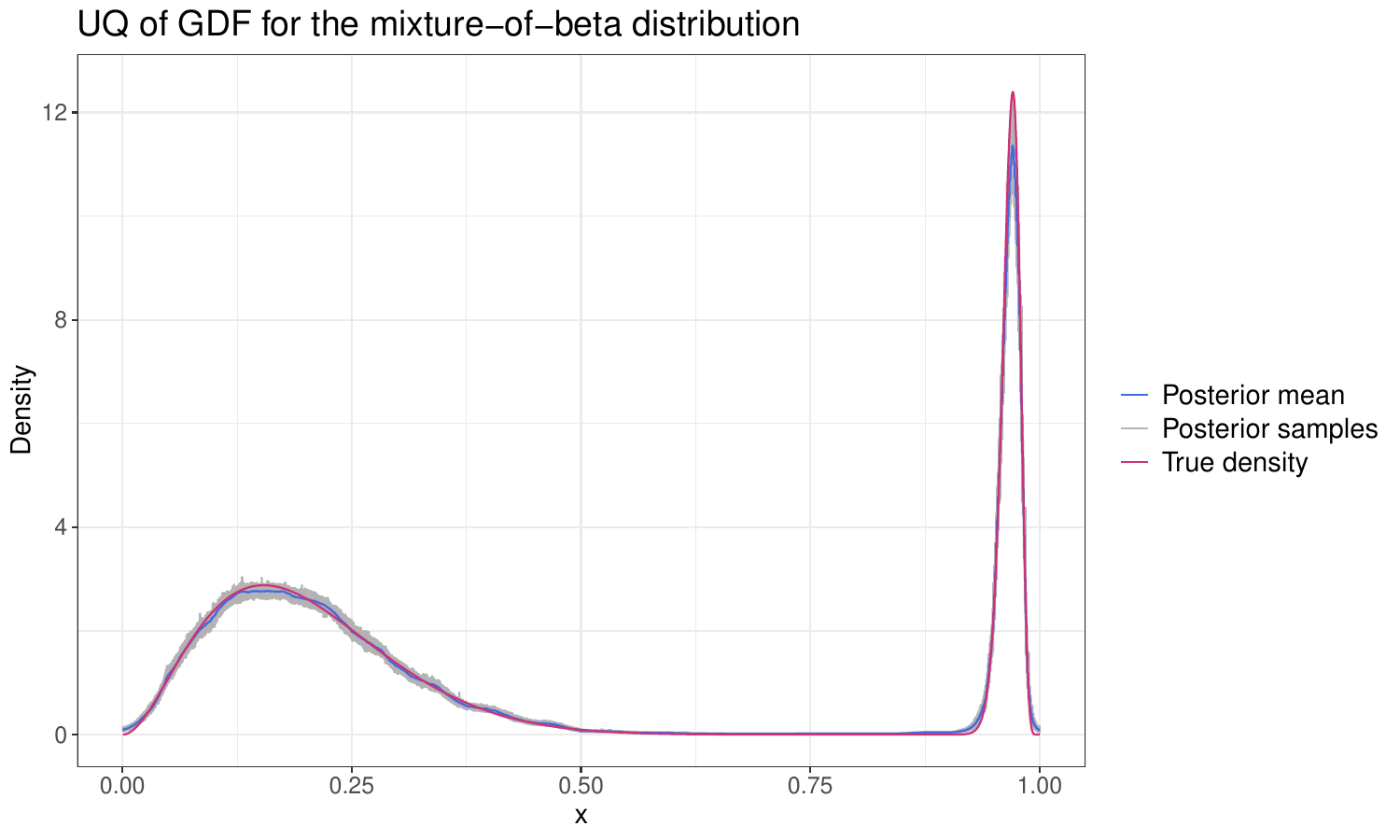}
\caption[betaUQ]{{\rev{Uncertainty quantification for the mixture-of-beta distribution.}}}
\label{fig:beta}
\end{figure}

For the following two distributions, the sample size and the choice of
forest size for each scenario are exactly the same as those for
Example \ref{eg:beta}. \rev{When the sample size is $5 \times 10^2$,
  $1\times 10^3$, and $1 \times 10^4$, we set the inverse temperature
parameter $\psi$ to be 0.5, 0.2, and 0.1, respectively.} We still run \rev{50} replicates for each case and report the average KL and the standard deviation over the \rev{50} replicates. 

\begin{example} (Trigonometric distribution)
  \label{eg:trig}
 Assume the density is
  \begin{eqnarray*}
    f(x) = 1+ \sin(2 \pi x - \pi/2).
  \end{eqnarray*}
\end{example}
The numerical results for this example are summarized in Table
\ref{tab:trig}.

\begin{table}[h]
  \centering
  \rev{
  \begin{tabular}{l l c c c}
    \toprule
   & & $n = 5 \times 10^2$ & $n = 1 \times 10^3$ & $n = 1 \times 10^4$
    \\
    \midrule
    {\footnotesize Optional P\'olya tree} 
        & 
            & 0.038 (0.019)
                &0.0265 (0.0091)
                    &0.0044 (0.0012) \\
    \midrule
    \multirow{5}{*}{{\footnotesize Bayesian density forest}}
        &{\footnotesize $\iota_\beta =2^{-6}$}
            &0.029 (0.011) 
                & 0.016 (0.0062) 
                    &0.0028 (0.0007)\\
    &{\footnotesize $\iota_\beta =2^{-10}$}
        & 0.023 (0.0086)
            & 0.013 (0.0040)
                &0.00195 (0.0005)\\
     &{\footnotesize $\iota_\beta =2^{-14}$} 
        &0.027 (0.010) 
            & 0.015 (0.0044)
                &0.0025 (0.0005)\\
    &{\footnotesize $\iota_\beta =0$} 
        &0.033 (0.015) 
            & 0.0185 (0.0077)
                &0.0032 (0.0008)\\
    \midrule
    {\footnotesize {Predictive GDF}} 
        & \multirow{2}{*}{\footnotesize $n_0/n =0.5$} 
            & 0.0185 (0.0070)
                & 0.013 (0.0050) 
                    & 0.0020 (0.0006) \\
    {\footnotesize {MAP GDF}} 
        & 
           & 0.020 (0.0071)
                & 0.012 (0.0035) 
                    & 0.0019 (0.0006) \\
    \addlinespace
   {\footnotesize {Predictive GDF}} 
        & \multirow{2}{*}{\footnotesize $n_0/n =0.6$} 
            & 0.0165 (0.00645)
                & 0.011 (0.00475) 
                    & 0.0018 (0.0005) \\
    {\footnotesize {MAP GDF}} 
        & 
            &0.020 (0.0081)
                & 0.012 (0.0043)
                    & 0.0018 (0.0006) \\
    \addlinespace
   {\footnotesize {Predictive GDF}} 
        & \multirow{2}{*}{\footnotesize $n_0/n =0.7$} 
            & 0.016 (0.0065)
                & 0.00955 (0.0029) 
                    & 0.0016 (0.0005) \\
    {\footnotesize {MAP GDF}} 
        & 
            &0.0185 (0.00855)
                & 0.011 (0.0038)
                    & 0.0017 (0.0005) \\
     \addlinespace
   {\footnotesize {Predictive GDF}} 
        & \multirow{2}{*}{\footnotesize $n_0/n =0.8$} 
            & 0.016 (0.0064)
                & 0.0095 (0.0035) 
                    & 0.0015 (0.0003) \\
    {\footnotesize {MAP GDF}} 
        & 
            &0.019 (0.0080)
                & 0.011 (0.0035)
                    & 0.0020 (0.0010) \\
    \midrule
    {\footnotesize KDE }&& 0.0079 (0.0042)& 0.00565 (0.0023) &0.0009 (0.0003)
     \\
    \midrule
    \multirow{5}{*}{{\footnotesize Unsupervised RF}}
        &{\footnotesize depth 4}
            &0.032 (0.0044) 
                & 0.0255 (0.0021) 
                    &0.23 (0.0002)\\
        &{\footnotesize depth 5}
            & 0.022 (0.0068)
                & 0.011 (0.0031)
                    &0.0063 (0.0004)\\
         &{\footnotesize depth 6} 
            & 0.026 (0.010)
                & 0.0145 (0.0052)
                    & 0.0029 (0.0006) \\
        &{\footnotesize depth 8} 
            &0.0655 (0.020) 
                & 0.046 (0.010)
                    &0.0058(0.0011)\\
    \bottomrule
  \end{tabular}
  }
  \caption[mixtrig]{Simulation results for the trigonometric example. For each sample size $n$, \rev{50} replicates are performed to calculate the average KL divergence. The standard deviation of the KL over \rev{50} replicates is reported in parentheses.}
  \label{tab:trig}
\end{table}

\begin{figure}[h]
\centering
\includegraphics[width=0.9\linewidth]{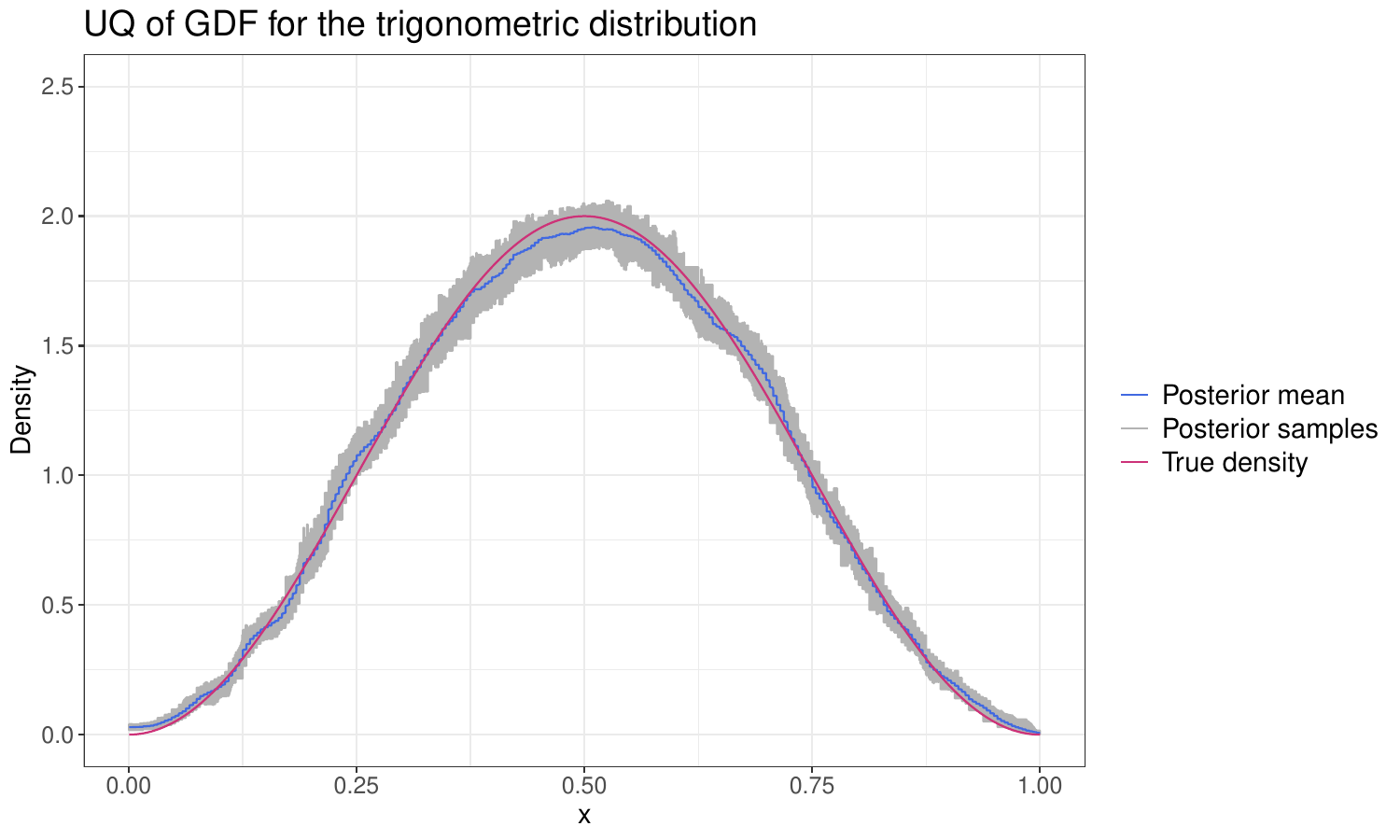}
\caption[trigUQ]{{\rev{Uncertainty quantification for the
      trigonometric distribution.}}}
\label{fig:trig}
\end{figure}

\rev{\begin{example}(Semi-circle distribution)
  \label{eg:circle}
 Assume the true density function is
  \begin{eqnarray*}
   f_0(x) =
    \begin{cases}
       0 & \mbox{for } 0 \leq x < 0.5, \\
      \frac{8}{\pi} \sqrt{(1 - (4x -3)^2 )} & \mbox{for } 0.5
      \leq x \leq 1.
    \end{cases}
  \end{eqnarray*}
\end{example}}

\rev{The numerical results for this example are summarized in Table \ref{tab:circle}.}

\begin{table}[h]
  \centering
  \rev{
  \begin{tabular}{l l c c c}
    \toprule
   & & $n = 5 \times 10^2$ & $n = 1 \times 10^3$ & $n = 1 \times 10^4$
    \\
    \midrule
    {\footnotesize Optional P\'olya tree} 
        & 
            & 0.038 (0.018)
                &0.025 (0.0084)
                    &0.0036 (0.0009) \\
    \midrule
    \multirow{5}{*}{{\footnotesize Bayesian density forest}}
        &{\footnotesize $\iota_\beta =2^{-6}$}
            &0.024 (0.011) 
                & 0.0125 (0.0053) 
                    &0.0022 (0.0005)\\
    &{\footnotesize $\iota_\beta =2^{-10}$}
        & 0.019 (0.0071)
            & 0.011 (0.0038)
                &0.0015 (0.0004)\\
     &{\footnotesize $\iota_\beta =2^{-14}$} 
        &0.021 (0.0069) 
            & 0.012 (0.0039)
                &0.00155 (0.0004)\\
    &{\footnotesize $\iota_\beta =0$} 
        &0.030 (0.013) 
                & 0.014 (0.0060) 
                    &0.0023 (0.0006)\\
    \midrule
    {\footnotesize {Predictive GDF}} 
        & \multirow{2}{*}{\footnotesize $n_0/n =0.5$} 
            & 0.0175 (0.0061)
                & 0.0105 (0.0036) 
                    & 0.0020 (0.0006) \\
    {\footnotesize {MAP GDF}} 
        & 
           & 0.019 (0.0070)
                & 0.011 (0.0045) 
                    & 0.0019 (0.0005) \\
    \addlinespace
   {\footnotesize {Predictive GDF}} 
        & \multirow{2}{*}{\footnotesize $n_0/n =0.6$} 
            & 0.015 (0.0055)
                & 0.0097 (0.0032) 
                    & 0.0019 (0.0006) \\
    {\footnotesize {MAP GDF}} 
        & 
            &0.017 (0.0061)
                & 0.010 (0.0035)
                    & 0.0019 (0.0005) \\
    \addlinespace
   {\footnotesize {Predictive GDF}} 
        & \multirow{2}{*}{\footnotesize $n_0/n =0.7$} 
            & 0.015 (0.0046)
                & 0.0093 (0.0029) 
                    & 0.0019 (0.0005) \\
    {\footnotesize {MAP GDF}} 
        & 
            &0.018 (0.0066)
                & 0.010 (0.0030)
                    & 0.0018 (0.0005) \\
     \addlinespace
   {\footnotesize {Predictive GDF}} 
        & \multirow{2}{*}{\footnotesize $n_0/n =0.8$} 
            & 0.016 (0.0049)
                & 0.010  (0.0033) 
                    & 0.0021 (0.0012) \\
    {\footnotesize {MAP GDF}} 
        & 
            &0.019 (0.0069)
                & 0.011 (0.0041)
                    & 0.0019 (0.0005) \\
    \midrule
    {\footnotesize KDE }&& 0.018 (0.0044)& 0.012 (0.0036) &0.0029 (0.0004)
     \\
    \midrule
    \multirow{5}{*}{{\footnotesize Unsupervised RF}}
        &{\footnotesize depth 4}
            &0.033 (0.0019) 
                & 0.027 (0.0015) 
                    &0.026 (0.0002)\\
        &{\footnotesize depth 5}
            & 0.020 (0.00415)
                & 0.013 (0.0020)
                    &0.011 (0.0002)\\
         &{\footnotesize depth 6} 
            & 0.022 (0.0053)
                & 0.012 (0.0031)
                    & 0.0050 (0.0003) \\
        &{\footnotesize depth 8} 
            &0.0635 (0.0145) 
                & 0.033 (0.0073)
                    &0.0038(0.0007)\\
    \bottomrule
  \end{tabular}
  }
  \caption[mixtrig]{\rev{Simulation results for the semi-circle example. For each sample size $n$, \rev{50} replicates are performed to calculate the average KL divergence. The standard deviation of the KL over \rev{50} replicates is reported in parentheses.}}
  \label{tab:circle}
\end{table}

\begin{figure}[h]
\centering
\includegraphics[width=0.9\linewidth]{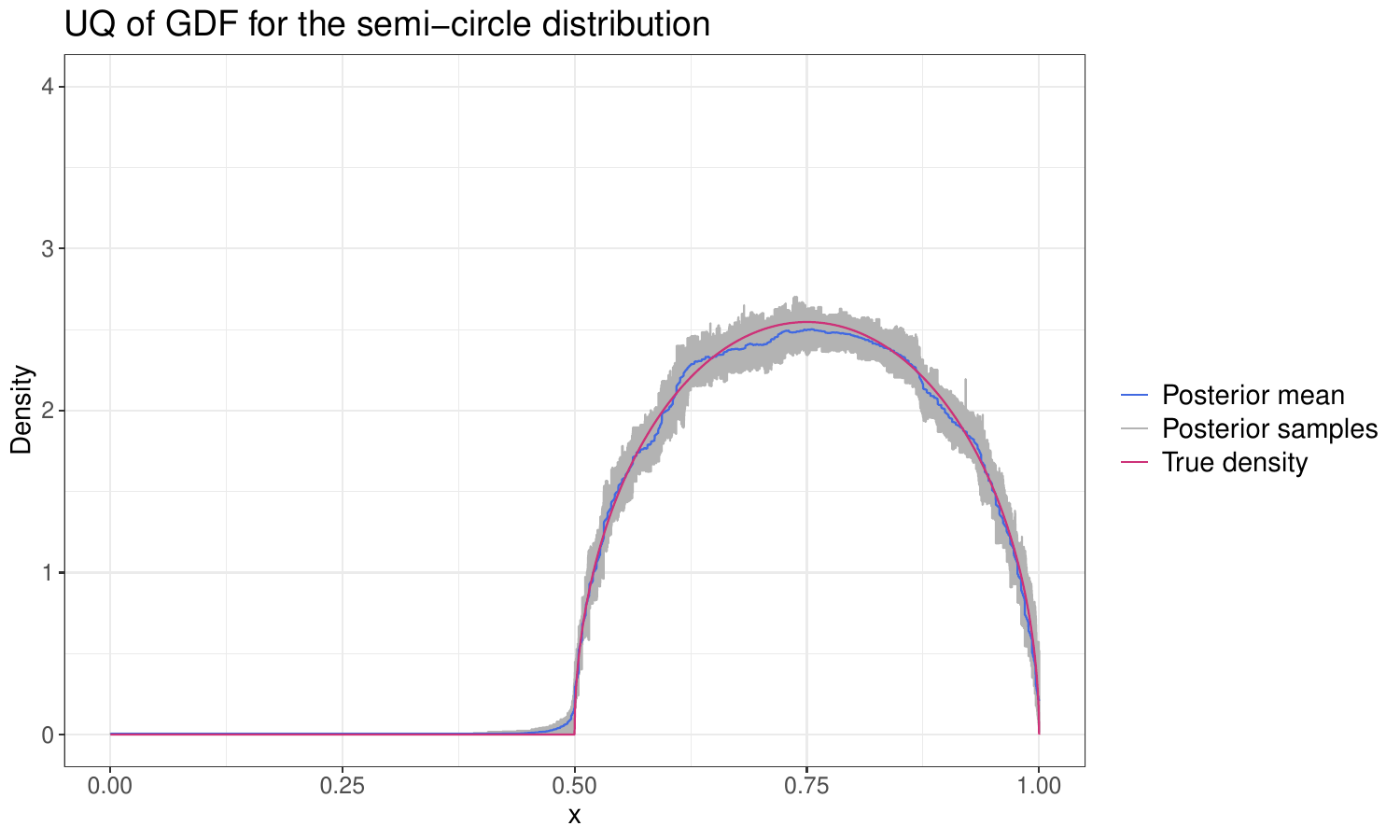}
\caption[circleUQ]{{\rev{Uncertainty quantification for the
      semi-circle distribution.}}}
\label{fig:circle}
\end{figure}

From these \rev{three} one-dimensional examples, we clearly see that when the true density has sharp variations, tree-based methods generally perform better than KDE; while when the variation of the true density is mild, KDE becomes better as it can better take advantage of the
smoothness. Both versions of the forest estimators can outperform the
single-tree method and tend to have a smaller variance.

\subsection{Uncertainty Quantification}
\label{subsec:UQ}
\rev{Under the Bayesian framework, 
  one can naturally quantify the
  uncertainty level through the full posterior distribution over
  unknown parameters. This built-in uncertainty quantification (UQ)
  also extends to the generalized Bayesian framework, in which the
  likelihood-based posterior is replaced by the risk-based Gibbs
  posterior. In this subsection, we construct a credible set
  for the density function under the generalized Bayesian framework
  through sampling.}

\rev{The learning of the posterior distribution is still based on the sample-splitting strategy. We split the full dataset into two subsets
  $\Y_n^{\en}$ and $\Y_n^{\text{leaf}}$. The subset $\Y_n^{\en}$ is
  used to learn the Gibbs posterior on forest topologies, and
  $\Y_n^{\text{leaf}}$ is used to learn the posterior distribution on
  tree-supported probability measures given the forest topology.}

  \rev{We learn the Gibbs posterior on forest topologies using the method
  described in Section~\ref{sec:gibbs}. Note
  that in this step, the
  subset $\Y_n^{\en}$ needs to be further split into two subsets.}

\rev{Given any forest topology $\newforest = \forest$, we learn the posterior distribution on
tree-supported probability measures under the OPT prior. To avoid the
difficulty of investigating correlations among trees, we further split
the dataset $\Y_n^{\text{leaf}}$. In specific,
$\Y_n^{\text{leaf}}$ is randomly split into $B$ folds, where $B$ is
the number of trees within the forest, with each fold being used for
inference for one tree under the OPT prior.  }

\rev{With these two components, we can draw samples of forest-supported
densities from the ``posterior distribution''. We first sample a
forest topology from the Gibbs posterior learned in the first
part. Then, for a sampled forest topology $\newforest$, we further draw a density function
supported by it from the distribution learned in the second
part. Specifically, the sampled density can be written as
$f_\newforest = \frac{1}{B} \sum_{b=1}^B f_{\T_b}$, where each
$f_{\T_b}$ is a probability measure sampled from the posterior distribution under
the constrained OPT prior on the tree structure $\T_b$ given the $b$-th
fold of the dataset $\Y_n^{\text{leaf}}$.}

\rev{Admittedly, this is a very simple approach to UQ. In this
  approach, we may not use data efficiently, as we split the
  dataset multiple times. Moreover, we have not yet investigated the
  frequentist coverage properties of the UQ approach from a 
  theoretical perspective. In particular, it remains an open problem
  how the temperature parameter and the sample-splitting ratio
  influence the frequentist coverage. It is unclear whether the resulting credible sets are optimal in terms of their width, either. We still
  apply this method here as a proof of concept to illustrate the
  feasibility of UQ under the generalized Bayesian framework through
  numerical results, while acknowledging that several open problems
  remain and that the approach could be further improved in multiple
  ways.}

\rev{For the three one-dimensional examples studied in
  Section~\ref{subsec:1Dsimulation} (Example~\ref{eg:beta} to
  Example~\ref{eg:circle}), we apply this method for uncertainty quantification. In our simulations, $\Y^\en_n $  is of size $1 \times
  10^4$, and $\Y^{\text{leaf}}$ is of size $4 \times 10^4$. In the
  forest, there are $B= 200$ trees. The inference data is further
  divided into $200$ folds. This results in $200$ data for
  each tree for inference. We draw $200$ sample densities from
  the posterior via the two-step procedure described above. 
  For the three distributions, the credible sets are
  visualized in Figure~\ref{fig:beta}, Figure~\ref{fig:trig}, and
  Figure~\ref{fig:circle}, respectively.}

\rev{Based on the numerical results, the credible set constructed
  under the generalized Bayesian framework has
reasonably good coverage. With inherited spatial adaptivity from 
tree-based methods, the credible set can capture sharpe changes and
local features of
the underlying distribution. However, when $f_0$ is close to zero, the coverage
is less satisfying. We think this is 
acceptable, as it is fundamentally difficult to quantify the
uncertainty level when only limited amount of data are available.}

\bibliographystyle{plainnat}
\bibliography{ref_tree}       

@ARTICLE{amit1997,
  author={Amit, Yali and Geman, Donald},
  journal={Neural Computation}, 
  title={Shape Quantization and Recognition with Randomized Trees}, 
  year={1997},
  volume={9},
  number={7},
  pages={1545--1588},
  journaldoi={10.1162/neco.1997.9.7.1545}}

@article{arlot2014analysis,
      title={Analysis of purely random forests bias}, 
      author={Sylvain Arlot and Robin Genuer},
      year={2014},
      eprint={1407.3939},
      archivePrefix={arXiv},
      journal={arXiv preprint: arXiv:1407.3939},
      primaryClass={math.ST}
}

@article{athey2019,
author = {Susan Athey and Julie Tibshirani and Stefan Wager},
title = {{Generalized random forests}},
volume = {47},
journal = {The Annals of Statistics},
number = {2},
publisher = {Institute of Mathematical Statistics},
pages = {1148--1178},
year = {2019},
journaldoi = {10.1214/18-AOS1709},
journalURL = {https://doi.org/10.1214/18-AOS1709}
}

@article{awaya2021,
  author  = {Naoki Awaya and Li Ma},
  title   = {Unsupervised Tree Boosting for Learning Probability Distributions},
  journal = {Journal of Machine Learning Research},
  year    = {2024},
  volume  = {25},
  number  = {198},
  pages   = {1--52},
  journalurl     = {http://jmlr.org/papers/v25/22-0980.html}
}

@article{barron1999,
year={1999},
journal={Probability Theory and Related Fields},
volume={113},
number={3},
title={Risk bounds for model selection via penalization},
author={Barron, Andrew and Birg{\'e}, Lucien and Massart, Pascal},
pages={301-413},
language={English}
}

@incollection{bernardo1999bayesian,
  title={Bayesian model averaging and model search strategies},
  editor={Bernardo, Jos{\'e} M and Berger, James O and Dawid, AP and Smith, Adrian FM},
  author={Clyde, Merlise A.},
  booktitle={Bayesian statistics},
  volume={6},
  pages={157--185},
  year={1999},
  publisher={Oxford University Press},
  address={Oxford, United Kingdom}
}

@article{biau2008,
author = {Biau, G\'{e}rard and Devroye, Luc and Lugosi, G\'{a}bor},
title = {Consistency of Random Forests and Other Averaging Classifiers},
year = {2008},
issue_date = {6/1/2008},
publisher = {JMLR.org},
volume = {9},
journalissn = {1532-4435},
journal = {J. Mach. Learn. Res.},
jounralmonth = {jun},
pages = {2015--2033},
numpages = {19} }

@article{biau2012analysis,
  title={Analysis of a random forests model},
  author={Biau, G{\'e}rard},
  journal={The Journal of Machine Learning Research},
  volume={13},
  pages={1063--1095},
  year={2012},
  journalpublisher={JMLR. org}
}

@book{breiman1984classification,
  title={Classification and Regression Trees},
  author={Breiman, L. and Friedman, J. and Stone, C.J. and Olshen, R.A.},
  bookisbn={9780412048418},
  booklccn={83019708},
  bookurl={https://books.google.com/books?id=JwQx-WOmSyQC},
  year={1984},
  publisher={Taylor \& Francis}
}

@article{breiman2001random,
  title={Random forests},
  author={Breiman, Leo},
  journal={Machine Learning},
  volume={45},
  number={1},
  pages={5--32},
  year={2001},
  publisher={Springer}
}

@article{castillo2017,
author = {Isma{\"e}l Castillo},
title = {{P{\'o}lya tree posterior distributions on densities}},
volume = {53},
journal = {{Annales de l'Institut Henri Poincar{\'e}, Probabilit{\'e}s et Statistiques}},
number = {4},
publisher = {Institut Henri Poincar{\'e}},
pages = {2074--2102},
year = {2017},
journaldoi = {10.1214/16-AIHP784},
journalURL = {https://doi.org/10.1214/16-AIHP784}
}

@article{castillo2022spikeslab,
author = {Isma{\"e}l Castillo and Romain Mismer},
title = {{Spike and slab P{\'o}lya tree posterior densities: Adaptive inference}},
volume = {57},
journal = {Annales de l'Institut Henri Poincar{\'e}, Probabilit{\'e}s et Statistiques},
number = {3},
publisher = {Institut Henri Poincar{\'e}},
pages = {1521--1548},
year = {2021},
journaldoi = {10.1214/20-AIHP1132},
journalURL = {https://doi.org/10.1214/20-AIHP1132}
}

@article{castillo2022,
author = {Isma{\"e}l Castillo and Thibault Randrianarisoa},
title = {{Optional P{\'o}lya trees: {P}osterior rates and uncertainty quantification}},
volume = {16},
journal = {Electronic Journal of Statistics},
number = {2},
publisher = {Institute of Mathematical Statistics and Bernoulli Society},
pages = {6267--6312},
year = {2022},
journaldoi = {10.1214/22-EJS2086},
journalURL = {https://doi.org/10.1214/22-EJS2086}
}

@article{castillo2021uncertainty,
author = {Isma{\"e}l Castillo and Veronika Ro{\v c}kov{\'a}},
title = {{Uncertainty quantification for {B}ayesian CART}},
volume = {49},
journal = {The Annals of Statistics},
number = {6},
journalpublisher = {Institute of Mathematical Statistics},
pages = {3482--3509},
year = {2021},
journaldoi = {10.1214/21-AOS2093},
journalURL = {https://doi.org/10.1214/21-AOS2093}
}

@article{cattaneo2023inference,
      title={Inference with {M}ondrian Random Forests}, 
      author = {Cattaneo, Matias Damian and Klusowski, Jason Matthew and Underwood, William George},
      journal = {Journal of the Royal Statistical Society Series B: Statistical Methodology},
    pages = {qkaf077},
    year = {2025},
    month = {12},
    journalissn = {1369-7412},
    journaldoi = {10.1093/jrsssb/qkaf077},
    journalurl = {https://doi.org/10.1093/jrsssb/qkaf077},
    eprint = {https://academic.oup.com/jrsssb/advance-article-pdf/doi/10.1093/jrsssb/qkaf077/65824218/qkaf077.pdf},
}

@article{chi2022,
author = {Chien-Ming Chi and Patrick Vossler and Yingying Fan and Jinchi Lv},
title = {{Asymptotic properties of high-dimensional random forests}},
volume = {50},
journal = {The Annals of Statistics},
number = {6},
publisher = {Institute of Mathematical Statistics},
pages = {3415--3438},
year = {2022},
journaldoi = {10.1214/22-AOS2234},
journalURL = {https://doi.org/10.1214/22-AOS2234}
}

@Incollection{chipman1997,
  title = 	 {A {B}ayesian approach to {CART}},
  author =       {Chipman, Hugh and George, Edward I. and McCulloch, Robert E.},
  booktitle = 	 {Proceedings of the Sixth International Workshop on Artificial Intelligence and Statistics},
  pages = 	 {91--102},
  year = 	 {1997},
  bookeditor = 	 {Madigan, David and Smyth, Padhraic},
  volume = 	 {R1},
  series = 	 {Proceedings of Machine Learning Research},
  bookmonth = 	 {04--07 Jan},
  bookpublisher =    {PMLR},
  bookurl = 	 {https://proceedings.mlr.press/r1/chipman97a.html},
  booknote =         {Reissued by PMLR on 30 March 2021.},
  address ={Fort Lauderdale, FL, USA}
}

@article{chipman2010,
author = {Hugh A. Chipman and Edward I. George and Robert E. McCulloch},
title = {{BART: Bayesian additive regression trees}},
volume = {4},
journal = {The Annals of Applied Statistics},
number = {1},
publisher = {Institute of Mathematical Statistics},
pages = {266--298},
year = {2010},
journaldoi = {10.1214/09-AOAS285},
journalURL = {https://doi.org/10.1214/09-AOAS285}
}

@Inbook{coifman1995,
author="Coifman, R. R.
and Donoho, D. L.",
bookeditor="Antoniadis, Anestis
and Oppenheim, Georges",
chapter="Translation-Invariant De-Noising",
title="Wavelets and Statistics",
year="1995",
publisher="Springer New York",
address="New York, NY",
pages="125--150",
bookisbn="978-1-4612-2544-7",
bookdoi="10.1007/978-1-4612-2544-7_9",
bookurl="https://doi.org/10.1007/978-1-4612-2544-7_9"
}

@incollection{pmlrcui21,
  title = 	 {{G}{B}{H}{T}: Gradient Boosting Histogram Transform for Density Estimation},
  author =       {Cui, Jingyi and Hang, Hanyuan and Wang, Yisen and Lin, Zhouchen},
  booktitle = 	 {Proceedings of the Thirty-Eighth International Conference on Machine Learning},
  pages = 	 {2233--2243},
  year = 	 {2021},
  journaleditor = 	 {Meila, Marina and Zhang, Tong},
  journalvolume = 	 {139},
  journalseries = 	 {Proceedings of Machine Learning Research},
  journalmonth = 	 {18--24 Jul},
  journalpublisher =    {PMLR},
  journalurl = 	 {https://proceedings.mlr.press/v139/cui21c.html}
}

@incollection{denil2014,
author = {Denil, Misha and Matheson, David and De Freitas, Nando},
title = {Narrowing the Gap: Random Forests in Theory and in Practice},
year = {2014},
journalpublisher = {JMLR.org},
booktitle = {Proceedings of the 31st International Conference on International Conference on Machine Learning},
pages = {665--673},
address = {Beijing, China},
journalseries = {ICML'14} }

@article{donoho1996,
ajournal = "Ann. Statist.",
author = "Donoho, David L. and Johnstone, Iain M. and Kerkyacharian, G{\'e}rard and Picard, Dominique",
journaldoi = "10.1214/aos/1032894451",
journal = "The Annals of Statistics",
journalmonth = "04",
number = "2",
pages = "508--539",
publisher = "The Institute of Mathematical Statistics",
title = "Density estimation by wavelet thresholding",
journalurl = "http://dx.doi.org/10.1214/aos/1032894451",
volume = "24",
year = "1996"
}

@article {ferguson1974,
    AUTHOR = {Ferguson, Thomas S.},
     TITLE = {Prior distributions on spaces of probability measures},
   FJOURNAL = {Ann. Statist.},
  JOURNAL = {The Annals of Statistics},
    VOLUME = {2},
      YEAR = {1974},
     PAGES = {615--629},
      jounralISSN = {0090-5364},
   MRCLASS = {62F15},
MRREVIEWER = {I. J. Good},
}

@article{freund1997,
title = {A Decision-Theoretic Generalization of On-Line Learning and an Application to Boosting},
journal = {Journal of Computer and System Sciences},
volume = {55},
number = {1},
pages = {119--139},
year = {1997},
journalissn = {0022-0000},
journaldoi = {https://doi.org/10.1006/jcss.1997.1504},
journalurl = {https://www.sciencedirect.com/science/article/pii/S002200009791504X},
author = {Yoav Freund and Robert E Schapire}
}

@article{friedman2001,
author = {Jerome H. Friedman},
title = {{Greedy function approximation: {A} gradient boosting machine.}},
volume = {29},
journal = {The Annals of Statistics},
number = {5},
publisher = {Institute of Mathematical Statistics},
pages = {1189--1232},
year = {2001},
journaldoi = {10.1214/aos/1013203451},
journalURL = {https://doi.org/10.1214/aos/1013203451}
}

@ARTICLE{gey2005,
  author={Gey, S. and Nedelec, E.},
  journal={IEEE Transactions on Information Theory}, 
  title={Model selection for CART regression trees}, 
  year={2005},
  volume={51},
  number={2},
  pages={658-670},
  journaldoi={10.1109/TIT.2004.840903}}

@article{ghosal00,
author = "Ghosal, Subhashis and Ghosh, Jayanta K. and van der Vaart, Aad W.",
journaldoi = "10.1214/aos/1016218228",
journal = "The Annals of Statistics",
journalmonth = "04",
number = "2",
pages = "500--531",
publisher = "The Institute of Mathematical Statistics",
title = "Convergence rates of posterior distributions",
journalurl = "http://dx.doi.org/10.1214/aos/1016218228",
volume = "28",
year = "2000"
}

@article {ghosal2001Bern,
    AUTHOR = {Ghosal, Subhashis},
     TITLE = {Convergence rates for density estimation with {B}ernstein
              polynomials},
   JOURNAL = {The Annals of Statistics},
    VOLUME = {29},
      YEAR = {2001},
    NUMBER = {5},
     PAGES = {1264--1280},
      journalISSN = {0090-5364},
     CODEN = {ASTSC7},
   jounralMRCLASS = {62G07 (62G20)},
  journalMRNUMBER = {1873330 (2002i:62079)},
       journalDOI = {10.1214/aos/1013203453},
       journalURL = {http://dx.doi.org/10.1214/aos/1013203453},
}

@article{hanson2006,
author = {Timothy E Hanson},
title = {Inference for Mixtures of Finite Polya Tree Models},
journal = {Journal of the American Statistical Association},
volume = {101},
number = {476},
pages = {1548--1565},
year = {2006},
publisher = {ASA Website},
journaldoi = {10.1198/016214506000000384},
journalURL = {    
        https://doi.org/10.1198/016214506000000384
},
eprint = {    
        https://doi.org/10.1198/016214506000000384
}
}

@ARTICLE{ho1998,
  author={Tin Kam Ho},
  journal={IEEE Transactions on Pattern Analysis and Machine Intelligence}, 
  title={The random subspace method for constructing decision forests}, 
  year={1998},
  volume={20},
  number={8},
  pages={832--844},
  journaldoi={10.1109/34.709601}}

@article{hoeting1999,
 journalISSN = {08834237},
 journalURL = {http://www.jstor.org/stable/2676803},
 author = {Jennifer A. Hoeting and David Madigan and Adrian E. Raftery and Chris T. Volinsky},
 journal = {Statistical Science},
 number = {4},
 pages = {382--401},
 journalpublisher = {Institute of Mathematical Statistics},
 title = {Bayesian Model Averaging: A Tutorial},
 urldate = {2023-10-20},
 volume = {14},
 year = {1999}
}

@article{jeong2020art,
  author  = {Seonghyun Jeong and Veronika Ro{\v c}kov{\'a}},
  title   = {The Art of {BART}: Minimax Optimality over Nonhomogeneous Smoothness in High Dimension},
  journal = {Journal of Machine Learning Research},
  year    = {2023},
  volume  = {24},
  number  = {337},
  pages   = {1--65},
  journalurl     = {http://jmlr.org/papers/v24/22-0382.html}
}

@article{klusowski2024,
author = {Jason M. Klusowski and Peter M. Tian},
title = {Large Scale Prediction with Decision Trees},
journal = {Journal of the American Statistical Association},
volume = {119},
number = {545},
pages = {525--537},
year = {2024},
journalpublisher = {ASA Website},
journaldoi = {10.1080/01621459.2022.2126782},
journalURL = {     
        https://doi.org/10.1080/01621459.2022.2126782
},
journaleprint = {     
        https://doi.org/10.1080/01621459.2022.2126782
}
}

@article {lavine1992,
    AUTHOR = {Lavine, Michael},
     TITLE = {Some aspects of {P}{\'o}lya tree distributions for statistical
              modelling},
   JOURNAL = {The Annals of Statistics},
    VOLUME = {20},
      YEAR = {1992},
    NUMBER = {3},
     PAGES = {1222--1235},
      journalISSN = {0090-5364},
     CODEN = {ASTSC7},
   MRCLASS = {62A15 (62G99)},
  journalMRNUMBER = {1186248 (93k:62006b)},
MRREVIEWER = {A. L. Rukhin},
       journalDOI = {10.1214/aos/1176348767},
       journalURL = {http://dx.doi.org/10.1214/aos/1176348767},
}

@article{Linero2018BayesianRT,
    author = {Linero, Antonio R. and Yang, Yun},
    title = "{Bayesian Regression Tree Ensembles that Adapt to Smoothness and Sparsity}",
    journal = {Journal of the Royal Statistical Society Series B: Statistical Methodology},
    volume = {80},
    number = {5},
    pages = {1087--1110},
    year = {2018},
    journalmonth = {09},
    journalissn = {1369-7412},
    journaldoi = {10.1111/rssb.12293},
    journalurl = {https://doi.org/10.1111/rssb.12293},
    journaleprint = {https://academic.oup.com/jrsssb/article-pdf/80/5/1087/49270333/jrsssb\_80\_5\_1087.pdf},
}

@article{liu2023,
  author  = {Linxi Liu and Dangna Li and Wing Hung Wong},
  title   = {Convergence Rates of a Class of Multivariate Density Estimation Methods Based on Adaptive Partitioning},
  journal = {Journal of Machine Learning Research},
  year    = {2023},
  volume  = {24},
  number  = {50},
  pages   = {1--64},
  journalurl     = {http://jmlr.org/papers/v24/20-060.html}
}

@InProceedings{liu2024,
  title = 	 {Spatial properties of {B}ayesian unsupervised trees},
  author =       {Liu, Linxi and Ma, Li},
  booktitle = 	 {Proceedings of Thirty Seventh Conference on Learning Theory},
  pages = 	 {3556--3581},
  year = 	 {2024},
  editor = 	 {Agrawal, Shipra and Roth, Aaron},
  volume = 	 {247},
  series = 	 {Proceedings of Machine Learning Research},
  month = 	 {30 Jun--03 Jul},
  journalpublisher =    {PMLR},
  journalurl = 	 {https://proceedings.mlr.press/v247/liu24b.html}
}

@article{liu2026b,
  author  = {Linxi Liu and Li Ma},
  title   = {Supplement to ``{O}n forest-type tree ensemble approaches to density learning under a
generalized Bayesian framework''},
  year    = {2026},
}

@article{LJW,
author = {Lu, Luo and Jiang, Hui and Wong, Wing H.},
title = {Multivariate Density Estimation by {B}ayesian Sequential Partitioning},
journal = {Journal of the American Statistical Association},
volume = {108},
number = {504},
pages = {1402-1410},
year = {2013},
journaldoi = {10.1080/01621459.2013.813389},
eprint = {http://dx.doi.org/10.1080/01621459.2013.813389}
}

@article{ma2017,
author = {Li Ma},
title = {{Adaptive Shrinkage in {P}{\'o}lya Tree Type Models}},
volume = {12},
journal = {Bayesian Analysis},
number = {3},
publisher = {International Society for Bayesian Analysis},
pages = {779 -- 805},
year = {2017},
journaldoi = {10.1214/16-BA1021},
journalURL = {https://doi.org/10.1214/16-BA1021}
}

@article{mentch2016,
author = {Mentch, Lucas and Hooker, Giles},
title = {Quantifying Uncertainty in Random Forests via Confidence Intervals and Hypothesis Tests},
year = {2016},
issue_date = {January 2016},
journalpublisher = {JMLR.org},
volume = {17},
number = {1},
journalissn = {1532-4435},
journal = {J. Mach. Learn. Res.},
journalmonth = {jan},
pages = {841--881},
numpages = {41} }

@article{mourtada2020,
author = {Jaouad Mourtada and St{\'e}phane Ga{\"i}ffas and Erwan Scornet},
title = {{Minimax optimal rates for Mondrian trees and forests}},
volume = {48},
journal = {The Annals of Statistics},
number = {4},
journalpublisher = {Institute of Mathematical Statistics},
pages = {2253--2276},
year = {2020},
journaldoi = {10.1214/19-AOS1886},
journalURL = {https://doi.org/10.1214/19-AOS1886}
}

@article{nieto2012,
 journalISSN = {03036898, 14679469},
 journalURL = {http://www.jstor.org/stable/41411163},
 author = {Luis E. Nieto-Barajas and Peter M{\"u}ller },
 journal = {Scandinavian Journal of Statistics},
 number = {1},
 pages = {166--184},
 publisher = {Wiley},
 title = {Rubbery Polya Tree},
 urldate = {2025-02-25},
 volume = {39},
 year = {2012}
}

@article{ommen2017,
author = {Peter Gr{\"u}nwald and Thijs van Ommen},
title = {{Inconsistency of Bayesian Inference for Misspecified Linear Models, and a Proposal for Repairing It}},
volume = {12},
journal = {Bayesian Analysis},
number = {4},
publisher = {International Society for Bayesian Analysis},
pages = {1069--1103},
year = {2017},
journaldoi = {10.1214/17-BA1085},
journalURL = {https://doi.org/10.1214/17-BA1085}
}

@article{paddock2003,
title ={Randomised Polya Tree Models for Nonparametric Bayesian Inference},
author={Susan M. Paddock and Fabrizio Ruggeri and Michael Lavine and Mike West},
year ={2003},
journal ={Statistica Sinica},
volume = {13},
pages = {443--460},
}

@article{randrianarisoa2021smoothing,
      author = {Thibault Randrianarisoa},
title = {Smoothing and adaptation of shifted {P}{\'o}lya tree ensembles},
volume = {28},
journal = {Bernoulli},
number = {4},
publisher = {Bernoulli Society for Mathematical Statistics and Probability},
pages = {2492--2517},
year = {2022},
journaldoi = {10.3150/21-BEJ1426},
journalURL = {https://doi.org/10.3150/21-BEJ1426}
}

@article{scornet2015,
author = {Erwan Scornet and G{\'e}rard Biau and Jean-Philippe Vert},
title = {{Consistency of random forests}},
volume = {43},
journal = {The Annals of Statistics},
number = {4},
publisher = {Institute of Mathematical Statistics},
pages = {1716 -- 1741},
year = {2015},
journaldoi = {10.1214/15-AOS1321},
journalURL = {https://doi.org/10.1214/15-AOS1321}
}

@ARTICLE{scornet2016,
  author={Scornet, Erwan},
  journal={IEEE Transactions on Information Theory}, 
  title={Random Forests and Kernel Methods}, 
  year={2016},
  volume={62},
  number={3},
  pages={1485--1500},
  jouranldoi={10.1109/TIT.2016.2514489}}

@article{scott1979,
     jstor_articletype = {research-article},
     title = {On Optimal and Data-Based Histograms},
     author = {Scott, David W.},
     journal = {Biometrika},
     jstor_issuetitle = {},
     volume = {66},
     number = {3},
     jstor_formatteddate = {Dec., 1979},
     pages = {605-610},
     journalurl = {http://www.jstor.org/stable/2335182},
     journalISSN = {00063444},
     language = {English},
     year = {1979},
     publisher = {Biometrika Trust},
    }

@article{scott1985ASH,
author = {David W. Scott},
title = {{Averaged shifted histograms: effective nonparametric density estimators in several dimensions}},
volume = {13},
journal = {The Annals of Statistics},
number = {3},
publisher = {Institute of Mathematical Statistics},
pages = {1024--1040},
year = {1985},
journaldoi = {10.1214/aos/1176349654},
journalURL = {https://doi.org/10.1214/aos/1176349654}
}

@book{scott1992,
  author    = {Scott, David W.},
  title     = {Multivariate Density Estimation: Theory, Practice, and Visualization},
  publisher = {Wiley},
  year      = {1992}
}

@article{shen1994,
     jstor_articletype = {research-article},
     title = {Convergence Rate of Sieve Estimates},
     author = {Shen, Xiaotong and Wong, Wing Hung},
     journal = {The Annals of Statistics},
     jstor_issuetitle = {},
     volume = {22},
     number = {2},
     jstor_formatteddate = {Jun., 1994},
     pages = {580--615},
     eprint = {http://www.jstor.org/stable/2242281},
     journalISSN = {00905364},
     language = {English},
     year = {1994},
     publisher = {Institute of Mathematical Statistics},
     copyright = {Copyright © 1994 Institute of Mathematical Statistics},
    }

@article{shen2001,
author = "Shen, Xiaotong and Wasserman, Larry",
journaldoi = "10.1214/aos/1009210686",
journal = "The Annals of Statistics",
journalmonth = "06",
number = "3",
pages = "687--714",
publisher = "The Institute of Mathematical Statistics",
title = "Rates of convergence of posterior distributions",
journalurl = "http://dx.doi.org/10.1214/aos/1009210686",
volume = "29",
year = "2001"
}

@article{shi2024,
      title={Multivariate {G}aussian Approximation for Random Forest via Region-based Stabilization}, 
      author={Zhaoyang Shi and Chinmoy Bhattacharjee and Krishnakumar Balasubramanian and Wolfgang Polonik},
      year={2024},
      eprint={2403.09960},
      journal={arXiv preprint: arXiv2403.09960},
      archivePrefix={arXiv},
      primaryClass={math.ST},
      journalurl={https://arxiv.org/abs/2403.09960}, 
}

@InProceedings{syrgkanis2020,
  title = 	 {Estimation and Inference with Trees and Forests in High Dimensions},
  author =       {Syrgkanis, Vasilis and Zampetakis, Manolis},
  booktitle = 	 {Proceedings of Thirty Third Conference on Learning Theory},
  pages = 	 {3453---3454},
  year = 	 {2020},
  editor = 	 {Abernethy, Jacob and Agarwal, Shivani},
  volume = 	 {125},
  series = 	 {Proceedings of Machine Learning Research},
  month = 	 {09-12 Jul},
  journalpublisher =    {PMLR},
  journalurl = 	 {https://proceedings.mlr.press/v125/syrgkanis20a.html}
}

@article{syring2023,
author = {Nicholas Syring and Ryan Martin},
title = {{Gibbs posterior concentration rates under sub-exponential type losses}},
volume = {29},
journal = {Bernoulli},
number = {2},
publisher = {Bernoulli Society for Mathematical Statistics and Probability},
pages = {1080 -- 1108},
year = {2023},
journaldoi = {10.3150/22-BEJ1491},
jounralURL = {https://doi.org/10.3150/22-BEJ1491}
}

@article {tribouley1995,
author = {Tribouley, Karine},
title = {Practical estimation of multivariate densities using wavelet methods},
journal = {Statistica Neerlandica},
volume = {49},
number = {1},
publisher = {Blackwell Publishing Ltd},
journalissn = {1467-9574},
journalurl = {http://dx.doi.org/10.1111/j.1467-9574.1995.tb01454.x},
journaldoi = {10.1111/j.1467-9574.1995.tb01454.x},
pages = {41--62},
year = {1995},
}

@article{veronika2020,
author = {Veronika Ro{\v c}kov{\'a} and St{\'e}phanie van der Pas},
title = {{Posterior concentration for Bayesian regression trees and forests}},
volume = {48},
journal = {The Annals of Statistics},
number = {4},
publisher = {Institute of Mathematical Statistics},
pages = {2108--2131},
year = {2020},
journaldoi = {10.1214/19-AOS1879},
journalURL = {https://doi.org/10.1214/19-AOS1879}
}

@article{rockova2021ideal,
author = {Veronika Ro{\v c}kov{\'a} and Judith Rousseau},
title = {Ideal Bayesian Spatial Adaptation},
journal = {Journal of the American Statistical Association},
journalvolume = {0},
journalnumber = {ja},
journalpages = {1--27},
year = {forthcoming},
publisher = {Taylor & Francis},
journaldoi = {10.1080/01621459.2023.2241705},

journalURL = { 
    
        https://doi.org/10.1080/01621459.2023.2241705
    
    

},
journaleprint = { 
    
        https://doi.org/10.1080/01621459.2023.2241705
    
    

}

}

@incollection{pmlr-wen22c,
  title = 	 {Random Forest Density Estimation},
  author =       {Wen, Hongwei and Hang, Hanyuan},
  booktitle = 	 {Proceedings of the 39th International Conference on Machine Learning},
  pages = 	 {23701--23722},
  year = 	 {2022},
  bookeditor = 	 {Chaudhuri, Kamalika and Jegelka, Stefanie and Song, Le and Szepesvari, Csaba and Niu, Gang and Sabato, Sivan},
  bookvolume = 	 {162},
  bookseries = 	 {Proceedings of Machine Learning Research},
  journalmonth = 	 {17--23 Jul},
  bookpublisher =    {PMLR},
  journalurl = 	 {https://proceedings.mlr.press/v162/wen22c.html},
  address={Baltimore, Maryland}
}

@article{wong1995,
author = "Wong, Wing Hung and Shen, Xiaotong",
journaldoi = "10.1214/aos/1176324524",
journal = "The Annals of Statistics",
journalmonth = "04",
number = "2",
pages = "339--362",
publisher = "The Institute of Mathematical Statistics",
title = "Probability Inequalities for Likelihood Ratios and Convergence Rates of Sieve {M}{L}{E}S",
journalurl = "http://dx.doi.org/10.1214/aos/1176324524",
volume = "23",
year = "1995"
}

@article{wong2010,
author = "Wong, Wing Hung and Ma, Li",
journaldoi = "10.1214/09-AOS755",
journal = "The Annals of Statistics",
journalmonth = "06",
number = "3",
pages = "1433--1459",
publisher = "The Institute of Mathematical Statistics",
title = "Optional {P}{\'o}lya tree and {B}ayesian inference",
journalurl = "http://dx.doi.org/10.1214/09-AOS755",
volume = "38",
year = "2010"
}

@article{zhang2006a,
 journalISSN = {00905364},
 journalURL = {http://www.jstor.org/stable/25463505},
 author = {Tong Zhang},
 journal = {The Annals of Statistics},
 number = {5},
 pages = {2180--2210},
 publisher = {Institute of Mathematical Statistics},
 title = {From $\epsilon$-Entropy to {KL}-Entropy: Analysis of Minimum Information Complexity Density Estimation},
 urldate = {2022-11-28},
 volume = {34},
 year = {2006}
}

@ARTICLE{zhang2006b,
  author={Tong Zhang},
  journal={IEEE Transactions on Information Theory}, 
  title={Information-theoretic upper and lower bounds for statistical estimation}, 
  year={2006},
  volume={52},
  number={4},
  pages={1307-1321},
  journaldoi={10.1109/TIT.2005.864439}}


\end{appendix}
\end{document}